\documentclass[phd,lfcs,logo,normalheadings,notimes,twoside]{infthesis-nesta}

\usepackage[british]{babel}

\usepackage{mlmodern} 
\usepackage[T1]{fontenc}

\usepackage{titlesec}
\titleformat{\chapter}[display]
{\normalfont\Large\bfseries}
{\chaptername\ \thechapter}
{-10pt}
{\Huge}
\titlespacing*{\chapter}
{0pt}{1ex}{1ex}  
\titleformat{\section}
{\normalfont\large\bfseries}{\thesection}{1em}{}
\titlespacing*{\section}
{0pt}{3ex}{.75ex}  
\titleformat{\subsection}
{\normalfont\fontsize{12}{15}\bfseries}{\thesubsection}{1em}{}
\titlespacing*{\subsection}
{0pt}{2.5ex}{.2ex}  

\renewcommand{\thesection}{\thechapter.\arabic{section}}
\renewcommand{\thesubsection}{\thesection.\arabic{subsection}}

\renewcommand{\thechapter}{\arabic{chapter}}

\counterwithout{figure}{chapter}
\makeatletter
\renewcommand{\p@subsection}{}
\renewcommand{\p@subsubsection}{}
\makeatother
\counterwithout{footnote}{chapter}

\usepackage[leqno]{amsmath}
\usepackage{amssymb}
\usepackage{amsthm}
\usepackage{mathtools}
\usepackage{stmaryrd}
\usepackage[mathcal]{euscript}

\usepackage{caption}
\usepackage{subcaption}
\usepackage{graphicx}
\usepackage{epigraph}
\usepackage[dvipsnames]{xcolor}
\definecolor{LinkColor}{HTML}{697B8C}
\definecolor{CiteColor}{HTML}{8C6B69}
\definecolor{URLColor}{HTML}{898C69}
\usepackage[colorlinks=true,linkcolor=LinkColor
,citecolor=CiteColor
,urlcolor=URLColor
]{hyperref}
\usepackage[hyphenbreaks]{breakurl} 
\usepackage[capitalise,nameinlink,noabbrev]{cleveref}

\usepackage[most]{tcolorbox}
\definecolor{HighlightColor}{HTML}{DDFFCD}
\tcbset{on line, toprule=-2pt, bottomrule=-2pt,
	boxsep=2pt, left=0pt,right=0pt,top=0pt,bottom=0pt,
	colframe = white, colback=HighlightColor,  
	highlight math style={enhanced}
}
\newcommand{\highlighttikzcd}[1]{\tcbhighmath{#1}}

\makeatletter
\newcommand{\dummylabel}[2]{\def\@currentlabel{#2}\label{#1}}
\makeatother

\usepackage[backend=bibtex,style=alphabetic,defernumbers=true,doi=false,isbn=false,url=false,useprefix=true,sorting =anyt,backref=true,giveninits=true,maxnames=99,maxalphanames=99]{biblatex} 
\bibliography{bibliography.bib}
	
\usepackage{enumitem} 
\setlist[enumerate]{topsep=0.5ex,itemsep=-0.5ex,partopsep=1ex,font=\normalfont,leftmargin=2\parindent}

\theoremstyle{plain}
\newtheorem{theorem}{Theorem}[chapter]   
\newtheorem*{theorem*}{Theorem}
\newtheorem{proposition}[theorem]{Proposition}
\newtheorem{corollary}[theorem]{Corollary}
\newtheorem{lemma}[theorem]{Lemma}
\newtheorem{conjecture}[theorem]{Conjecture}
\theoremstyle{definition}
\newtheorem{definition}[theorem]{Definition}
\newtheorem{example}[theorem]{Example}
\newtheorem{construction}[theorem]{Construction}
\theoremstyle{remark}
\newtheorem{remarknumbered}[theorem]{Remark}

\makeatletter
\newtheoremstyle{remarknamedstyle}%
{\topsep}{\topsep}
{}{}
{\itshape}{.}
{0.5em}
{\thmname{\@ifempty{#3}{#1}\@ifnotempty{#3}{#3}}}
\makeatother
\theoremstyle{remarknamedstyle}
\newtheorem*{remark}{Remark}
\makeatletter
\newtheoremstyle{theoremnamedstyle}%
{\topsep}{\topsep}
{\itshape}{}
{\bfseries}{.}
{0.5em}
{\thmname{\@ifempty{#3}{#1}\@ifnotempty{#3}{#3}}}
\makeatother
\theoremstyle{theoremnamedstyle}
\newtheorem*{theoremnamed}{Theorem}

\crefname{theoremnamed}{Theorem}{Theorems}
\crefname{construction}{Construction}{Constructions}

\usepackage{wrapfig}
\usepackage{tikz,tikz-cd,quiver}
\usetikzlibrary{backgrounds}
\usetikzlibrary{arrows,arrows.meta}
\usetikzlibrary{shapes,shapes.geometric,shapes.misc}
\usetikzlibrary{intersections}
\usetikzlibrary{decorations.pathmorphing,decorations.pathreplacing}
\usetikzlibrary{decorations.markings}
\usetikzlibrary{patterns,patterns.meta}
\pgfdeclarelayer{edgelayer}
\pgfdeclarelayer{nodelayer}
\pgfsetlayers{background,edgelayer,nodelayer,main}
\tikzstyle{none}=[inner sep=0mm]
\tikzstyle{tikzfig}=[baseline=-0.25em,scale=0.5]
\tikzstyle{RRed_text}=[text={rgb,255: red,190; green,49; blue,26}]
\tikzstyle{dot_point}=[fill=black, draw=none, shape=circle, inner sep=1pt]
\tikzstyle{thick_dotted}=[-, dotted, thick]
\tikzstyle{thin_dotted}=[-, dotted]
\tikzstyle{green_open}=[-, fill={rgb,255: red,225; green,255; blue,235}, very thick]
\tikzstyle{green_cone}=[-, draw=none, fill={rgb,255: red,176; green,255; blue,176}, fill opacity=.5]
\tikzstyle{blue_cone}=[-, draw=none, fill={rgb,255: red,205; green,249; blue,255}, fill opacity=.34]
\tikzstyle{blue_open}=[-, fill={rgb,255: red,239; green,255; blue,255}, very thick]
\tikzstyle{yellow_opencone}=[fill={rgb,255: red,255; green,255; blue,220}, draw=none]
\tikzstyle{yellow_open}=[-, fill={rgb,255: red,255; green,255; blue,220}, fill opacity=1, very thick]
\tikzstyle{yellow_cone}=[-, draw=none, fill={rgb,255: red,255; green,244; blue,185}, fill opacity=.5]
\tikzstyle{yellow_cone_light}=[-, draw=none, fill={rgb,255: red,255; green,244; blue,185}, fill opacity=.2]
\tikzstyle{red_open}=[-, fill={rgb,255: red,255; green,228; blue,228}, very thick]
\tikzstyle{red_cone}=[-, fill={rgb,255: red,255; green,230; blue,230}, fill opacity=.5, draw=none]
\tikzstyle{red_dashed}=[-, fill={rgb,255: red,255; green,228; blue,228}, very thick, draw={rgb,255: red,71; green,26; blue,26}, dash pattern={{on 10pt off 2pt on 1pt off 2pt}}]
\tikzstyle{grey_open}=[-, fill={rgb,255: red,240; green,240; blue,240}, very thick]
\tikzstyle{grey_cone}=[-, draw=none, fill={rgb,255: red,218; green,218; blue,218}, fill opacity=.5]
\tikzstyle{real_line}=[-, ultra thick]
\tikzstyle{real_line_blue}=[ultra thick, draw={rgb,255: red,135; green,197; blue,255}, {|-|}]
\tikzstyle{real_line_red}=[ultra thick, draw={rgb,255: red,191; green,0; blue,3}, {|-|}]
\tikzstyle{real_line_yellow}=[ultra thick, draw={rgb,255: red,217; green,145; blue,0}, {|-|}]
\tikzstyle{real_line_grey}=[ultra thick, draw={rgb,255: red,171; green,171; blue,171}, {|-|}]
\tikzstyle{preimage_frames}=[->]
\newcommand{\tikzfig}[1]{{\tikzstyle{every picture}=[tikzfig]\input{#1.tikz}}}

\newcommand{\Opens}{\mathop{\mathcal{O}\mspace{-4mu}}}

\newcommand{\Powerset}{\mathop{\mathcal{P}\mspace{-1mu}}}
\newcommand{\Sl}{\mathop{\mathrm{S}\ell}\mspace{-1mu}}
\DeclareMathOperator{\graph}{graph}
\DeclareMathOperator{\im}{im}
\DeclareMathOperator{\id}{id}
\newcommand{\pr}{\mathrm{pr}}
\newcommand{\sqleq}{\sqsubseteq} 
\newcommand{\sqgeq}{\sqsupseteq}
\DeclareMathOperator{\dom}{dom}
\newcommand{\op}{{\mathrm{op}}}
\DeclareMathOperator{\Ideals}{\mathfrak{I}}
\DeclareMathOperator{\ideal}{\mathfrak{i}}

\newcommand{\aleq}{\mathrel{\triangleleft}}
\DeclareMathOperator{\IP}{IP}
\DeclareMathOperator{\IF}{IF}
\newcommand{\hull}[1]{\widehat{#1}}
\DeclareMathOperator{\MSub}{\calM Sub}

\DeclareMathOperator{\OL}{OL}

\newcommand{\refines}{\Subset}

\def\dashedboxtikz{\tikz\node[draw=black,thick,densely dashed] {\phantom{x}};}
\newcommand{\inlinedashedbox}[1]{\tikz[baseline]{\node[draw=black,thick,densely dashed,anchor=base] {\phantom{#1}};}}
\let\rightarrowtail\relax\newcommand{\rightarrowtail}{\begin{tikzcd}[ampersand replacement=\&,cramped, column sep=1.2em]  \phantom{}\ar[r,tail] \& \phantom{} \end{tikzcd}}

\DeclareMathOperator{\Clopup}{Clop^{\up}}

\renewcommand{\leq}{\leqslant}
\renewcommand{\geq}{\geqslant}
\renewcommand{\epsilon}{\varepsilon}

\newcommand{\calB}{\mathcal{B}}
\newcommand{\calF}{\mathcal{F}}
\newcommand{\calG}{\mathcal{G}}
\newcommand{\calH}{\mathcal{H}}

\newcommand{\calE}{\mathcal{E}}

\newcommand{\calM}{\mathcal{M}}

\newcommand{\calA}{\mathcal{A}}

\newcommand{\simS}{\sim_\mathrm{S}}

\newcommand{\Leq}{\mathrel{\trianglelefteqslant}}
\newcommand{\LeqU}{\mathrel{\trianglelefteqslant_\mathrm{U}}}
\newcommand{\LeqL}{\mathrel{\trianglelefteqslant_\mathrm{L}}}
\newcommand{\LeqEM}{\mathrel{\trianglelefteqslant_\mathrm{EM}}}
\newcommand{\leqU}{\leq_\mathrm{U}}
\newcommand{\leqL}{\leq_\mathrm{L}}
\newcommand{\up}{\mathop{\uparrow}\!}
\newcommand{\down}{\mathop{\downarrow}\!}
\newcommand{\triup}{\text{\rotatebox[origin=c]{-90}{$\triangleleft$}}}
\newcommand{\tridown}{\text{\rotatebox[origin=c]{90}{$\triangleleft$}}}

\usepackage{trimclip}
\makeatletter
\newcommand{\tworightarrow}{\mathrel{\text{\two@rightarrow}}}
\newcommand{\two@rightarrow}{%
	\sbox0{$\m@th\rightarrow$}%
	\smash{\rlap{\kern0.1\wd0\clipbox{{.3\width} {-\height} 0pt {-\height}}{$\m@th\rightarrow$}}}%
	$\m@th\rightarrow$%
}
\makeatletter
\newcommand{\twoleftarrow}{\mathrel{\text{\two@leftarrow}}}
\newcommand{\two@leftarrow}{%
	\sbox0{$\m@th\leftarrow$}%
	\smash{\rlap{\kern0.1\wd0\clipbox{{-.1\width} {-\height} 3pt {-\height}}{$\m@th\leftarrow$}}}%
	$\m@th\leftarrow$%
}

\newcommand{\preUp}{\scalebox{1.2}[1.05]{\rotatebox[origin=c]{90}{$\rightarrowtriangle$}}}
\newcommand{\preDown}{\scalebox{1.2}[1.05]{\rotatebox[origin=c]{-90}{$\rightarrowtriangle$}}}

\newcommand{\Up}{\text{\preUp}}
\newcommand{\Down}{\text{\preDown}}
\newcommand{\Upsub}[1]{\Up_{\mspace{-3.5mu}#1}}
\newcommand{\Downsub}[1]{\Down_{\mspace{-1.5mu}#1}}
\newcommand{\Upi}{\Upsub{i}}
\newcommand{\Downi}{\Downsub{i}}

\newcommand{\Cov}{\mathrm{Cov}}
\newcommand{\CovLeq}{\Cov_{\Leq}}
\newcommand{\Covcaus}{\mathrm{Caus}}
\newcommand{\Covchron}{\mathrm{Chron}}

\newcommand{\sqleqdown}{\mathrel{\sqsubset_\tridown}}
\newcommand{\Dcaus}{D_{\caus}}
\newcommand{\Dchron}{D_{\chron}}
\newcommand{\DCaus}{D_\Covcaus}
\newcommand{\DChron}{D_\Covchron}
\newcommand{\DLeq}{D_{\Leq}}

\newcommand{\cat}[1]{\mathbf{#1}}
\newcommand{\Set}{\cat{Set}}
\newcommand{\Top}{\cat{Top}}
\newcommand{\Frm}{\cat{Frm}}
\newcommand{\Loc}{\cat{Loc}}

\newcommand{\Ord}{\cat{Ord}}
\newcommand{\OrdLoc}{\cat{OrdLoc}}
\newcommand{\OrdTop}{\cat{OrdTop}}

\newcommand{\OC}{\mathrm{OC}}
\newcommand{\Rel}{\cat{Rel}}
\newcommand{\OrdC}{\mathrm{Ord}(\cat{C})}
\newcommand{\OrdD}{\mathrm{Ord}(\cat{D})}
\newcommand{\FBT}{\cat{P}}

\newcommand{\ord}{\mathrm{Ord}}

\newcommand{\Heyt}{\cat{Heyt}}

\newcommand{\coHeyt}{\cat{coHeyt}}
\newcommand{\biHeyt}{\cat{biHeyt}}
\newcommand{\Esakia}{\cat{Esak}}
\newcommand{\coEsakia}{\cat{coEsak}}
\newcommand{\biEsakia}{\cat{biEsak}}

\DeclareMathOperator{\pt}{pt}
\DeclareMathOperator{\pf}{pf}
\DeclareMathOperator{\loc}{Loc}
\DeclareMathOperator{\Kleisli}{K\mspace{-1mu}\ell}
\DeclareMathOperator{\Sh}{Sh}

\newcommand{\caus}{\mathrel{\preccurlyeq}}

\DeclareFontFamily{U}{matha}{\hyphenchar\font45}
\DeclareFontShape{U}{matha}{m}{n}{
	<5> <6> <7> <8> <9> <10> gen * matha
	<10.95> matha10 <12> <14.4> <17.28> <20.74> <24.88> matha12
}{}
\DeclareSymbolFont{matha}{U}{matha}{m}{n}
\DeclareMathSymbol{\prechron}{3}{matha}{'316} 
\newcommand{\chron}{\mathrel{\prechron}}



\makeatletter
\newcommand{\newparallel}{\mathrel{\mathpalette\new@parallel\relax}}
\newcommand{\new@parallel}[2]{%
	\begingroup
	\sbox\z@{$#1T$}
	\resizebox{!}{\ht\z@}{\raisebox{\depth}{$\m@th#1/\mkern-5mu/$}}%
	\endgroup
}
\makeatother

\title{Towards\\Point-Free\\Spacetimes}
\author{Nesta van der Schaaf}
\submityear{2024}\date{\today}

\abstract{In this thesis we propose and study a theory of \emph{ordered locales}, a type of point-free space equipped with a preorder structure on its frame of opens. It is proved that the Stone-type duality between topological spaces and locales lifts to a new adjunction between a certain category of ordered topological spaces and the newly introduced category of ordered locales.

As an application, we use these techniques to develop point-free analogues of some common aspects from the causality theory of Lorentzian manifolds. In particular, we show that so-called \emph{indecomposable past sets} in a spacetime can be viewed as the points of the \emph{locale of futures}. This builds towards a point-free causal boundary construction. Furthermore, we introduce a notion of \emph{causal coverage} that leads naturally to a generalised notion of Grothendieck topology incorporating the order structure. From this naturally emerges a localic notion of \emph{domain of dependence}, which is generally distinct from the traditional notion in spacetimes.}

\begin{document}
\begin{preliminary}
\maketitle
\begin{acknowledgements}
	First, I want to thank my supervisor, Chris Heunen, for all the guidance over the past years. It was a pleasure, and I feel especially grateful for being allowed the freedom to let the project wander from monoidal categories and Hilbert spaces, to spacetimes and point-free topology. One of the things I appreciated most about our weekly meetings is that you kept me on my toes by showing genuine curiosity, not letting me get away with shallow explanations. This taught me not to rush so much (a lesson I am admittedly still learning). Thank you Chris!
	
	On the academic front I also want to thank Prakash Panangaden. The meetings we had in Edinburgh were always fun and inspiring, and I hope that our collaboration will continue in the future.
	
	Thanks also to my thesis examiners, Antony Maciocia and Rui Soares Barbosa, for all their time invested, and for their help with improving the manuscript by spotting typos and clumsy exposition.
	
	Thanks to all my office mates and colleagues from the Forum and Quantum Tunnel. Admittedly, I worked from home a lot, but you made my time in the office (and during Meadows picnics) over the years fun and memorable. Thanks in particular to Malin Altenm\"uller, Pablo Andrés-Martínez, Andrew Beckett, Robert Booth, Carmen Constantin, Matthew Di Meglio, Nuiok Dicaire, Robert Furber, Kengo Hirata, Robin Kaarsgaard, Steinar Laenen, Jean-Simon Pacaud Lemay, Jesse Sigal, and Adithya Sireesh. 
	
	Next, I want to thank the people at the Sejny Summer Institute 2022. This was a formative experience for me (and I'm not just talking about the sauna sessions). Special thanks to Jan Głowacki and Hamed Mohammady, and also to Patrick Fraser, Ryszard Paweł Kostecki and Nicetu Tibau Vidal.
	
	Similarly, thanks also to the people who made me feel welcome at my first-ever conference QPL 2023, and to the people from the Adjoint School 2023 who made my first-ever time in the USA a fun one. Special thanks to Clémence Chanavat, Nathan Haydon, Tein van der Lugt, and Ari Rosenfield.
	
	Finally, I should mention Klaas Landsman and Walter van Suijlekom, who opened my world from mostly physics to include also mathematics, where everything seems to make just a bit more sense.
	
	I arrived in Edinburgh in the middle of September, 2020, during the coronavirus pandemic. After the mandatory two-week self-quarantine I was allowed to leave my room, and got to meet my fellow G/6 flatmates, and later the other inhabitants of Pentland House. It would be difficult to overstate what a great time that was. I could probably fill a whole thesis worth of stories from just that single year, but will spare the reader (if you know you know). Thank you, Abhinav, Adam, Cristian, Nikos, Rahul (honorary Pentland House member), Saba, and Shikha (honorary G/6 member). \emph{``Does she love home, or, Edinburgh?''}

	I need to also acknowledge the long lasting support from Jelle, Jelle and Jonathan. We've known each other for decades, and I'm sure this friendship shall continue for decades more, no matter the distance. 
	
	Lastly, Saba. It would be impossible to put into just a few words the gratitude I have for you. Your unconditional love and support make life bearable. Thank you, for everything.
		
	\hspace*{\fill}Swansea, 10 May 2024
\end{acknowledgements}
\mydeclaration
\tableofcontents
\end{preliminary}

\markboth{Introduction}{}
\chapter{Introduction}\label{section:introduction}
\setlength{\epigraphwidth}{.5\textwidth}
\epigraph{``Is anybody working here\\ in chronotopology?''}{Ursula Le Guin, \emph{The Dispossessed}, p.~88}

The main contribution of this thesis is the introduction of \emph{ordered locales}, and the subsequent development of their theory. This chiefly has two purposes. First, it can be seen as a study in pure mathematics, most notably as a generalisation of \emph{Stone duality} to the setting of ordered topological spaces. Second, we propose ordered locales as a new mathematical framework in which to study \emph{spacetimes} using techniques from lattice theory and logic. In this thesis we combine both these viewpoints, giving a detailed development of the theory of ordered locales (\cref{chapter:ordered locales}), but simultaneously pointing out its implications for the study of point-free spacetimes (\cref{chapter:topftir}). The basic idea is as follows.

\begin{figure}[t]\centering
	\input{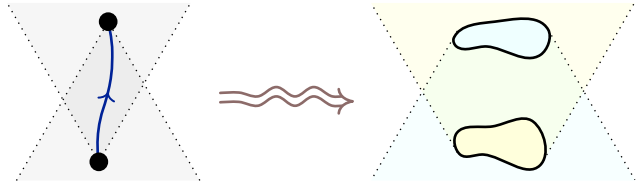}
	\caption{Transition from causal order on points to causal order on regions.}
	\label{figure:idea of order on regions}
\end{figure}

Stone-type dualities fit into the following (incomplete) list. The central theme is to give a \emph{topological representation} to common algebraic or lattice-theoretic objects. Amongst the most famous dualities are:
	\begin{enumerate}[label = \textbullet]
		\item Finite Boolean algebras as finite sets;
		\item \emph{Birkhoff duality:} finite distributive lattices as finite posets \cite{birkhoff1937RingsSets};
		\item \emph{Tarski duality:} complete atomic Boolean algebras as sets \cite{tarski1935ZurGrundlegungBoole};
		\item \emph{Stone duality:} Boolean algebras as Stone spaces \cite{stone1936TheoryRepresentationBooleana};
		\item \emph{Priestley duality:} bounded distributive lattices as Priestley spaces \cite{priestley1972OrderedTopologicalSpaces};
		\item \emph{Esakia duality:} Heyting algebras as Esakia spaces \cite{esakia1974topologicalKripkemodels}.
	\end{enumerate}
Our motivation is, however, in the converse direction. We are not after a topological representation of a certain algebraic object. What we are after is an \emph{algebraic}, \emph{point-free} representation of ordered spaces. Thus we need to fill in the gap of the following list:
	\begin{enumerate}[label = \textbullet]
		\item \emph{Generalised Stone duality:} topological spaces as frames/locales (\cref{section:locales and spaces});
		\item \emph{Ordered Stone duality:} \underline{ordered} topological spaces as $\inlinedashedbox{ordered locales}$.
	\end{enumerate}
In this thesis we provide one possible answer: \emph{ordered locales}, introduced in \cref{section:ordered locales}. Our \emph{modus operandi} for arriving at this notion is explained briefly as follows. The starting point is the definition of an ordered topological space, the object to be dualised, which is simply a set $S$ simultaneously equipped with a topology $\Opens S$ and a preorder $\leq$. Intuitively, locales are what you get when leaving out the underlying point-set in the definition of a topological space (see \cref{section:locales and spaces} for a precise definition). Therefore, the only means of describing what goes on in a locale $X$ is by using its collection of (abstract) \emph{open subsets}~$\Opens X$. Locales form a genuine generalisation of topological spaces in the sense that there exist locales that do not admit any ``points'' whatsoever. It therefore does not immediately make sense to define a preorder $\leq$ on the ``points'' of a locale. Rather, we need to be able to describe the order solely in terms of its open regions. Perhaps the most immediate way to do this is to simply replace the order $\leq$ on \emph{points} by an order~$\Leq$ on \emph{regions}. The visual intuition of this idea is illustrated in \cref{figure:idea of order on regions}. Thus, an ordered locale shall consist of a pair $(X,\Leq)$, where $X$ is a locale, and~$\Leq$ a suitable preorder on its collection of regions. One of our main results shows that this simple procedure indeed leads to a suitable point-free analogue of ordered spaces (\cref{theorem:adjunction ordtopOC bullet and ordloc bullet}).

\section{Why spacetimes as ordered locales?}
The above gives a motivation from duality theory to study ordered locales. But ordered locales were in no small part developed to provide an alternative mathematical framework in which to study spacetimes. The very definition of them is deeply inspired by philosophical considerations on what mathematically (should) constitute spacetime. In this section we outline these ideas. It starts with an argument for the fundamental status of the causal structure (\cref{section:arguments for causal structure}). We then turn to an argument for point-free spaces in physics (\cref{section:arguments for topology}), where in the classical case one argues against the empirical verifiability of points, and in the quantum case one resorts to a measure-theoretic argument. 

Once these arguments are accepted, the question becomes how to define point-free analogues of the mathematical structures typically associated to a spacetime. At first thought this seems like a daunting task: how do we develop a theory of differential geometry on locales? This we leave to future research(ers).\footnote{Although with a recent point-free fundamental theorem of calculus \cite{vickers2023FundamentalTheoremCalculus}, localic differential geometry seems more within reach.} Fortunately, thanks to the work of mathematical relativity theorists of the past century, we know that spacetimes are to a large extent characterised by their underlying \emph{causal order}, which is a much more elementary structure: a \emph{preorder} (\cref{definition:preorders}). So instead we are to answer the much more elementary question: how do we develop a theory of \emph{preordered point-free spaces?} This is precisely what ordered locales provide.

\subsection{Mathematical structure of spacetimes}\label{section:mathematical structure of spacetimes}
\begin{wrapfigure}{R}{.4\textwidth}
	\begin{tikzpicture}[
		> = Stealth,
		box/.style = {draw, semithick, font=\footnotesize,
			minimum height=2.5ex, text width=25mm, align=flush center},
		every edge/.append style = {draw, ->},
		scale=1.1   
		]
		\draw (0,3) node [box](metric){metric $g$};
		\draw (-1.5,1) node[box](causality){causal order $\caus$};
		\draw (1,2) node[box](smooth){smooth structure};
		\draw (1.5,1) node[box](topology){topology $\Opens M$};
		\draw (0,0) node[box](events){set of events $M$};
		\draw[->] (topology) to (events);
		\draw[->] (smooth) to (topology);
		\draw[->] (causality) to (events);
		\draw[->] (metric) to (causality);
		\draw[->] (metric) to (smooth);
		\draw[dashed] (causality) to (topology);
	\end{tikzpicture}
\end{wrapfigure}
We will formally define spacetimes in \cref{section:spacetimes}, but for now provide a brief conceptual overview of their mathematical structure. Modern general relativity is based on the framework of \emph{Lorentzian manifolds}. These are smooth manifolds~$M$, providing differentiable structure, equipped with a \emph{metric}~$g$, providing geometric structure of distance, angle, curvature and geodesic. Each of these ingredients further unpack to more elementary mathematical structures. For instance, the metric is defined on the tangent bundle of $M$, which in turn is defined using the smooth structure (not to speak of the real numbers~$\mathbb{R}$). The smooth structure itself is (traditionally\footnote{As an independent side-note: the theory of \emph{diffeology} provides an elegant notion of smooth structure on a bare set, without presupposing topology, that faithfully subsumes the theory of smooth manifolds \cite{iglesias-zemmour2013Diffeology}. In particular, there exist topologically trivial spaces that nonetheless carry non-trivial smooth information. An example is the irrational Kronecker foliation, whose diffeological structure recovers the same characterisation of the Morita-equivalence classes of irrational rotation algebras~\cite{donato1983exemple} (cf.~\cite{rieffel1981AlgebrasAssociatedIrrational}). Diffeology has also been used to study the foundations of relativity theory \cite{blohmann2013groupoid,glowacki2019groupoidSymmetryGR}, but that is a different story.}) based on a given topology. The metric also induces \emph{lightcone} structure, a type of conformal structure, that in turn induces the \emph{causal order}. And at the bottom of it all lies the underlying point-set $M$ of spacetime events. A more fine-grained exposition on the mathematical structure of spacetimes is found in e.g.~\cite{ehlers2012RepublicationGeometryFree,panangaden2014CausalityPhysicsComputation}. For us the coarser diagram above suffices.

Putting aside the tremendous empirical success of general relativity, it is not straightforward to provide motivation from physical first-principles to support all this mathematical baggage. To put it stronger:~there are arguments \emph{against} using these types of mathematical structures. We provide a brief overview of these arguments from the literature below. On top of that, there are several philosophical problems plaguing the contemporary theory of spacetimes, such as the \emph{hole argument} \cite{landsman2023ReopeningHoleArgumenta}, the problem of \emph{local extendibility} \cite{beem1980MinkowskiSpacetimeLocally,norton2011ObservationallyIndistinguishableSpacetimes}, and the problem of \emph{hole-freeness} \cite{krasnikov2009EvenMinkowskiSpace,manchak2009SpacetimeHolefree,robertsNotesOnHoles}.

To resolve these problems, the general approach is thus to develop a theory of spacetimes by starting with the more elementary structures at the bottom of this diagram, the ontological commitment to which we are more comfortable with. In particular we shall argue for the \emph{causal order}~$\caus$ and the \emph{topology}~$\Opens M$.

\begin{figure}[b]\centering
	\input{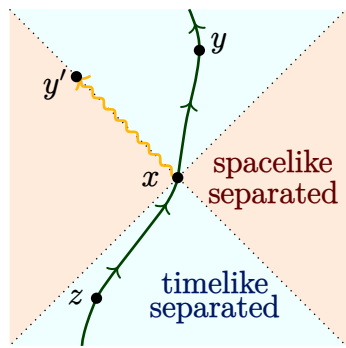}
	\caption{A point $x\in M$ in spacetime with its local causal structure. Here $z\chron x\chron y$ and $x\caus y'$, but $x\not\chron y'$. The interpretation is that a massive particle can travel from $z$ to $y$ through $x$, but only a lightray can travel from $x$ to $y'$.}
	\label{figure:lightcones}
\end{figure}

The causal order is part of what we call the \emph{causal structure} of a spacetime~$M$. It is comprised of two relations on the set of events of the spacetime: the \emph{chronology} $\chron$ and the \emph{causality} $\caus$. Their informal definition is as follows:
	\begin{enumerate}[label = \textbullet]
		\item $x\caus y$ iff $y$ is reachable from $x$ \emph{without} going \emph{faster} than light;
		\item $x\chron y$ iff $y$ is reachable from $x$ going \emph{slower} than light.
	\end{enumerate}
There is also a third relation, called the \emph{horismos:}
	\begin{enumerate}[label = \textbullet]
		\item $x\to y$ iff $y$ is reachable from $x$ going \emph{exactly} at the speed of light;
	\end{enumerate}
but it is fully determined by the causality and chronology, so will play no role in this thesis. Again, the technical details of these definitions are made sense of in \cref{section:spacetimes}, but for now the intuitive picture in \cref{figure:lightcones} suffices.

\subsection{Arguments for causal structure}\label{section:arguments for causal structure}
In this section we review the remarkable result that shows the differentiable and conformal structure of a spacetime is classified by its causal structure.

We briefly define some terminology, whose technical details can be made precise using \cref{section:spacetimes}. A \emph{pseudo-Riemannian manifold} is a smooth manifold $M$ equipped with a \emph{metric (tensor)} $g$, determining a symmetric nondegenerate bilinear form $g_x$ on each tangent space $T_xM$. A \emph{smooth isometry} between pseudo-Riemannian manifolds $h\colon M\to M'$ is a diffeomorphism such that $h_\ast(g) = g'$. This is the appropriate notion of isomorphism between pseudo-Riemannian manifolds: it preserves the smooth structure and the entirety of the metric structure. More laxly, $h$ is called a \emph{smooth conformal isometry} if $h$ is a diffeomorphism, and there exists a smooth non-vanishing map $\Omega\colon M'\to \mathbb{R}$ such that $h_\ast(g) = \Omega^2 g'$. This means $h$ preserves all geometric structure, except possibly for a scaling factor:~$\Omega$.

\emph{Spacetimes} form a particular class of pseudo-Riemannian manifolds. A \emph{continuous timelike curve} is a continuous curve $\gamma\colon I\to M$, defined on an interval $I\subseteq \mathbb{R}$, that suitably respects the chronology relation $\chron$ of the spacetime. (Since $\gamma$ is only continuous, the definition is slightly more subtle than in the smooth case.) We think of these as continuous worldlines of particles moving strictly slower than the speed of light. Similarly, a \emph{continuous null geodesic} is a curve such as $\gamma$ that respects the horismos relation. We think of these as continuous lightrays. A function $h\colon M\to M'$ between spacetimes is said to preserve timelike/null geodesics if $h\circ \gamma$ is timelike/null whenever $\gamma$ is. We now have all the ingredients needed to state the theorems. In the following, all spacetimes are assumed to be of arbitrary but equal dimension.

Hawking's theorem, appearing first in published form as \cite[Theorem~5]{hawking1976NewTopologyCurved}, says the following. (The slightly generalised statement here is from \cite[\S 3]{malament1977ClassContinuousTimelike}.)

\begin{theoremnamed}[Hawking's Theorem]\label{theorem:hawking}
	If $h\colon M\to M'$ is a homeomorphism where both $h$ and the inverse $h^{-1}$ preserve future directed continuous null geodesics, then $h$ is a smooth conformal isometry.
\end{theoremnamed}

We interpret this as saying that the (manifold) topology of a spacetime $M$, together with the data of its lightrays, determines all structure of $M$ up to conformal scaling. This was later substantially generalised by Malament in \cite{malament1977ClassContinuousTimelike}. First, it is noted in their Lemma~3.1 that the condition of preserving continuous null geodesics is actually implied by the preservation of continuous timelike curves. Moreover, it turns out the topological structure is preserved as soon as the continuous timelike curves are preserved. 
\begin{theoremnamed}[Malament's Theorem]\label{theorem:malament}
	If $h\colon M\to M'$ is a bijection between spacetimes, where both $h$ and the inverse $h^{-1}$ preserve future directed continuous timelike curves, then $h$ is a homeomorphism. By \nameref{theorem:hawking} it follows that in this case~$h$ must be a smooth conformal isometry. 
\end{theoremnamed}
This remarkably shows that the class of continuous timelike curves of a spacetime determines almost all of its structure, leaving room only for the conformal scaling factor $\Omega$ to vary.

These results are stated in terms of continuous timelike curves. Malament further elaborates that in a reasonable class of spacetimes, the same characterisation goes through restricting simply to the chronological order $\chron$. Let $I^\pm(x)$ denote the set of those $y\in M$ chronologically in the future/past of $x$ (see \cref{section:causal relations of spacetime}). A spacetime~$M$ is called \emph{distinguishing} if $I^+(x) = I^+(y)$ or $I^-(x)= I^-(y)$ implies $x=y$. A function $h\colon M \to M'$ between spacetimes is called a \emph{causal isomorphism} if $h$ is a bijection and $x\chron y$ iff $h(x)\chron h(y)$. Malament's Lemma~4.3 then says that if $h$ is a causal isomorphism between distinguishing spacetimes, both $h$ and $h^{-1}$ preserve future directed continuous timelike curves. (They furthermore give a counterexample showing the result does not hold outside of the class of distinguishing spacetimes.) This leads us to the conclusion.
\begin{theoremnamed}[Corollary]
	A causal isomorphism between distinguishing spacetimes is a smooth conformal isometry.
\end{theoremnamed}
This shows that the chronology relation of a spacetime is incredibly densely packed with topological and geometric information. With an eye towards abstraction, taking pairs $(M,\chron)$ is hereby very much justified. 

It does raise the question whether to prefer the chronology~$\chron$ or the causality~$\caus$ as the more fundamental notion. This can be debated. The work of Harris on \emph{chronological sets} \cite{harris1998UniversalityFutureChronological,harris2000TopologyFutureChronological} takes $\chron$ as fundamental, while the \emph{causal set} programme \cite{bombelli1987SpacetimeCausalSet} takes $\caus$ as primitive. The paper \cite{kronheimer1967StructureCausalSpaces} combines both into the notion of a \emph{causal space}. A topological view favours causality, since the chronology is recovered as the \emph{interior} of the causal order (\cref{proposition:chronology is interior of causality}). As it happens, in the theory of ordered locales both $\chron$ and $\caus$ determine the same amount of information (\cref{example:localic cones in spacetime}). For that reason, since the theory of preorders is simply more well known, we choose to work with the causal order~$\caus$ as the primitive notion. See however the remarks in \cref{section:causal sites vs ordered locales}.\vspace*{-.5\baselineskip}

\subsection{Arguments for topology}\label{section:arguments for topology}
Before we even resort to arguments from physics, we need to point out that it is widely accepted by some mathematicians that \emph{locales} provide the correct notion of ``space'' over the traditional definition of a topological space. The deep history of this conclusion is sketched in \cite{johnstone1983PointPointlessTopology,johnstone2001ElementsHistoryLocale,vickers2007LocalesToposesSpaces}, which we can only briefly outline here. Here are some ways in which locales behave better than spaces:
	\begin{enumerate}
		\item Any intersection of dense sublocales is again dense \cite[Theorem~1.5]{isbell1972AtomlessPartsSpaces}. Here \emph{sublocale} is the appropriate point-free analogue of sub\emph{set}. A particularly illustrative example of this occurs in the familiar real line $\mathbb{R}$. There we obtain two dense subsets: the rationals $\mathbb{Q}$ and the irrationals $\mathbb{I}$. As subsets we get $\mathbb{Q}\cap \mathbb{I}=\varnothing$, but as sublocales $\widetilde{\mathbb{Q}}\wedge \widetilde{\mathbb{I}}\neq\varnothing$. In fact, it can be shown that $\widetilde{\mathbb{Q}}\wedge \widetilde{\mathbb{I}}=\mathbb{R}_{\neg\neg}$, the sublocale of \emph{regular opens}. From our classical topological intuition we might think at first glance that, while non-empty, the sublocale~$\mathbb{R}_{\neg\neg}$ must be quite close to being trivial. But this is not true:~it is again dense, and what is more, we will see later (\cref{corollary:IPs from regular opens}) that the sublocale $M_{\neg\neg}$ of a spacetime contains all the so-called indecomposable past sets (\cite{geroch1972IdealPointsSpaceTime}), so is in fact rather \emph{non}-trivial.
		
		Another interesting example of such phenomena is as follows \cite[Example~C1.2.8]{johnstone2002SketchesElephantTopos1} (see also \cite{henry2018wonderfulConsequencesMinimalDense}, whose presentation we follow). Proceeding somewhat informally: take any set $A$, and consider the space $A^\mathbb{N}$ of sequences with values in $A$, equipped with the product topology. For every $x\in A$ we can define an open subset $D_x\subseteq A^\mathbb{N}$ consisting of those sequences that attain the value $x$ at least once. This is dense, since the basic opens of $A^\mathbb{N}$ are products $\prod_{n\in\mathbb{N}}U_n$, where for only finitely many $m$ we have $U_m\neq A$. Thus for infinitely many $n$ we get $x\in U_n = A$. The intersection of the induced sublocales is therefore non-empty:
			\[
				\bigwedge_{x\in A}\widetilde{D_x}\neq\varnothing.
 			\]
 		The interpretation of this intersection is the space of \emph{surjective} functions $\mathbb{N}\to A$, which, irrespective of $A$, is thus non-trivial. In particular, we may set $A=\mathbb{R}$, or any other uncountable set. \enlargethispage*{\baselineskip}Of course such functions do not actually exist set-theoretically, which here just means that $\bigwedge_{x\in A}\widetilde{D_x}$ has \emph{no points}. But as a locale it is \emph{non-trivial}.
 		
 		\item Relatedly, the \emph{Banach-Tarski paradox} is a theorem in measure theory that states the unit ball in three-dimensional Euclidean space can be decomposed into finitely many disjoint pieces, reassembled, only to return two identical copies of the original ball \cite{banach1924DecompositionEnsemblesPoints}. This paradox is resolved in the setting of locales \cite{leroy2013ThEorieMesure,simpson2012MeasureRandomnessSublocales}.
 		 
		\item \emph{Tychonoff's theorem} states that the product $\prod_{i\in I} S_i$ of topological spaces is compact if and only if each individual $S_i$ is compact. It is well-known that this statement is equivalent to the Axiom of Choice. But the definition of compactness readily generalises to locales (it can be stated purely in terms of open subsets), and the \emph{localic} Tychonoff theorem holds constructively~\cite{johnstone1981TychonoffTheoremAxiom}.
		\item In the theory of topological groups it is natural to study \emph{closed} subgroups. This is perhaps most apparent in Lie theory, where by \emph{Cartan's theorem} every closed subgroup of a Lie group is necessarily a \emph{Lie} subgroup. In point-free topology it is redundant to restrict to closed subgroups: any subgroup of a localic group is closed~\cite{isbell1988RemarksLocalicGroups}.
		
		\item In classical topology we have the distinct notions of \emph{connectedness} and \emph{path-connectedness}. The topologist's sine curve provides an example of a space that is connected but not path-connected. More generally, there are spaces that are connected and locally connected but not path-connected~\cite[p.~73]{steen1995CounterexamplesTopology}. However, in point-free topology it holds that connected locally connected locales are necessarily path-connected~\cite{moerdijk1986ConnectedLocallyConnected}.
	\end{enumerate}
For more on \citetitle{johnstone1983PointPointlessTopology} we refer to \cite{johnstone1983PointPointlessTopology}.

But point-free spaces can also be motivated from a physical perspective. The idea that the topological structure of physical spacetime is more fundamental than its set of points is by no means new, and has been pointed out by numerous other authors. Appealing to authority, we collect a few quotes. In his lecture notes on general topology and quantum topology, Isham writes:
	\begin{quotation}\small
		``However, reflection~[\ldots]~suggests a more radical possibility. It is clear that what is really being asserted is that the most important feature of space is not the points which it contains but rather the open subsets and the lattice relations between them. But then, since, physically speaking, a ``point'' is a most peculiar concept anyway, why not drop it altogether and deal directly with frames/locales?''\hfill\cite[p.~59]{isham1990IntroductionGeneralTopology}
	\end{quotation}
Sorkin, one of the founders of \emph{causal set theory} \cite{bombelli1987SpacetimeCausalSet}, writes similarly:
	\begin{quotation}\small
		``From an ``operational'' perspective, an individual point of~[a topological space]~$S$ is a very ideal limit of what we can directly measure. A much better correlate of a single ``position-determination'' would probably be an open subset of $S$. Moreover, even for continuum physics, the individual points (or ``events'') of $S$ exist only as carriers for the topology, and thereby also for higher-level constructs such as the differentiable structure and the metric and ``matter'' fields: not the points per se, but only this kind of relation involving them has physical meaning.''
	\hfill\cite[p.~927]{sorkin1991FinitarySubstituteContinuous}
	\end{quotation}
Indeed, intuitively, if we model spacetime locally by the four-dimensional space~$\mathbb{R}^4$, then to determine the precise coordinates of some point $x\in \mathbb{R}^4$ would require an infinite amount of measurement precision, which seems quite unreasonable from an operational point of view.

From a quantum angle, consider the position wave function of a particle inhabiting spacetime. The amplitude of this wave function determines the probability distribution of finding the particle at a certain position upon measurement. However, since singletons in spacetime have measure zero, the probability of finding the particle at a single point will always be zero. Rather, it only makes sense to assign probabilities to \emph{regions}. A more detailed argument along these lines is in~\cite{forrest1996OntologyTopologyTheory}. On a technical level, this corresponds to the fact that the complete Boolean algebra of measurable subsets of the spacetime, modulo null sets, has no~atoms, cf.~\cite[\S 2]{arntzenius2003QuantumMechanicsPointless}.

The existence of spacetime points has also been debated in quantum field theory~(QFT):
	\begin{quotation}\small
		``There is another problem which has to be faced. The quantum field~$\Phi$ at a point cannot be an honest observable. Physically this appears evident because a measurement at a point would necessitate infinite energy.''
		\hspace*{\fill}\cite[\S I.5.2]{haag1996LocalQuantumPhysics}
	\end{quotation}
An overview of several no-go theorems about spacetime points in \emph{algebraic} quantum field theory (AQFT) in particular are summarised in \cite[\S 6]{halvorson2007AlgebraicQuantumField}. For instance, it can be shown that under translation covariance the canonical position and momentum fields \emph{commute} at any spacetime point \cite[\S 6.1.1]{halvorson2007AlgebraicQuantumField}. In \cref{section:aqft on ordered locales} we sketch how ordered locales provide the minimal ingredients on which to define an~AQFT. Philosophical arguments against the existence of points related to quantum field theory are also presented in \cite[\S 6]{arntzenius2003QuantumMechanicsPointless}.


We also briefly address the question of \emph{which} topology to consider on spacetime. Most straightforwardly, we shall always take the \emph{manifold} topology. This is the topology characterised by the causal structure of spacetimes via \nameref{theorem:malament}. This is also argued for in \cite{heathcote1988ZeemanGobelTopologies}, and backed by recent results as \cite{sorkin2019ManifoldtopologyCausalOrder}. We shall not consider alternative topologies such as those in \cite{zeeman1967TopologyMinkowskiSpace,gobel1976ZeemanTopologiesSpacetimes,hawking1976NewTopologyCurved}, since the development of our theory relies heavily on the nice relation between the manifold topology and causal structure of the spacetime (\cref{section:causal relations of spacetime}).

Finally, even if one is not convinced by the philosophical arguments presented here, we think that the resulting mathematical structures are nevertheless of interest. The treatment of spacetimes in point-free topology opens up several new toolboxes from already deeply developed mathematical theories. For instance, we will now have access to compactification techniques from \emph{biframe} theory to help construct causal boundaries of spacetimes. Similarly, we shall see domains of dependence emerge from a canonical \emph{Grothendieck topology} on ordered locales. This in turn opens doors to more categorical treatments of spacetime, or to logical axiomatisations. See further \cref{section:motivation for ordered locales}.


\section{Locales and spaces}\label{section:locales and spaces}
In this thesis we study the interplay between two types of spaces: the well-known and pervasive \emph{topological spaces}, and their littler-known point-free counterpart \emph{locales}. In this section we describe their basic definitions and relation.

The modern definition of a topological space in terms of open subsets originates from the early 1900s, with Hausdorff's definition in 1914 \cite{hausdorff1914GrundzuegeMengenlehre} in terms of neighbourhoods and Kuratowski's definition in 1922 \cite{kuratowski1922OperationAnalysisSitus} (at age 26) in terms of closure operators being essentially equivalent the one we use today.

\begin{definition}\label{definition:topology}
	A topological space is a set $S$ equipped with a \emph{topology}:~a collection $\Opens S$ of so-called \emph{open} subsets of $S$, subject to the following axioms:
		\begin{enumerate}[label = (\roman*)]
			\item $\varnothing \in \Opens S$ and $S\in\Opens S$;
			\item if $U,V\in \Opens S$ then $U\cap V\in \Opens S$;
			\item if $U_i\in \Opens S$ then $\bigcup U_i \in \Opens S$.
		\end{enumerate}
	In other words, $\Opens S$ is a collection of subsets closed under finite intersections and arbitrary unions. A function $g\colon S\to T$ between topological spaces is called \emph{continuous} if the preimage $g^{-1}(V)$ is open in $S$ whenever $V\subseteq T$ is open. Together, they form the category\footnote{The basic definitions of category theory that are used in this thesis are outlined in \cref{section:categories}. A textbook reference is \cite{maclane1998CategoriesWorkingMathematician}.} $\Top$ of topological spaces and continuous functions.
\end{definition}

The precise history is of course more complex than can just be attributed to Hausdorff and Kuratowski, and than what we can do justice to here. A more in-depth modern historical overview of the definition of topological spaces can be found in e.g.~\cite{moore2008EmergenceOpenSets}.

\begin{example}
	The \emph{discrete topology} on a set $S$ is just $\Opens S = \Powerset (S)$, the powerset of $S$, which is the largest topology on $S$. Dually, the \emph{codiscrete topology} on $S$ is $\Opens S = \{\varnothing, S\}$, which is the smallest possible topology.
\end{example}

\begin{example}
	Any metric space $(S,d)$ can be viewed as a topological space, by generating a topology out of the open $d$-balls.
\end{example}

\begin{example}
	The real number line $\mathbb{R}$ has a topology generated by the open intervals $(a,b)$. This generalises to higher-dimensional spaces $\mathbb{R}^n$, and we call the topology here the \emph{Euclidean topology}. Even more generally, any manifold $M$ has a canonical topology generated by its local Euclidean structure.
\end{example}

\begin{example}\label{example:subset topology}
	If $S$ is a topological space, the \emph{subset topology} on a subset $A\subseteq S$ is defined as $\Opens A :=\{U\cap A: U\in \Opens S\}$. This is the smallest topology on $A$ that makes the inclusion function $i\colon A\hookrightarrow S$ continuous.
\end{example}

The definition of a locale is essentially that of a topological space, except that it is stated without the existence of the underlying point-set $S$. The way this is made precise is by viewing $\Opens S$ as an entity on its own, independent of the powerset~$\Powerset(S)$. The fundamental structure here is the inclusion partial order relation~($\subseteq$), together with the fact that this admits a finite intersection operation~($\cap$) and arbitrary union operation ($\bigcup$). It is exactly these types of structures that are studied in \emph{lattice theory}, which was developed into some maturity in the 1930s and resulting in a textbook by Birkhoff \cite{birkhoff1948LatticeTheory} (with first edition published in 1940). Already at that time, it was Stone who laid bare the connection between lattice theory and topology \cite{stone1936TheoryRepresentationBooleana,stone1937ApplicationsTheoryBoolean} by proving the now famous \emph{Stone representation theorem} of Boolean algebras (the algebraic systems that describe classical logic). As Johnstone emphasises in his historical overview \cite{johnstone2001ElementsHistoryLocale}, this development
\begin{quotation}\small
	``\ldots forced topologists to take the algebraic aspects of the structure of a topological space seriously, and in so doing provided a powerful impetus towards freeing topologists from the preconception that (classical Euclidean) space could and should be described in terms of the totality of its points.''
\end{quotation}
For more on the history of locale theory we refer to \cite{johnstone2001ElementsHistoryLocale}, and the introduction of \cite{johnstone1982StoneSpacesa}. 

To state the definition of a locale, we therefore first need to describe the lattice theoretic underpinnings. The structures appearing here are heuristically motivated in \cite[\S 1.2]{picado2021NotesPointFreeTopology}. There are two relevant ways to describe the notion of a lattice. Firstly, as a partial order structure $(L,\sqleq)$ (recall \cref{definition:preorders}) that admits certain greatest lower bounds, called \emph{meets}, or admitting certain least upper bounds, called \emph{joins}. Alternatively, one can take the meets and joins as primary, and think of a lattice as a set $L$ with certain (binary) operations $\wedge$ and $\vee$ (from which the order $\sqleq$ is then derived). To make things precise, an appropriate starting point is the following.

\begin{definition}
	A \emph{lattice} $(L,\wedge,\vee)$ is a set $L$ together with binary operations~$\wedge$ and~$\vee$ that are commutative, associative and idempotent:
		\begin{align*}
			x\wedge y = y\wedge x, \qquad  x\wedge (y\wedge z) = (x\wedge y)\wedge z, \qquad  x\wedge x = x;\\
			x\vee y = y\vee x,\qquad  x\vee (y\vee z) = (x\vee y)\vee z, \qquad  x\vee x = x.
		\end{align*}
	A \emph{distributive lattice} is a lattice where the \emph{distributive law} holds:
		\[
		x\wedge (y\vee z) = (x\wedge y)\vee (x\wedge z),
		\quad\text{or equivalently:}\quad
		x\vee(y\wedge z) = (x\vee y)\wedge (x\vee z).
		\]
	Given any lattice $L$ we obtain a partial order $\sqleq$ as follows:
		\[
			x\sqleq y
			\quad\text{if and only if}\quad
			x\wedge y = x
			\quad\text{or equivalently}\quad
			x\vee y = y.
		\]
	We say $L$ is \emph{bounded} if it has a greatest element $\top$ and a least element $\bot$ with respect to this order. This order admits binary meets and joins that precisely recover $\wedge$ and $\vee$, respectively.
	
	A partial order $(L,\sqleq)$ is called a \emph{complete lattice} if it admits arbitrary greatest lower bounds, or equivalently (\cref{theorem:lattice has joins iff it has meets}), arbitrary least upper bounds. For subsets $A\subseteq L$ these are denoted by $\bigwedge A$ and $\bigvee A$, respectively. If $A$ is of the form ${\{ x_i\in L:i\in I\}}$, then we also write these as $\bigwedge_{i\in I}x_i$ and $\bigvee_{i\in I}x_i$. Note that complete lattices are automatically bounded, since $\bot = \bigvee\varnothing$ and $\top = \bigwedge\varnothing$. In a lattice $(L,\wedge,\vee)$, with respect to the the canonical induced order~$\sqleq$, the elements $x\wedge y$ and $x\vee y$ provide the meets and joins of the subset $\{x,y\}$, respectively.
\end{definition}

See \cref{section:lattice theory} for more technical details on lattices. For us, the most important type of lattice is the following.

\begin{definition}
	A \emph{frame} is a complete lattice $L$ in which the \emph{infinite distributive law} holds:
		\[
			x\wedge \bigvee_{i\in I} y_i = \bigvee_{i\in I} (x\wedge y_i).
		\]
	A \emph{morphism of frames} is a function $h\colon L\to M$ that preserves finite meets and arbitrary joins. Explicitly:
		\[
			h(\bot)=\bot,
			\quad
			h(\top) = \top,
			\quad
			h(x\wedge y) = h(x)\wedge h(y),
			\quad\text{and~}
			h\left(\bigvee x_i\right)=\bigvee h(x_i).
		\]
	Together, they form the category $\Frm$ of frames and frame maps.
\end{definition}

\begin{remark}
	While any frame $L$ admits \emph{arbitrary} meets, their morphisms are only required to preserve \emph{finite} ones. Similarly, for abstract reasons (\cref{section:heyting algebras}), any frame~$L$ is canonically equipped with the structure of a complete \emph{Heyting algebra}. But again, the morphisms of frames do not preserve this structure. So while at the object-level these notions are equivalent, they are distinguished at the categorical level. In later chapters we study how the Heyting algebra structure interacts with order structure.
\end{remark}

\begin{example}\label{example:OS is a frame}
	Importantly, if $S$ is a topological space, the collection of open subsets $\Opens S$ is a frame, with
		\[
			\bot = \varnothing,
			\quad
			\top = S,
			\quad
			U\wedge V = U\cap V,
			\quad\text{and~}
			\bigvee U_i = \bigcup U_i.
		\]
	That the infinite distributive law holds follows from the elementary properties of set operations.
\end{example}

Note that whilst $\Opens S$ lives in the ambient lattice $\Powerset (S)$ of arbitrary subsets of~$S$, the general definition of a frame $L$ can be stated completely without such representation. Frames now provide a purely algebraic, order-theoretic means of doing topology. In fact, expanding on \cref{example:OS is a frame}, there exists a functor
	\[
		\Top \longrightarrow \Frm^\op;
		\qquad
		S\longmapsto \Opens S,
	\]
but it lands in the \emph{opposite}\footnote{If $\cat{C}$ is a category, then $\cat{C}^\op$ is the category with the same objects, but whose morphism are formally reversed; see \cref{definition:opposite category}.} category of frames, since if $g\colon S\to T$ is continuous, the induced map of frames $g^{-1}\colon \Opens T\to \Opens S$ goes in the other direction. Thus we may think of $\Frm$ as a category of algebraic representations of a kind of geometric objects living in $\Frm^\op$. This is a type of duality between algebra and geometry, found commonly in mathematics. The definition of a locale is now simply a matter of formality.

\begin{definition}\label{definition:locales}
	The \emph{category of locales} is defined as $\Loc:= \Frm^\op$. Therefore, a \emph{locale} $X$ is entirely represented by its underlying \emph{frame of opens} $\Opens X\in\Frm$. A \emph{morphism of locales} $f\colon X\to Y$ is entirely defined by its underlying map of frames $f^{-1}\colon \Opens Y\to \Opens X$. 
\end{definition}

\vspace*{-.5em}
The precise relation between the categories $\Top$ and $\Loc$ will be worked out in detail in \cref{section:adjunction top and loc}.


\vspace*{-.8em}
\section{The definition of an ordered locale}\label{section:definition ordered locales}
\vspace*{-.4em}
Now that we have set the stage, a quick remark on the mathematical definition of an ``ordered locale,'' for which there is \emph{a priori} no unique approach. Here we list at least three feasible avenues.
\begin{enumerate}[label=(\roman*)]
	\item As we outlined above (\cref{figure:idea of order on regions}): lifting the order $\leq$ on points of a space~$S$ to an order $\Leq$ on the (open) subsets of~$S$. A canonical way to do this is with the so-called \emph{Egli-Milner order}, which is a preorder on the powerset $\Powerset (S)$ defined as follows:\vspace*{-.6em}
		\[
			A\LeqEM B
			\qquad\text{if and only if}\qquad 
			\begin{array}{l}
				\text{$\forall a\in A\exists b\in B: a\leq b$, and}\\
				\text{$\forall b\in B\exists a\in A: a\leq b$.}
			\end{array}
		\]
	See also \cref{section:ordered spaces as ordered locales}. This preorder can then in turn be restricted to the open subsets $\Opens S$. One definition of ``ordered locale'' could follow from a sufficient axiomatisation of this preorder: we take pairs $(X,\Leq)$, where $X$ is a locale, and $\Leq$ is a preorder on the frame of open regions $\Opens X$, subject to suitable axioms. This is the task we undertook in \cite{heunenSchaaf2024OrderedLocales}, and is presented here in \cref{section:ordered locales}.
	\item The axioms of a preorder (reflexivity and transitivity, see \cref{definition:preorders}) can be captured categorically entirely with the use of finite limits. In any category with finite limits (actually pullbacks suffice) we get a notion of \emph{internal preorder}. A natural definition of ``ordered locale'' is then: an internal preorder in the category $\Loc$ of locales.\footnote{We thank Steve Vickers for suggesting to pursue this direction more seriously.} We investigate this approach slightly more generally in \cref{section:internal preorders}.
	\item Dually to the previous approach, in certain categories it is possible to define the notion of an \emph{internal locale} \cite[Section~C1.6]{johnstone2002SketchesElephantTopos2}. An ``ordered locale'' could then be defined as an internal locale in a suitable category $\Ord$ of preordered sets.\footnote{As suggested to us by Walter Tholen. We have not yet had time to investigate this approach, but only remark that it is not immediately obvious how to define internal locales in the category $\Ord$ of preordered sets and monotone functions, since it apparently does not have power objects.} \enlargethispage{2\baselineskip}
\end{enumerate}
How do we know any of these approaches are satisfactory? Our sanity check for a good definition of ``ordered locale'' shall be as follows: they must form a category $\OrdLoc$, in order so that the classic adjunction between topological spaces and locales \emph{``lifts''} to the ordered setting:
\[
\begin{tikzcd}[cramped]
	\OrdTop & \OrdLoc \\
	\Top & {\Loc.}
	\arrow[""{name=0, anchor=center, inner sep=0}, shift left=2, from=2-1, to=2-2]
	\arrow[""{name=1, anchor=center, inner sep=0}, shift left=2, from=2-2, to=2-1]
	\arrow[""{name=2, anchor=center, inner sep=0}, shift left=2, from=1-1, to=1-2]
	\arrow[""{name=3, anchor=center, inner sep=0}, shift left=2, from=1-2, to=1-1]
	\arrow[hook, from=2-1, to=1-1]
	\arrow[hook, from=2-2, to=1-2]
	\arrow["\dashv"{anchor=center, rotate=-90}, draw=none, from=0, to=1]
	\arrow["\dashv"{anchor=center, rotate=-90}, draw=none, from=2, to=3]
\end{tikzcd}
\]
Here $\OrdTop$ is (some) category of ordered topological spaces (\cref{section:ordered spaces}). In this thesis we focus mainly on the first avenue, as set out above. This gives a definition (\cref{definition:ordered locale}) that is easy to handle, but leads to an adjunction (\cref{theorem:adjunction ordtopOC bullet and ordloc bullet}) only after some technical restrictions on the class of ordered spaces (\cref{section:open cone condition}). Fortunately, this class includes smooth spacetimes. The second approach is technically more difficult to work with, but apparently much more general. We provide a sketch of a proof for a lifted adjunction in \cref{section:internal preorders}, and conjecture how it encapsulates the first approach. We also leave investigating the universality of these lifts to future work.



\section{Motivation for ordered locales}\label{section:motivation for ordered locales}

Reiterating what was said at the start of this \nameref{section:introduction}, one motivation for the development of the theory of ordered locales is to be able to state a natural extension of the Stone-type duality between topological spaces and locales to the setting of ordered spaces. \cref{section:adjunction} shall settle this aspect.

But in \cref{section:arguments for topology,section:arguments for causal structure} we have also seen arguments from both mathematics and physics in favour of locales equipped with a causal structure to model spacetimes. Does the framework of ordered locales actually provide a beneficial framework in which to study relativistic causality? To discuss this, we

\subsection{Causality in ordered locales}

Following \cref{section:definition ordered locales}, our main notion of ordered locale shall be a pair $(X,\Leq)$, where $X$ is a locale and $\Leq$ is a preorder on the frame of opens $\Opens X$ that is supposed to be an abstraction of the Egli-Milner order. What kind of typical notions from relativity theory can we investigate with this structure? In this section we summarise how several notions from causality theory can actually be treated quite naturally using ordered locales. 

\begin{enumerate}[label = $\star$]
	\item \emph{Ideal points} \cite{geroch1972IdealPointsSpaceTime} in a spacetime $M$ are precisely the \emph{coprime} elements in the frames $\im(I^+)$ and $\im(I^-)$ of future and past sets (\cref{section:ideal points}). Dually, we prove that ideal points can be recovered precisely as the \emph{points} of the locales $M^\triup$ and $M^\tridown$ of \emph{futures} and \emph{pasts} (\cref{corollary:IPs as primes in spacetime}). This is strong evidence for a localic \emph{causal boundary} construction, the precise definition of which is still an unsettled question, cf.~\cite{flores2011FinalDefinitionCausal}. We also have the remarkable result that ideal points in $M$ are recovered as ideal points in the strictly \emph{pointless} double-negation locale $M_{\neg\neg}$ (\cref{corollary:IPs from regular opens}). 
	
	\item \emph{Strongly causal} spacetimes \cite[Definition~4.74]{minguzzi2019LorentzianCausalityTheory} are classified as those spacetimes that are \emph{convex} as ordered locales (\cref{corollary:strongly causal iff convex}). We obtain a category theoretic characterisation of strong causality in \cref{theorem:strongly causal iff sublocale}: a spacetime $M$ is strongly causal precisely if the canonical map ${M\to M^\triup\times M^\tridown}$ defines a \emph{sublocale} (\cref{section:sublocales}).
	
	\item \emph{Domains of dependence} \cite{geroch1970DomainofDependence} emerge naturally (\cref{section:domains of dependence}) from a newly defined \emph{causal coverage} relation (\cref{section:causal coverage}) of an ordered locale. This structure is naturally associated to a type of \emph{Grothendieck topology} (\cref{theorem:causal grothendieck topology}), which opens up the toolbox of sheaf theory to analyse spacetimes. \emph{Deterministic sheaves} with respect to a causal coverage encode deterministic evolution of initial data (\cref{section:deterministic sheaves}). We further conjecture that the localic domains of dependence have implications for the problem of ``holes'' \cite{krasnikov2009EvenMinkowskiSpace} in spacetime (\cref{section:holes in spacetime}).
\end{enumerate}

Combined, ordered locales provide us a new framework with which we can use lattice-theoretic, topological, sheaf-theoretic, and even logical techniques to analyse spacetimes.

\begin{remark}[Disclaimer]
	The work in this thesis is \emph{classical}. Despite working with locales, at present we have not attempted to define the theory fully constructively. Similarly, we leave generalisations to \emph{toposes} (an even more generalised notion of point-free space \cite{vickers2007LocalesToposesSpaces}) up to future work (but see \cref{section:internal preorders}). Interspersed throughout the thesis are brief speculative remarks aimed in these directions. One might take our work as a first small step towards a more mature, constructive theory of spacetimes.
\end{remark}

\subsection{Ordered locales versus causal sites}\label{section:causal sites vs ordered locales}\enlargethispage{1\baselineskip}
The early stages of the work that led to this thesis were heavily inspired by the definition of a \emph{causal site} in \cite{christensen2005CausalSitesQuantum}. In particular, the founding idea for \cref{section:causal coverage} comes from the notion of ``coverage'' from \cite[Definition~2.12]{christensen2005CausalSitesQuantum}.

That paper contains, to our knowledge, one of the only attempts in the physics literature at a definition of ``spacetime'' in the spirit of point-free topology.%
\footnote{There is a body of literature on spacetimes through the lens of \emph{mereology}, see e.g.~\cite{wake2011SpacetimeMereology}. We are not familiar enough with this field to comment on the relevance to the present work.}
The two-page paper \cite{szabo1986QuantumCausalStructures} also deserves a mention, containing a similar idea. Lastly, we mention also \cite{vakarelov2020PointfreeTheoriesSpace} and references therein, which also contain ideas that are similar in spirit.

The full definition of a causal site is given in \cite[Definition~2.2]{christensen2005CausalSitesQuantum}, but the core idea is that it consists of a set of abstract regions, equipped with an inclusion partial order $\subseteq$ to model subregion inclusions, and a \emph{strict} partial order~$\prec$ on non-empty regions to model the causal relations. Recall that strictness here means \emph{irreflexivity:} $A\not\prec A$ for all $A$ non-empty. If $M$ is a smooth spacetime where there exist no closed timelike curves, i.e.~where the chronology $\chron$ is irreflexive, then the relation
\[
A\prec B
\qquad\text{if and only if}
\qquad
\forall a\in A\forall b\in B: a\chron b
\]
defined on non-empty subsets $A,B\subseteq M$ defines the motivating example of a causal site. 

In a previous version of the current work, we developed a locale-theoretic definition of the definition of a causal site, where the set of regions forms a genuine frame and not just a join-semilattice, and proved that this gives rise to an adjunction with a certain category of \emph{strictly} ordered spaces. While strict orders and partial orders are equivalent in the point-wise setting, one of the downsides strictness brings in the context of regions is that $A\prec B$ implies $A$ and $B$ must be disjoint. Similarly, strictness excludes behaviour of closed causal loops, which is undesirable for a general mathematical theory of point-free ordered spaces. From a practical perspective, preorders are simply more well known.

Ordered locales are based on \emph{preorders}, and thus admit causal relations $A\Leq B$ between overlapping regions, and also admit non-trivial closed loops. One of the main advantages of ordered locales is that the canonical causal coverage is ``pullback stable'' (in the sense of \cref{lemma:properties of coverage from locale}(d)), which is not the case for causal sites \cite[Remark~2.13]{christensen2005CausalSitesQuantum}. The consequence being that ordered locales give rise to a canonical, generalised Grothendieck topology (\cref{theorem:causal grothendieck topology}), whereas causal sites do not.


\subsection{Tensor topology and quantum causality}\label{section:tensor topology}
\emph{Monoidal categories} are categories $\cat{C}$ equipped with a \emph{tensor product} functor $\otimes\colon \cat{C}\times\cat{C}\to \cat{C}$, together with a \emph{tensor unit} $I\in \cat{C}$. Thinking of the objects of $\cat{C}$ as ``systems,'' and arrows as ``processes'' between them, the tensor product~$\otimes$ allows us to think of processes happening in \emph{parallel}, and the tensor unit $I$ is often interpreted as the ``empty'' system. These types of structures form the mathematical foundation of the process theoretic view of (quantum) physics and computation, where the objects of $\cat{C}$ represent physical systems or state spaces, and the arrows represent physical processes or computations between them. For instance, the category $\Set$ of sets and functions can be interpreted as a theory of classical deterministic computation. The tensor product is just the Cartesian product of sets, and the tensor unit is a singleton. In \emph{categorical quantum mechanics} (CQM) one often studies categories arising from Hilbert spaces: for instance, the category
of finite-dimensional Hilbert spaces and unitary maps can be interpreted as the process theory of reversible quantum computation. Here the tensor product is the usual tensor product of Hilbert spaces, and the tensor unit is the one-dimensional Hilbert space. We refer to~\cite{heunen2019CategoriesQuantumTheory} for a more detailed introduction to these ideas.

The framework of \emph{tensor topology} for monoidal categories was developed recently in \cite{enriquemoliner2018SpaceMonoidalCategories,enriquemoliner2020TensorTopology,soaresbarbosa2023SheafRepresentationMonoidal}. One of its main upshots is \cite[Theorem~9.4]{soaresbarbosa2023SheafRepresentationMonoidal}, which tells us that any (suitable) monoidal category $\cat{C}$ can be viewed to live on a `space' $\mathrm{ZI}(\cat{C})$. This allows us to interpret the systems $C\in\cat{C}$ as inhabiting a definite `location' in $\mathrm{ZI}(\cat{C})$. One of the initial questions of this PhD project was: \emph{``what if this space $\mathrm{ZI}(\cat{C})$ was not just a space, but something like a \emph{spacetime}?''} In particular, the goal was to use an abstract causal structure on $\mathrm{ZI}(\cat{C})$ to provide a rigorous way to study the causal structure of abstract process theories $\cat{C}$. The hope is that for categories $\cat{C}$ that model classical computation this would give a new perspective on \emph{concurrency} \cite{lamport1978TimeClocksOrdering,panangaden2014CausalityPhysicsComputation}, while for categories used in CQM this might provide new insights on \emph{categorical causality} \cite{coecke2013CausalCategoriesRelativistically,coecke2016TerminalityImpliesNosignalling,kissinger2017CategoricalSemanticsCausal,pinzani2019CategoricalSemanticsTime,gogioso2021FunctorialEvolutionQuantum,hefford2023PrePromonoidalStructure} or even \emph{quantum causality} and \emph{causal indefiniteness}~\cite{brukner2014QuantumCausality,goswami2018IndefiniteCausalOrder}.

What type of space is $\mathrm{ZI}(\cat{C})$? Two instructive examples are that, for any (locally compact Hausdorff) topological space $S$, we get (see \cite[\S 3]{enriquemoliner2020TensorTopology} for details):
	\[
		\mathrm{ZI}(\Sh(S))\cong \Opens S,
		\qquad\text{and}\qquad
		\mathrm{ZI}(\cat{Hilb}_{C_0(S)})\cong \Opens S.
	\]
Here, $\Sh(S)$ is the Cartesian monoidal category of sheaves on $S$ \cite{maclane1994SheavesGeometryLogic}, and $\cat{Hilb}_{C_0(S)}$ is the monoidal category of Hilbert-modules over the commutative C$^\ast$-algebra of $C_0(S)$ \cite{heunen2018FrobeniusStructuresHilbert}. Observe here that in the first case, where the topological space $S$ is allowed to be arbitrary, the `space' $\mathrm{ZI}(\cat{C})$ does not recover the underlying point-set of~$S$, but rather only its topology $\Opens S$. Thus the space that $\cat{C}$ lives on is inherently of a \emph{point-free} type. In fact, tensor topology allows us to think of monoidal categories as a categorification of locales. Before answering the question of causal structures on $\mathrm{ZI}(\cat{C})$, we therefore first had to answer the decategorified question:~\emph{``what is an appropriate notion of causality on point-free spaces?''} This naturally led to the theory of ordered locales. The categorification of ordered locales to the setting of monoidal categories is left to future work.

\section{Outline of the thesis}
\begin{enumerate}[label = \textbullet,align=left,itemsep=.5ex]
	\item[\cref{chapter:ordered spaces}.] We start with a brief introduction on the theory of ordered topological spaces (\cref{section:ordered spaces}) and spacetimes (\cref{section:spacetimes}). There are no new results here. Most important is the introduction of a class of ordered spaces with \emph{open cones} (\cref{definition:open cones}), and the observation that the causality relation of any smooth spacetime satisfies the open cone condition (\cref{corollary:spacetimes have open cones}).
	
	\item[\cref{chapter:ordered locales}.] The core of this part of the thesis is adopted from the paper \cite{heunenSchaaf2024OrderedLocales}:
		\begin{quote}
			\fullcite{heunenSchaaf2024OrderedLocales},
		\end{quote}
	introducing the main definition of \emph{ordered locales} (\cref{section:ordered locales}).
	
	\cref{section:adjunction} is devoted to proving one of our main results: there is an adjunction between suitable categories of ordered spaces and ordered locales (\cref{theorem:adjunction ordtopOC bullet and ordloc bullet}). 
	
	This chapter also contains some new work, not found in \cite{heunenSchaaf2024OrderedLocales}. Mainly this is the newfound connection between ordered locales and biframes (\cref{corollary:isomorphism reflective biframes}) and the introduction of \emph{parallel ordered} locales (\cref{section:parallel ordered locales}), which play a crucial role in the development of later chapters.
	
	\item[\cref{chapter:topftir}.] In this part we explore to what extent we can develop common notions of causality in the framework of ordered locales.
	
	In \cref{section:secondary structure} we introduce \emph{convexity} (\cref{section:convexity}) and characterise the strong causality condition of spacetimes locale theoretically (\cref{theorem:strongly causal iff sublocale}).
	
	In \cref{section:causal boundaries} we propose the use of ordered locales for a generalised \emph{causal boundary} construction. One of our main observations is that the \emph{indecomposable past} and \emph{future sets} of a spacetime can be recovered as the points of the induced \emph{locales of futures} and \emph{pasts} (\cref{corollary:IPs as primes,corollary:IPs as primes in spacetime}). Remarkably, they are equivalently recovered via the truly point-free locale of regular opens of the spacetime (\cref{corollary:IPs from regular opens}).
	
	Lastly, in \cref{section:causal coverage} we introduce and study the notion of \emph{causal coverages}. Our main result here is that this gives rise to a canonical type of \emph{Grothendieck topology} (\cref{theorem:causal grothendieck topology}), and \emph{domains of dependence} (\cref{section:domains of dependence}). In \cref{section:domains of dependence in spactime} we compare the traditional definition of domains of dependence to the localic version, and find that they are generally distinct.
	
	In the final \cref{chapter:applications} we outline smaller ideas and directions for future research that could not fill a whole section by themselves. In particular we discuss AQFTs on ordered locales (\cref{section:aqft on ordered locales}), \emph{causal Heyting implication} (\cref{section:caucal heyting implication}), sheaves for causal coverages (\cref{section:deterministic sheaves}), and implications for the problem of \emph{hole-freeness} (\cref{section:holes in spacetime}).
	
	\item[\nameref{chapter:appendix}.] Here, we provide several technical details on category theory (\cref{section:categories}), lattice theory (\cref{section:lattice theory}) and locale theory (\cref{section:locale theory}) that are necessary throughout the main text.
	
	In \cref{section:internal preorders} we provide a proof that for certain adjunctions between categories there exists a lift to their induced categories of \emph{internally preordered objects} and monotone morphisms (\cref{theorem:internal ordered adjunction}). We sketch of a proof that in the setting of spaces and locales this generalises the adjunction from \cref{theorem:adjunction ordtopOC bullet and ordloc bullet} in \cref{section:adjunction}. We conjecture that there will be a correspondence between our main definition of ordered locales (\cref{section:ordered locales}) and internal ordered locales with a suitable localic open cone condition (\cref{conjecture:internally ordered locales with OC are externally ordered locales}).
\end{enumerate}

All work was conducted in collaboration with the author's PhD supervisor Chris Heunen. The work in \cref{section:convexity} and the suggestion of finding a localic \emph{causal ladder} (\cref{section:causal ladder}) was additionally in collaboration with Prakash Panangaden. The early stages of the work in \cref{section:causal boundaries} also benefited from discussions with Prakash Panangaden. 
\part{Ordered spaces and spacetimes}\label{chapter:ordered spaces}
\chapter{Ordered spaces}
\label{section:ordered spaces}

\section{Preordered sets}
In this section we introduce the mathematical bedrock of this thesis: the notion of a \emph{preorder}. Preorders are most typically defined within a chosen framework of sets and relations. Recall that a \emph{relation} $R$ between sets $X$ and $Y$ is (defined by) a subset $R\subseteq X\times Y$ of their cartesian product. Those pairs $(x,y)$ in $X\times Y$ that are contained in $R$ are said to be \emph{related} by $R$. As is common, we write $xRy$ whenever $(x,y)\in R$. A relation \emph{on} $X$ is just a relation from $X$ to itself. A preorder on $X$ is just a special type of relation on $X$.

\begin{definition}\label{definition:preorders}
	A \emph{preorder} (or \emph{preordered set}) is a set $P$, equipped with a relation~$\leq$, satisfying the following two axioms:
	\begin{enumerate}[label=(\roman*)]
		\item $x\leq x$ for all $x\in P$; \hfill(\emph{reflexivity})
		\item if $x\leq y$ and $y\leq z$, then $x\leq z$. \hfill (\emph{transitivity})
	\end{enumerate}
	A \emph{partial order} (or \emph{poset}) is a preordered set $P$ satisfying additionally:
	\begin{enumerate}[label=(\roman*)]\setcounter{enumi}{2}
		\item $x\leq y$ and $y\leq x$ implies $x=y$. \hfill(\emph{anti-symmetry})
	\end{enumerate}
\end{definition}

\begin{remarknumbered}\label{remark:preorders are monads}
	A preorder is nothing but a \emph{monad} (cf.~\cref{definition:monads}) in the 2-category $\Rel$ of sets and relations. Indeed, a monad there is just a relation $R\colon P\to P$, and the existence of the unit and multiplication of the monad correspond precisely to the reflexivity and transitivity conditions of \cref{definition:preorders}. Such insights are exploited in the framework of \emph{Monoidal Topology}~\cite{hofmann2014MonoidalTopologyCategorical}. In \cref{proposition:properties of cones} we see another way to interpret order via monads.
\end{remarknumbered}

\begin{definition}\label{definition:monotone function between preorders}
	A function $g\colon P\to Q$ between preordered sets is called \emph{monotone} if $x\leq y$ implies $g(x)\leq g(y)$. We denote by $\Ord$ the category of preordered sets and monotone functions.
\end{definition}

Most of this thesis is dedicated to studying objects that are equipped with \emph{two} orders: an `intrinsic' partial order, describing a type of topology; and an `auxiliary' preorder, capturing what we generically call ``causality.'' The first one, usually denoted $\sqleq$, will be thought of as an abstract \emph{inclusion} relation. Hence we interpret $x\sqleq y$ as ``$x$ is contained in $y$,'' or ``$y$ contains $x$,'' and so forth. On the other hand, the preorder relation $\leq$ (and later also~$\Leq$) will be interpreted as encoding some type of causal information. Hence, whenever $x\leq y$, we say that $y$ is in the \emph{future} of $x$, and similarly that $x$ is in the \emph{past} of~$y$. The intrinsic order~$\sqleq$ will often play a more subdued role in our notation.

A central theme of this thesis is the interplay between order ($\leq$) and regions~($\sqleq$). In a preordered set this can be studied using the following concepts.

\begin{definition}\label{definition:cones}
	For a preordered set $(P,\leq)$, we define the \emph{up cone} and \emph{down cone} of a subset $A\subseteq P$, respectively, as:
	\[
	\up A:=\{x\in P:\exists a\in A: a\leq x\},
	\quad\text{and}\quad
	\down A:=\{y\in P:\exists a\in A: y\leq a \}.
	\]
	For singletons, we sometimes write $\up x:= \up\{x\}$ and $\down y := \down\{y\}$. With this notation, we can alternatively express the future and past cones as $\up A = \bigcup_{x\in A}\up x$ and $\down A= \bigcup_{y\in A}\down y$.
\end{definition}

\begin{remark}
	In the literature, the sets $\up A$ and $\down A$ are also typically called the \emph{upset} and \emph{downset} of $A$, respectively. Keeping the intuition of relativistic causality in mind, we will often also call $\up A$ and $\down A$ the \emph{future} and \emph{past cones} of $A$. The phrase ``cone'' comes from the notion of \emph{light cone} in relativity theory, which we discuss in-depth in \cref{section:causal relations of spacetime}.
\end{remark}

The up and down cones enjoy several nice properties.

\begin{lemma}\label{proposition:properties of cones}\label{lemma:properties of cones}
	If $(P,\leq)$ is a preordered set, then for any subsets $A,B\subseteq P$:
	\begin{enumerate}[label=(\alph*)]
		\item $\up A\subseteq \up B$ and $\down A\subseteq \down B$ whenever $A\subseteq B$;
		\item $A\subseteq \up A$ and $A\subseteq \down A$;
		\item $\up\up A\subseteq \up A$ and $\down\down A\subseteq \down A$.
	\end{enumerate}
	Further, for any family $(A_i)_{i\in I}$ of subsets of $P$:
	\begin{enumerate}[label=(\alph*)]\setcounter{enumi}{3}
		\item $\bigcup_{i\in I}\up A_i = \up\left(\bigcup_{i\in I}A_i\right)$ and $\bigcup_{i\in I}\down A_i = \down\left(\bigcup_{i\in I}A_i\right)$.
	\end{enumerate}
\end{lemma}
\begin{proof}
	The first property follows by construction. The second and third follow by reflexivity and transitivity, respectively. For the final property, note first that $\up \varnothing = \varnothing = \down \varnothing$, so we may assume the index $I$ is non-empty. In that case, for every $j\in I$ we get $A_j\subseteq \bigcup_{i\in I}A_i$, so property~(a) already implies $\bigcup_{i\in I}\up A_i\subseteq \up\left(\bigcup_{i\in I} A_i\right)$. For the converse, pick $x\in \up\left(\bigcup_{i\in I}A_i\right)$. Then we can find $j\in I$ and $a\in A_j$ such that $a\leq x$. But this just means $x\in \up A_j$, and hence $x\in \bigcup_{i\in I}\up A_i$. The argument for the down cones is analogous.
\end{proof}

\begin{remark}
	In other words, we get join-preserving monads $\up(-)\colon\Powerset(P)\to \Powerset(P)$ and $\down(-)\colon \Powerset(P)\to\Powerset(P)$ on the powerset of $P$. Either one of these monads carry all the information of $\leq$, since:
	\[
	x\leq y
	\quad\text{if and only if}\quad
	x\in \down \{y\}
	\quad\text{if and only if}\quad
	y\in\up \{x\}.
	\]
	Note also that properties~(b) and~(c) of course combine to give strict equalities $\up \up A = \up A$ and $\down \down A=\down A$.
\end{remark}

A map of locales $f\colon X\to Y$ (\cref{definition:locales}) is only defined in terms of its \emph{preimage} map $f^{-1}$, rather than something that takes the ``elements'' of $X$ to~$Y$. The \cref{definition:monotone function between preorders} of monotonicity therefore does not readily adopt to this setting. For our purposes, the following lemma provides a suitable region-based alternative. 

\begin{lemma}\label{lemma:monotonicity in terms of cones}
	The following conditions for a function $g\colon P\to Q$ between preordered sets are equivalent:
	\begin{enumerate}[label=(\alph*)]
		\item $g$ is monotone (\cref{definition:monotone function between preorders});
		\item $\up g^{-1}(A)\subseteq g^{-1}(\up A)$ for all $A\subseteq Q$;
		\item $\down g^{-1}(A)\subseteq g^{-1}(\down A)$ for all $A\subseteq Q$.
	\end{enumerate}
\end{lemma}
\begin{proof}
	First we show that (a) implies (b). For that, take an element $y\in \up g^{-1}(A)$, meaning that there exists $x\in g^{-1}(A)$ with $x\leq y$. Monotonicity of $g$ implies $g(x)\leq g(y)$. But $g(x)\in A$, and hence $g(y)\in \up A$. Or, in other words: ${y\in g^{-1}(\up A)}$.
	
	For the converse, suppose that $x\leq y$ in $P$. This translates equivalently to $y\in \up \{ x\}$; moreover note that $x\in g^{-1}(\{g(x)\})$. Hence (b) implies:
	\[
	y\in \up \{x\}\subseteq \up g^{-1}(\{g(x)\})\subseteq g^{-1}(\up\{g(x)\}).
	\]
	But this just means $g(y)\in \up\{g(x)\}$, that is: $g(x)\leq g(y)$. The equivalence between (a) and (c) is proved dually.
\end{proof}
%
%
\section{Ordered spaces with open cones}\label{section:open cone condition}
Before we start with a theory of point-free ordered spaces, we first need to clarify what we mean by ordinary ordered spaces.

\begin{definition}\label{definition:ordered space}
	An \emph{ordered space} is a pair $(S,\leq)$, where $S$ is a topological space, and $\leq$ is a preorder defined on the underlying set of elements of $S$.
\end{definition}

\begin{figure}[t]\centering
	\input{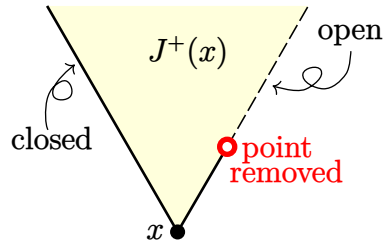}
	\caption{Minkowski space with one point removed. Here $J^+(x)$ is not closed.}
	\label{figure:minkowski space point removed}
\end{figure}

Ordered spaces are important in topology, with particularly seminal work being that of Nachbin; see \cite{nachbin1965TopologyOrder}. Much of ordinary topology generalises to this setting: for instance, one can define \emph{order separation} conditions that (when $\leq$ is equality) recover the familiar separation conditions \cite{mccartan1968SeparationAxiomsTopological}. We will see the \emph{$T_0$-ordered} separation condition later in \cref{definition:T0-ordered}. In particular, the \emph{$T_2$-ordered} separation condition generalises Hausdorffness, and simply states that the graph of $\leq$ is closed in $S\times S$. Such ordered spaces are particularly well-behaved, and have also been used in the study of spacetime causal structure \cite{minguzzi2013ConvexityQuasiuniformizabilityClosed,sorkin2019ManifoldtopologyCausalOrder}. Note however that the causal order $\caus$ of a spacetime generally does not enjoy these properties: \cref{figure:minkowski space point removed} gives a simple example of a non-$T_2$-ordered spacetime. In fact, we will see in \cref{example:minkowski space point removed not T0} that this is not even $T_0$-ordered.


To treat these spacetimes in our framework we therefore do not assume any order separation properties, and indeed our definition of ordered spaces treats the topology and order as totally separate.\footnote{Here we diverge from some literature, such as \cite{nachbin1965TopologyOrder}, where the term ``ordered space'' is reserved for the $T_2$-ordered spaces.} However, it turns out that to be able to define a point-free analogue we need to impose at least some compatibility between them: the \emph{open cone condition}. The precise condition can be motivated in several ways. As we shall see later in \cref{proposition:loc well defined iff OC}, this condition will be necessary to obtain a well-defined adjunction with ordered locales. A more heuristic explanation is as follows. We have seen in \cref{proposition:properties of cones} that the future and past cone operators of a preordered set define monads on its powerset. In this case, we get monads on the powerset $\Powerset (S)$ of a space $S$. However, now there is some topological structure involved, so we should rather look at the lattice $\Opens S$ of \emph{open} subsets of $S$. It then becomes very natural to ask for the up and down cones to restrict to functions $\Opens S\subseteq \Powerset (S)$. This amounts exactly to the following.

\begin{definition}\label{definition:open cones}
	We say an ordered space $(S,\leq)$ has \emph{open cones} if the up and downsets~$\up U$ and~$\down U$ are open whenever $U\subseteq S$ is open.
\end{definition}

\begin{remark}
	It is possible to write down a laxer version of this definition, requiring only the up cones, or only the down cones, to preserve open subsets. Below we encounter some important examples that do not generally satisfy the full open cone condition, but may satisfy one of these laxer conditions. We have a brief \cref{section:upper-lower adjunctions} that discusses this further.
\end{remark}

We discuss some examples.

\begin{example}\label{example:trivial OC}
	Let $(P,\leq)$ be any preordered set. If we equip $P$ with the discrete topology, it is clear that the open cone condition holds, since any subset is open. Similarly, if we equip $P$ with the codiscrete topology the resulting ordered space also has open cones, since $\up P = P = \down P$ and $\up \varnothing = \varnothing =\down \varnothing$ are always open.
	
	Approaching the question from the other extreme, any topological space $S$ becomes an ordered space with equality as the preorder. This has open cones since $\up U = U=\down U$ is open whenever $U$ is open. Similarly, ${\leq} = S\times S$ has open cones, since $\up U = S = \down U$.
\end{example}

\begin{example}
	Consider the real line $\mathbb{R}$ with the Euclidean topology, and standard ordering $\leq$. Using \cref{lemma:open cones via basis} below it is easy to see that this space has open cones. If we instead topologise $\mathbb{R}$ by taking, say, the basis consisting of intervals of the form $(a,\infty)$, for $a\in\mathbb{R}$, we still get a space with open cones, since $\up (a,\infty) = (a,\infty)$ and $\down(a,\infty) = \mathbb{R}$ are both open.
\end{example}

\begin{example}\label{example:interval topology}
	If $(L,\leq)$ is a lattice, its \emph{interval topology} is the one generated by the subbase consisting of subsets of the form $\up x$ and $\down x$, for $x\in L$. The basic opens are the \emph{intervals} $\up x \cap \down y$, containing those $z\in L$ with $x\leq z \leq y$. The resulting ordered space always has open cones. 
\end{example}

\begin{figure}[t]\centering
	\input{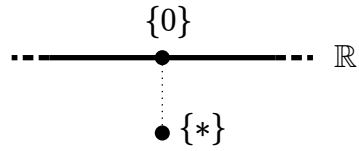}
	\caption{Illustration of the space in \cref{example:non-OC space}.}
	\label{figure:non-OC space}
\end{figure}

\begin{example}\label{example:non-OC space}
	Here is an ordered space that does not have open cones. Consider the real line $\mathbb{R}$ with the Euclidean topology, imagined horizontally. Then adjoin a point ``below'' the origin of this line, see \cref{figure:non-OC space}. Formally, we get the space $\{\ast\}\sqcup \mathbb{R}$, with the preorder generated by $\ast\leq 0$. The future $\up\{\ast\}$ of the adjoined point is the disjoint union $\{\ast\}\sqcup \{0\}$, which is not open since the singleton $\{0\}$ is not open in the Euclidean topology.
\end{example}

There are also important examples from order theory where either only the future or past cones preserve open subsets.

\begin{example}\label{example:upper topology}
	Any preordered set $(P,\leq)$ induces the \emph{upper topology} on $P$, where exactly the sets of the form $\up A$ are open, for $A\subseteq P$. Obviously this has open upper cones, but it will not typically have open lower cones. In fact, this is the case if and only if $\down \up A = \up\down\up A$ for every subset~$A\subseteq P$. This equation clearly fails in the preorder $\{a,b,c,d\}$ where
		\[
			a\leq b\geq c\leq d.
		\]
	Setting $A = \{a\}$ we get $\down \up A = \{a,b,c\}$, but $\up\down\up A$ equals the whole set.
	
	This also shows that the \emph{Scott topology} on $(P,\leq)$, which is the topology where the open subsets are precisely the upsets in $P$ that are \emph{inaccessible by directed joins} \cite[\S 7.3]{vickers1989TopologyLogic}, forms an example of a class of ordered spaces that have upper open cones but not necessarily upper lower cones.
\end{example}

\begin{example}\label{example:specialisation order}
	Let $S$ be any topological space. The \emph{specialisation order} of $S$ is defined as
	\[
	x\leq y
	\qquad\text{if and only if}\qquad
	x\in\overline{\{y\}}.
	\]
	Equivalently, $x\leq y$ if and only if $y\in U$ for every open neighbourhood $U$ of~$x$. The specialisation order is always a preorder, and it is a partial order if and only if $S$ is $T_0$. Similarly, $S$ is a $T_1$-space if and only if the specialisation order is equality.
	
	From the characterisation above it is easy to see that $U = \up U$ for every $U\in \Opens S$, showing that the specialisation order satisfies the upper open cone condition. On the other hand, the lower open cone condition is not always satisfied. Take for instance the preorder $(P,\leq)$ from \cref{example:upper topology} with the upper topology. The specialisation order induced by this topology returns the original preorder~$\leq$~\cite[Proposition~II.1.8]{johnstone1982StoneSpacesa}, but we saw that this ordered space does not have open lower cones. If $S$ is a $T_1$-space then the specialisation order is equality, which trivially has open cones by \cref{example:trivial OC}.
\end{example}

The following lemma shows that it suffices to check the open cone condition on basic opens.

\begin{lemma}\label{lemma:open cones via basis}
	Let $S$ be an ordered space, and let $\mathcal{B}$ be a basis for its topology. Then $S$ has open cones if and only if $\up B$ and $\down B$ are open in $S$ for every $B\in\mathcal{B}$. 
\end{lemma}
\begin{proof}
	The ``only if'' direction is immediate. Conversely, any open $U\subseteq S$ can be written as a union $\bigcup_{i\in I} B_i$ of basic opens $B_i\in\mathcal{B}$. Using \cref{proposition:properties of cones}(d) we then calculate
		\[
			\up U = \bigcup_{i\in I}\up B_i,
		\]
	showing that $\up U$ is a union of opens, and hence itself open. Dually, $\down U$ will be open, and hence $S$ has open cones.
\end{proof}

\begin{example}\label{example:vertical R2}
	The following is an important example of a space with open cones that will later provide instructive counterexamples for generalising definitions from relativity to the localic setting. As a space, we consider $\mathbb{R}^2$ with the standard Euclidean topology. The order structure is defined via
		\[
			(x,y)\leq (a,b)
			\quad\text{if and only if}\quad 
			x=a \text{ and } y\leq b.
		\]
	Hence two points in $\mathbb{R}^2$ are related when they occupy the same position on the first, horizontal axis. We claim that this space has open cones. To see this, recall that a basis for the Euclidean topology on $\mathbb{R}^2$ is given by sets of the form $U\times V$, where $U,V\in\Opens \mathbb{R}$. In turn, we may assume that $V= (b,t)$ is a basic open interval in $\mathbb{R}$. Now it is easy to see that
		\[
			\up\left(U\times (b,t)\right) = U\times (b,\infty)
			\quad\text{and}\quad
			\down\left(U\times (b,t)\right)= U\times (-\infty,t),
		\]
	so the open cone condition follows by \cref{lemma:open cones via basis}.
\end{example}


An important class of examples has been left out so far: spacetimes. We shall prove in \cref{section:OC as pushup} that the pair $(M,\caus)$ of a smooth spacetime equipped with its causal relation always satisfies the open cone condition.

Ending this section, we provide another characterisation of the open cone condition.

\begin{proposition}\label{proposition:open cones iff interiors of open neighbourhoods}
	An ordered space $S$ has open cones if and only if $x\leq y$ implies that $y\in (\up U)^\circ$ and $x\in (\down V)^\circ$ for every open $U\ni x$ and $V\ni y$.
\end{proposition}
\begin{proof}
	Suppose first that $S$ has open cones, and take $x\leq y$. In that case $y\in\up\{x\}$, so if $U$ is an open neighbourhood of $x$, then $y\in\up U = (\up U)^\circ$. 
	
	Conversely, to show $S$ has open cones, we need $\up U\subseteq (\up U)^\circ$ for every open set $U\subseteq S$. If $y\in \up U$, then there exists $x\in U$ with $x\leq y$. The hypothesis then gives $y\in (\up U)^\circ$, as desired. 
\end{proof}

\begin{remark}
	Put another way, in an ordered space with open cones we have for every element $x$ that $\up \{x\}\subseteq \bigcap_{U\ni x}(\up U)^\circ$. In that sense, the open subsets of the space have ``enough room'' to contain the future of $x$. This property fails for instance for $x=\ast$ in the space described in \cref{example:non-OC space}.
\end{remark}

\chapter{Causal structure of spacetimes}
\label{section:spacetimes}
Our main references for the traditional theory of spacetimes are the classic piece \cite{penrose1972TechniquesDifferentialTopology}, the modern review article \cite{minguzzi2019LorentzianCausalityTheory}, and the recent book \cite{landsman2021FoundationsGeneralRelativity}. The latter two have been indispensable for us. We should also mention the classic book \cite{hawking1973LargeScaleStructure}. None of the results stated in this section are new.

A \emph{smooth spacetime} is a type of differential geometric object that models the physical notion that is often described as the \emph{arena} in which physical objects and processes take place. We assume the reader is at least intuitively familiar with the notion of a \emph{(real) smooth manifold}. The full differential geometric machinery of spacetimes is developed throughout \cite{landsman2021FoundationsGeneralRelativity}, and especially in Chapters~2 and~3.

For a smooth manifold $M$, denote by $TM$ and $T^\ast M$ its tangent and cotangent bundles, respectively. 
A \emph{pseudo-Riemannian metric} (or just \emph{metric}) on $M$ is a smooth family of functions ${g_x\colon T_xM\times T_x M\to\mathbb{R}}$, indexed by $x\in M$, that are:
\begin{enumerate}[label = \textbullet]
	\item \emph{bilinear:} for all $v,w\in T_xM$ the maps $g_x(v,-)$ and $g_x(-,w)$ are linear;
	\item \emph{symmetric:} $g_x(v,w)=g_x(w,v)$ for all $v,w\in T_xM$;
	\item \emph{nondegenerate:} $g_x(v,w) = 0$ for all $w\in T_xM$ if and only if $v=0$.
\end{enumerate}
Note importantly that $g_x$ does not have to be \emph{positive definite}, meaning there can exist non-zero tangent vectors $v\in T_x$ such that nevertheless $g_x(v,v)=0$.

The \emph{signature} of $g$ is defined by the number of positive and negative numbers occurring in the diagonal form of $g_x$ in an orthonormal basis of $T_xM$. By continuity, if $M$ is connected, this signature will be independent on $x\in M$. Hence the signature of $g$ is unique, up to sign change. We call a metric $g$ on an $n+1$-dimensional manifold \emph{Lorentzian} if its signature is
	\[
		(-\underbrace{+\cdots +}_\text{$n$ times}).
	\]
The first dimension is interpreted as a \emph{time dimension}, while the other $n$ dimensions are interpreted as \emph{spatial dimensions}. In the four-diensional case, the signature would be denoted $(-+++)$. For the purpose of examples it will be convenient to allow spacetimes that have dimension other than four (but still at least two). Notably, this allows two-dimensional Minkowski space, having signature $(-+)$. Explicitly, $n+1$-dimensional \emph{Minkowski space} is the Euclidean space~$\mathbb{R}^{n+1}$ equipped with the form
	\[
		g(v,w) = -v_0w_0 + \sum_{i=1}^{n} v_iw_i,
	\]
for $v,w\in\mathbb{R}^{n+1}$, written in components as $v=(v_0,v_1,\ldots,v_n)$. Here we conflate the manifold~$\mathbb{R}^{n+1}$ with its tangent spaces.

A \emph{Lorentzian manifold} $(M,g)$ is a smooth manifold $M$ of dimension at least two, equipped with a Lorentzian metric $g$. We say a Lorentzian manifold $(M,g)$ is \emph{time orientable} if there exists a smooth vector field $T$ on $M$ such that ${g_x(T_x,T_x)<0}$ for every $x\in M$. This ensures that there is a globally consistent choice of which direction is the ``future,'' and this notion is again unique up to sign change. Assuming the existence of a time orientability is not much of a restriction, since any Lorentzian manifold has a time orientable double covering \cite[\S 1.6]{penrose1972TechniquesDifferentialTopology}.

We can now finally give the formal definition of a spacetime.

\begin{definition}\label{definition:spacetime}
	A \emph{smooth spacetime} is a connected Lorentzian manifold with time orientation.
\end{definition}

\begin{remark}
	Any non-compact smooth manifold can be equipped with the structure of a spacetime, or if it is compact then as soon as its Euler characteristic is zero \cite[Theorem~1.8]{minguzzi2019LorentzianCausalityTheory}.
\end{remark}


\section{Causality relations in spacetime}\label{section:causal relations of spacetime}
Now that the stage is set, we can finally define the causal relations of spacetime. Our main reference here is \cite[Section~1.11]{minguzzi2019LorentzianCausalityTheory} and \cite[Section~5.3]{landsman2021FoundationsGeneralRelativity}, whose conventions we adopt. For the rest of this section, we fix a smooth spacetime~$(M,g)$. 

A tangent vector $v\in T_xM$ at some point $x\in M$ in the spacetime is called:
	\begin{enumerate}[label = \textbullet]
		\item \emph{timelike} if $g_x(v,v)<0$;
		\item \emph{spacelike} if $v=0$ or $g_x(v,v)>0$;
		\item \emph{lightlike} if $v\neq 0$ and $g_x(v,v)=0$;
		\item \emph{null} if $g_x(v,v)=0$;
		\item \emph{causal} if it is timelike or lightlike, i.e. $g_x(v,v)\leq 0$ and $v\neq 0$.
	\end{enumerate}
If $T$ is a global time orientation on $M$, we further say $v$ is:
	\begin{enumerate}[label = \textbullet]
		\item \emph{future directed} (fd) if $g_x(T_x,v)<0$;
		\item \emph{past directed} (pd) if $g_x(T_x,v)>0$.
	\end{enumerate}
For visual intuition, recall \cref{figure:lightcones}.
Similarly, the terminology generalises to curves $\gamma$ in spacetime by classifying its tangent vectors. Following \cite[\S 1.11]{minguzzi2019LorentzianCausalityTheory}, a piecewise continuously differentiable curve $\gamma\colon I\to M$, defined on any interval $I\subseteq\mathbb{R}$, is called \emph{causal} or \emph{timelike} if all its tangent vectors are causal or timelike, respectively. For simplicity, we shall assume that all causal and timelike curves are \emph{future directed}, which means that each of their tangent vectors are future directed.

The \emph{causal} and \emph{chronological relation} of $M$ are now defined, respectively, as:
	\begin{align*}
		x\caus y
		&\quad\text{if and only if}\quad
		\text{there exists a causal curve from $x$ to $y$, or $x=y$};\\
		x\chron y
		&\quad\text{if and only if}\quad
		\text{there exists a timelike curve from $x$ to $y$}.
	\end{align*}
Since causal and timelike curves compose, the following result is immediate.

\begin{proposition}\label{proposition:causality relations}
	For any smooth spacetime, the causal relation $\caus$ is a preorder, and the chronological relation $\chron$ is transitive. 
\end{proposition}

Further, we define for $x\in M$ its \emph{causal future} and \emph{past} as
	\[
		J^+(x):=\{y\in M: x\caus y\}
		\qquad\text{and}\qquad 
		J^-(x) := \{z\in M: z\caus x\}.
	\]
Similarly, the \emph{chronological future} and \emph{past} are defined as
	\[
		I^+(x):=\{y\in M: x\chron y\}
		\qquad\text{and}\qquad 
		I^-(x) := \{z\in M: z\chron x\}.
	\]
For intuition, it helps to unpack this definition in Minkowski space. For ${x,y\in\mathbb{R}^{n+1}}$ we get that $x\caus y$ precisely when $y-x$ is a causal tangent vector, and this gives:
	\[
		J^+(x) = \{y\in \mathbb{R}^{n+1}: y_0 \geq x_0 + \left\lVert \vec{y}-\vec{x}\right\rVert \},
	\]
where $\vec{x}= (x_1,\ldots, x_n)$ denotes the spatial part of $x$. The formula for $J^-(x)$ is obtained dually. The cones $I^\pm$ are obtained analogously by replacing the inequality~$\geq$ by a strict inequality $>$:
	\[
	I^+(x) = \{y\in \mathbb{R}^{n+1}: y_0 > x_0 + \left\lVert \vec{y}-\vec{x}\right\rVert \},
	\]
so we see that the chronological cones are in fact obtained as the topological interior of the causal cones.

It turns out that generally there is a remarkably strong connection between the manifold topology of a spacetime $M$ and its causal structure. The reason for this is ultimately that any spacetime locally looks like Minkowski space. This is made precise using \emph{normal neighbourhoods}. To explain these, we recall some technology from \cite[\S 1.10]{penrose1972TechniquesDifferentialTopology} and \cite[\S 5.2]{landsman2021FoundationsGeneralRelativity}. The \emph{exponential map} at $x\in M$ is the function $\exp_x\colon V_x\to M$, defined on a suitable open subset $V_x\subseteq T_xM$ of the tangent space, defined by $\exp_x(v):=\gamma(1)$, where $\gamma$ is the (affinely parametrised) geodesic that starts at $x$ with tangent vector $v$. This region contains an open subset $U_x\subseteq V_x$ such that $\exp_x|_{U_x}$ defines a diffeomorphism onto its image in $M$. This image is called a \emph{normal neighbourhood} of $x$.

\begin{proposition}[{\cite[\S2.8]{penrose1972TechniquesDifferentialTopology}}]\label{proposition:chronological cones are open}
	For every $x\in M$, the chronological cones $I^\pm(x)$ are open. Consequently, $I^\pm(A)$ are open for any $A\subseteq M$.
\end{proposition}
\begin{proof}[Proof idea]
	Let $y\in I^+(x)$, so there exists a timelike curve $\gamma\colon [0,1]\to M$ with $\gamma(0)=x$ and $\gamma(1)=y$. We can find some $z$ on the image of $\gamma$ whose normal neighbourhood contains $y$. In particular, we can then find a timelike vector $v\in T_zM$ such that $y= \exp_z(v)$. Since $T_zM$ is homeomorphic to Minkowski space, we can find some open neighbourhood $W$ around $v$ consisting only of timelike vectors. Then $V=\exp_z(W)$ is an open neighbourhood of $y$ contained in $I^+(z)$. This implies $y\in V\subseteq I^+(x)$, showing $I^+(x)$ is open.
\end{proof}

The proofs of the following two results can be found in \cite[Proposition~5.4]{landsman2021FoundationsGeneralRelativity} or Theorems~2.24 and~2.27 in \cite{minguzzi2019LorentzianCausalityTheory}.

\begin{proposition}\label{proposition:push-up}
	In any smooth spacetime:
		\[
			x\chron y\caus z\quad\text{or}\quad x\caus y\chron z
			\quad\text{implies}\quad
			x\chron z.
		\]
\end{proposition}

\begin{remark}
	The above property is also called the \emph{push-up} property. We revisit it in \cref{section:OC as pushup} in relation to the open cone condition. In axiomatic approaches to causality, the relations $\caus$ and $\chron$, together with laws like this, are taken to be fundamental \cite{kronheimer1967StructureCausalSpaces,carter1971CausalStructureSpacetime}.
\end{remark}

\begin{proposition}\label{proposition:chronology is interior of causality}
	For every $A\subseteq M$ in a 	spacetime:
		\begin{align*}
			I^+(A) = \left(J^+(A)\right)^\circ
			&\qquad\text{and}\qquad
			I^-(A) = \left(J^-(A)\right)^\circ;\\
			\overline{I^+(A)} =\overline{J^+(A)}
			&\qquad\text{and}\qquad
			\overline{I^-(A)} =\overline{J^-(A)}.
		\end{align*}
\end{proposition}

\begin{remark}
	The last pair of relations between the closures of the causal and chronological cones show that the chronological cone $I^\pm(x)$ of any point $x\in M$ must be non-empty. Namely, otherwise $x\in J^\pm(x)\subseteq \overline{J^\pm(x)}= \overline{I^\pm(x)}=\varnothing$, a contradiction. This is quite special, and it is not to be expected that $(\up x)^\circ\neq\varnothing$ in an arbitrary ordered space, even with open cones. For instance, in the space with open cones from \cref{example:vertical R2} the set $(\up x)^\circ$ is always empty. 
\end{remark}

\begin{lemma}\label{lemma:chronological cones intersect opens}
	Let $x\in M$ be a point in a smooth spacetime. For every open neighbourhood $U$ of $x$, the set $U\cap I^\pm(x)$ is non-empty.
\end{lemma}
\begin{proof}[Proof idea]
	This property clearly holds in Minkowski space. It generalises to arbitrary smooth spacetimes using \cite[Theorem~5.5]{landsman2021FoundationsGeneralRelativity}, which says that in a normal neighbourhood $W$ of $x$ the chronological future of $x$ can equivalently be computed in $\exp_x^{-1}(W)\subseteq T_xM$, i.e.~in Minkowski space. In particular, for any open neighbourhood $U$ of $x$ there exists $v\in \exp_x^{-1}(U\cap W)$ such that $0\chron v$, and hence $x\chron \exp_x(v)$, proving that $U\cap W\cap I^+(x)$ is non-empty.
\end{proof}

\cref{proposition:causality relations} translates to the following properties about causal and chronological cones.

\begin{proposition}\label{proposition:properties of cuasal and chronological cones}
	If $M$ is a smooth spacetime, then for any subsets $A,B\subseteq M$:
	\begin{enumerate}[label=(\alph*)]
		\item $J^\pm(A)\subseteq J^\pm(B)$ whenever $A\subseteq B$;
		\item $A\subseteq J^\pm(A)$;
		\item $J^\pm(J^\pm(A))\subseteq J^\pm(A)$.
	\end{enumerate}
	Similarly (but note the second property):
	\begin{enumerate}[label=(\alph*)]\setcounter{enumi}{3}
		\item $I^\pm(A)\subseteq I^\pm(B)$ whenever $A\subseteq B$;
		\item $A^\circ \subseteq I^\pm(A)$;
		\item $I^\pm(I^\pm(A))\subseteq I^\pm(A)$.
	\end{enumerate}
\end{proposition}
\begin{proof}
	Statements (a)--(c) are just \cref{proposition:properties of cones} for the causality preorder of a spacetime. Claim~(d) is elementary. The last two claims follow from~(b) and~(c) using \cref{proposition:chronology is interior of causality,proposition:push-up}.
\end{proof}

\begin{remark}
	Since the chronological cones are open (\cref{proposition:chronological cones are open}), properties~(e) and~(f) together imply that, just as for the causal cones, we actually get an equality: $I^\pm(I^\pm(A))=I^\pm(A)$.
	
	Note that the inclusion $A^\circ \subseteq I^\pm(A)$ cannot in general be improved to include the entirety of $A$. For instance, in many spacetimes (e.g.~Minkowski space) it is the case that $x\notin I^\pm(x)$, so it fails for $A=\{x\}$.
\end{remark}

\begin{proposition}[{\cite[\S2.11]{penrose1972TechniquesDifferentialTopology}}]\label{proposition:chronological cone of closure}
	For any subset $A\subseteq M$: $I^\pm(A) = I^\pm(\overline{A})$.
\end{proposition}
\begin{proof}
	The inclusion ``$\subseteq$'' clearly holds by monotonicity. Conversely, take $x\in \overline{A}$ and $x\chron y$. Recall that $x\in \overline{A}$ if and only if every open neighbourhood of $x$ intersects $A$. Now $x\in I^-(y)$ is open by \cref{proposition:chronological cones are open}, so there exists some $z\in A$ such that $z\chron y$. Hence $y\in I^+(A)$. The proof for $I^-$ is dual.
\end{proof}


\begin{definition}\label{definition:past set}
	A subset $F\subseteq M$ of spacetime is called a \emph{future set} if $I^+(F)\subseteq F$. Dually, a \emph{past set} is a subset $P\subseteq M$ such that $I^-(P)\subseteq P$. 
\end{definition}

\begin{remark}\label{remark:examples of future sets}
	For an arbitrary subset $A\subseteq M$, by \cref{proposition:properties of cuasal and chronological cones}(f) $I^+(A)$ is a future set; by \cref{proposition:push-up} $J^+(A)$ is a future set, and by \cref{proposition:chronological cone of closure} even $\overline{J^+(A)}$ is a future set. Similarly, of course, $I^-(A)$, $J^-(A)$, and $\overline{J^-(A)}$ are past sets. 
\end{remark}

The following result from \cite[Proposition~2.84]{minguzzi2019LorentzianCausalityTheory} shows that the interior of future and past sets can equivalently be calculated using chronological cones.

\begin{proposition}\label{proposition:interior of future/past set}
	For any future set $F$ and any past set $P$ in a smooth spacetime: $F^\circ = I^+(F)$ and $P^\circ = I^-(P)$. 
\end{proposition}
\begin{proof}
	By definition, $I^+(F)\subseteq F$, so since chronological cones are open (\cref{proposition:chronological cones are open}) we get $I^+(F)\subseteq F^\circ$. Conversely, if $x\in F^\circ$ then by \cref{lemma:chronological cones intersect opens} there exists $y\in F^\circ$ such that $y\chron x$, and hence $x\in I^+(F^\circ)\subseteq I^+(F)$, as desired. The proof for past sets is dual.
\end{proof}

\begin{remark}
	Together with \cref{proposition:chronological cones are open}, this shows that the open future and past sets are precisely the images of $I^+$ and $I^-$, respectively.
\end{remark}

\subsection{Regularity of the chronology}
The above string of results provide us with a strong connection between the topology and causal structure of a spacetime. The next few results are some elementary, yet remarkable corollaries. These topological properties will later become essential to relate our localic definition of causal boundary (\cref{section:causal boundaries}) to the ones from the spacetime literature \cite{geroch1972IdealPointsSpaceTime}. Briefly, we recall that the \emph{Heyting negation} (\cref{definition:heyting algebra}) in the lattice $\Opens S$ of opens of a topological space $S$ can be calculated as:
	\[
		\neg U = \bigcup\left\{V\in \Opens S: U \cap V = \varnothing \right\} = \left(S\setminus U\right)^\circ.
	\]
In turn, the double Heyting negation of $U\in \Opens S$ may be calculated as the interior of the closure of an open subset $U$:
	\[
		\neg\neg U = \neg (S\setminus U)^\circ 
		=
		\left(S\setminus \left(S\setminus U\right)^\circ\right)^\circ
		=
		(\overline{U})^\circ.
	\]
The next two lemmas show that the chronological cones $I^\pm$ of a spacetime $M$ cancel out the double Heyting negation in $\Opens M$.

\begin{lemma}\label{lemma:chronological cone of double negation}
	For any open $U\subseteq M$ we have $I^\pm(U)=I^\pm(\neg\neg U)$. 
\end{lemma}
\begin{proof}
	Double negation is a monad on open subsets (\cref{lemma:double negation is monad}), and hence we get the inclusion $I^\pm(U)\subseteq I^\pm(\neg\neg U)$. The converse follows using \cref{proposition:chronological cone of closure} for the last equality:
	\[
	I^\pm(\neg\neg U) = I^\pm \left((\overline{U})^\circ \right) \subseteq I^\pm(\overline{U}) = I^\pm(U).
	\qedhere
	\]
\end{proof}

\begin{lemma}\label{lemma:chronological cones are regular opens}
	Let $A\subseteq M$ be a subset of spacetime. Then $\neg\neg I^\pm(A)= I^\pm(A)$. 
\end{lemma}
\begin{proof}
	This can be seen from the following string of calculations:
		\[
			\neg\neg I^\pm(A)
			=
			\left(\overline{I^\pm(A)}\right)^\circ
			=
			\left(\overline{J^\pm(A)}\right)^\circ
			=
			I^\pm \left(\overline{J^\pm(A)}\right) 
			=
			I^\pm (J^\pm(A))
			=
			I^\pm(A).
		\]
	The first step is by definition, the second step is \cref{proposition:chronology is interior of causality}, the third step is \cref{proposition:interior of future/past set}, the fourth step is \cref{proposition:chronological cone of closure}, and the final step is the push-up principle in \cref{proposition:push-up}.
\end{proof}

In summary:
	\[
		\neg\neg \circ I^\pm = I^\pm = I^\pm\circ\neg\neg.
	\]
The significance of these two lemmas will be further explored in \cref{section:causal boundaries}. In brief, they tell us that the chronological cone structure of a spacetime are fully captured on the so-called \emph{double negation sublocale} of $M$ (\cref{section:IPs in double-negation sublocale}). A special feature of this locale is that it is \emph{Boolean}. 

\subsection{Open cones as a push-up principle}\label{section:OC as pushup}
In this section we highlight a connection between the open cone condition and the \emph{push-up} principle from \cref{proposition:push-up}, which we recall states that from either $x\caus y\chron z$ or $x\chron y\caus z$ we may conclude that $x\chron z$. Intuitively this means that once you occupy the interior of $I^+(x)$, you can never again reach its boundary, even going at the speed of light. Alternatively, it says that a causal curve that arises from $x\caus y \chron z$ can be ``shortened'' to a chronological curve $x\chron z$.

The push-up principle can be translated into conditions that relate the cones~$I^\pm$ and $J^\pm$, analogous to how transitivity of a preorder translates to idempotence the cones. It is straightforward to verify that if $x\caus y \chron z$ implies $x\chron z$, this gives rise to the following two inclusions (which are actually equalities by \cref{proposition:properties of cuasal and chronological cones}):
	\[
		I^+(J^+(A))\subseteq I^+(A)
		\qquad\text{and}\qquad
		J^-(I^-(A))\subseteq I^-(A),
	\]
holding for all $A\subseteq M$. Dually, if $x\chron y\caus z$ implies $x\chron z$ we get
	\[
		J^+(I^+(A))\subseteq I^+(A)
		\qquad\text{and}\qquad
		I^-(J^-(A))\subseteq I^-(A).
	\]
Conversely, letting $A$ be a singleton, each of these inclusions gives rise to one half of the push-up property. However, in smooth spacetimes, two of these four inclusions carry more information than the others. Namely, from $I^\pm(A) = (J^\pm(A))^\circ$ (\cref{proposition:chronology is interior of causality}) it already follows that the inclusions
	\[
		I^\pm (J^\pm(A))
		=
		\left(
		J^\pm(J^\pm(A))
		\right)^\circ
		= J^\pm(A)^\circ = I^\pm(A)
	\]
always hold. Thus, when thinking in terms of cones and granting that the chronology can be obtained as the interior of the causality, the essence of the push-up principle seems to be rather captured by the inclusions
	\[
		J^+(I^+(A))\subseteq I^+(A)
		\qquad\text{and}\qquad
		J^-(I^-(A))\subseteq I^-(A).
	\]
The following result shows that abstracting this condition to the setting of ordered spaces and restricting to open subsets gives a new characterisation of the open cone condition.

\begin{lemma}\label{lemma:open cones from push-up}
	An ordered space has open cones if and only if $\up(\up U)^\circ\subseteq (\up U)^\circ$ and $\down(\down U)^\circ\subseteq(\down U)^\circ$ for all open subsets $U$.
\end{lemma}
\begin{proof}
	If the open cone condition holds, we get from a simple calculation using \cref{proposition:properties of cones}(c) that
	\[
	\up (\up U)^\circ = \up\up U = \up U = (\up U)^\circ,
	\]
	and analogously for downsets. For the converse, simply note that ${U=U^\circ\subseteq (\up U)^\circ}$ by \cref{proposition:properties of cones}(b), and hence $\up U\subseteq \up (\up U)^\circ\subseteq (\up U)^\circ$ follows from the assumption, proving that $\up$~preserves opens. Together with the dual condition this gives open cones.
\end{proof}

\newpage
\begin{theorem}\label{corollary:spacetimes have open cones}
	Any smooth spacetime has open cones.
\end{theorem}
\begin{proof}
	For spacetimes we have $\up U = J^+(U)$ and $\down U = J^-(U)$. By \cref{proposition:chronology is interior of causality} we know further that $(J^\pm(U))^\circ = I^\pm(U)$. Hence we see that the push-up principle of \cref{proposition:push-up} implies
		\[
			\up (\up U)^\circ = J^+(I^+(U)) \subseteq I^+(U) = (\up U)^\circ,
		\]
	and dually for the past cones, which is exactly the equivalent condition for open cones outlined in \cref{lemma:open cones from push-up}.
\end{proof}

\begin{corollary}
	For any open subset $U\subseteq M$ of a smooth spacetime:
		\[
			J^+(U) = I^+(U)
			\qquad\text{and}\qquad
			J^-(U) = I^-(U).
		\]
\end{corollary}

Because of this we see that the causal and chronological cones determine the same amount of information on the locale induced by the spacetime, as alluded to in \cref{section:arguments for causal structure}.

\begin{example}
	There exist non-smooth spacetimes that do not satisfy the push-up principle, and hence whose causality relation will not satisfy the open cone condition \cite{grant2020FutureNotAlways}. 
\end{example}

\section{The causal ladder}\label{section:causal ladder}
The \emph{causal ladder} of spacetimes is a hierarchy of properties that makes their casual structure particularly well-behaved from a physical point of view. We adopt these definitions and characterisations from \cite{minguzzi2019LorentzianCausalityTheory}, but see also e.g.~\cite{minguzzi2008CausalHierarchySpacetimesc}. In order from weak to strong, a spacetime $M$ is called:
	\begin{enumerate}[label = \textbullet]
		\item \emph{chronological} if there are no closed timelike loops, i.e.~when $\chron$ is \emph{irreflexive};
		\item \emph{causal} if there are no closed causal loops, i.e.~when $\caus$ is a \emph{partial} order;
		\item \emph{strongly causal} if every point admits an arbitrarily small \emph{causally convex} neighbourhood;
		\item \emph{causally continuous} if it is \emph{weakly distinguishing}, meaning points are equal as soon as their chronological cones are equal, and the chronological cones are \emph{outer continuous};
		\item \emph{causally simple} if it is causal and the graph of $\caus$ is \emph{closed};
		\item \emph{globally hyperbolic} if the \emph{causal diamonds} $J^+(x)\cap J^-(y)$ are \emph{compact}.
	\end{enumerate}
It is apparent from these definitions that some of these properties are \emph{topological} in nature. It is therefore natural to ask whether these classes of spacetime have analogues in the setting of ordered locales,\footnote{This problem was suggested to us by Prakash Panangaden, who collaborated with us on this problem, and resulted in the work of localic convexity in \cref{section:convex locales}.} and if they can be classified fully in terms of localic language. While we currently do not know the answer, besides a localic characterisation of strong causality in \cref{section:convex locales}, the techniques developed in \cref{chapter:topftir} were in part intended to help us answer these questions.

%
%
%
%
\part{Ordered locales}\label{chapter:ordered locales}
\chapter{Ordered locales and their axioms}\label{section:ordered locales}
In this chapter we introduce the main definition of our point-free analogue of ordered spaces. The central idea is to translate all definitions made in terms of an underlying point-set into a definition that can be stated just with the open subsets. Therefore, instead of studying preorders on points, we study preorders on open regions. And this gives us the basic idea behind ordered locales.

\section{The main definition of an ordered locale}
\begin{figure}[b]
	\centering
	\begin{subfigure}[b]{0.4\textwidth}\centering
		\input{tikz/cones0.tikz}
		\caption{Archetypal example of $U\Leq V$.}
	\end{subfigure}\hfil
	\begin{subfigure}[b]{0.45\textwidth}\centering
		\input{tikz/inkscape-export/fig-12}
		\caption{Intuition of the future localic cone.}
	\end{subfigure}
	\caption{Illustrations of $\Leq$ and $\Up$.}
	\label{figure:archetypal Leq}
\end{figure}
\begin{definition}\label{definition:ordered locale}
	An \emph{ordered locale} $(X,\Leq)$ is a locale $X$, equipped with an additional preorder $\Leq$ on its frame of opens $\Opens X$, such that the following law holds:
		\[\tag{$\vee$}\label{axiom:V}
			\forall i\in I: U_i\Leq V_i 
			\quad \text{implies} \quad 
			\bigvee_{i\in I}U_i\Leq \bigvee_{i\in I}V_i.
		\]
\end{definition}
%

The visual intuition we have in mind is that of the Egli-Milner order (recall \cref{section:definition ordered locales}) in Minkowski space, see \cref{figure:archetypal Leq}(a). Axiom~\eqref{axiom:V} says that the graph of $\Leq$ is closed under all joins in the lattice $\Opens X\times \Opens X$. We will see in \cref{lemma:properties of localic cones} that it has something to do with the order-preservation of the localic analogues of the up and down cones. 

Even though $\Opens X$ is a frame, we emphasise further that the graph of $\Leq$ is \emph{not} assumed to be a sub\emph{frame}. That is to say, we do not assume $\Leq$ respects finite meets. Intuitively, the reason for this is that the up and down cones of an ordered space do not preserve intersections. A straightforward illustration in Minkowski space is as follows. In \cref{figure:meets not preserved} we have open sets $U$, $V$ and $W$ that are related by the Egli-Milner order as $U\Leq W$ and $V\Leq W$. The preservation of finite meets would imply $\varnothing = U\cap V \Leq W\cap W = W$, which in turn implies $W=\varnothing$, a contradiction. In \cref{section:parallel ordered locales} we discuss a laxer, more appropriate axiom that relates the order to the meet structure.

\vspace*{10pt}
\noindent\begin{minipage}{0.62\textwidth}\setlength{\parindent}{1.5em}
\indent Before moving on, we introduce an important concept, alluded to in the previous remark. We have seen that the up and down cones of an ordered space (\cref{definition:cones}) give us tools to study how the order structure interacts with the topological structure. We can make analogous definitions in the localic setting. This will make a handy tool in studying the properties of ordered locales.
\end{minipage}%
\hfill%
\begin{minipage}{.33\textwidth}\centering
	\input{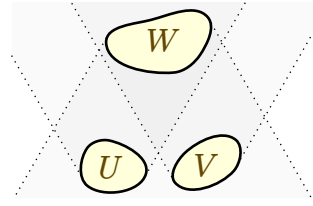}
	\captionof{figure}{Example of $\Leq$ not preserving meets in a Minkowski-type space.}
	\label{figure:meets not preserved}
\end{minipage}
\vspace*{5pt}

\begin{definition}\label{definition:localic cones}
	The \emph{localic cones} of an ordered locale $(X,\Leq)$ are defined as
	\[
	\Up U:=\bigvee \left\{ V\in\Opens X: U\Leq V\right\}
	\qquad\text{and}\qquad
	\Down U:=\bigvee \left\{ W\in\Opens X: W\Leq U\right\}.
	\]
\end{definition}
For visual intuition, see \cref{figure:archetypal Leq}(b). Wherever more care is required, we may write $\Upsub{X}$ and $\Downsub{X}$ to distinguish localic cones belonging to different ordered locales.

Just like the cones of a preordered space, we shall find the localic cones define monads on the frame of opens $\Opens X$. To prove this, we need the following simple lemma. Recall that $\sqleq$ denotes the abstract region inclusion relation of a locale.

\begin{lemma}\label{lemma:ordered locales satisfy L}
	In any ordered locale $(X,\Leq)$ the following two laws hold:
	\begin{center}\vspace{-4mm}
		\dummylabel{axiom:L}{\textsc{l}$^\pm$}%
		\begin{minipage}{.3\linewidth}
			\[\tag{\textsc{l}$^+$}\label{axiom:L+}
			\begin{tikzcd}[every label/.append style = {font = \normalsize},column sep=0.25cm, row sep=0.2cm]
				{U} & {U'} \\
				{V} & \highlighttikzcd{\exists V'}
				\arrow["\sqleq"{anchor=center, rotate=-90}, draw=none, from=2-1, to=1-1]
				\arrow["\sqleq"{anchor=center, rotate=-90}, draw=none, from=2-2, to=1-2]
				\arrow["{}"{description}, "\Leq"{anchor=center}, draw=none, from=1-1, to=1-2]
				\arrow["{}"{description}, "\Leq"{anchor=center}, draw=none, from=2-1, to=2-2]
			\end{tikzcd}
			\]
		\end{minipage}%
		\hfil
		\begin{minipage}{.3\linewidth}
			\[\tag{\textsc{l}$^-$}\label{axiom:L-}
			\begin{tikzcd}[every label/.append style = {font = \normalsize},column sep=0.25cm, row sep=0.2cm]
				{U} & {U'} \\
				\highlighttikzcd{\exists V} & {V'.}
				\arrow["\sqleq"{anchor=center, rotate=-90}, draw=none, from=2-1, to=1-1]
				\arrow["\sqleq"{anchor=center, rotate=-90}, draw=none, from=2-2, to=1-2]
				\arrow["{}"{description}, "\Leq"{anchor=center}, draw=none, from=1-1, to=1-2]
				\arrow["{}"{description}, "\Leq"{anchor=center}, draw=none, from=2-1, to=2-2]
			\end{tikzcd}
			\]
		\end{minipage}
	\end{center}
\end{lemma}
\begin{remark}
	Before giving the elementary proof, here is how these diagrams (and similar ones below) are parsed. Axiom~\eqref{axiom:L+} says that: \emph{if} $U\Leq U'$ and $U\sqleq V$, \emph{then} there exists $V'$ such that $U'\sqleq V'$ and $V\Leq V'$. Axiom~\eqref{axiom:L-} is interpreted analogously. The visual intuition is in \cref{figure:axiom L+}.
\end{remark}
\begin{proof}
	For the first square, note that $U\Leq U'$ by assumption, and $V\Leq V$ by reflexivity, so that through axiom~\eqref{axiom:V} we get $U\vee V\Leq U'\vee V$. Hence we can set $V':= U'\vee V$. The proof of the second square is dual.
\end{proof}

\begin{figure}[b]\centering
	\input{tikz/cones3.tikz}
	\caption{Illustration of~\eqref{axiom:L+}.}
	\label{figure:axiom L+}
\end{figure}

\begin{remark}
	We will see in \cref{corollary:identity maps are monotone} that the laws~\eqref{axiom:L} are really an instance of the more general notion of a \emph{monotone map}.
\end{remark}

Analogously to the properties of up and down cones in an ordered space (recall \cref{proposition:properties of cones}), we find the localic cones behave similarly nicely. In particular, they form monads:

\begin{lemma}\label{lemma:properties of localic cones}
	For a locale $X$ and a preorder $\Leq$ on its frame of opens:
	\begin{enumerate}[label = (\alph*)]
		\item if $U\Leq V$ then $U\sqleq \Down V$ and $V \sqleq \Up U$;
		\item $U\sqleq \Up U$ and $U \sqleq \Down U$.
	\end{enumerate}
	If axiom~\eqref{axiom:V} is satisfied, then furthermore:
	\begin{enumerate}[label = (\alph*)]\setcounter{enumi}{2}
		\item $U \Leq \Up U$ and $\Down U \Leq U$; 
		\item $\Up \Up U = \Up U$ and $\Down \Down U = \Down U$;
		\item if $U\sqleq V$ then $\Up U\sqleq \Up V$ and $\Down U\sqleq \Down V$. 
	\end{enumerate}
	
	\begin{proof}
		Properties~(a) and~(c) follow by definition of the localic cones, and~(b) follows similarly by reflexivity of $\Leq$. From~(c) it also follows that $U\Leq \Up U\Leq \Up\Up U$, so by transitivity of $\Leq$ and property~(a) we get $\Up\Up U\sqleq \Up U$. Lastly, if $U\sqleq V$ it follows from~(c) and~\eqref{axiom:L+} that there exists $V'\in \Opens X$ such that $V\Leq V'\sqgeq \Up U$. By~(a) this implies $\Up U\sqleq V'\sqleq \Up V$, as desired. The proof for the down cones follows similarly using~\eqref{axiom:L-}. This finishes the proof.
	\end{proof}
\end{lemma}

\begin{remark}
	Similar to how a preorder $\leq$ on a set is determined fully by its cones, here too it is somewhat a matter of taste whether we want to work with either the additional relation $\Leq$ or the localic cones $\Up$ and $\Down$ as the primitive notion. We will make this more precise in \cref{section:ordered locales from monads}.
\end{remark}

\begin{remarknumbered}\label{remark:axiom L vs V}
	Note that the proof of \cref{lemma:properties of localic cones} relies only on the axioms~\eqref{axiom:L}, rather than the full force of~\eqref{axiom:V}. If one wanted to develop a more general theory of ``ordered lattices'' in which joins do not necessarily exist, the axioms~\eqref{axiom:L} could be taken as primary instead of~\eqref{axiom:V}. Note that they only depend on the inclusion order $\sqleq$, so they could be stated in a setting as general as that of preorders.
\end{remarknumbered}

Before moving on to examples of ordered locales, we finish with some remarks on how the localic cones interact with meets and joins. It is convenient to state explicitly what \cref{lemma:monotone laxly respects meets and joins} implies about the localic cones.

\begin{lemma}\label{lemma:localic cones lax inclusions}
	Let $(U_i)_{i\in I}$ be any family of opens in an ordered locale. Then the localic up cones satisfy the following inclusions:
		\[
			\bigvee_{i\in I}\Up U_i \sqleq \Up \bigvee_{i\in I} U_i
			\qquad\text{and}\qquad
			\Up \bigwedge_{i\in I} U_i \sqleq \bigwedge_{i\in I}\Up U_i.
		\]
\end{lemma}
\begin{proof}
	\cref{lemma:properties of localic cones}(e) allows us to apply \cref{lemma:monotone laxly respects meets and joins}.
\end{proof}

\begin{remarknumbered}\label{remark:join preservation not necessary}
	Even though we do not expect the localic cones to preserve meets, given the fact that the cones of an ordered space preserve joins (\cref{lemma:properties of cones}(d)), it is natural to expect this in the localic setting too. However, we will see this is not always the case. It fails for instance for quite trivial reasons in \cref{example:inclusion is ordered locale} below. In \cref{definition:localic cones preserve joins} we introduce a further axiom on ordered locales that requires the localic cones to preserve all joins. We believe this to be a natural condition; it is necessary for several constructions (see e.g. \cref{section:ordered locales from biframes}), and it also arises automatically in the approach of internally ordered locales with open cones (\cref{section:internal preorders}).
\end{remarknumbered}

Indeed, the localic cones will not generally preserve meets; look for instance at the situation in \cref{figure:meets not preserved}. However, we do get the following weaker property.

\begin{lemma}\label{lemma:localic cones preserve meets in image}
	For any $U,V\in\Opens X$ in an ordered locale:
		\[
			\Up U \wedge \Up V = \Up \left(\Up U \wedge \Up V\right)
			\qquad\text{and}\qquad
			\Down U \wedge \Down V = \Down \left(\Down U \wedge \Down V\right).
		\]
\end{lemma}
\begin{proof}
	This is a standard property of closure operators on lattices: the inclusions from left to right follow by the unit \cref{lemma:properties of localic cones}(b). The converse follows by applying \cref{lemma:localic cones lax inclusions}, and then using the idempotence from \cref{lemma:properties of localic cones}(d).
\end{proof}

\section{First examples of ordered locales}\label{section:examples ordered locales}
We are finally ready to introduce some examples of ordered locales, and we start with some elementary ones.
\begin{example}\label{example:equality is ordered locale}
	For any locale $X$ we get an ordered locale $(X,=)$. It is clear that~\eqref{axiom:V} is satisfied. Its localic cones are the identities: $\Up U = U = \Down U$.
\end{example}

\begin{example}\label{example:inclusion is ordered locale}
	Less trivially, given any locale $X$, we naturally get the inclusion partial order $\sqleq$ on $\Opens X$. From elementary properties of complete lattices it follows that axiom~\eqref{axiom:V} is satisfied (\cref{lemma:join respects inclusion}):
	\[
	\forall i\in I: U_i\sqleq V_i
	\qquad\text{implies}\qquad
	\bigvee_{i\in I} U_i\sqleq \bigvee_{i\in I} V_i.
	\]
	We will see later that this ordered locale is a localic generalisation of the \emph{specialisation order} (\cref{example:ordered locale with inclusion order gives specialisation order}). Its localic cones behave somewhat strangely: we find $\Down U = U$, but $\Up U = X$ always equals the top element. Thus the localic up cone does not generally preserve joins, as alluded to in \cref{remark:join preservation not necessary}, since that would imply in particular that the empty join is preserved: $\varnothing = \Up \varnothing = X$.
\end{example}

\begin{example}
	Given any ordered locale $(X,\Leq)$, the pair $(X,\trianglerighteqslant)$ is also an ordered locale. The localic cones of $\trianglerighteqslant$ are obtained as the dual of those of $\Leq$.
\end{example}


\subsection{Ordered spaces as ordered locales}
\label{example:ordered spaces as ordered locales}\label{section:ordered spaces as ordered locales}
For us, one of the most important examples of ordered locales are the ones coming from ordered \emph{spaces}. Given an ordered space $(S,\leq)$, what kind of order can we put on $\Opens S$? There are three canonical choices for this (see for instance \cite[Section~11.1]{vickers1989TopologyLogic}):
\begin{enumerate}[label = \textbullet]
	\item the \emph{upper order:} $U\Leq_{\mathrm{U}}V$ if and only if $V\subseteq \up U$;
	\item the \emph{lower order:} $U\Leq_{\mathrm{L}}V$ if and only if $U\subseteq \down V$;
	\item the \emph{Egli-Milner order:} $U\Leq_{\mathrm{EM}} V$ if and only if $U\Leq_{\mathrm{U}}V$ and $U\Leq_{\mathrm{L}}V$.
\end{enumerate}
Note that this definition of $\LeqEM$ is equivalent to the one in \cref{section:definition ordered locales}. Each of these relations ${\Leq}\in\{\Leq_{\mathrm{U}},\Leq_{\mathrm{L}},\Leq_{\mathrm{EM}}\}$ turns $\loc(S)$ into an ordered locale. This can be seen directly: first note that using the properties of cones \cref{lemma:properties of cones} it is straightforward to prove that $\Leq$ is a preorder. To see that axiom~\eqref{axiom:V} is also satisfied, fix two families of open subsets of $S$, indexed by some set $I$, that satisfy $U_i\Leq V_i$ for all $i\in I$. In the case of the lower order, this means in particular that $U_i\subseteq \down V_i$, so taking joins gives
\[
\bigcup_{i\in I} U_i \subseteq \bigcup_{i\in I} \down V_i \subseteq \down\bigcup_{i\in I} V_i,
\]
where the last inclusion is a simple argument using monotonicity (\cref{lemma:monotone laxly respects meets and joins}). Dually, for the upper order we would get $\bigcup_{i\in I}V_i\subseteq \up \bigcup_{i\in I} U_i$. So for each $\Leq$ it follows that $\bigcup_{i\in I}U_i\Leq \bigcup_{i\in I}V_i$. In conclusion:

\begin{proposition}\label{proposition:space with EM order is ordered locale}
	If $(S,\leq)$ is an ordered space, then $\loc(S)$ is an ordered locale when equipped with the upper, lower, or Egli-Milner order.
\end{proposition}

\begin{remark}
	Throughout this thesis we mainly work with the Egli-Milner order. For this reason, if the symbol $\Leq$ appears in the context of ordered spaces it will always be assumed to be the Egli-Milner order, unless otherwise specified.
\end{remark}

\begin{lemma}\label{lemma:localic cones in a space}
	The localic cones of an open subset $U\subseteq S$ of an ordered space are calculated as:
	\[
	\Up U = (\up U)^\circ
	\qquad\text{and}\qquad
	\Down U = (\down U)^\circ.
	\]
\end{lemma}
\begin{proof}
	To prove this, first we claim that
	\[
	U \Leq (\up U)^\circ
	\qquad\text{and}\qquad
	(\down U)^\circ \Leq U.			
	\]
	Looking at the definition of the Egli-Milner order, we see that $U\Leq (\up U)^\circ$ if and only if $(\up U)^\circ\subseteq \up U$ and $U\subseteq \down (\up U)^\circ$. The first inclusion holds by construction of the interior operator, and the latter follows since $U= U^\circ \subseteq (\up U)^\circ$, and hence $U\subseteq \down U \subseteq \down (\up U)^\circ$. 
	
	It now follows from \cref{lemma:properties of localic cones}(a) that $(\up U)^\circ \subseteq \Up U$, so we are left to prove the converse. For that, take an element $x\in \Up U$. By construction of the localic cones, there exists an open neighbourhood $V\in \Opens S$ of $x$ such that $U\Leq V$. By definition of the Egli-Milner order we thus get in particular $x\in V\subseteq \up U$, from which it immediately follows that $x\in (\up U)^\circ$. The proof for down cones is dual.
\end{proof}

\begin{corollary}\label{corollary:localic cones in space with open cones}
	In an ordered space with open cones: $\Up U =\up U$ and $\Down U = \down U$.
\end{corollary}

\begin{example}\label{example:localic cones in spacetime}
	For any open region $U$ in a spacetime $M$, we have
		\[
			\Up U = I^+(U) = J^+(U)
			\qquad\text{and}\qquad
			\Down U = I^-(U)=J^-(U).
		\]
\end{example}


\begin{figure}
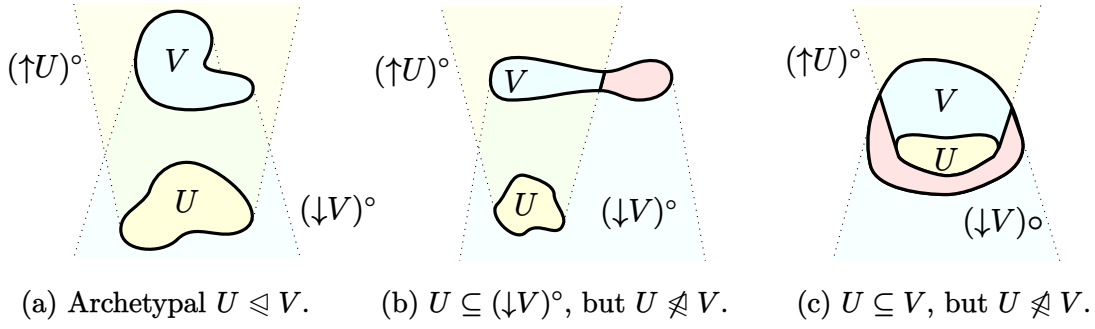

	\centering
	\begin{subfigure}[b]{0.29\textwidth}
		\input{tikz/cones0.tikz}
		\caption{Archetypal $U\Leq V$.}
	\end{subfigure}\hfill
	\begin{subfigure}[b]{0.33\textwidth}
		\input{tikz/cones0b.tikz}
		\caption{$U\subseteq (\down V)^\circ$, but $U\not\Leq V$.}
	\end{subfigure}\hfill
	\begin{subfigure}[b]{0.29\textwidth}
		\input{tikz/cones0c.tikz}
		\caption{$U\subseteq V$, but $U\not\Leq V$.}
	\end{subfigure}
	\caption{Illustrations of $\Leq$ in a Minkowski-like ordered space.}
	\label{figure:illustration of causal order on opens}
\end{figure}

\subsection{Ordered locales from pairs of monads}
\label{section:ordered locales from monads}
\label{example:order from pair of monads}
This section is dedicated to an important class of ordered locales: those that arise from pairs of monads. We shall see in \cref{lemma:pair of monads gives ordered locale} that these are precisely the ordered locales whose causal order is fully determined by their localic cones. The proofs here are a direct generalisation of those in \cref{example:ordered spaces as ordered locales}.

Let $X$ be any locale, and consider two monads $u$ and $d$ on its frame of opens~$\Opens X$, which we think of as abstract localic cones. Recall this means that (\cref{definition:monads}):
	\begin{enumerate}[label = (\roman*)]
		\item if $U\sqleq V$ then $u(U)\sqleq u(V)$ and $d(U)\sqleq d(V)$;\hfill(\emph{functoriality})
		\item $U\sqleq u(U)$ and $U\sqleq d(U)$;\hfill(\emph{unit})
		\item $u(u(U))\sqleq u(U)$ and $d(d(U))\sqleq d(U)$.\hfill(\emph{multiplication})
	\end{enumerate}
We define a relation on $\Opens X$ by:
	\[
		U\Leq V
		\qquad\text{if and only if}\qquad
		U\sqleq d(V)\text{~and~}V\sqleq u(U).
	\]
	
\begin{lemma}\label{lemma:pair of monads gives ordered locale}
	For any two monads $u$ and $d$ defined on the frame of opens of a locale $X$, the induced pair $(X,\Leq)$ forms an ordered locale. Moreover, the localic cones of this order return the original monads.
\end{lemma}

\begin{proof}
	We can see that $\Leq$ is reflexive and transitive precisely because of the unit and multiplication of the monads. Namely $U\Leq U$ if and only if $U\sqleq u(U)$ and $U\sqleq d(U)$. Furthermore, if $U\Leq V\Leq W$, then $U\sqleq d(V)\sqleq d(d(W))\sqleq d(W)$, and similarly $W\sqleq u(V)\sqleq u(u(U))\sqleq u(U)$, so it follows that $U\Leq W$.
	
	It also satisfies~\eqref{axiom:V}: suppose $U_i\Leq V_i$ for all $i\in I$. In particular we get inclusions $U_i\sqleq d(V_i)$, so taking the join and using \cref{lemma:monotone laxly respects meets and joins} gives
	\[
	\bigvee_{i\in I} U_i \sqleq \bigvee_{i\in I}d(V_i) \sqleq d\left(\bigvee_{i\in I} V_i\right).
	\]
	Dually we find $\bigvee_{i\in I}V_i \sqleq u \left(\bigvee_{i\in I} U_i\right)$, so we conclude $\bigvee_{i\in I} U_i \Leq \bigvee_{i\in I} V_i$.
	
	We are left to show that the localic cones equal $u$ and $d$. To see this, first note that by \cref{definition:localic cones}, the localic future $\Up U$ is defined as the join over all $V\in\Opens X$ such that $U\Leq V$. By definition, any such $V$ is contained in $u(U)$, so we get $\Up U \sqleq u(U)$. For the converse inclusion, we claim that $U\Leq u(U)$. Again by definition, this holds if and only if $U\sqleq d(u(U))$ and $u(U)\sqleq u(U)$. The second inclusion is trivial, and the first holds by composing the units of both monads. Thus we see $\Up U = u(U)$, and with a dual argument we find $\Down U = d(U)$. 
\end{proof}

\begin{remark}
	Note in this proof that axiom~\eqref{axiom:V} arises from the functoriality of the monads, and that the monads themselves do not need to preserve any joins.
\end{remark}

The ordered locales that arise from a pair $(u,d)$ of monads through \cref{lemma:pair of monads gives ordered locale} are special in the sense that the condition $U\Leq V$ can be described fully in terms of the localic cones. We capture this property in a new axiom.

\begin{definition}\label{definition:cones give order}
	In an ordered locale $(X,\Leq)$, we say \emph{order is determined by cones} if the following axiom holds:
		\[\tag{\textsc{c}-$\Leq$}\label{axiom:cones give order}
		U\Leq V \quad \text{if and only if} \quad U\sqleq \Down V \text{ and } V\sqleq \Up U.
		\]
\end{definition}

\begin{remark}
	Note that, by \cref{lemma:properties of localic cones}(a), the ``only if'' implication always holds. Therefore an ordered locale $(X,\Leq)$ satisfies~\eqref{axiom:cones give order} precisely when $U\sqleq \Down V$ and $V\sqleq \Up U$ implies $U\Leq V$.
\end{remark}

\begin{proposition}\label{proposition:ordered locales vs monads}
	Given a locale $X$, there is a bijective correspondence between preorders $\Leq$ on $\Opens X$ such that $(X,\Leq)$ is an ordered locale satisfying~\eqref{axiom:cones give order}, and pairs of monads $(u,d)$ on $\Opens X$.
\end{proposition}
\begin{proof}
	We have seen in \cref{lemma:properties of localic cones} that any ordered locale $(X,\Leq)$ induces monads $u:= \Up(-)$ and $d:=\Down(-)$ on its frame of opens $\Opens X$. Conversely, given a pair of monads $(u,d)$ on $\Opens X$, we obtain an ordered locale via \cref{lemma:pair of monads gives ordered locale}, the localic cones of which return precisely the original monads. Therefore, this ordered locale satisfies~\eqref{axiom:cones give order}. It is clear that these constructions are mutually inverse.
\end{proof}

\begin{remark}
	Most ordered locales we are interested in will have the property that is showcased in \cref{definition:cones give order}. In fact, the special class of \emph{parallel ordered} locales, introduced in \cref{section:parallel ordered locales}, automatically has this property (\cref{lemma:wedge implies cones determine order}). Ordered spaces with open cones belong to this class (\cref{proposition:space with OC is parallel ordered}).
\end{remark}

\begin{example}\label{example:ordered space comes from monads}
	The upper, lower, and Egli-Milner orders on an ordered space $(S,\leq)$ from \cref{section:ordered spaces as ordered locales} are recovered through this more general construction. The Egli-Milner order specifically is recovered from the pair of monads $U\mapsto (\up U)^\circ$ and $U\mapsto (\down U)^\circ$. The upper and lower orders are recovered similarly by replacing one of these monads by the constant map $U\mapsto S$. It follows automatically that the ordered locale coming from an ordered space satisfies~\eqref{axiom:cones give order}. In the case that~$S$ has open cones, the monads can be simplified to $U\mapsto \up U$ and $U \mapsto \down U$ (which do not define monads on $\Opens S$ otherwise), so that also \cref{corollary:localic cones in space with open cones} follows from this more general setting.
\end{example}

\subsection{Ordered locales whose cones preserve joins}
Recall from \cref{lemma:properties of cones}(d) that the up and down cones of a preorder preserve arbitrary unions. In the list of properties of localic cones \cref{lemma:properties of localic cones} we have not mentioned an analogous property. This is because, with our current \cref{definition:ordered locale} of ordered locales, this is not automatic. We therefore introduce it as an additional axiom.

\begin{definition}\label{definition:localic cones preserve joins}
	We say the localic cones of an ordered locale $(X,\Leq)$ \emph{preserve joins} if the following law holds:
	\[\tag{\textsc{c}-$\vee$}\label{axiom:LV}
	\Up\bigvee_{i\in I} U_i = \bigvee_{i\in I}\Up U_i
	\qquad\text{and}\qquad
	\Down \bigvee_{i\in I} U_i = \bigvee_{i\in I}\Down U_i.
	\]
\end{definition}

Any ordered locale coming from the Egli-Milner order of an ordered space with open cones is an example.

\begin{proposition}
	If $(S,\leq)$ has open cones, then $\loc(S)$ satisfies~\eqref{axiom:LV}.
\end{proposition}
\begin{proof}
	This follows simply from the fact that the localic cones equal the point-wise cones (\cref{corollary:localic cones in space with open cones}), and the latter preserve joins by \cref{lemma:properties of cones}(d).
\end{proof}

However, even for ordered locales based on spaces $S$, axiom~\eqref{axiom:LV} can fail for quite trivial reasons. Consider for instance the lower order from \cref{section:ordered spaces as ordered locales}, which we recall is defined by $U\Leq_{\mathrm{L}} V$ if and only if $U\subseteq \down V$. Thus we find that $\Down U = (\down U)^\circ$, but $\Up U = S$. In particular, the empty subset $\varnothing$ is the empty join, however $\Up \varnothing = S$ is non-empty, generally. The upper order clearly suffers from a similar problem, and we take this as a strong argument to prefer the Egli-Milner order over the upper and lower orders. But~\eqref{axiom:LV} can fail in less trivial ways, even for finite joins, as the next example shows.

\begin{figure}[t]\centering
	\input{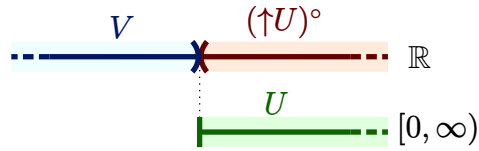}
	\caption{Illustration of the space in \cref{example:LV can fail in spaces}.}
	\label{figure:LV can fail in spaces}
\end{figure}

\begin{example}\label{example:LV can fail in spaces}
	This example is similar in spirit to the ordered space in \cref{example:non-OC space} that does not have open cones. We consider the ordered space $(S,\leq)$ whose underlying space is the disjoint union ${(\{0\}\times [0,\infty))\cup (\{1\}\times \mathbb{R})}$, and whose order is generated by setting $(0,x)\leq (1,y)$ if and only if $x=y$. This space does not have open cones either. We construct two open subsets of $S$ such that $\Up(U\cup V)\neq \Up U \cup \Up V$. Namely, take:
	\[
	U:= \{0\}\times [0,\infty)
	\qquad\text{and}\qquad
	V:= \{1\}\times (-\infty,0).
	\]
	Then $\up V = V$, and $\up U = U \cup \left(\{1\}\times [0,\infty)\right)$. We therefore find that
	\[
	\Up(U\cup V) = \left(\up U \cup \up V\right)^\circ 
	=
	S,
	\]
	while
	\[
	\Up U \cup \Up V = \left(\up U\right)^\circ \cup V
	=
	\left(\{0\}\times [0,\infty)\right)\cup\left(\{1\}\times \left(\mathbb{R}\setminus \{0\}\right)\right).
	\]
	In other words, $\Up U \cup \Up V$ misses precisely the point $(1,0)$. See \cref{figure:LV can fail in spaces}. We can see that~\eqref{axiom:LV} fails essentially due to the interior operator not preserving joins.
\end{example}

\subsection{Ordered locales from biframes}
\label{section:ordered locales from biframes}
Topological spaces generalise to \emph{bitopological spaces} (or \emph{bispaces}), which were introduced in \cite{kelly1963BitopologicalSpaces} to study the topological aspects of quasi-metric spaces (i.e.~sets equipped with a non-symmetric distance function), and are defined simply as tuples $(S,\tau_1,\tau_2)$ where $\tau_1$ and $\tau_2$ are topologies on a set $S$. Bispaces have become an important tool in point-set topology, and have been used to study ordered spaces.

The notion of a \emph{biframe} is the corresponding point-free analogue, introduced in \cite{banaschewski1983BiframesBispaces}. Formally, a \emph{biframe} is a tuple $L=(L_0,L_1,L_2)$ of frames, where $L_1,L_2\hookrightarrow L_0$ are \emph{subframes}, i.e.~subsets that are closed under finite meets and arbitrary joins of $L_0$, that together generate $L_0$ in the sense that $L_1\cup L_2$ is a subbasis for $L_0$. Explicitly, this means that every $x\in L_0$ can be written as a join of the form $x=\bigvee_{i\in I} (y_i\wedge z_i)$, for $y_i\in L_1$ and $z_i\in L_2$. Refer to \cite{schauerte1992Biframes} for more details. 

To connect the theory of biframes to that of ordered locales, in this section we shall think of $L_1$ as an abstract frame of up sets, and $L_2$ as an abstract frame of down sets. Then $L_0$ is the frame generated by convex combinations of the two. In other words, we are trying to think of tuples $(\Opens X,\im(\Up),\im(\Down))$ as biframes, where:
\[
\im(\Up):= \left\{\Up U: U\in \Opens X\right\}
\qquad\text{and}\qquad
\im(\Down):= \left\{\Down V:V\in\Opens X\right\}.
\]


First, we characterise when, starting with an ordered locale $(X,\Leq)$, the tuple $(\Opens X,\im(\Up),\im(\Down))$ defines a biframe. To start, the images of the localic cones need to be subframes of $\Opens X$, and the following two lemmas tell us that this is the case precisely if~\eqref{axiom:LV} holds.

\begin{lemma}\label{lemma:localic cone of join}
	In any ordered locale:
	\[
	\Up \bigvee_{i\in I} U_i = \Up \bigvee_{i\in I}\Up U_i
	\qquad\text{and}\qquad
	\Down \bigvee_{i\in I} U_i = \Down \bigvee_{i\in I} \Down U_i.
	\]
\end{lemma}
\begin{proof}
	First, for all $i\in I$ we have $U_i\sqleq \Down U_i$, so taking joins: $\bigvee U_i \sqleq \bigvee \Down U_i$. By monotonicity it follows that $\Down \bigvee U_i \sqleq \Down \bigvee \Down U_i$. For the converse direction, \cref{lemma:localic cones lax inclusions} gives $\bigvee \Down U_i \sqleq \Down \bigvee U_i$, and again applying $\Down$ to both sides and using idempotence gives the result.
\end{proof}

\begin{lemma}\label{lemma:past sets form a frame}
	If $X$ is an ordered locale, then $\im(\Down)$ is a subframe of $\Opens X$ if and only if the localic past cone $\Down$ preserves joins. Dually for future cones.
\end{lemma}
\begin{proof}
	It is always the case that $\im(\Down)$ is closed under finite meets in $\Opens X$ since, using \cref{lemma:localic cones preserve meets in image}, we get $\Down U\wedge \Down V=\Down(\Down U\wedge \Down V)\in\im(\Down)$, and it clearly contains the top element $X=\Down X \in \im(\Down)$. Hence $\im(\Down)$ is a sub meet-semilattice. 
	
	For it to be a subframe we additionally need it to be closed under arbitrary joins. This means that if $\Down U_i$ is a family of elements in $\im(\Down)$, then we must have $\bigvee \Down U_i\in \im(\Down)$. Hence $\im(\Down)$ is a subframe if and only if we have the equations $\bigvee \Down U_i = \Down \bigvee \Down U_i$, and by \cref{lemma:localic cone of join} this is equivalent to $\Down$ preserving joins.
\end{proof}

\begin{remark}
	Therefore, any ordered locale $X$ with~\eqref{axiom:LV} defines two frames $\im(\Up)$ and $\im(\Down)$, which in turn can be interpreted as \emph{locales} $X^{\triup}$ and $X^\tridown$, respectively. We study these more in \cref{section:locale of pasts and futures}, and turn out to form an important ingredient in a \emph{causal boundary} construction (\cref{section:causal boundaries}).
\end{remark}

\begin{corollary}\label{corollary:ordered locale forms biframe iff}
	Let $(X,\Leq)$ be an ordered locale. The tuple $(\Opens X,\im(\Up),\im(\Down))$ defines a biframe if and only if $X$ satisfies~\eqref{axiom:LV}, and every $U\in \Opens X$ is of the form $\bigvee (\Up V_i \wedge \Down W_i)$.
\end{corollary}
\begin{proof}
	This follows from \cref{lemma:past sets form a frame} and by unpacking what it means for the set ${\im(\Up)\cup \im(\Down)}$ to generate $\Opens X$.
\end{proof}

The class of ordered locales whose frame of opens is generated by $\im(\Up)\cup \im(\Down)$ are called \emph{convex}. These generalise strongly causal spacetimes. We shall formally define and study these in \cref{section:convex locales}, and now turn our attention to the converse question: given a biframe $(\Opens X,L_1,L_2)$, can we turn $X$ into an ordered locale such that $\im(\Up) = L_1$ and $\im(\Down) = L_2$? It is not to be expected that such a procedure works in general, since $\im(\Up)$ and $\im(\Down)$ are subframes of quite a specific form. Namely, they are the images of join-preserving monads. The following lemma helps us state this more abstractly.

\begin{lemma}\label{lemma:monads are adjoints}
	Let $t\colon \Opens X\to \Opens X$ be a monotone function. Then $t$ is a monad if and only if it is left adjoint to the inclusion of its image:
	\[
	\begin{tikzcd}[cramped,column sep=large]
		{\im(t)} & {\Opens X.}
		\arrow[""{name=0, anchor=center, inner sep=0}, "i"', shift right=2, hook, from=1-1, to=1-2]
		\arrow[""{name=1, anchor=center, inner sep=0}, "t"', shift right=2, from=1-2, to=1-1]
		\arrow["\dashv"{anchor=center, rotate=-90}, draw=none, from=1, to=0]
	\end{tikzcd}
	\]
\end{lemma}
\begin{proof}
	That $t\dashv i$ means that for every $U\in \Opens X$ and $t(V)\in \im(t)$ we have
	\[
	t(U)\sqleq t(V) 
	\qquad\text{if and only if}\qquad
	U\sqleq i\circ t(V) = t(V).
	\]
	It is easy to see that this holds if $t$ is a monad: namely from $t(U)\sqleq t(V)$ the unit gives $U\sqleq t(U)\sqleq t(V)$, and conversely if $U\sqleq t(V)$ then the idempotence gives $t(U)\sqleq t\circ t(V) = t(V)$.
	
	Conversely, if $t\dashv i$ then $t(U)\sqleq t(U)$ implies $U\sqleq t(U)$, providing the unit for~$t$. Similarly, $t(U)\sqleq t(U)$ implies $t\circ t(U)\sqleq t(U)$, giving the idempotence.
\end{proof}
Note that, as a left adjoint, $t$ necessarily preserves joins by \cref{theorem:join preserving implies left adjoint}.

\begin{remark}
	This is effectively a heavily de-categorified statement of the fact that every monad comes from an adjunction. In particular, the result holds much more generally replacing $\Opens X$ by an arbitrary preordered set. See e.g.~\cite[\S 7.27]{davey2002IntroductionLatticesOrder}.
\end{remark}

Combining \cref{lemma:past sets form a frame,lemma:monads are adjoints} tells us that the biframes $(\Opens X,L_1,L_2)$ arising from (convex) ordered locales are those of the form where the subframe inclusion maps $L_1,L_2\hookrightarrow \Opens X$ admit left adjoints. These can then directly be used to construct the localic cones, which automatically preserve joins. We capture these types of biframes into a definition:

\begin{definition}\label{definition:reflective biframe}
	A biframe $(L_0,L_1,L_2)$ is called \emph{reflective} if the subframe inclusions $L_1,L_2\hookrightarrow L_0$ have left Galois adjoints.
\end{definition}

\begin{remark}
	The terminology comes from category theory, where a (full) subcategory is called \emph{reflective} if the inclusion functor has a left adjoint. We are not aware of any literature on reflective biframes, but the theory of reflective subcategories is well known. (Note subframes are automatically full subcategories.)
\end{remark}

\begin{proposition}\label{proposition:convex locales are biframes}
	There is a bijective correspondence between convex ordered locales satisfying~\eqref{axiom:cones give order} and~\eqref{axiom:LV}, and reflective biframes.
\end{proposition}
\begin{proof}
	Starting with such an ordered locale $(X,\Leq)$, \cref{corollary:ordered locale forms biframe iff,lemma:monads are adjoints} give us a biframe $(\Opens X, \im(\Up),\im(\Down))$ with the desired properties.
	
	Conversely, if $(\Opens X, L_1,L_2)$ is a reflective biframe, then we get
	\[
	\begin{tikzcd}[cramped,column sep=large]
		{L_1} & {\Opens X}
		\arrow[""{name=0, anchor=center, inner sep=0}, shift right=2, hook, from=1-1, to=1-2]
		\arrow[""{name=1, anchor=center, inner sep=0}, "u"', shift right=2, from=1-2, to=1-1]
		\arrow["\dashv"{anchor=center, rotate=-90}, draw=none, from=1, to=0]
	\end{tikzcd}
	\qquad\text{and}\qquad
	\begin{tikzcd}[cramped,column sep=large]
		{L_2} & {\Opens X.}
		\arrow[""{name=0, anchor=center, inner sep=0}, shift right=2, hook, from=1-1, to=1-2]
		\arrow[""{name=1, anchor=center, inner sep=0}, "d"', shift right=2, from=1-2, to=1-1]
		\arrow["\dashv"{anchor=center, rotate=-90}, draw=none, from=1, to=0]
	\end{tikzcd}
	\]
	From the general theory of Galois connections it follows that, since the subframe inclusion maps are injections, the maps $u$ and $d$ must be surjections (\cref{proposition:galois adjunction implies surjection iff injection}). In other words: $L_1 = \im(u)$ and $L_2 = \im(d)$. Subsequently, by \cref{lemma:monads are adjoints} the maps $u$ and $d$ define monads on $\Opens X$, and as left adjoints they preserve all joins (\cref{theorem:adjoints preserve limits}). Via the general procedure of \cref{section:ordered locales from monads} we therefore obtain an ordered locale $(X,\Leq)$ satisfying~\eqref{axiom:cones give order} and~\eqref{axiom:LV}, and whose localic cones return~$u$ and~$d$. In particular $L_1 = \im(\Up)$ and $L_2 = \im(\Down)$. That these two constructions are mutually inverse is clear. 
\end{proof}

\begin{remark}
	This correspondence is object-level. To make it categorical, we first need to define what we mean by a morphism of ordered locales, which we do in the next section. We find in \cref{corollary:isomorphism reflective biframes} that the correspondence above indeed has a categorical interpretation.
	
	Further note that if we take triples $(L_0,L_1,L_2)$ of subframes $L_1,L_2\subseteq L_0$ that are reflective (in the sense of \cref{definition:reflective biframe}), but drop the biframe condition, we obtain a more general correspondence where on the locale side the convexity condition is dropped. 
\end{remark}

\begin{remark}
	One reason the biframe framework is appealing is that the theory of their compactifications is well studied \cite{schauerte1992Biframes}. With \cref{proposition:convex locales are biframes} we have therefore opened a new toolbox to construct ``order compactifications'' of ordered locales, and in turn, of spacetimes. Although our work on this is not definitive, we explore this idea further in \cref{section:causal boundaries}.
\end{remark}

\section{Morphisms of ordered locales}\label{section:morphisms ordered locales}
In this section we define morphisms of ordered locales, and thereby obtain a category $\OrdLoc$. A morphism $(X,\Leq)\to (Y,\Leq)$ of ordered locales will consist of a morphism $f\colon X\to Y$ between the underlying locales that additionally respects the order structure, in a sense we specify below. The problem lies in finding a suitable region-based alternative to the definition of a monotone function between preordered sets (\cref{definition:monotone function between preorders}). Here we consider two approaches in formulating such a definition. The first is directly in terms of $\Leq$, and is motivated by the perspective, outlined in \cref{remark:preorders are monads}, that $\Leq$ can be seen as a monad internal to $\Rel$. Monotonicity then comes down to the fact that $f$ induces a certain morphism of monads on $\Leq$, but we spell this out explicitly in \cref{definition:monotone map of ordered locales}.

The second is in terms of the localic cones $\Up$ and $\Down$, and is motivated by the characterisation of monotone functions in terms of cones in \cref{lemma:monotonicity in terms of cones}. We prove in \cref{proposition:monotonicity in terms of localic cones} that both approaches are equivalent.

We proceed as follows. Given a map of locales $f\colon X\to Y$, define the following relation:
\[
R_f:=
\left\{ (U,V): U\sqleq f^{-1}(V)\right\}\subseteq \Opens X\times \Opens Y.
\]
In general, for a relation $R$ we denote by $R^\dagger$ its opposite relation. We can then phrase monotonicity of $f$ in terms of how $R_f$ relates to the causal orders on $X$ and $Y$.

\begin{definition}\label{definition:monotone map of ordered locales}
	A map of locales $f\colon X\to Y$ between ordered locales $(X,\Leq_X)$ and $(Y,\Leq_Y)$ is called \emph{monotone} if
	\[
	{\Leq_X}\circ {R_f^\dagger} \subseteq {R_f^\dagger}\circ {\Leq_Y} \qquad\text{and}\qquad
	{R_f}\circ {\Leq_X} \subseteq {\Leq_Y}\circ {R_f}.
	\]
	Separately, these two inclusions may also be called \emph{upper monotonicity} and \emph{lower monotonicity}, respectively.
\end{definition}

\begin{remark}
	Beware that the word ``monotone'' carries two meanings here. In the first place, we have maps of \emph{frames} $f^{-1}$, which are monotone in the sense that $V\sqleq W$ implies $f^{-1}(V)\sqleq f^{-1}(W)$. However, when viewed as a map of \emph{locales} $f$, monotonicity is meant in the ``causal'' sense of \cref{definition:monotone map of ordered locales}. In the context of ordered locales, the first sense of monotonicity is almost always taken for granted, so the word should be understood as in the causal sense, unless explicitly stated otherwise. We hope this causes no confusion.
\end{remark}

\begin{remarknumbered}\label{remark:monotonicity in terms of maps of monads}
	The monotonicity of a localic map can be phrased in terms of the 2-categorical structure of $\Rel$. Namely, it is equivalence to the existence of the following 2-cells, respectively:
	\[
	\begin{tikzcd}[ampersand replacement=\&]
		{\Opens X} \& {\Opens X} \\
		{\Opens Y} \& {\Opens Y}
		\arrow["{R_f^\dagger}", from=2-1, to=1-1]
		\arrow["\Leq"', from=2-1, to=2-2]
		\arrow["\Leq", from=1-1, to=1-2]
		\arrow["\subseteq"{marking, allow upside down}, draw=none, from=1-1, to=2-2]
		\arrow["{R_f^\dagger}"', from=2-2, to=1-2]
	\end{tikzcd}
	\qquad\text{and}\qquad
	\begin{tikzcd}[ampersand replacement=\&]
		{\Opens X} \& {\Opens X} \\
		{\Opens Y} \& {\Opens Y.}
		\arrow["\Leq"', from=2-1, to=2-2]
		\arrow["\Leq", from=1-1, to=1-2]
		\arrow["{R_f}"', from=1-1, to=2-1]
		\arrow["{R_f}", from=1-2, to=2-2]
		\arrow["\supseteq"{marking, allow upside down}, draw=none, from=2-1, to=1-2]
	\end{tikzcd}		
	\]
	In monadic language, it means that the relation $R_f^\dagger$ defines a morphism of monads ${\Leq_Y}\to {\Leq_X}$ and ${\trianglerighteqslant_Y}\to {\trianglerighteqslant_X}$ internal to $\Rel$. Using this definition, it should be straightforward to generalise our notion of an ordered locale to the setting of locales internal to any \emph{topos} (in which case $\Rel$ would be replaced by the category of relations internal to that topos).
\end{remarknumbered}

It is convenient to unpack the admittedly opaque monotonicity condition above into more concrete terms.

\begin{lemma}\label{lemma:monotonicity unpacked}
	A map $f\colon X\to Y$ between ordered locales is monotone if and only~if: 
	\[
	\begin{tikzcd}[every label/.append style = {font = \normalsize},column sep=0.3cm, row sep=0.35cm]
		{U} & {U'} \\
		{V} & {\highlighttikzcd{\exists V'}}
		\arrow["R_f"{anchor=center}, draw=none, from=2-1, to=1-1]
		\arrow["R_f"{anchor=center}, draw=none, from=2-2, to=1-2]
		\arrow["{}"{description}, "\Leq"{anchor=center}, draw=none, from=1-1, to=1-2]
		\arrow["{}"{description}, "\Leq"{anchor=center}, draw=none, from=2-1, to=2-2]
	\end{tikzcd}
	\qquad\text{and}\qquad 
	\begin{tikzcd}[every label/.append style = {font = \normalsize},column sep=0.3cm, row sep=0.35cm]
		{U} & {U'} \\
		\highlighttikzcd{\exists V} & {V'.}
		\arrow["R_f"{anchor=center}, draw=none, from=2-1, to=1-1]
		\arrow["R_f"{anchor=center}, draw=none, from=2-2, to=1-2]
		\arrow["{}"{description}, "\Leq"{anchor=center}, draw=none, from=1-1, to=1-2]
		\arrow["{}"{description}, "\Leq"{anchor=center}, draw=none, from=2-1, to=2-2]
	\end{tikzcd}		
	\]
\end{lemma}
\begin{proof}
	The upper monotonicity condition for $f$ says that we have an inclusion of graphs: ${\Leq}\circ R_f^\dagger\subseteq R_f^\dagger\circ{\Leq}$. Now, $(V,U')$ is contained in the left hand side if and only if there exists $U\in \Opens X$ such that $(U,V)\in R_f$ and $U\Leq U'$. This is exactly the initial data given in the first square. Upper monotonicity then says $(V,U')\in R_f^\dagger\circ {\Leq}$, which holds if and only if there exists $V'\in\Opens Y$ such that $V\Leq V'$ and $(U',V')\in R_f$, precisely as desired. That the second square corresponds to lower monotonicity is proved dually.
\end{proof}

\begin{figure}[t]\centering
	\input{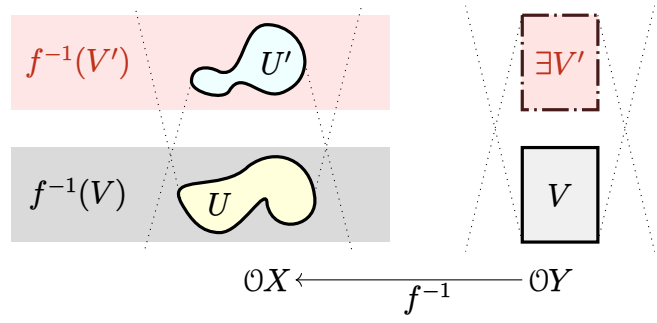}
	\caption{Archetypal example of a monotone map of locales $f$.}
	\label{figure:monotone maps}
\end{figure}

\begin{remark}
	Note the resemblance of the diagrams in the previous lemma to the diagrams in \cref{remark:monotonicity in terms of maps of monads}. The existential quantifiers are pointed in the same direction as the inclusion 2-cells. The visual intuition behind monotonicity is not so straightforward, but we have attempted an illustration in \cref{figure:monotone maps}.
\end{remark}

\begin{corollary}\label{corollary:identity maps are monotone}
	Let $X$ be a locale. A preorder $\Leq$ on $\Opens X$ satisfies~\eqref{axiom:L} if and only if the localic identity map $\id_X$ is monotone.
\end{corollary}

The next observation shows that ordered locales, together with monotone maps of locales, form a category.


\begin{proposition}\label{proposition:ordered locales form a category}
	Ordered locales and monotone maps form a category $\OrdLoc$.
\end{proposition}
\begin{proof}
	\cref{lemma:ordered locales satisfy L,corollary:identity maps are monotone} tell us identity morphisms are monotone, so we are left to show that if $f\colon X\to Y$ and $h\colon Y\to Z$ are monotone maps of locales, then their composition $h\circ f$ is again monotone. Observe that $R_{h\circ f}= R_h\circ R_f$, and using lower monotonicity calculate:
	\[
	{R_{h\circ f}}\circ {\Leq} 
	=
	{R_h}\circ {R_f}\circ {\Leq} 
	\subseteq 
	{R_h}\circ {\Leq} \circ {R_f}
	\subseteq
	{\Leq}\circ {R_h}\circ {R_f}
	=
	{\Leq} \circ {R_{h\circ f}}.
	\]
	That upper monotonicity is preserved is proved dually.
\end{proof}

The next results shows that monotonicity can be expressed using localic cones in a condition that is analogous to \cref{lemma:monotonicity in terms of cones}.

\begin{proposition}\label{proposition:monotonicity in terms of localic cones}
	A localic map $f\colon X\to Y$ between ordered locales is monotone if and only if for all $V\in\Opens Y$:
	\[
	\Up f^{-1}(V)\sqleq f^{-1}(\Up V)
	\qquad\text{and}\qquad
	\Down f^{-1}(V)\sqleq f^{-1}(\Down V).
	\]
\end{proposition}
\begin{proof}
	For this proof, we use the characterisation of monotonicity in \cref{lemma:monotonicity unpacked}. Assume first that $f$ is upper monotone. Using \cref{lemma:properties of localic cones}(c) we get the relation $f^{-1}(V)\Leq \Up f^{-1}(V)$, and $(f^{-1}(V),V)\in R_f$ trivially holds, so there exists $V'\in \Opens Y$ such that $V\Leq V'$ and $(\Up f^{-1}(V),V')\in R_f$. Using \cref{lemma:properties of localic cones}(a) to unpack this, we get 
	\[
	\Up f^{-1}(V)\sqleq f^{-1}(V')\sqleq f^{-1}\left(\Up V\right).
	\]
	
	Conversely, suppose that the stated inclusion holds. Take $U\Leq U'$ and ${(U,V)\in R_f}$. Then $V\Leq \Up V$ by \cref{lemma:properties of localic cones}(c), and further
	\[
	U'\sqleq \Up U \sqleq \Up f^{-1}(V) \sqleq f^{-1}\left(\Up V\right)
	\]
	follows by \cref{lemma:properties of localic cones}(a) and (e), together with the hypothesis, giving ${(U',\Up V)\in R_f}$. This proves $f$ is upper monotone. The proof for lower monotonicity is analogous.
\end{proof}

\begin{remark}
	This result suggests that, especially from the categorical point of view, the ordered locales that do not satisfy~\eqref{axiom:cones give order} can mostly be ignored, since even if $\Leq$ is not determined by its cones, monotonicity is always completely encoded solely by $\Up$ and $\Down$.
	
	We can make this precise. For any $\Leq$, use the procedure of \cref{section:ordered locales from monads} to define a new order $\overline{\Leq}$ from the localic cones $\Up$ and $\Down$:
		\[
			U\mathrel{\overline{\Leq}}V
			\qquad\text{if and only if}\qquad
			V\sqleq \Up U\text{~and~} U\sqleq \Down V.
		\]
	By \cref{lemma:pair of monads gives ordered locale} $(X,\overline{\Leq})$ is an ordered locale whose cones are the same as that of~$\Leq$. Moreover, by \cref{lemma:properties of localic cones}(a) we get that $U\Leq V$ implies $U\mathrel{\overline{\Leq}}V$. The characterisation in \cref{proposition:monotonicity in terms of localic cones} of monotonicity in terms of localic cones then proves that the identity map $\id_X$ defines an isomorphism $(X,\Leq)\cong (X,\overline{\Leq})$ of ordered locales. The categories $\OrdLoc$ and that of ordered locales satisfying~\eqref{axiom:cones give order} are therefore equivalent. At present, we do not know of an example where $\overline{\Leq}$ does not strictly equal $\Leq$.
\end{remark}

\cref{proposition:monotonicity in terms of localic cones} also says that the monotonicity of $f$ is equivalent to the frame map $f^{-1}$ defining morphisms of monads. This elevates the correspondence in \cref{proposition:ordered locales vs monads} to an isomorphism of categories, made precise as follows. If $s$ is a monad on a preordered set $P$, and $t$ is a monad on a preordered set $Q$, then a \emph{morphism of monads} $(P,s)\to (Q,t)$ consists of a monotone map $g\colon P\to Q$ such that $t\circ g(x)\sqleq g\circ s(x)$ for all $x\in P$, cf.~the categorical definition in \cref{definition:monads}. Thus we see that the monotonicity condition for $f$ in \cref{proposition:monotonicity in terms of localic cones} says precisely that $f^{-1}$ defines morphisms of monads
	\[
		(\Opens Y,\Upsub{Y})\longrightarrow (\Opens X,\Upsub{X})
		\qquad\text{and}\qquad
		(\Opens Y,\Downsub{Y})\longrightarrow(\Opens X,\Downsub{X}).
	\]
Motivated by this, define a new category of \emph{bimonads} on frames, which are simply frames $L$ equipped with two monads: $(L,u,d)$, and where the morphisms $(L,u,d)\to (L',u',d')$ consist of morphisms of frames $h\colon L \to L'$ that satisfy $u'\circ h(x)\sqleq h\circ u(x)$ and $d'\circ h(x)\sqleq h\circ d(x)$ for all $x\in L$. 

\begin{proposition}\label{proposition:ordered locales are bimonads}
	There is an isomorphism of categories
	\[
		\left\{\begin{array}{c}
			\text{\normalfont ordered locales} \\
			\text{\normalfont with~\eqref{axiom:cones give order}}
		\end{array}\right\}
		\cong
		\left\{\begin{array}{c}
			\text{\normalfont bimonads} \\
			\text{\normalfont on frames}
		\end{array}\right\}^\op
		.
	\]
\end{proposition}

In a similar spirit to the previous result, the next lemma will show that the monotone maps between ordered locales are actually just biframe morphisms. This means the correspondence in \cref{proposition:convex locales are biframes} actually lifts to an isomorphism of categories. A \emph{biframe map} $h\colon (L_0,L_1,L_2)\to (M_0,M_1,M_2)$ is a map of frames $h\colon L_0\to M_0$ for which the restrictions $h|_{L_i}$ define maps of frames $L_i\to M_i$ \cite[Definition~1.1]{schauerte1992Biframes}. It is equivalent to say that $h(L_i)\subseteq M_i$. The next result shows that monotonicity is captured by precisely such a condition.

\begin{lemma}\label{lemma:monotonicity in terms of biframe maps}
	A map $f\colon X\to Y$ between ordered locales is monotone if and only if:
	\[
	f^{-1}\left(\im(\Upsub{Y})\right)\subseteq \im(\Upsub{X})
	\qquad\text{and}\qquad
	f^{-1}\left(\im(\Downsub{Y})\right) \subseteq \im(\Downsub{X}).
	\]
\end{lemma}
\begin{proof}
	If $f$ is monotone we get the inclusions from \cref{proposition:monotonicity in terms of localic cones}, since then $\Up f^{-1}(\Up V) \sqleq f^{-1}(\Up\Up V) = f^{-1}(\Up V)$, so $f^{-1}(\Up V)\in \im(\Up)$. Dually for past cones.
	
	Conversely, suppose the stated inclusions hold. Then
	\[
	\Up f^{-1}(V)\sqleq \Up f^{-1}(\Up V) \sqleq f^{-1}(\Up V),
	\]
	where the last inclusion is precisely due to the hypothesis. Similarly for down cones, so $f$ is monotone by \cref{proposition:monotonicity in terms of localic cones}.
\end{proof}

Recall from \cref{section:ordered locales from biframes} that an ordered locale is called \emph{convex} if it is generated by $\im(\Up)\cup \im(\Down)$. 
\begin{corollary}\label{corollary:isomorphism reflective biframes}
	The correspondence of \cref{proposition:convex locales are biframes} extends to an isomorphism of categories\vspace{-1ex}
		\[
		\left\{\begin{array}{c}
			\text{\normalfont convex ordered locales} \\
			\text{\normalfont with~\eqref{axiom:cones give order} and~\eqref{axiom:LV}}
		\end{array}\right\}
		\cong
		\left\{\begin{array}{c}
			\text{\normalfont reflective} \\
			\text{\normalfont biframes}
		\end{array}\right\}^\op
		.
		\]
\end{corollary}

\section{More examples of ordered locales}
\label{section:more examples of ordered locales}
\subsection{The lattice of ordered locale orders}\label{section:lattice of ordered locale orders}
In this section we briefly study the structure of the set $\OL(X)$ of all preorders~$\Leq$ that make $X$ into an ordered locale. This should be seen as analogous to, for instance, the lattice of topologies or preorders on a fixed set. We can define an order on $\OL(X)$ through graph inclusions. So ${\Leq} \subseteq {\Leq'}$ if and only if $U\Leq V$ implies $U \Leq' V$. In fact, $\OL(X)$ is itself a complete lattice.

\begin{construction}\label{example:intersection of orders}\label{construction:intersection of orders}
	Let $(\Leq_i)_{i\in I}$ be any family of preorders on the frame of opens of a locale $X$ that satisfy property~\eqref{axiom:V}. We define a new preorder $\Leq$ by:
	\[
	U\Leq V
	\qquad\text{if and only if}\qquad
	\forall i\in I: U\Leq_i V.
	\]
	Using the commutativity of universal quantification we find that the relation $\Leq$ again satisfies axiom~\eqref{axiom:V}. Hence we obtain a new ordered locale $(X,\Leq)$. Note that if the index $I$ is empty, the resulting relation $\Leq$ vacuously relates every pair of opens. In other words: the collection of preorders on $\Opens X$ that satisfy~\eqref{axiom:V} forms a complete lattice. The new order $\Leq$ is the meet $\bigwedge_{i\in I}\Leq_i$.
\end{construction}

The following lemmas shows that the property~\eqref{axiom:cones give order} is preserved under arbitrary meets in $\OL(X)$. This also allows us to calculate the cones of meets explicitly.

\begin{lemma}
	If ${\Leq}\subseteq {\Leq'}$ then $\Up \sqleq \Up'$ and $\Down \sqleq \Down'$. If $\Leq'$ satisfies~\eqref{axiom:cones give order}, then the converse implication also holds.
\end{lemma}
\begin{proof}
	Take ${\Leq}\subseteq {\Leq'}$, so that $U\Leq \Up U$ implies $U\Leq'\Up U$, from which it immediately follows that $\Up U \sqleq \Up' U$. Analogously for down cones.
	
	If $\Leq'$ is determined by its cones, which additionally contain the cones of $\Leq$, then $U\Leq V$ implies $U\sqleq \Down' V$ and $V\sqleq \Up'U$, which in turn implies $U\Leq' V$. Thus we get ${\Leq}\subseteq {\Leq'}$.
\end{proof}

\begin{lemma}
	Consider a family $({\Leq_i})_{i\in I}$ in $\OL(X)$ of orders satisfying~\eqref{axiom:cones give order}, whose cones we denote by $\Upi$ and $\Downi$. Then their meet ${\Leq} = \bigwedge_{i\in I} {\Leq_i}$ also satisfies~\eqref{axiom:cones give order}, and the cones are calculated as $\Up U = \bigwedge_{i\in I}\Upi U$ and $\Down U = \bigwedge_{i\in I}\Downi U$.
\end{lemma}
\begin{proof}
	We start by proving the stated equation for the cone $\Up$, the proof of which is done in two steps. First we prove that $\Up U \sqleq \bigwedge_{i\in I} \Upi U$, for which we need to show that $\Up U \sqleq \Upi U$ for every~${i\in I}$. In turn, it suffices to show that if $U\Leq V$ then $V\sqleq \Upi U$. However this follows easily from the definition of $\Leq$ in \cref{construction:intersection of orders} and \cref{lemma:properties of localic cones}(a).
	
	Secondly, we need to prove $\bigwedge_{i\in I} \Upi U \sqleq \Up U$, for which it suffices to prove $U\Leq \bigwedge_{i\in I}\Upi U$. Since each $\Leq_j$ satisfies~\eqref{axiom:cones give order}, to prove that $U\Leq_j \bigwedge_{i\in I}\Upi U$ it suffices to obtain the following two inclusions:
		\[
			U\sqleq \Downsub{j}\bigwedge_{i\in I}\Upi U
			\qquad\text{and}\qquad
			\bigwedge_{i\in I}\Upi U \sqleq \Upsub{j} U.
		\]
	The first inclusion follows by combining $U\sqleq \Downsub{j} U$ and $U\sqleq \bigwedge_{i\in I}\Upi U$, and the second inclusion is obvious. Since $j\in I$ was arbitrary, it follows that $U\Leq \bigwedge_{i\in I} \Upi U$, so by \cref{lemma:properties of localic cones}(a) the desired inclusion follows. The proof of the equation for the down cone is dual.
	
	That $\Leq$ satisfies~\eqref{axiom:cones give order} now follows straightforwardly. Namely, if $V\sqleq \Up U$ and $U\sqleq \Down V$, by the previous part of the proof it follows that for every $j\in I$ we get $V\sqleq \bigwedge_{i\in I}\Upi U \sqleq \Upsub{j} U$, and similarly $U\sqleq \Downsub{j} V$. Since $\Leq_j$ satisfies~\eqref{axiom:cones give order} this implies $U\Leq_j V$, and since this argument works for every $j\in I$ we get $U\Leq V$, as desired.
\end{proof}


\subsection{Orders from maps}
\label{example:ordered locale on domain}
	Consider a map of locales $f\colon X\to (Y,\Leq_Y)$. We describe the largest order $\Leq_f$ on $\Opens X$ that turns $X$ into an ordered locale and makes $f$ monotone. This generalises \cite[Example~3.10]{heunenSchaaf2024OrderedLocales} from sublocale inclusions to arbitrary morphisms.
	
	First, fixing $U\in \Opens X$, we write \[R_f(U):= \{ V\in\Opens Y: (U,V)\in R_f\}= \{V\in\Opens Y: U\sqleq f^{-1}(V)\}.\] We define the order on $\Opens X$ by:
		\[
			U\Leq_f U'
			\qquad\text{if and only if}\qquad
			\begin{array}{l}
				\forall V\in R_f(U)\colon \exists V'\in R_f(U')\colon V\Leq_YV' \text{~and}\\
				\forall V'\in R_f(U')\colon \exists V\in R_f(U)\colon V\Leq_YV'.
			\end{array}
		\]
	This is the precise characterisation of monotonicity in the diagrams of \cref{lemma:monotonicity unpacked}, so we can state this equivalently (albeit more informally) as:
		\[
		U \Leq_f U'
		\qquad\text{if and only if}\qquad
		\begin{tikzcd}[every label/.append style = {font = \normalsize},column sep=0.3cm, row sep=0.35cm]
			{U} & {U'} \\
			{V} & \highlighttikzcd{\exists V'}
			\arrow["R_f"{anchor=center}, draw=none, from=2-1, to=1-1]
			\arrow["R_f"{anchor=center}, draw=none, from=2-2, to=1-2]
			\arrow["{}"{description}, "\Leq_Y"{anchor=center}, draw=none, from=2-1, to=2-2]
		\end{tikzcd}
		\qquad\text{and}\qquad 
		\begin{tikzcd}[every label/.append style = {font = \normalsize},column sep=0.3cm, row sep=0.35cm]
			{U} & {U'} \\
			\highlighttikzcd{\exists V} & {V'.}
			\arrow["R_f"{anchor=center}, draw=none, from=2-1, to=1-1]
			\arrow["R_f"{anchor=center}, draw=none, from=2-2, to=1-2]
			\arrow["{}"{description}, "\Leq_Y"{anchor=center}, draw=none, from=2-1, to=2-2]
		\end{tikzcd}		
		\]
	The condition also has a simplified form using localic cones:
		\[
			U\Leq_f U'
			\qquad\text{if and only if}\qquad
			\begin{array}{l}
				\forall V\in R_f(U)\colon \Up V \in R_f(U') \text{~and}\\
				\forall V'\in R_f(U')\colon \Down V'\in R_f(U).
			\end{array}
		\]
	We claim that $\Leq_f$ defines an ordered locale structure on $X$. First, it is easy to see that $\Leq_f$ inherits reflexivity from $\Leq_Y$. For transitivity, use $U\Leq_f U'$ and $U'\Leq_f U''$ to calculate as follows:
		\[
			\begin{tikzcd}[every label/.append style = {font = \normalsize},column sep=0.3cm, row sep=0.35cm]
				{U} & {U'} & {U''}\\
				{V} & |[color=OliveGreen]|{\dashedboxtikz} & |[color=OliveGreen]|{\dashedboxtikz}
				\arrow["R_f"{anchor=center}, draw=none, from=2-1, to=1-1]
				\arrow["R_f"{anchor=center}, draw=none, from=2-2, to=1-2]
				\arrow["R_f"{anchor=center}, draw=none, from=2-3, to=1-3]
				\arrow["{}"{description}, "\Leq_Y"{anchor=center}, draw=none, from=2-1, to=2-2]
				\arrow["{}"{description}, "\Leq_Y"{anchor=center}, draw=none, from=2-2, to=2-3]
			\end{tikzcd}
		\qquad\text{so}\qquad
			\begin{tikzcd}[every label/.append style = {font = \normalsize},column sep=0.3cm, row sep=0.35cm]
				{U} & {U'} & {U''}\\
				{V} & \highlighttikzcd{\exists V'} & {\dashedboxtikz}
				\arrow["R_f"{anchor=center}, draw=none, from=2-1, to=1-1]
				\arrow["R_f"{anchor=center}, draw=none, from=2-2, to=1-2]
				\arrow["R_f"{anchor=center}, draw=none, from=2-3, to=1-3]
				\arrow["{}"{description}, "\Leq_Y"{anchor=center}, draw=none, from=2-1, to=2-2]
				\arrow["{}"{description}, "\Leq_Y"{anchor=center}, draw=none, from=2-2, to=2-3]
			\end{tikzcd}
		\qquad\text{so}\qquad
			\begin{tikzcd}[every label/.append style = {font = \normalsize},column sep=0.3cm, row sep=0.35cm]
				{U} & {U'} & {U''}\\
				{V} & {V'} & \highlighttikzcd{\exists V''}{\color{Black}.}
				\arrow["R_f"{anchor=center}, draw=none, from=2-1, to=1-1]
				\arrow["R_f"{anchor=center}, draw=none, from=2-2, to=1-2]
				\arrow["R_f"{anchor=center}, draw=none, from=2-3, to=1-3]
				\arrow["{}"{description}, "\Leq_Y"{anchor=center}, draw=none, from=2-1, to=2-2]
				\arrow["{}"{description}, "\Leq_Y"{anchor=center}, draw=none, from=2-2, to=2-3]
			\end{tikzcd}
		\]
	An analogous calculation holds for the dual condition, so it follows that $U\Leq_f U''$. Lastly, we need to check for~\eqref{axiom:V}. Take a family of related opens $U_i\Leq_f U_i'$. The case that the index is empty amounts to $\varnothing \Leq_f \varnothing$, which is covered by reflexivity, so we may assume the family is non-empty. We need to find $V'\in\Opens Y$ to complete the following diagram:
		\[
			\begin{tikzcd}[every label/.append style = {font = \normalsize},column sep=0.3cm, row sep=0.35cm]
				{\bigvee U_i} & {\bigvee U_i'} \\
				{V} & \highlighttikzcd{\exists V'}{\color{Black}.}
				\arrow["R_f"{anchor=center}, draw=none, from=2-1, to=1-1]
				\arrow["R_f"{anchor=center}, draw=none, from=2-2, to=1-2]
				\arrow["{}"{description}, "\Leq_Y"{anchor=center}, draw=none, from=2-1, to=2-2]
			\end{tikzcd}
		\]
	If $\left(\bigvee U_i,V\right)\in R_f$, then in particular the individual opens satisfy $U_i \sqleq f^{-1}(V)$. Since $U_i\Leq_f U_i'$, we can thus find $V_i'\in \Opens Y$ such that $U_i'\sqleq f^{-1}(V_i')$ and $V\Leq_Y V_i'$. Then $V'= \bigvee V_i'$ gives a desired open. Namely, by~\eqref{axiom:V} we get $V\Leq_Y \bigvee V_i'$, and moreover $\bigvee U_i'\sqleq \bigvee f^{-1}(V_i') = f^{-1}\left(\bigvee V_i'\right)$. The dual part of the proof is analogous.
	
	The conclusion is that $(X,\Leq_f)$ is an ordered locale, and by construction it makes $f$ monotone. It is in fact the supremum in $\OL(X)$ of all orders that make $f$ monotone. Namely, if ${\Leq'}\in \OL(X)$ also makes $f$ monotone, it follows directly from \cref{lemma:monotonicity unpacked} that ${\Leq'}\subseteq {\Leq}_f$. And since $\Leq_f$ itself makes $f$ monotone, it must be the largest such.
	
Now consider the case that $f$ defines a \emph{sublocale} (\cref{section:sublocales}), which for the purposes of this discussion means that $f^{-1}$ is a surjection, or equivalently that $f^{-1}\circ f_\ast = \id_{\Opens X}$ \cite[Lemma~IX.4.2]{maclane1994SheavesGeometryLogic}. In particular, this allows $U\sqleq f^{-1}\circ f_\ast(U)$ for every $U\in\Opens X$, i.e.~$f_\ast(U)\in R_f(U)$. Denoting the localic cones of $(X,\Leq_f)$ by $\Upsub{f}$ and $\Downsub{f}$, the monotonicity of $f$ then implies through \cref{proposition:monotonicity in terms of localic cones}:
	\[
		\Upsub{f} U = \Upsub{f} f^{-1}(f_\ast(U))\sqleq f^{-1}(\Upsub{Y}f_\ast(U)).
	\]
We do not think that this inclusion can be simplified in general to an explicit equation for $\Upsub{f}$ in terms of $\Upsub{Y}$, $f^{-1}$ and $f_\ast$. We only get:
	\[
		f_\ast\circ\Upsub{f}\sqleq \Upsub{Y}\circ f_\ast
		\qquad\text{and}\qquad
		f_\ast \circ\Downsub{f}\sqleq \Downsub{Y}\circ f_\ast.
	\]
More could be said in the case there exists of a further left adjoint $f_!\dashv f^{-1}$, e.g.~when $f$ is \emph{open} \cite[\S IX.7]{maclane1994SheavesGeometryLogic}, but we will not work this out here. 

We apply these ideas in \cref{section:IPs in double-negation sublocale}.

\subsection{Even more examples of ordered locales}
\begin{example}\label{example:ordered subframe}
	Consider an ordered locale $X$, and take a subframe $\Opens L \subseteq \Opens X$. We obtain an ordered locale $(L,\Leq)$, where the order is just that of $\Opens X$ restricted to $\Opens L$. Then the induced quotient map of locales $e\colon X\tworightarrow L$ is monotone if and only if $\Opens L$ is closed under the cones of $X$.
	
	To see this, denote for clarity the cones of $X$ by $\Upsub{X}$ and $\Downsub{X}$, and those of $L$ by $\Upsub{L}$ and $\Downsub{L}$, meaning that for $V\in \Opens L$:
		\[
			\Upsub{L} V = \bigvee\left\{W\in\Opens L: V\Leq W\right\}
			\qquad\text{and}\qquad
			\Downsub{L} V = \bigvee\left\{W\in\Opens L:W\Leq V\right\}.
		\]
	For any $V\in\Opens L$ we clearly get $\Upsub{L} V\sqleq \Upsub{X} V$, since the left-hand-side is simply a join over a smaller set than on the right-hand-side. Dually, we have the inclusion $\Downsub{L} V \sqleq \Downsub{X} V$. Since the frame map $e^{-1}\colon \Opens L \hookrightarrow \Opens X$ is just the inclusion function, from \cref{proposition:monotonicity in terms of localic cones} it can be seen that the monotonicity of $e$ reduces precisely to the converse inclusions: $\Upsub{X} V\sqleq \Upsub{L} V$ and $\Down _X V\sqleq \Downsub{L}V$. Hence $e$ is monotone if and only if the cones of $L$ coincide with the cones in $X$. 
	
	This, in turn, is equivalent to $L$ being closed under $\Upsub{X}$ and $\Downsub{X}$. To see this, note that this condition is now clearly implied by $e$ being monotone. Conversely, if $L$ is closed under the cones of $X$, then $V\Leq \Upsub{X} V$ holds in $\Opens L$, and implies $\Upsub{X} V\sqleq \Upsub{L} V$. Together with the dual condition, this implies $e$ is monotone.
\end{example}

\begin{remark}
	If $(X,\Leq)$ is an ordered locale with~\eqref{axiom:LV}, then $\im(\Up)$ and $\im(\Down)$ are subframes of $\Opens X$ (\cref{lemma:past sets form a frame}). However, in general, $\im(\Up)$ is not closed under~$\Down$, and dually $\im(\Down)$ is not closed under $\Up$. Thus the previous example does not produce monotone maps of ordered locales in this case.
\end{remark}

\begin{example}\label{example:comonads give ordered locale}
	Consider a pair $(a,b)$ of \emph{co}monads on $\Opens X$, which are just monads on the opposite poset $\Opens X^\op$, and define the following relation:
		\[
			U\Leq V
			\qquad\text{if and only if}\qquad
			a(U)\sqleq V\text{~and~} b(V)\sqleq U.
		\]
	That $\Leq$ is a preorder follows from a similar argument to that in \cref{section:ordered locales from monads} for monads, but using the \emph{co}unit and \emph{co}multiplication instead. In order for $\Leq$ to satisfy~\eqref{axiom:V}, the \cref{lemma:monotone laxly respects meets and joins} can no longer be used, so we need to assume that $a$ and $b$ themselves preserve joins. In that case, if $U_i\Leq V_i$, we get $a(U_i)\sqleq V_i$, and so
		\[
			a\left(\bigvee U_i\right)\sqleq \bigvee a(U_i)\sqleq \bigvee V_i,
		\]
	and similarly $b\left(\bigvee V_i\right)\sqleq \bigvee U_i$, so $\bigvee U_i\Leq \bigvee V_i$. Thus we obtain an ordered locale. We know of no natural examples of join-preserving comonads, but point out that \cite{akbartabatabai2021ImplicationSpacetime} uses a join-preserving map with a counit to model a ``temporal modal'' operator. We study the connection to this work more in \cref{section:caucal heyting implication}.
\end{example}

\section{Parallel ordered locales}\label{section:parallel ordered locales}
\begin{figure}[t]
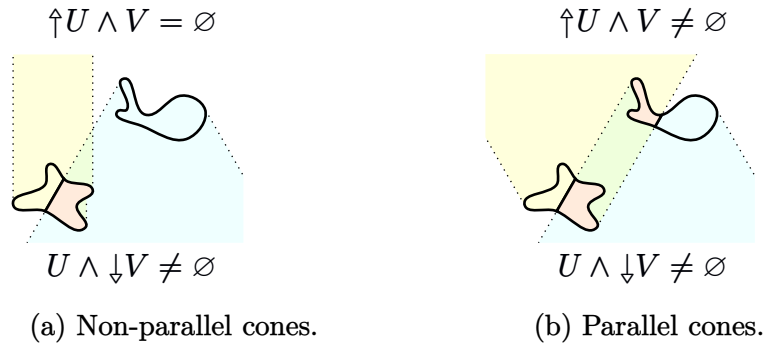
\centering
	\begin{subfigure}[b]{0.3\textwidth}
		\input{tikz/inkscape-export/fig-24-counter}
		\caption{Non-parallel cones.}
	\end{subfigure}\hfil
	\begin{subfigure}[b]{0.3\textwidth}
		\input{tikz/inkscape-export/fig-24-parallel}
		\caption{Parallel cones.}
	\end{subfigure}
	\caption{Intuition of non-parallel vs.~parallel cones.}
	\label{figure:parallel vs nonparallel ordered}
\end{figure}
We have seen in \cref{example:order from pair of monads} that any pair of monads on the frame of opens of a locale determine an order structure on it, where the monads serve as the localic cones. In fact, by \cref{proposition:ordered locales are bimonads} we know that virtually all ordered locales are of this form. But we have imposed no relation between these monads whatsoever. This contradicts the behaviour of a preorder $\leq$ on a space, whose cones run ``parallel'' in the sense that
	\[
		x\in \down\{y\}
		\quad\text{if and only if}\quad	
		x\leq y
		\quad\text{if and only if}\quad	
		y\in\up\{x\}.
	\]
In other words: $x$ is in the past of $y$ if and only if $y$ is in the future of $x$. In the localic world this can fail. An explicit counterexample to this is in \cref{example:minkowski different speeds of light}, where we define localic cones on $\mathbb{R}^2$ such that $x$ is in the past of~$y$, but $y$ is \emph{not} in the future of~$x$. See \cref{figure:parallel vs nonparallel ordered} for visual intuition of what goes wrong. To remedy this, here we introduce a special class of ordered locales where the localic cones have to be related in such a way that excludes this behaviour. They are called the \emph{parallel ordered locales}.

\begin{definition}\label{definition:parallel ordered locales}
	If $(X,\Leq)$ is an ordered locale, we say $\Leq$ \emph{respects meets} if the following two laws hold:
	\begin{center}\vspace{-4mm}
		\dummylabel{axiom:wedge}{$\wedge^\pm$}%
		\begin{minipage}{.3\linewidth}
			\[\tag{$\wedge^+$}\label{axiom:wedge+}
			\begin{tikzcd}[every label/.append style = {font = \normalsize},column sep=0.25cm, row sep=0.2cm]
				{U} & \highlighttikzcd{\exists U'} \\
				{V} & {V'}
				\arrow["\sqleq"{anchor=center, rotate=-90}, draw=none, from=2-1, to=1-1]
				\arrow["\sqleq"{anchor=center, rotate=-90}, draw=none, from=2-2, to=1-2]
				\arrow["{}"{description}, "\Leq"{anchor=center}, draw=none, from=1-1, to=1-2]
				\arrow["{}"{description}, "\Leq"{anchor=center}, draw=none, from=2-1, to=2-2]
			\end{tikzcd}
			\]
		\end{minipage}%
		\hfil
		\begin{minipage}{.3\linewidth}
			\[\tag{$\wedge^-$}\label{axiom:wedge-}
			\begin{tikzcd}[every label/.append style = {font = \normalsize},column sep=0.25cm, row sep=0.2cm]
				\highlighttikzcd{\exists U} & {U'} \\
				{V} & {V'.}
				\arrow["\sqleq"{anchor=center, rotate=-90}, draw=none, from=2-1, to=1-1]
				\arrow["\sqleq"{anchor=center, rotate=-90}, draw=none, from=2-2, to=1-2]
				\arrow["{}"{description}, "\Leq"{anchor=center}, draw=none, from=1-1, to=1-2]
				\arrow["{}"{description}, "\Leq"{anchor=center}, draw=none, from=2-1, to=2-2]
			\end{tikzcd}
			\]
		\end{minipage}
	\end{center}
	We say $\Leq$ \emph{preserves bottom} if the bottom element of $\Opens X$ is only related to itself:
	\[\tag{$\varnothing$}\label{axiom:empty}
	U\Leq \varnothing \Leq V \quad\text{implies}\quad U=\varnothing=V.
	\]
	Lastly, we say $(X,\Leq)$ is \emph{parallel ordered} if all three axioms~\eqref{axiom:wedge} and~\eqref{axiom:empty} hold. Denote by $\OrdLoc_{\newparallel}$ the full subcategory of ordered locales consisting of parallel ordered locales.
\end{definition}

\begin{figure}[t]
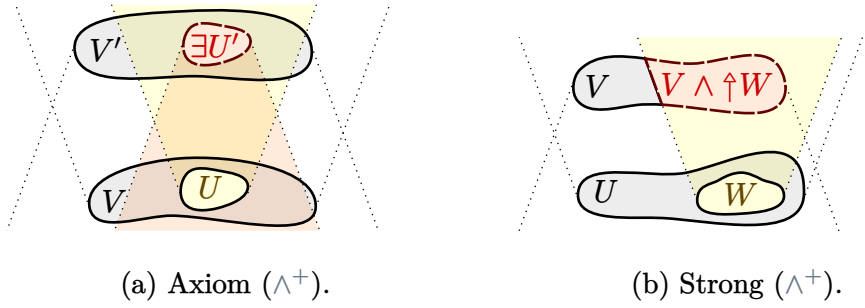
\centering
	\begin{subfigure}[b]{0.4\textwidth}
		\input{tikz/inkscape-export/fig-23}
		\caption{Axiom~\eqref{axiom:wedge+}.}
	\end{subfigure}\hfil
	\begin{subfigure}[b]{0.4\textwidth}
		\input{tikz/inkscape-export/fig-23-strong}
		\caption{Strong~\eqref{axiom:wedge+}.}
	\end{subfigure}
	\caption{Illustration of~\eqref{axiom:wedge+} and \cref{lemma:wedge iff strong wedge} in a Minkowski-like space.}
	\label{figure:axiom wedge+}
\end{figure}

\begin{remark}
	What is meant by ``respects meets'' is explained further in \cref{lemma:wedge iff frobenius,lemma:wedge iff strong wedge}. The intended reading of axiom~\eqref{axiom:wedge+} is analogous to how we have interpreted other diagrams of this form: if $U\sqleq V$ and $V\Leq V'$, \emph{then} there exists $U'$ such that $U\Leq U'$ and $U'\sqleq V'$. Axiom~\eqref{axiom:wedge-} is parsed analogously. The visual intuition behind these axioms is portrayed in \cref{figure:axiom wedge+}(a).
\end{remark}

\begin{lemma}\label{lemma:wedge implies cones determine order}
	In any ordered locale satisfying~\eqref{axiom:wedge}, axiom~\eqref{axiom:cones give order} holds.
\end{lemma}
\begin{proof}
	The ``only if'' direction in~\eqref{axiom:cones give order} holds in any ordered locale by \cref{lemma:properties of localic cones}(a). For the other direction, suppose that $U\sqleq \Down V$ and $V\sqleq \Up U$. Using \cref{lemma:properties of localic cones}(c) we hence get $U\Leq \Up U\sqgeq V$ and $U\sqleq \Down V\Leq V$. Applying axioms~\eqref{axiom:wedge-} and~\eqref{axiom:wedge+}, respectively, we get:
	\[
	\begin{tikzcd}[every label/.append style = {font = \normalsize},column sep=0.25cm, row sep=0.2cm]
		\highlighttikzcd{\exists U'} & {V} \\
		{U} & {\Up U}
		\arrow["\sqleq"{anchor=center, rotate=-90}, draw=none, from=2-1, to=1-1]
		\arrow["\sqleq"{anchor=center, rotate=-90}, draw=none, from=2-2, to=1-2]
		\arrow["{}"{description}, "\Leq"{anchor=center}, draw=none, from=1-1, to=1-2]
		\arrow["{}"{description}, "\Leq"{anchor=center}, draw=none, from=2-1, to=2-2]
	\end{tikzcd}
	\qquad\text{and}\qquad
	\begin{tikzcd}[every label/.append style = {font = \normalsize},column sep=0.25cm, row sep=0.2cm]
		{U} & \highlighttikzcd{\exists V'} \\
		{\Down V} & {V.}
		\arrow["\sqleq"{anchor=center, rotate=-90}, draw=none, from=2-1, to=1-1]
		\arrow["\sqleq"{anchor=center, rotate=-90}, draw=none, from=2-2, to=1-2]
		\arrow["{}"{description}, "\Leq"{anchor=center}, draw=none, from=1-1, to=1-2]
		\arrow["{}"{description}, "\Leq"{anchor=center}, draw=none, from=2-1, to=2-2]
	\end{tikzcd}
	\]
	Applying~\eqref{axiom:V} to the top two rows then gives $U\vee U'\Leq V\vee V'$, which reduces to $U\Leq V$, as desired.
\end{proof}

\begin{remark}\label{lemma:LV implies empty}
	Note that in any ordered locale satisfying~\eqref{axiom:cones give order}, the axiom~\eqref{axiom:empty} holds if and only if $\Up\varnothing = \varnothing = \Down\varnothing$. So an ordered locale is parallel ordered if and only if~\eqref{axiom:wedge} holds and the localic cones preserve the empty region. We note in passing that axiom~\eqref{axiom:LV} implies axiom~\eqref{axiom:empty}, since the empty join is the empty region.
\end{remark}

To explain the terminology ``parallel ordered,'' and to justify the intuition outlined at the start of this section and in \cref{figure:parallel vs nonparallel ordered}, we need to show how the axioms~\eqref{axiom:wedge} can be reformulated in terms of statements about the localic cones.

\begin{lemma}\label{lemma:wedge iff frobenius}
	An ordered locale $(X,\Leq)$ satisfies~\eqref{axiom:wedge} if and only if~\eqref{axiom:cones give order} holds and the localic cones satisfy the following two conditions:
	\begin{center}\vspace{-4mm}
		\dummylabel{axiom:frobenius}{\textsc{f}$^\pm$}
		\begin{minipage}{.4\linewidth}
			\[\tag{\textsc{f}$^-$}\label{axiom:frobenius+}
			\Up U\wedge V\sqleq \Up\left(U\wedge \Down V\right)
			\]
		\end{minipage}%
		\hfil
		\begin{minipage}{.4\linewidth}
			\[\tag{\textsc{f}$^+$}\label{axiom:frobenius-}
			\Down U\wedge V\sqleq \Down\left(U\wedge \Up V\right).
			\]
		\end{minipage}
	\end{center}
\end{lemma}
\begin{proof}
	We already know by \cref{lemma:wedge implies cones determine order} that~\eqref{axiom:wedge} imply~\eqref{axiom:cones give order}. To show that~\eqref{axiom:wedge+} implies~\eqref{axiom:frobenius-}, note that by \cref{lemma:properties of localic cones}(c) we get $\Down U\wedge V\sqleq \Down U \Leq U$, so~\eqref{axiom:wedge+} says there exists an open region $W$ such that $\Down U\wedge V\Leq W\sqleq U$. In turn, \cref{lemma:properties of localic cones}(a) and (e) give $W\sqleq \Up(\Down U \wedge V)\sqleq \Up\Down U\wedge \Up V$, so we find $W\sqleq U\wedge \Up V$. Applying \cref{lemma:properties of localic cones}(a) to $\Down U\wedge V\Leq W$ then gives the desired inclusion:
	\[
	\Down U\wedge V\sqleq \Down W \sqleq \Down \left(U\wedge \Up V\right).
	\]

	Conversely, we need to prove that~\eqref{axiom:frobenius-} implies~\eqref{axiom:wedge+} under~\eqref{axiom:cones give order}. Given ${U\sqleq V\Leq V'}$, define $U':= \Up U \wedge V'\sqleq V'$. Since $U\sqleq V\sqleq \Down V'$ by \cref{lemma:properties of localic cones}(a), we find using~\eqref{axiom:frobenius-} that
	\[
	U=U\wedge \Down V' \sqleq \Down \left(V'\wedge \Up U\right) = \Down U'.
	\]
	That $U'\sqleq \Up U$ follows by construction, so~\eqref{axiom:cones give order} gives $U\Leq U'$, as desired, hence showing that~\eqref{axiom:wedge+} holds. That~\eqref{axiom:wedge-} follows from~\eqref{axiom:frobenius+} is proved similarly.
\end{proof}

\begin{figure}[b]\centering
	\input{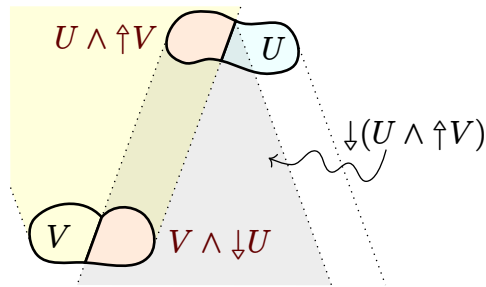}
	\caption{Illustration of~\eqref{axiom:frobenius-} in a Minkowski-like space.}
	\label{figure:axiom:frobenius-}
\end{figure}

\begin{remark}
	These new conditions are visualised in \cref{figure:axiom:frobenius-}. Note the resemblance between these laws and the so-called ``Frobenius reciprocity condition'' of open locale maps \cite[\S IX.7]{maclane1994SheavesGeometryLogic}.
\end{remark}

\begin{proposition}
	Given a locale $X$, there is a bijective correspondence between preorders ${\Leq}$ on $\Opens X$ such that $(X,\Leq)$ is an ordered locale satisfying~\eqref{axiom:wedge}, and pairs of monads $(u,d)$ on $\Opens X$ satisfying~\eqref{axiom:frobenius}.
\end{proposition}
\begin{proof}
	This follows immediately by combining \cref{proposition:ordered locales vs monads} and \cref{lemma:wedge implies cones determine order,lemma:wedge iff frobenius}.
\end{proof}

The following is a very useful reformulation of axioms~\eqref{axiom:wedge}, and is the one we most often use in computations.
\begin{lemma}\label{lemma:wedge iff strong wedge}
	An ordered locale $(X,\Leq)$ satisfies~\eqref{axiom:wedge} if and only if the following two laws are satisfied:
	\[
	\begin{tikzcd}[every label/.append style = {font = \normalsize},column sep=0.25cm, row sep=0.2cm]
		{W} & \highlighttikzcd{\Up W\wedge V} \\
		{U} & {V}
		\arrow["\sqleq"{anchor=center, rotate=-90}, draw=none, from=2-1, to=1-1]
		\arrow["\sqleq"{anchor=center, rotate=-90}, draw=none, from=2-2, to=1-2]
		\arrow["{}"{description}, "\Leq"{anchor=center}, draw=none, from=1-1, to=1-2]
		\arrow["{}"{description}, "\Leq"{anchor=center}, draw=none, from=2-1, to=2-2]
	\end{tikzcd}
	\qquad\text{and}\qquad
	\begin{tikzcd}[every label/.append style = {font = \normalsize},column sep=0.25cm, row sep=0.2cm]
		\highlighttikzcd{\Down W\wedge U} & {W} \\
		{U} & {V.}
		\arrow["\sqleq"{anchor=center, rotate=-90}, draw=none, from=2-1, to=1-1]
		\arrow["\sqleq"{anchor=center, rotate=-90}, draw=none, from=2-2, to=1-2]
		\arrow["{}"{description}, "\Leq"{anchor=center}, draw=none, from=1-1, to=1-2]
		\arrow["{}"{description}, "\Leq"{anchor=center}, draw=none, from=2-1, to=2-2]
	\end{tikzcd}
	\]
	
	\begin{proof}
		Clearly these conditions imply~\eqref{axiom:wedge}. To show the converse, note that in the second square we get from~\eqref{axiom:wedge-} that there exists $W'\in\Opens X$ such that ${U\sqgeq W'\Leq W}$. Hence, using \cref{lemma:properties of localic cones} and~\eqref{axiom:frobenius+} we get
		\[
		W =\Up W'\wedge W\sqleq \Up U \wedge W \sqleq \Up \left(U \wedge \Down W\right),
		\]
		so that by~\eqref{axiom:cones give order} we get $\Down W\wedge U\Leq W$. The proof for the first square is similar.
	\end{proof}
\end{lemma}

\begin{remark}
	From this point on we use~\eqref{axiom:wedge} to refer interchangeably to the squares in either \cref{lemma:wedge iff strong wedge} or \cref{definition:parallel ordered locales}. The visual intuition of these stronger axioms is in \cref{figure:axiom wedge+}(b).
	
	A remark analogous to \cref{remark:axiom L vs V} applies here: for a more general setting of ``ordered lattices'' where meets do not exist, axioms~\eqref{axiom:wedge} can still be validly interpreted, since they only depend on the inclusion order $\sqleq$. The previous lemma just says that, in the case of frames, the axioms take a particularly nice form.
\end{remark}

\begin{lemma}\label{lemma:parallel order gives order between arbitrary opens}
	In a parallel ordered locale $(X,\Leq)$, for arbitrary $U,V\in\Opens X$:
	\[
	U\wedge \Down V \Leq V\wedge \Up\left( U \wedge \Down V \right)
	\qquad\text{and}\qquad
	U\wedge \Down \left(\Up U \wedge V\right) \Leq \Up U \wedge V.
	\]
\end{lemma}
\begin{proof}
	By \cref{lemma:wedge implies cones determine order,lemma:wedge iff frobenius} it suffices to give an argument using the localic cones utilising~\eqref{axiom:frobenius}. We prove $U\wedge \Down V \Leq V\wedge \Up\left( U \wedge \Down V \right)$. Trivially, the right hand term is contained in the future of the left hand term. For the other inclusion, use~\eqref{axiom:frobenius-} to calculate:
		\[
			U\wedge \Down V
			=
			\Down V \wedge (U\wedge \Down V)
			\sqleq
			\Down(V\wedge \Up(U\wedge\Down V)),
		\]
	which is precisely the desired inclusion. The other case is proved dually using~\eqref{axiom:frobenius+}.
\end{proof}

Of course, generally, the intersection $U\wedge \Down V$ may well be empty, so that the lemma tells us simply that $\varnothing\Leq \varnothing$. While seemingly trivial, this is actually an important case (see the next proposition), and captures the essence of what we mean by parallel orderedness. Later, it will tell us (\cref{lemma:heyting negation of past and future set}) that Heyting negation respects the order structure. The visual idea is captured in \cref{figure:parallel vs nonparallel ordered}.

\begin{proposition}\label{proposition:cones are parallel}
	For any $U,V\in\Opens X$ in a parallel ordered locale $(X,\Leq)$ we have:
	\[
	U\wedge \Down V =\varnothing
	\qquad\text{if and only if}\qquad
	\Up U \wedge V =\varnothing.
	\]
\end{proposition}
\begin{proof}
	Assume that $U\wedge\Down V =\varnothing$. Using~\eqref{axiom:empty} and~\eqref{axiom:frobenius+} immediately gives that $\Up U\wedge V \sqleq \Up (U\wedge\Down V)= \Up \varnothing =\varnothing$. The other implication follows similarly using~\eqref{axiom:frobenius-}.
\end{proof}

To finish this section, we discuss some examples.

\begin{proposition}\label{proposition:space with OC is parallel ordered}
	If $(S,\leq)$ has open cones, then $\loc(S)$ is parallel ordered.
\end{proposition}
\begin{proof}
	From the characterisation of the localic cones in a space, it follows that $\Up \varnothing = (\up \varnothing)^\circ = \varnothing$, and similarly $\Down \varnothing = \varnothing$, so $\loc(S)$ satisfies~\eqref{axiom:empty} (even without open cones).
	
	In a space with open cones, axiom~\eqref{axiom:frobenius+} says that $\up U \cap V \subseteq \up\left(U\cap \down V\right)$. This is almost trivially true: if $y\in \up U\cap V$, then there exists $x\in U$ such that $x\leq y$. But then $x\in \down \{y\}\subseteq \down V$, so $y \in \up\left(U\cap \down V\right)$. The dual axiom~\eqref{axiom:frobenius-} follows similarly, and since ordered spaces satisfy~\eqref{axiom:cones give order}, the result follows by \cref{lemma:wedge iff frobenius}.
\end{proof}

\begin{example}\label{example:non-OC is not parallel}
	Locales induced by spaces without open cones are not always parallel ordered. Consider the space $\{\ast\}\sqcup \mathbb{R}$ from \cref{example:non-OC space,figure:non-OC space}, where $\ast\leq 0$. We can then see that~\eqref{axiom:frobenius-} is violated, since it would imply a contradiction:
		\[
			\{\ast\} = \{\ast\}\cap \Down \mathbb{R} \subseteq \Down\left( \Up \{\ast\}\cap \mathbb{R}\right) = \Down \varnothing = \varnothing.
		\]
\end{example}

\begin{example}\label{example:lower order is not parallel}
	Starting with an ordered space $(S,\leq)$ with open cones, and equipping its frame of opens $\Opens S$ with the lower order ($U\Leq_\mathrm{L} V$ if and only if $U\subseteq \down V$) yields an ordered locale that satisfies~\eqref{axiom:wedge}, but does \emph{not} satisfy~\eqref{axiom:empty}. Namely, we find that $\Up U = S$ for every $U\in\Opens S$. Hence $\Up\varnothing \wedge V = V$, but $\varnothing\wedge\Down V =\varnothing$. This is further evidence in favour of the Egli-Milner order.
\end{example}

\begin{figure}[t]
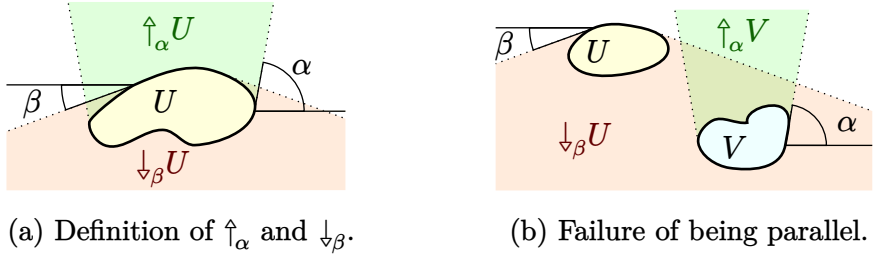
\centering
	\begin{subfigure}[b]{0.33\textwidth}\centering
		\input{tikz/inkscape-export/fig-25}
		\caption{Definition of $\Upsub{\alpha}$ and $\Downsub{\beta}$.}
	\end{subfigure}\hfil
	\begin{subfigure}[b]{0.45\textwidth}\centering
		\input{tikz/inkscape-export/fig-25-counter}
		\caption{Failure of being parallel.}
		\label{figure:minkowski different cones not parallel}
	\end{subfigure}
	\caption{Illustration of the cones in \cref{example:minkowski different speeds of light}.}
	\label{figure:minkowski diferent cones}
\end{figure}

\begin{example}\label{example:minkowski different speeds of light}
	We construct a modified example of Minkowski space where the speed of light into the future is independent from the speed of light into the past. To do this, we consider the underlying space $\mathbb{R}^2$, and equip it with a family of preorders $\leq_\alpha$ defined by the Minkowski metric with speed of light $1/\alpha$ (\cref{section:spacetimes}). Explicitly:
		\[
			x\leq_\alpha y
			\qquad\text{if and only if}\qquad
			\parbox{.35\textwidth}{the straight line connecting $x$ and $y$ makes an angle with the horizontal axis of at least $\alpha$ radians.}
		\]
	For $\alpha\in (0,\pi/2)$ the pair $(\mathbb{R}^2,\leq_\alpha)$ comes from a spacetime, and so must have open cones (\cref{corollary:spacetimes have open cones}). In particular we get monads $\Upsub{\alpha}$ and $\Downsub{\alpha}$ on $\Opens\mathbb{R}^2$.
	
	For any two numbers $\alpha,\beta\in (0,\pi/2)$ we therefore get an ordered locale $\loc(\mathbb{R}^2)$ whose order $\Leq$ is defined by the pair of cones $\Upsub{\alpha}$ and $\Downsub{\beta}$. See \cref{figure:minkowski diferent cones}(a). We claim that this locale is parallel ordered if and only if $\alpha = \beta$. In the latter case, where the slopes are indeed equal, the ordered locale $\loc(\mathbb{R}^2)$ is just the one coming from Minkowski space, which is parallel ordered by \cref{proposition:space with OC is parallel ordered}. Conversely, suppose for the sake of contradiction that $\alpha\neq \beta$. It is then straightforward to construct an example, such as in \cref{figure:minkowski diferent cones}(b), that contradicts the conclusion of \cref{proposition:cones are parallel}, showing that $\loc(\mathbb{R}^2)$ is not parallel ordered.
\end{example}

\begin{remark}[Summary]
	For any space $(S,\leq)$ with open cones, $(\loc(S),\Leq)$ is an ordered locale that satisfies all of the axioms~\eqref{axiom:cones give order},~\eqref{axiom:LV}, and~\eqref{axiom:wedge}. We therefore propose this as an abstract setting in which to study point-free spacetimes. Throughout the rest of the thesis, we will see these axioms play a crucial role in several abstract constructions inspired by relativity theory, such as \emph{causal boundaries} (\cref{section:causal boundaries}) and \emph{domains of dependence} (\cref{section:causal coverage}). 
\end{remark}


\chapter{An adjunction}\label{section:adjunction}
We now turn to a more mathematically based justification for the \cref{definition:ordered locale} of ordered locales. Namely, we prove that the adjunction between topological spaces and locales lifts, after certain technical restrictions, to an adjunction between the categories of ordered topological spaces and of ordered locales:
	\[
	\begin{tikzcd}[cramped]
		\OrdTop & \OrdLoc \\
		\Top & {\Loc.}
		\arrow[""{name=0, anchor=center, inner sep=0}, "\loc", shift left=2, from=2-1, to=2-2]
		\arrow[""{name=1, anchor=center, inner sep=0}, "\pt", shift left=2, from=2-2, to=2-1]
		\arrow[""{name=2, anchor=center, inner sep=0}, shift left=2, from=1-1, to=1-2]
		\arrow[""{name=3, anchor=center, inner sep=0}, shift left=2, from=1-2, to=1-1]
		\arrow[shift left=2, hook, from=2-1, to=1-1]
		\arrow[shift left=2, hook, from=2-2, to=1-2]
		\arrow[shift left=2, from=1-1, to=2-1]
		\arrow[shift left=2, from=1-2, to=2-2]
		\arrow["\dashv"{anchor=center, rotate=-90}, draw=none, from=0, to=1]
		\arrow["\dashv"{anchor=center, rotate=-90}, draw=none, from=2, to=3]
	\end{tikzcd}
	\]
The idea here is that we start with the functors $\loc$ and $\pt$ and extend them to the ordered setting. To be precise, we require that they are compatible with the obvious forgetful functors:
	\[
		\begin{tikzcd}[ampersand replacement=\&]
			\OrdTop \& \OrdLoc \\
			\Top \& {\Loc,}
			\arrow["\loc", shift left=2, from=2-1, to=2-2]
			\arrow["\pt", shift left=2, from=2-2, to=2-1]
			\arrow["\loc", shift left=2, from=1-1, to=1-2]
			\arrow["\pt", shift left=2, from=1-2, to=1-1]
			\arrow[from=1-1, to=2-1]
			\arrow[from=1-2, to=2-2]
		\end{tikzcd}
		\qquad\text{or explicitly:}\qquad
		\begin{tikzcd}[ampersand replacement=\&,row sep=0ex]
			{(S,\leq)\longmapsto \left(\loc(S),\dashedboxtikz\right),} \\
			{(X,\Leq)\longmapsto \left(\pt(X),\dashedboxtikz\right).}
		\end{tikzcd}
	\]
We need to fill in the dashed boxes. In other words, given an ordered space $(S,\leq)$ we need to define a suitable order $\Leq$ on $\Opens \loc(S)=\Opens S$, and given an ordered locale $(X,\Leq)$ we need to define a suitable preorder $\leq$ on the space $\pt(X)$. In terms of morphisms, we need to show that if $g\colon S\to T$ is a continuous monotone function between ordered spaces, then $\loc(g)\colon \loc(S)\to \loc(T)$ is a monotone map of locales, and that if $f\colon X\to Y$ is a monotone map of locales, then the function $\pt(f)\colon \pt(X)\to \pt(Y)$ is continuous and monotone. 

This section is based on the main result of~\cite{heunenSchaaf2024OrderedLocales}, and is adopted here largely unchanged, except for extended exposition and more detailed proofs.

\section{Recap of the adjunction between spaces and locales}
\label{section:adjunction top and loc}
To establish notation, we first recall the ordinary adjunction between $\Top$ and $\Loc$. We give a detailed construction of this adjunction (\cref{theorem:adjunction top and loc}), since the proof of the new adjunction (\cref{theorem:adjunction ordtopOC bullet and ordloc bullet}) relies on it directly. More details can be found in \cite{johnstone1982StoneSpacesa,vickers1989TopologyLogic}, and our notation sticks particularly closely to~\cite[Chapter~IX]{maclane1994SheavesGeometryLogic}.

We have actually already encountered all the ingredients that constitute the functor $\loc\colon \Top\to\Loc$ previously. Now we put them together formally:

\begin{definition}
	There is a functor:\vspace{-1em}
	\[
		\parbox{.2\textwidth}{\begin{align*}
				\loc: \Top&\longrightarrow \Loc;\\
				S & \longmapsto \loc(S);\\
				g & \longmapsto \loc(g);
		\end{align*}}
		\quad\text{where}\quad
		\begin{aligned}
			\text{$\loc(S)$ is the locale defined by $\Opens \loc(S):=\Opens S$;}\\
			\text{$\loc(g)$ is the map defined by $\loc(g)^{-1}:= g^{-1}$.}
		\end{aligned}
	\]
\end{definition}

In more words, $\loc(S)$ is the locale whose frame of opens is defined by the concrete lattice of open subsets of $S$. This lattice is a frame precisely by definition of a topology. Similarly, if $g\colon S\to T$ is a continuous function, then $\loc(g)\colon \loc(S)\to \loc(T)$ is the map of locales whose preimage map is $\loc(g)^{-1}= g^{-1}\colon \Opens T\to \Opens S$, which is of course well defined precisely because $g$ is continuous.

\begin{lemma}\label{lemma:loc(g) frame map}
	If $g\colon S\to T$ is a continuous function, then $g^{-1}\colon \Opens T\to \Opens S$ is a map of frames.
\end{lemma}
\begin{proof}
	That the preimage map $g^{-1}$ sends open subsets $V\in \Opens T$ to an open subset $g^{-1}(V)\in \Opens S$ is equivalent to continuity. That this defines a map of frames follows from elementary properties of preimages.
\end{proof}

Next we describe the functor $\pt\colon \Loc\to \Top$. On the level of objects, given any locale $X$, we need to construct a \emph{space of points} $\pt(X)$. Recall from \cref{section:points of locales} that a \emph{point} of a locale $X$ is a completely prime filter $\calF\subseteq \Opens X$, representing a system of open regions ``zooming in'' on an ideal ``point.'' We remind ourselves that this means $\calF$ is a filter:
\begin{enumerate}[label = (F\arabic*)]\setcounter{enumi}{-1}
	\item $X \in \calF$; \hfill(\emph{non-empty})
	\item if $U,V\in\calF$ then $U\wedge V\in\calF$; \hfill(\emph{downwards directedness})
	\item if $U\in\calF$ and $U\sqleq V$, then $V\in \calF$; \hfill(\emph{upwards closure})
\end{enumerate}
that is moreover \emph{completely prime}, meaning that
\begin{enumerate}[label = (F\arabic*)]\setcounter{enumi}{2}
	\item $\varnothing\notin \calF$; \hfill(\emph{proper})
	\item if $\bigvee_{i\in I}U_i\in\calF$, then $U_i\in\calF$ for some $i\in I$. \hfill(\emph{inaccessibility by joins})
\end{enumerate}

\begin{remark}[Explanation]
	It is helpful to unpack the intuition behind this definition. We think of a completely prime filter $\calF\subseteq \Opens X$ as a system of open regions that contain an ``imaginary'' point. If we call this point $x$, then we get the informal rule:
		\[
			U\in \calF
			\qquad\text{if and only if}\qquad
			 x\in U,
		\]
	and the axioms take a more familiar form:
	\begin{enumerate}[label = (F\arabic*)]\setcounter{enumi}{-1}
		\item $x\in X$;
		\item if $x\in U$ and $x\in V$, then $x\in U\wedge V$;
		\item if $x\in U$ and $U\sqleq V$, then $x\in V$;
		\item $x\notin \varnothing$;
		\item if $x\in \bigvee_{i\in I} U_i$, then there is $i\in I$ such that $x\in U_i$.
	\end{enumerate}
\end{remark}

This intuitive rule can be made literal inside of a space.

\begin{lemma}\label{lemma:neighbourhood filter is completely prime}
	If $x\in S$ is an element of a space, then its neighbourhood filter defines a completely prime filter in $\Opens S$:
		\[
			\calF_x:=\{U\in\Opens S: x\in U \}.
		\]
\end{lemma}
\begin{proof}
	This follows by elementary properties of sets and elements.
\end{proof}

\begin{definition}\label{definition:space of points}
	For a locale $X$, we denote by $\pt(X)$ the set of completely prime filters $\calF\subseteq \Opens X$. We define a topology on this set by declaring the opens to be exactly those subsets of the form
		\[
			\pt(U):= \{\calF\in\pt(X): U\in \calF \},
		\]
	where $U\in\Opens X$. That this defines a topology follows from the next lemma.
\end{definition}

\begin{lemma}\label{lemma:properties of pt(U)}
	In the space of points of any locale:
	\begin{enumerate}[label = (\alph*)]
		\item $\pt(\varnothing) = \varnothing$ and $\pt(X)= \pt(X)$;
		\item $\pt(U\wedge V) = \pt(U)\cap \pt(V)$;
		\item $\pt\left(\bigvee U_i\right) =\bigcup \pt(U_i)$.
	\end{enumerate}
\end{lemma}
\begin{proof}
	For (a), by (F3) there is no $\calF\in\pt(X)$ such that $\varnothing\in\calF$, and the second (seemingly trivial) equation follows by (F0). The equation in (b) follows from (F1) and (F2), and similarly (c) follows from (F4) and (F2).
\end{proof}

\begin{definition}\label{definition:pt(f)}
	Let $f\colon X\to Y$ be a morphism of locales. We define the function
		\[
			\pt(f)\colon \pt(X)\longrightarrow \pt(Y);
			\qquad \calF \longmapsto \pt(f)(\calF):=\left\{V\in\Opens Y: f^{-1}(V)\in \calF\right\}.
		\]
\end{definition}

\begin{lemma}\label{lemma:pt(f) continuity}
	For any map of locales $f\colon X\to Y$, the map $\pt(f)$ is a well-defined continuous function. Moreover, $\pt(f)^{-1}(\pt(V)) = \pt(f^{-1}(V))$ for every $V\in\Opens Y$.
\end{lemma}
\begin{proof}
	We need to prove that $\pt(f)(\calF)$ is a completely prime filter in $\Opens Y$. This follows from the fact that $f^{-1}$ is a map of frames, and that $\calF$ itself is a completely prime filter in $\Opens X$. Explicitly, (F0) and (F3) hold because $f^{-1}$ preserves top and bottom elements, and (F2) holds since $f^{-1}$ is monotone. For (F1), if ${V,W\in\pt(f)(\calF)}$ then $f^{-1}(V),f^{-1}(W)\in \calF$, and so $f^{-1}(V\wedge W) = f^{-1}(V)\wedge f^{-1}(W)\in \calF$, which just means $V\wedge W\in \pt(f)(\calF)$. Property (F4) is proved similarly.
	
	To prove continuity it suffices to prove the given equation, which in turn can be derived as follows:
		\begin{align*}
			\pt(f)^{-1}\left(\pt(V)\right)
			&=
			\left\{\calF\in \pt(X):\pt(f)(\calF)\in\pt(V) \right\}
			\\&=
			\left\{\calF\in\pt(X): V\in\pt(f)(\calF)\right\}&&\text{(definition $\pt(V)$)}
			\\&=
			\left\{\calF\in\pt(X):f^{-1}(V)\in\calF \right\}&&\text{(definition $\pt(f)$)}
			\\&=
			\pt\left(f^{-1}(V)\right).
			&&\qedhere
		\end{align*}
\end{proof}

\begin{corollary}
	There is a functor:
	\begin{align*}
		\pt: \Loc&\longrightarrow \Top;\\
		X & \longmapsto \pt(X);\\
		f & \longmapsto \pt(f).
	\end{align*}
\end{corollary}
\begin{proof}
	This assignment is well-defined on objects and arrows, so we are left to show functoriality. Since as frame maps $\id_X^{-1} = \id_{\Opens X}$, it is easy to see that $\pt(\id_X)=\id_{\pt(X)}$. For composition, take $f\colon X\to Y$ and $h\colon Y\to Z$. Then for $\calF\in\pt(X)$:
		\begin{align*}
			\pt(h)\circ\pt(f)(\calF)
			&=
			\left\{W\in\Opens Z: h^{-1}\in \pt(f)(\calF)\right\}
			\\&=
			\left\{W\in \Opens Z: f^{-1}(h^{-1}(W))\in \calF \right\}
			\\&=
			\pt(h\circ f)(\calF). \qedhere
		\end{align*}
\end{proof}

With the functors established, we can now work towards the adjunction. To do this, we shall build the unit $\eta$ and counit $\epsilon$ explicitly.

\begin{lemma}\label{lemma:unit}
	For any topological space $S$, there is a continuous function
		\[
		\eta_S\colon S\longrightarrow \pt(\loc(S));
		\qquad
		x\longmapsto\calF_x:=\left\{U\in\Opens S: x\in U\right\}.
		\]
	This constitutes a natural transformation $\eta\colon \id_{\Top}\to \pt\circ \loc$.
\end{lemma}
\begin{proof}
	That $\eta_S$ is well defined follows by \cref{lemma:neighbourhood filter is completely prime}. For continuity, recall that (\cref{definition:space of points}) a generic open in $\pt(\loc(S))$ is of the form $\pt(U)$, and so it follows immediately that its preimage is open:
		\[
			\eta_S^{-1}(\pt(U))= \left\{x\in S: \calF_x \in \pt(U)\right\}
			=
			\left\{x\in S: x\in U\right\} = U.
		\]
	Naturality amounts to the fact that any continuous function $g\colon S\to T$ induces a commutative square
		\[
			\begin{tikzcd}[ampersand replacement=\&]
				S \& {\pt(\loc(S))} \\
				T \& {\pt(\loc(T)).}
				\arrow["{\pt(\loc(g))}", from=1-2, to=2-2]
				\arrow["g"', from=1-1, to=2-1]
				\arrow["{\eta_S}", from=1-1, to=1-2]
				\arrow["{\eta_T}"', from=2-1, to=2-2]
			\end{tikzcd}
		\]
	To see this, calculate for $x\in S$:
		\begin{align*}
			\pt(\loc(g))\circ \eta_S(x) 
			&=
			\pt(\loc(g))(\calF_x) &&\text{(definition $\eta_S$)}
			\\&=
			\left\{V\in\Opens T: g^{-1}(V)\in\calF_x\right\} &&\text{(\cref{definition:pt(f)})}
			\\&=
			\left\{V\in\Opens T: x\in g^{-1}(V)\right\}&&\text{(definition $\calF_x$)}
			\\&=
			\left\{V\in\Opens T: g(x)\in V\right\}
			\\&= \calF_{g(x)} =
			\eta_T\circ g(x). &&\text{(definition $\eta_T$)}\qedhere
		\end{align*}
\end{proof}

\begin{remarknumbered}\label{remark:completely prime filters vs elements}
	The function $\eta_S$ takes an element $x\in S$ of the space, which we traditionally think of as a ``point,'' and constructs a \emph{localic} point $\calF_x$, which is how the topology of $S$ sees $x$. How the topology sees the points (i.e.~elements) is not always accurate. Two things can happen: either the topology is too coarse to distinguish different elements of $S$. For instance, in a codiscrete space, where $\Opens S = \{\varnothing, S\}$, there exists only one completely prime filter: $\calF = \{S\}$, while there can of course be many distinct elements.
	
	On the other hand, the topology can ``hallucinate'' points of $S$. A clear example of this is the space $\mathbb{R}$, with topology generated by the intervals of the form $(a,\infty)$, for $a\in\mathbb{R}$. Distinct elements generate distinct filters, but there exists a completely prime filter $\calF=\Opens \mathbb{R}\setminus\{\varnothing\}$ that does not come from a point. It zooms in on an imaginary point at infinity.
	
	A space where the topology and the underlying set of elements see the same points is called \emph{sober} (\cref{definition:sober spaces}). We discuss these further in \cref{section:fixed points of old adjunction}.	
\end{remarknumbered}

\begin{lemma}\label{lemma:counit}
	For any locale $X$, there is a map of locales
		\[
		\epsilon_X:\loc(\pt(X))\longrightarrow X
		\quad\text{defined by}\quad
		\epsilon_X^{-1}\colon \Opens X\longrightarrow \Opens\pt(X);
		\quad
		U\longmapsto \pt(U).
		\]
	This constitutes a natural transformation $\epsilon\colon \loc\circ\pt\to \id_{\Loc}$.
\end{lemma}
\begin{proof}
	It follows immediately from \cref{lemma:properties of pt(U)} that $\epsilon_X^{-1}$ is a map of frames. For naturality, we need to show that if $f\colon X\to Y$ is a map of locales, we get a commuting square:
		\[
			\begin{tikzcd}[ampersand replacement=\&]
				{\loc(\pt(X))} \& X \\
				{\loc(\pt(Y))} \& {Y,}
				\arrow["f", from=1-2, to=2-2]
				\arrow["{\loc(\pt(f))}"', from=1-1, to=2-1]
				\arrow["{\epsilon_X}", from=1-1, to=1-2]
				\arrow["{\epsilon_Y}"', from=2-1, to=2-2]
			\end{tikzcd}
		\qquad\text{meaning:}\qquad
			\begin{tikzcd}[ampersand replacement=\&]
				{\Opens\pt(X)} \& {\Opens X} \\
				{\Opens\pt(Y)} \& {\Opens Y.}
				\arrow["{f^{-1}}"', from=2-2, to=1-2]
				\arrow["{\epsilon_X^{-1}}"', from=1-2, to=1-1]
				\arrow["{\epsilon_Y^{-1}}", from=2-2, to=2-1]
				\arrow["{\pt(f)^{-1}}", from=2-1, to=1-1]
			\end{tikzcd}
		\]
	To see this, take $V\in \Opens Y$ and calculate:
		\begin{align*}
			\pt(f)^{-1}\circ \epsilon_Y^{-1}(V)
			&=
			\pt(f)^{-1}(\pt(V)) &&\text{(definition $\epsilon_Y$)}
			\\&=
			\pt\left(f^{-1}(V)\right) &&\text{(\cref{lemma:pt(f) continuity})}
			\\&=
			\epsilon_X^{-1}\circ f^{-1}(V). &&\text{(definition $\epsilon_X$)}\qedhere
		\end{align*}
\end{proof}

Putting it all together gives the following celebrated theorem.

\begin{theorem}\label{theorem:adjunction top and loc}
	There is an adjunction:
		\[
			\begin{tikzcd}[cramped,column sep=large]
				\Top & {\Loc.}
				\arrow[""{name=0, anchor=center, inner sep=0}, "\loc", shift left=2, from=1-1, to=1-2]
				\arrow[""{name=1, anchor=center, inner sep=0}, "\pt", shift left=2, from=1-2, to=1-1]
				\arrow["\dashv"{anchor=center, rotate=-90}, draw=none, from=0, to=1]
			\end{tikzcd}
		\]
\end{theorem}
\begin{proof}
	The following is an expanded proof of that found in \cite[Theorem~IX.3.1]{maclane1994SheavesGeometryLogic}. We prove that there is a natural bijection:
		\begin{align*}
			\Loc(\loc(S),X)  &\xrightarrow{~ \sim ~}  \Top(S,\pt(X))\\
			f&\xmapsto{\hspace{3ex}} \pt(f)\circ \eta_S\\ \epsilon_X\circ \loc(g)&{~\reflectbox{\ensuremath{\xmapsto{\hspace{3ex}}}}} ~g	.				
		\end{align*}
	The first function is well defined by \cref{lemma:pt(f) continuity,lemma:unit}, and the second by \cref{lemma:loc(g) frame map,lemma:counit}. We are therefore just left to show that they are mutually inverse. This is done through the following calculations, which mostly involve unpacking definitions.
	
	First, we claim that $f = \epsilon_X\circ \loc(\pt(f))\circ \loc(\eta_S)$. To show this, first use naturality of the counit (\cref{lemma:counit}) to see $\epsilon_X\circ \loc(\pt(f)) = f\circ \epsilon_{\loc(S)}$, and then calculate for the underlying frame maps at $U\in \Opens X$:
		\begin{align*}
			\eta_S^{-1}\circ \epsilon_{\loc(S)}^{-1}\circ f^{-1}(U)
			&=
			\eta_S^{-1}\left(\pt\left(f^{-1}(U)\right)\right)&&\text{(definition $\epsilon_{\loc(S)}$)}
			\\&=
			f^{-1}(U), &&\text{(proof of \cref{lemma:unit})}
		\end{align*}
	giving the desired equality.
	
	For the other direction, we need to show that $g = \pt(\epsilon_X)\circ \pt(\loc(g))\circ \eta_S$. First use the naturality of the unit (\cref{lemma:unit}) to see $\pt(\loc(g))\circ \eta_S = \eta_{\pt(X)} \circ g$. Next, we describe the function $\pt(\epsilon_X)\colon \pt(\loc(\pt(X)))\to \pt(X)$. Namely, for any completely prime filter $\calF\subseteq \Opens \pt(X)$ we get:
		\begin{align*}
			\pt(\epsilon_X)(\calF)
			&= \left\{U\in \Opens X: \epsilon_X^{-1}(U)\in\calF \right\} &&\text{(\cref{definition:pt(f)})}
			\\&=
			\left\{U\in \Opens X: \pt(U)\in \calF\right\}. &&\text{(definition $\epsilon_X$)}
		\end{align*}
	Combined, if $x\in S$ then $g(x)\subseteq \Opens X$ is a completely prime filter, and we find:
		\begin{align*}
			\pt(\epsilon_X)\circ \eta_{\pt(X)}\circ g(x)
			&=
			\pt(\epsilon_X)(\calF_{g(x)}) &&\text{(definition $\eta_{\pt(X)}$)}
			\\&=
			\left\{U\in \Opens X: \pt(U)\in \calF_{g(x)}\right\} &&\text{(above)}
			\\&=
			\{U\in \Opens X: U\in g(x)\} = g(x). &&\text{(definition $\pt(U)$)}
		\end{align*}
	This establishes that the functions are mutually inverse, and hence form a bijection. Finally, that this is bijection is natural follows since composition and $\eta$ and~$\epsilon$ are natural. We therefore obtain the desired adjunction.
\end{proof}

\subsection{Sober spaces and spatial locales}
\label{section:fixed points of old adjunction}
In this section we discuss the \emph{fixed points} of the adjunction $\loc\dashv \pt$. These are the spaces $S$ and locales $X$ whose unit map $\eta_S$ and respectively counit map $\epsilon_X$ are isomorphisms. This section recalls the well-known characterisation of when this is the case.

We start with spaces. To characterise when $\eta_S$ is a homeomorphism, we need to study the relation between the elements $x\in S$ and the points of $\loc(S)$. The unit map itself associates the completely prime filter $\calF_x= \eta_S(x)$ to any $x\in S$. But in the development of this theory it is customary to introduce an alternative way of describing the element $x$ localically. While the singleton $\{x\}$ is not typically an element of $\Opens S$, there is another way to obtain an open subset: the closure of the singleton defines a closed subset $\overline{\{x\}}$, and so its complement satisfies
	\[
		P_x:= S\setminus \overline{\{x\}} \in \Opens S.
	\]
It can now be seen that $P_x$ satisfies the conditions in the following lattice-theoretic definition.
\begin{definition}\label{definition:prime element}
	Let $L$ be a bounded distributive lattice. A \emph{prime element} is an element $p\in L$ such that the following two conditions hold:
		\begin{enumerate}[label = (P\arabic*)]
		\item $p\neq \top$;
		\item if $x\wedge y\sqleq p$ then $x\sqleq p$ or $y\sqleq p$.
	\end{enumerate}
\end{definition}
\begin{remark}
	In a distributive lattice, prime elements can equivalently be defined by using equality in~(P2).
\end{remark}

It is readily verified that the open $P_x\in \Opens S$ defines a prime element in this sense. In \emph{sober} spaces there is a bijective correspondence between primes of $\Opens S$ and elements of $S$.

\begin{definition}\label{definition:sober spaces}
	A topological space $S$ is called \emph{sober} if for every prime element $P\in \Opens S$ there exists a unique element $x\in S$ such that $P = S\setminus\overline{\{x\}}$. Here the prime elements of $\Opens S$ are precisely the open subsets $P$ satisfying:
		\begin{enumerate}[label = (P\arabic*)]
			\item $P\neq S$;
			\item if $U\cap V\subseteq P$ then $U\subseteq P$ or $V\subseteq P$.
		\end{enumerate}
\end{definition}

The next few results set out to show that it is precisely the sober spaces that form fixed points for the adjunction $\loc\dashv \pt$. To recast sobriety in more familiar terms, we recall the following definition.

\begin{definition}\label{definition:T0 and enough points}
	A topological space $S$ is called a \emph{$T_0$-space} if:
	\[
	x\neq y
	\qquad\text{implies}\qquad
	\exists U\in \Opens X\colon x\in U, y\notin U
	\text{~or~}
	\exists V\in \Opens Y\colon y\in V, x\notin V.
	\]
	Expressed differently, a space is $T_0$ if and only if the unit $\eta_S$ is injective.
	
	Similarly, a space $S$ is said to \emph{have enough points} (or \emph{weakly sober}) if $\eta_S$ is surjective. This means that for every completely prime filter $\calF\subseteq \Opens S$ there exists $x\in S$ such that $\calF= \calF_x$.
\end{definition}

\begin{remark}
	Continuing on \cref{remark:completely prime filters vs elements}, we see that a codiscrete space $S$ will rarely be~$T_0$. On the other hand, the space $\mathbb{R}$ with topology generated by intervals $(a,\infty)$ is a $T_0$-space, but does not have enough points. The set $\calF = \{(a,\infty):a\in\mathbb{R}\}$ defines a completely prime filter, but is not in the image of $\eta$.
\end{remark}

The following shows that prime elements of $\Opens X$ provide an equivalent way to think about localic points.

\begin{lemma}\label{lemma:prime elements are cpf}
	For any locale $X$, there is a bijective correspondence between completely prime filters $\calF\subseteq \Opens X$ and prime elements $P\in\Opens X$, given explicitly by
		\[
			\calF\longmapsto P_\calF := \bigvee\left\{U\in\Opens X: U\notin \calF \right\}
			\quad\text{and}\quad
			P\longmapsto \calF_P := \left\{ U\in\Opens X: U\not\sqleq P\right\}.
		\]
\end{lemma}
\begin{proof}
	The (classical) proof is a relatively elementary translation between the axioms of prime elements and completely prime filters. First we prove that $P_\calF$ is prime and $\calF_P$ is a completely prime filter, and then prove the assignments are mutually inverse.
	
	Start by observing that $P_\calF\notin \calF$: otherwise~(F4) gives a $U\in\calF$ such that $U\notin\cal F$, a contradiction. Since $X\in\calF$ we thus get $P_\calF\neq X$. Using upwards closure~(F2) it follows more generally that:
		\[
			U\sqleq P_\calF\qquad\text{if and only if}\qquad
			U\notin \calF.
		\]
	Then if $U\wedge V\sqleq P_\calF$ implies $U\wedge V\notin \calF$, which by contrapositive of~(F1) gives $U\notin \calF$ or $V\notin \calF$. Equivalently: $U\sqleq P_\calF$ or $V\sqleq P_\calF$. Hence $P_\calF$ is a well-defined prime element of $\Opens X$.
	
	For the converse function, we check the axioms:
		\begin{enumerate}[label = (F\arabic*)]\setcounter{enumi}{-1}
			\item $X\in \calF_P$ since $P\neq X$;
			\item the contrapositive of primeness says $U,V\not\sqleq P$ implies $U\wedge V\not\sqleq P$, meaning here: $U,V\in\calF_P$ implies $U\wedge V\in \calF_P$;
			\item if $\calF\ni U\sqleq V$ then clearly $V\not\sqleq P$, so $V\in\calF_P$;
			\item $\varnothing\notin\calF_P$ since $\varnothing\sqleq P$;
			\item the contrapositive of~(F4) says if $U_i\notin \calF_P$ for all $i$, then $\bigvee U_i\notin\calF_P$. This holds since this means $U_i\sqleq P$ implies $\bigvee U_i\sqleq P$.
		\end{enumerate}
	Thus $\calF_P$ is a well-defined completely prime filter in $\Opens X$.
	
	We are left to show these functions are mutually inverse, but this follows from the following calculations:
		\[
			P_{\calF_P} = \bigvee\{ U\in\Opens X: U\notin \calF_P\} = \bigvee\{U\in\Opens X:U\sqleq P\} = P,
		\]
	and 
		\[
			\calF_{P_\calF} = \{U\in\Opens X: U\not\sqleq P_\calF\} = \{U\in\Opens X: U\in\calF\} = \calF.\qedhere
		\]
\end{proof}

\begin{remarknumbered}\label{remark:pt'}
	We could therefore have equivalently (at least classically) defined the space $\pt(X)$ of points in terms of prime elements $P\in\Opens X$. Since we will need this presentation of points in \cref{section:causal boundaries}, we give it a name: $\pt'(X)$. We equip this set with the unique topology so that the above bijection becomes a homeomorphism. Explicilty, open sets in $\pt'(X)$ are of the form
		\[
			\pt'(U) = \{P\in\pt'(X): U\not\sqleq P\}.
		\]
	For arrows $f\colon X\to Y$, the corresponding map $\pt'(f)$ is simply defined as the direct image map $f_\ast$, which preserves prime elements by \cite[\S II.3.4]{picado2012FramesLocalesTopology}.
\end{remarknumbered}

\begin{lemma}\label{lemma:sober iff unit bijective}
	A space $S$ is sober if and only if it is $T_0$ and has enough points.
\end{lemma}
\begin{proof}
	By definition, $S$ is sober if and only if there exists a bijective correspondence $x\mapsto S\setminus\overline{\{x\}}$ between the elements of $S$ and the prime elements of $\Opens S$. Composing this with the bijective correspondence in \cref{lemma:prime elements are cpf} gives
		\[
			x\longmapsto S\setminus\overline{\{x\}} \longmapsto \calF_{S\setminus\overline{\{x\}}}
			=
			\left\{U\in\Opens X: U\not\subseteq S\setminus\overline{\{x\}}\right\}
			=
			\left\{ U\in \Opens X: x\in U\right\},
		\]
	which is of course precisely the unit map $\eta_S$. Therefore $S$ is sober if and only if $\eta_S$ is bijective, which in turn is equivalent to the stated condition.
\end{proof}

\begin{lemma}\label{lemma:locales are sober}
	For any locale $X$, the space $\pt(X)$ is sober.
\end{lemma}
\begin{proof}
	Suppose that $\pt(P)\in \Opens\pt(X)$ is a prime element. We claim that
		\[
			\calF:= \{ U\in \Opens X: \pt(U)\not\subseteq \pt(P)\}
		\]
	defines a completely prime filter in $\Opens X$ such that $\pt(P)=\pt(X)\setminus\overline{\{\calF\}}$. We check the axioms of a completely prime filter. Throughout, we use the properties of \cref{lemma:properties of pt(U)}.
		\begin{enumerate}[label = (F\arabic*)]\setcounter{enumi}{-1}
			\item $X\in \calF$ iff $\pt(X)\not\subseteq \pt(P)$, which holds by (P1);
			\item if $U,V\in \calF$ then $\pt(U)\cap \pt(V)\not\subseteq \pt(P)$, otherwise (P2) would imply $\pt(U)$ or $\pt(V)\subseteq \pt(P)$, contradicting $U,V\in\calF$;
			\item if $\calF\ni U\sqleq V$, clearly $\pt(U)\not\subseteq\pt(P)$ implies $\pt(V) \not\subseteq\pt(P)$, so $V\in\calF$;
			\item $\varnothing\notin\calF$ is equivalent to $\varnothing = \pt(\varnothing)\subseteq\pt(P)$, which trivially holds;
			\item take $\bigvee U_i\in \calF$, meaning $\bigcup \pt(U_i)\not\subseteq\pt(P)$. Thus there exists $\calG\in\bigcup\pt(U_i)$ such that $\calG\notin\pt(P)$. In turn, there exists some index $i$ with $\calG\in\pt(U_i)$, and so $\pt(U_i)\not\subseteq \pt(P)$, or: $U_i\in\calF$.
		\end{enumerate}
	Hence $\calF$ defines a point of $X$. We are left to show $\pt(P)$ is the complement of this point:
		\begin{align*}
			\calG\in \pt(X)\setminus\overline{\{\calF\}}
			&\quad\text{iff}\quad
			\calG\notin \overline{\{\calF\}}
			\\&\quad\text{iff}\quad
			\exists\pt(U)\in\Opens \pt(X)\colon \calG\in\pt(U) \text{~and~} \calF\notin\pt(U)
			\\&\quad\text{iff}\quad
			\exists U\in \calG\colon U\notin\calF
			\\&\quad\text{iff}\quad
			\exists U\in \calG\colon \pt(U)\subseteq \pt(P)
			\\&\quad\text{iff}\quad
			\calG\in\pt(P).\qedhere
		\end{align*}
\end{proof}

\begin{proposition}\label{proposition:hausdorff is sober}
	Any Hausdorff space is sober.
\end{proposition}
\begin{proof}
	This is the proof of \cite[Theorem~IX.2.3]{maclane1994SheavesGeometryLogic}. Hausdorff spaces $S$ are~$T_0$, so by \cref{lemma:sober iff unit bijective} we are left to show that they have enough points. With \cref{lemma:prime elements are cpf} we thus need to show if $P\in\Opens S$ is prime, then $P=S\setminus\{x\}$ for some unique $x\in S$. In turn, it suffices to prove that $S\setminus P$ is a singleton. Suppose for the sake of contradiction that it instead contains two distinct points $x\neq y$. By Hausdorffness we can find disjoint $U,V\in \Opens S$ such that $x\in U$ and $y\in V$, while $S\setminus P = \left((S\setminus P)\cap U\right) \cup \left((S\setminus P)\cap V\right)$, and by taking complements one sees that this contradicts the primeness of $P$.
\end{proof}

\begin{proposition}\label{proposition:characterisation of sober spaces}
	For a space $S$, the following are equivalent:
		\begin{enumerate}[label = (\alph*)]
			\item $S$ is sober;
			\item the unit $\eta_S$ is a homeomorphism;
			\item $S\cong\pt(X)$ for some locale $X$.
		\end{enumerate}
\end{proposition}
\begin{proof}
	We know from (the proof of) \cref{lemma:sober iff unit bijective} that if $S$ is sober, then $\eta_S$ is bijective. In that case, we find that $\eta_S(U) = \pt(U)$, so that $\eta_S$ is not only continuous but open. It is therefore a homeomorphism, so (a) implies (b). Trivially (b) implies (c), and the latter implies (a) by \cref{lemma:locales are sober}.
\end{proof}

Next, we characterise the fixed points on the localic side of the adjunction.

\begin{definition}\label{definition:spatial locales}
	A locale $X$ is called \emph{spatial} if $\pt(U) = \pt(V)$ implies $U=V$. It is also said that $X$ \emph{has enough points}.
\end{definition}

\begin{remark}
	Note that the phrase ``having enough points'' means something different for locales than for spaces (\cref{definition:T0 and enough points}). We therefore prefer to use ``spatial,'' though both are common in the literature.
\end{remark}

\begin{lemma}\label{lemma:spaces are spatial}
	For any space $S$, the locale $\loc(S)$ is spatial.
\end{lemma}
\begin{proof}
	Contrary to the proof of \cref{lemma:locales are sober}, this one is easy.
	For a contrapositive argument, suppose $U\neq V$, so that without loss of generality there exists $x\in U$ but $x\notin V$. Then $\calF_x\in\pt(U)$ but $\calF_x\notin\pt(V)$, showing $\pt(U)\neq\pt(V)$.
\end{proof}

\begin{proposition}\label{proposition:characterisation of spatial locales}
	The following are equivalent, for any locale $X$:
		\begin{enumerate}[label = (\alph*)]
			\item $X$ is spatial;
			\item the counit $\epsilon_X$ is an isomorphism of locales;
			\item $X\cong \loc(S)$ for some space $S$.
		\end{enumerate}
\end{proposition}
\begin{proof}
	That (b) implies (c) is trivial, and (c) implies (a) by \cref{lemma:spaces are spatial}. Lastly, if $X$ is spatial then we can define a map of frames $\delta^{-1}\colon \Opens \pt(X) \to \Opens X$ by $\pt(U)\mapsto U$, and this clearly defines the localic inverse of the counit $\epsilon_X$.
\end{proof}

\begin{theorem}
	The adjunction $\loc\dashv\pt$ restricts to an equivalence between the full subcategories of $\Top$ and $\Loc$ of sober spaces and spatial~locales:
	\[
	\left\{\begin{array}{c}
		\text{\normalfont sober} \\
		\text{\normalfont spaces}
	\end{array}\right\}
	\simeq
	\left\{\begin{array}{c}
		\text{\normalfont spatial} \\
		\text{\normalfont locales}
	\end{array}\right\}
	.
	\]
\end{theorem}

\begin{remark}[Next steps]
	Now that we have fully exposed the classic adjunction $\loc\dashv\pt$ and calculated its fixed points, we finally move towards the extended adjunction between $\OrdTop$ and $\OrdLoc$. We reassure the reader that the additional work is not egregious. It comes down to the following. First: extending the functors $\loc$ and $\pt$ to the setting of $\OrdTop$ and $\OrdLoc$. Once successful, it is a simple matter of showing that the natural bijection from the proof of \cref{theorem:adjunction top and loc} restricts to a natural bijection between the hom-sets of $\OrdTop$ and $\OrdLoc$. This, in turn, can be accomplished by demonstrating that $\eta_S$ and $\epsilon_X$ are monotone in the sense of \cref{definition:monotone function between preorders,definition:monotone map of ordered locales}, respectively. As previously mentioned, unfortunately, these steps cannot be accomplished in full generality. Precisely what these restrictions are and how they arise will be explained in the next few sections.
\end{remark}

\section{Ordered locales from ordered spaces}
\label{section:functor ordtop to ordloc}
In this section we describe how the basic idea described in \cref{example:ordered spaces as ordered locales} can be used to construct a functor from $\OrdTop$ to $\OrdLoc$. Given a space, we obtain a locale $\loc(S)$ that is defined by the frame $\Opens \loc(S):= \Opens S$. We have seen in \cref{proposition:space with EM order is ordered locale} that the Egli-Milner order $\Leq$ on $\Opens S$ turns $\loc(S)$ into an ordered locale. This will be the object component of the functor.

Since we require the new functor to commute with the forgetful functor, we have no choice in how $\loc$ behaves on morphisms. A continuous monotone map $g$ is necessarily sent to the map of locales $\loc(g)$. The following result shows that this is well-defined precisely if we restrict to the subcategory of ordered spaces with open cones.

\begin{proposition}\label{proposition:loc well defined iff OC}
	The following are equivalent for an ordered space $T$:
	\begin{enumerate}[label = (\alph*)]
		\item $T$ has open cones;
		\item $\loc(g)$ is monotone for every continuous monotone $g\colon S\to T$;
		\item $\loc(g)$ is monotone for every monotone $g\colon \{0<1\}\to T$.
	\end{enumerate}
\end{proposition}
\begin{proof}
	Note that by \cref{corollary:localic cones in space with open cones} the localic cones of $T$ are just the ordinary cones, calculated pointwise, so the monotonicity of $\loc(g)$ is equivalent to that of $g$ itself by \cref{lemma:monotonicity in terms of cones}. Hence~(a) certainly implies~(b). 
	
	Trivially~(b) implies~(c).
	
	Lastly, we show that (c) implies (a). We need to show that if $U\in \Opens T$ then $\up U \subseteq (\up U)^\circ$. Let $y\in \up U$, so there exists $x\in U$ with $x\leq y$. Using this, define a function $g\colon \{0<1\}\to T$ by $g(0)=x$ and $g(1)=y$. We equip the domain of $g$ with the discrete topology, so $g$ is certainly monotone and continuous. By hypothesis, so is $\loc(g)$. Moreover, since $\{0<1\}$ has the discrete topology, it has open cones, so by \cref{corollary:localic cones in space with open cones} we find $\{0<1\}=\Up \{0\}$, and since $0\in g^{-1}(U)$ we thus get $\{0<1\}\subseteq \Up g^{-1}(U) \subseteq g^{-1}((\up U)^\circ)$.
	In particular $y=g(1)\in (\up U)^\circ$, as desired. A dual argument holds for down cones, so we conclude $T$ has open cones.
\end{proof}

With this, the functor $\loc\colon \Top\to \Loc$ therefore extends as follows:

\begin{corollary}\label{corollary:functor ordtopOC to ordloc}
	There is a functor:
	\begin{align*}
		\loc: \OrdTop_\OC &\longrightarrow \OrdLoc;\\
		(S,\leq) & \longmapsto (\loc(S),\Leq);\\
		g & \longmapsto \loc(g).
	\end{align*}
\end{corollary}

\begin{remark}
	\cref{proposition:loc well defined iff OC} says that, with the definition of an ordered locale as an axiomatisation of the Egli-Milner order, we cannot hope to obtain a more general adjunction between $\OrdTop$ and $\OrdLoc$, but have to restrict to spaces with open cones. From a mathematical point of view this can be seen as a drawback, given the existence of spaces without open cones. Nonetheless, the theory is still sufficiently general for us to consider physical applications to spacetimes.
\end{remark}

\section{Ordered spaces from ordered locales}\label{section:functor ordloc to ordtop}
In this section we extend the functor $\pt$ to a functor $\OrdLoc\to \OrdTop$. Starting with an ordered locale $(X,\Leq)$, we therefore need to define a preorder $\leq$ on the space of points $\pt(X)$. This is done as follows:

\begin{definition}\label{definition:order on points}
	For a locale $X$ and a preorder $\Leq$ on its frame of opens (not necessarily satisfying axiom~\eqref{axiom:V}), we define the following relation on its space of points. For $\calF,\calG\in\pt(X)$ we set:
		\[
			\calF\leq\calG
			\quad\text{if and only if}\quad
			\forall U\in\calF~\exists V\in\calG\colon U\Leq V
			\text{~and~}
			\forall V\in\calG~\exists U\in\calF\colon U\Leq V.
		\]
\end{definition}

\begin{proposition}
	The pair $(\pt(X),\leq)$ is an ordered space.
\end{proposition}
\begin{proof}
	We need to verify that $\leq$ is a preorder. For $\calF\leq\calF$, note that if $U\in\calF$ then $U\Leq U\in\cal F$, so we get reflexivity. For transitivity, take $\calF\leq\calG\leq \calH$. Then if $U\in \calF$, we can find $V\in\calG$ such that $U\Leq V$. In turn, there exists $W\in\calH$ such that $V\Leq W$, so by transitivity of $\Leq$ we get $U\Leq W$. A dual argument holds for the second condition, proving that $\leq$ is transitive.
\end{proof}

In a full-fledged ordered locale, the induced order on the points has a cleaner characterisation in terms of the localic cones.

\begin{lemma}\label{lemma:order on points in terms of localic cones}
	For points $\calF$ and $\calG$ in an ordered locale:
		\[
			\calF\leq \calG
			\qquad\text{if and only if}\qquad
			\forall U\in \calF\colon \Up U\in\calG \text{~and~}
			\forall V\in\calG\colon \Down V\in\calF.
		\]
\end{lemma}
\begin{proof}
	If $\calF\leq \calG$ then for any $U\in\calF$ we can find some open region $V\in\calG$ such that $U\leq V$, which with \cref{lemma:properties of localic cones}(a) gives $V\sqleq \Up U$. Since $\calG$ is upwards closed, this implies $\Up U\in\calG$. Dually we obtain that $\Down V\in\calF$ whenever $V\in\calG$. The converse implication follows from \cref{lemma:properties of localic cones}(c).
\end{proof}

\begin{remark}
	One is tempted to conclude from this that $\calF\leq \calG$ if and only if ${\bigwedge_{U\in\calF} \Up U \in \calG}$ and $\bigwedge_{V\in \calG}\Down V\in \calF$. While this provides a good visual intuition, see \cref{figure:order on cpf}, it is not technically correct. An easy counterexample is the vertically-ordered $\mathbb{R}^2$ from \cref{example:vertical R2}, in which for every point $x$ we get $\bigwedge_{U\ni x}\Up U = (\up x)^\circ = \varnothing$.
\end{remark}

\begin{figure}[b]\centering
	\input{tikz/inkscape-export/fig-26-big}
	\caption{Illustration of $\forall U\in \calF_x: y\in\Up U$.}
	\label{figure:order on cpf}
\end{figure}

\begin{remark}
	For completely prime filters induced by elements of a space, we get the following specific characterisation of this order. Let $x,y\in S$ be elements in an ordered space. Then in $\pt(\loc(S))$:
		\[
			\calF_x\leq \calF_y
			\qquad\text{if and only if}\qquad
			\forall U\in \calF_x\colon (\up U)^\circ\in\calF_y
			\text{~and~}
			\forall V\in \calF_y\colon (\down V)^\circ\in\calF_x.
		\]
	Compare this to the characterisation of the open cone condition in \cref{proposition:open cones iff interiors of open neighbourhoods}.
\end{remark}


\begin{remark}
	Recall from \cref{remark:pt'} that the points of a locale $X$ can equivalently be described as prime elements $P\in\Opens X$. In this language, the order $\calF\leq \calG$ on points from \cref{definition:order on points} translates to the following condition:
		\[
			P\leq Q
			\qquad\text{if and only if}\qquad
			\begin{array}{l}
				\text{$\Up U\sqleq Q$ implies $\Up U\sqleq P$, and}\\
				\text{$\Down V\sqleq P$ implies $\Down V\sqleq Q$.}
			\end{array}
		\]
	Intuitively, if $P=P_x$ and $Q=P_y$, this says $y\notin \Up U$ implies $x\notin \Up U$, and dually $x\notin \Down V$ implies $y\notin \Down V$. This contrapositively corresponds to $x\leq y$.
\end{remark}

The following shows that the point-wise cones $\up$ and $\down$ in $\pt(X)$ are bounded by the localic cones $\Up$ and $\Down$.

\begin{lemma}\label{lemma:cones of points contained in localic cones}
	Let $(X,\Leq)$ be an ordered locale, and pick $U\in\Opens X$. Then in $\pt(X)$:
		\[
			\up \pt(U)\subseteq \pt(\Up U)
			\qquad\text{and}\qquad
			\down \pt(U)\subseteq \pt(\Down U).
		\]
\end{lemma}
\begin{proof}
	For the first inclusion, take $\calG\in\up \pt(U)$. Then there exists a point $\calF\in \pt(U)$ such that $\calF\leq \calG$. In particular, since $U\in\calF$, this implies $\Up U\in\calG$ by \cref{lemma:order on points in terms of localic cones}, which translates to $\calG\in\pt(\Up U)$. The proof for past cones is dual.
\end{proof}

\begin{lemma}\label{lemma:pt(f) monotone}
	If $f$ is a monotone map of ordered locales, then $\pt(f)$ is monotone.
\end{lemma}
\begin{proof}
	Assume that $\calF \leq\calG$ for some points in the domain $X$. Let $V\in \pt(f)(\calF)$, meaning $f^{-1}(V)\in\calF$. Since $\calF\leq \calG$ this implies $\Up f^{-1}(V)\in\calG$. Using that $f$ is monotone and that $\calG$ is upwards closed, this gives $f^{-1}(\Up V)\in \calG$. In other words: $\Up V\in \pt(f)(\calG)$. A dual argument shows that if $W\in \pt(f)(\calG)$ then $\Down W\in \pt(f)(\calF)$, so we conclude that $\pt(f)(\calF)\leq \pt(f)(\calG)$ by \cref{lemma:order on points in terms of localic cones}. 
\end{proof}

\begin{corollary}\label{corollary:functor ordloc to ordtop}
	There is a functor:
		\begin{align*}
			\pt: \OrdLoc &\longrightarrow \OrdTop;\\
			(X,\Leq) & \longmapsto (\pt(X),\leq);\\
			f & \longmapsto \pt(f).
		\end{align*}
\end{corollary}

\begin{remark}
	Note that, unlike in the case of going from spaces to ordered locales (\cref{corollary:functor ordtopOC to ordloc}), we did not have to put any restrictions on the type of ordered locales we are considering. 
	
	However, in order to obtain an adjunction, the new functor $\pt$ needs to land in the subcategory of ordered spaces with open cones. Unfortunately, at present, we do not know if $\pt(X)$ always has open cones. To ensure this, we discuss a natural condition that guarantees the functor $\pt$ lands in the subcategory $\OrdTop_\OC$.
\end{remark}

\begin{definition}
	For ordered locales $(X,\Leq)$ we introduce the following axiom:
	\[\tag{$\bullet$}\label{axiom:P}
	U\Leq V \text{~in~} X 
	\qquad\text{implies}\qquad
	\pt(U)\Leq \pt(V)\text{~in~}\loc(\pt(X)).
	\]
\end{definition}

One interpretation of~\eqref{axiom:P} is that points must ``flow'' into the future and past: while $U\Leq V$ it is not possible for there to exists a point in $\pt(U)$ but that $\up\pt(U)$ be empty. Specifically, if~\eqref{axiom:P} holds then $U\Leq V$ implies $\pt(U)=\varnothing$ if and only if $\pt(V)=\varnothing$. This is made more precise in the following lemma.

\begin{lemma}\label{lemma:characterisation of axiom P}
	An ordered locale satisfies~\eqref{axiom:P} if and only if
	\[
	\up\pt(U) = \pt(\Up U)\qquad\text{and}\qquad \down\pt(U)=\pt(\Down U).
	\]
\end{lemma}
\begin{proof}
	For use throughout the the proof, we recall that by \cref{lemma:cones of points contained in localic cones} in an arbitrary ordered locale we get the inclusions $\up \pt(U)\subseteq \pt(\Up U)$ and ${\down\pt(U)\subseteq \pt(\Down U)}$. So in the following proof we only have to be concerned with the converse inclusions.
	
	Suppose first that~\eqref{axiom:P} holds. By \cref{lemma:properties of localic cones}(c) we get $U\Leq\Up U$, and hence $\pt(U)\Leq \pt(\Up U)$, which unpacks to give $\pt(\Up U)\subseteq \up\pt(U)$, as required. Equality follows by the remark in the first paragraph. The argument for downsets is analogous. 
	
	Now suppose the stated equations hold. Using \cref{lemma:properties of localic cones}(a), if $U\Leq V$ we get $U\sqleq \Down V$ and $V\sqleq \Up U$. These inclusions are respected by points, so we get $\pt(U)\subseteq \pt(\Down V)=\down \pt(V)$ and $\pt(V)\subseteq \pt(\Up U)=\up \pt(U)$, which is precisely what it means for $\pt(U)\Leq \pt(V)$ to hold. This shows~\eqref{axiom:P} holds.
\end{proof}

From these equations we can immediately see that, if $X$ satisfies~\eqref{axiom:P}, then $\pt(X)$ will have open cones, since $\pt(\Up U)$ and $\pt(\Down U)$ are open by construction. Using this, we can construct a functor on ordered locales that lands in $\OrdTop_\OC$.

\begin{definition}
	Let $\OrdLoc^\bullet$ denote the full subcategory of $\OrdLoc$ consisting of ordered locales satisfying axiom~\eqref{axiom:P}.
\end{definition}

\begin{corollary}\label{corollary:functor ordloc bullet to ordtopOC}
	There is a functor:
	\begin{align*}
		\pt: \OrdLoc^\bullet &\longrightarrow \OrdTop_\OC;\\
		(X,\Leq) & \longmapsto (\pt(X),\leq);\\
		f & \longmapsto \pt(f).
	\end{align*}
\end{corollary}

Having introduced the axiom~\eqref{axiom:P} on ordered locales, to obtain an adjunction we need the functor $\loc$ to land in the appropriate subcategory $\OrdLoc^\bullet$. Again, we currently do not know if this is automatically the case. The remedy turns out to be a relatively minor drawback.

\begin{lemma}\label{lemma:enough points and OC implies P}
	If $S$ is an ordered space with open cones and enough points, then $\loc(S)$ satisfies~\eqref{axiom:P}.
\end{lemma}
\begin{proof}
	Let $U$ and $V$ be open subsets of $S$ satisfying $U\Leq V$. We need to show that $\pt(U)\Leq \pt(V)$ with respect to the Egli-Milner order, meaning $\pt(U)\subseteq\down \pt(V)$ and $\pt(V)\subseteq\up \pt(U)$. For the first inclusion, take $\calF\in\pt(U)$. Since $S$ has enough points, there exists $x\in U$ such that $\calF=\calF_x$. This gives $x\in U\subseteq \down V$, so there exists an element $y\in V$ with $x\leq y$. This implies $\calF_x\leq\calF_y$ by \cref{lemma:open cones iff unit monotone}, and since $\calF_y\in\pt(V)$ this gives $\calF=\calF_x\in \down\pt(V)$, as desired. The inclusion $\pt(V)\subseteq\up\pt(U)$ is proved dually.
\end{proof}

\begin{definition}
	Let $\OrdTop^\bullet$ denote the full subcategory of $\OrdTop$ consisting of ordered spaces with enough points. Let $\OrdTop_\OC^\bullet$ be the full subcategory of $\OrdTop$ of ordered spaces with open cones and enough points.
\end{definition}

\begin{corollary}\label{corollary:functor ordtopOC bullet to ordloc bullet}
	There is a functor:
	\begin{align*}
		\loc: \OrdTop_\OC^\bullet &\longrightarrow \OrdLoc^\bullet;\\
		(S,\leq) & \longmapsto (\loc(S),\Leq);\\
		g & \longmapsto \loc(g).
	\end{align*}
\end{corollary}


\begin{example}
	It is not necessary for $S$ to have enough points for $\loc(S)$ to satisfy the property~\eqref{axiom:P}. Take for instance any topological space $S$ that does not have enough points, and equip it with the order that is equality of elements. The induced Egli-Milner order on $\loc(S)$ is just equality of open subsets, and in turn the induced order on $\pt(\loc(S))$ is equality of completely prime filters. The condition~\eqref{axiom:P} then says that $U=V$ has to imply $\pt(U)=\pt(V)$, which is true in any locale. This shows that the subcategory of $\OrdTop$ of ordered spaces whose induced locale satisfies~\eqref{axiom:P} is somewhat larger than $\OrdTop^\bullet$. But since we do not at present have a more intrinsic characterisation of such spaces, we prefer the current, slightly disadvantaged presentation. For applications to spacetimes this restriction is irrelevant, as the following example shows.
\end{example}

\begin{example}
	Any spacetime $M$ is Hausdorff, hence sober (\cref{proposition:hausdorff is sober}), and thus has enough points. Therefore the induced ordered locale $\loc(M)$ satisfies~\eqref{axiom:P} by \cref{lemma:enough points and OC implies P}, so by \cref{proposition:order in spatial locale from order on points} the order of $\loc(M)$ is completely determined by the causal relations between the points. (The converse is, however, not true:~see \cref{example:minkowski space point removed not T0}.)
\end{example}

\section{The adjunction}
In this section we finally prove that the adjunction $\loc\dashv\pt$ between $\Top$ and $\Loc$ lifts to the setting of ordered spaces and ordered locales. To summarise, we saw in \cref{corollary:functor ordtopOC to ordloc} that $\loc$ lifts to a functor $\OrdTop_\OC\to\OrdLoc$, where the open cone condition is necessary (and sufficient) to make the functor well-defined on morphisms. On the other hand, by \cref{corollary:functor ordloc to ordtop} the functor $\pt$ lifts to a functor $\OrdLoc\to \OrdTop$, but to ensure it lands in the subcategory of spaces with open cones, we restricted to the category $\OrdLoc^\bullet$ of ordered locales satisfying~\eqref{axiom:P}; see \cref{corollary:functor ordloc bullet to ordtopOC}. This, in turn, posed another restriction on the category of ordered spaces, where we restricted to those with enough points to ensure that the functor $\loc$ lands in $\OrdLoc^\bullet$; see \cref{corollary:functor ordtopOC bullet to ordloc bullet}. The process of restricting stabilises here, since $\loc$ lands in the category of \emph{sober} ordered spaces, which in particular have enough points.

To finish, we are left to prove the uni $\eta$ and counit $\epsilon$ of the old adjunction lift to the new functors. First, we see that monotonicity of the unit actually recovers precisely the open cone condition.

\begin{lemma}\label{lemma:open cones iff unit monotone}
	The unit $\eta_S\colon (S,\leq)\to (\pt(\loc(S)),\leq)$ is monotone if and only if $(S,\leq)$ has open cones.
\end{lemma}
\begin{proof}
	The monotonicity of $\eta_S$ means precisely that $x\leq y$ in $S$ implies ${\calF_x\leq \calF_y}$ in $\pt(\loc(S))$. Unpacking the definition of the order on points and using \cref{lemma:localic cones in a space}, this is equivalent to the open cone condition by \cref{proposition:open cones iff interiors of open neighbourhoods}.
\end{proof}

Next, we turn to the counit, which is always monotone.

\begin{lemma}\label{lemma:counit is monotone}
	The counit $\epsilon_X\colon (\loc(\pt(X)),\Leq)\to (X,\Leq)$ is monotone for every ordered locale $(X,\Leq)$.
\end{lemma}
\begin{proof}
	This follows from a calculation using \cref{lemma:localic cones in a space,lemma:cones of points contained in localic cones}:
	\[
	\Up \epsilon_X^{-1}(U)
	=
	\Up \pt(U)
	=
	\left(\up\pt(U)\right)^\circ
	\subseteq
	\up\pt(U)
	\subseteq
	\pt(\Up U)
	=
	\epsilon_X^{-1}(\Up U).
	\qedhere	
	\]
\end{proof}

This brings us to the main theorem of this chapter.

\begin{theorem}\label{theorem:adjunction ordtopOC bullet and ordloc bullet}
	There is an adjunction:
		\[
			\begin{tikzcd}[ampersand replacement=\&,column sep=large]
				{\OrdTop_\OC^\bullet} \& {\OrdLoc^\bullet.}
				\arrow[""{name=0, anchor=center, inner sep=0}, "\loc", shift left=2, from=1-1, to=1-2]
				\arrow[""{name=1, anchor=center, inner sep=0}, "\pt", shift left=2, from=1-2, to=1-1]
				\arrow["\dashv"{anchor=center, rotate=-90}, draw=none, from=0, to=1]
			\end{tikzcd}
		\]
\end{theorem}
\begin{proof}
	The functor $\loc$ is well defined by \cref{corollary:functor ordtopOC bullet to ordloc bullet}, and the functor $\pt$ is well defined by \cref{corollary:functor ordloc bullet to ordtopOC} together with the fact that $\pt(X)$ is always sober (\cref{lemma:locales are sober}). To obtain the adjunction, we claim that the natural bijection from the proof of the old adjunction (\cref{theorem:adjunction top and loc}) restricts to a natural bijection:
		\begin{align*}
			\OrdLoc^\bullet(\loc(S),X)  &\xrightarrow{~ \sim ~}  \OrdTop_\OC^\bullet(S,\pt(X))\\
			f&\xmapsto{\hspace{3ex}} \pt(f)\circ \eta_S\\ \epsilon_X\circ \loc(g)&{~\reflectbox{\ensuremath{\xmapsto{\hspace{3ex}}}}} ~g	.				
		\end{align*}
	That these functions are well-defined comes down to $\eta_S$ and $\epsilon_X$ being monotone, and this follows by \cref{lemma:open cones iff unit monotone,lemma:counit is monotone}, respectively. From the proof of \cref{theorem:adjunction top and loc} it now immediately follows that they define a natural bijection.
\end{proof}

\subsection{Fixed points}
Now that we have obtained an adjunction, we can calculate its fixed points. It turns out that, just as for the ordinary adjunction, the functors $\loc$ and $\pt$ will land in the subcategories of fixed points. 

For an ordered locale $(X,\Leq)$ to be a fixed point means that the counit $\epsilon_X$ is an isomorphism. This implies in particular that the un-ordered version of the counit is an isomorphism, which by \cref{proposition:characterisation of spatial locales} implies $X$ is a spatial locale. Similarly, by \cref{proposition:characterisation of sober spaces}, for $(S,\leq)$ to be a fixed point $S$ needs to be sober. Therefore, the fixed points of the new adjunction in \cref{theorem:adjunction ordtopOC bullet and ordloc bullet} will in particular be sober ordered spaces and spatial ordered locales. It turns out that the ordered locales do not actually pick up any additional conditions.

\begin{lemma}\label{lemma:counit inverse monotone iff bullet}
	For a spatial ordered locale $(X,\Leq)$, the inverse $\delta$ of the counit $\epsilon_X$ is monotone if and only if axiom~\eqref{axiom:P} holds.
\end{lemma}
\begin{proof}
	We get a map of frames $\delta^{-1}(\pt(U))=U$. The following calculation shows that if~\eqref{axiom:P} holds, then $\delta$ is monotone:
	\[
	\Up \delta^{-1}(\pt(U)) = \Up U = \delta^{-1}(\pt(\Up U)) \subseteq \delta^{-1}(\Up\pt(U)).
	\]
	The last inclusion is due to~\eqref{axiom:P}: by \cref{lemma:characterisation of axiom P,lemma:localic cones in a space} we get $\pt(\Up U)=(\up \pt(U))^\circ = \Up \pt(U)$.

	Conversely, suppose that $\delta$ is monotone. Given \cref{lemma:cones of points contained in localic cones}, it suffices to show that $\pt(\Up U)\subseteq \Up\pt(U)$. Using that $\delta$ is monotone gives an inclusion $\Up U = \Up \delta^{-1}(\pt(U))\subseteq \delta^{-1}(\Up \pt(U))$, and since $\delta$ is the inverse of $\epsilon_X$ this implies $\pt(\Up U)\subseteq \Up\pt(U)$, as desired. The proof involving down cones is dual.
\end{proof}

\begin{proposition}\label{proposition:order in spatial locale from order on points}
	In a spatial ordered locale satisfying~\eqref{axiom:P} and~\eqref{axiom:cones give order} we have
	\[
	U\Leq V \qquad\text{if and only if}\qquad \pt(U)\Leq \pt(V).
	\]
\end{proposition}
\begin{proof}
	One direction is provided by~\eqref{axiom:P} itself. For the converse, suppose that $\pt(U)\Leq \pt(V)$, meaning $\pt(U)\subseteq \down\pt(V)$ and $\pt(V)\subseteq\up\pt(U)$. Using the equations from \cref{lemma:characterisation of axiom P} and spatiality we hence get $U\sqleq \Down V$ and $V\sqleq \Up U$, which through~\eqref{axiom:cones give order} implies $U\Leq V$, as desired.
\end{proof}

For ordered spaces, the fixed points recover a condition that is well-known in the theory of ordered spaces. 

\begin{definition}\label{definition:T0-ordered}
	An ordered space $(S,\leq)$ is called \emph{$T_0$-ordered} if, whenever $x\not\leq y$, there either exists an open neighbourhood $U$ of $x$ such that $y\notin \up U$, or there exists an open neighbourhood $V$ of $y$ such that $x\notin\down V$.
\end{definition}

\begin{lemma}\label{lemma:T0 ordered iff inverse unit monotone}
	Let $(S,\leq)$ be a sober ordered space with open cones. Then the inverse of the unit $\eta_S^{-1}\colon \pt(\loc(S))\to S$ is monotone if and only if $S$ is $T_0$-ordered.
\end{lemma}
\begin{proof}
	Since the inverse of the unit sends a completely prime filter $\calF_x$ in $\Opens S$ to the unique element $x\in S$, we see that by contraposition $\eta_S^{-1}$ is monotone if and only if $x\not\leq y$ implies $\calF_x\not\leq\calF_y$. Using \cref{lemma:order on points in terms of localic cones}, the latter condition unpacks to:
		\[
			\calF_x\not\leq\calF_y
			\qquad\text{if and only if}\qquad
			\exists U\ni x\colon y\notin \up U \text{~or~} \exists V\ni y\colon x\notin \down V,
		\]
	so the $T_0$ condition is obtained.
\end{proof}

\begin{example}\label{example:minkowski space point removed not T0}
	Not every spacetime is $T_0$-ordered. An elementary example is Minkowski space with one point removed. To see that this fails the $T_0$-separation property, pick three distinct elements $x,y,p\in M$ in Minkowski space with ${x\caus p\caus y}$ that all lie on the same null geodesic. Passing to the spacetime $M\setminus \{p\}$, clearly $x\not\caus y$, since the only causal curve from $x$ to $y$ has now been discluded. However it remains the case that $\calF_x\leq \calF_y$, so by the proof of \cref{lemma:T0 ordered iff inverse unit monotone} the space $M\setminus\{p\}$ is not $T_0$-ordered.
\end{example}

The following shows that the points functor lands in the subcategory of fixed points. 

\begin{lemma}\label{lemma:pt(X) is T0-ordered}
	If $(X,\Leq)$ is an ordered locale, then $\pt(X)$ is $T_0$-ordered.
\end{lemma}
\begin{proof}
	By \cref{definition:order on points}, if $\calF\not\leq \calG$ we can either find $U\in \calF$ such that $\Up U\notin\calG$ or $V\in\calG$ such that $\Down V\notin\calF$. Using \cref{lemma:cones of points contained in localic cones} this implies, respectively, that $\calG\notin \up\pt(U)$ or $\calF\notin\down\pt(V)$.
\end{proof}

\begin{theorem}\label{theorem:equivalence ordtop and ordloc}
	The adjunction of \cref{theorem:adjunction ordtopOC bullet and ordloc bullet} restricts to an equivalence of categories between the full subcategory of $\OrdLoc^\bullet$ of spatial ordered locales, and the full subcategory of $\OrdTop_\OC^\bullet$ of sober $T_0$-ordered spaces:
		\[
			\left\{\begin{array}{c}
				\text{\normalfont sober $T_0$-ordered spaces} \\
				\text{\normalfont with open cones}
			\end{array}\right\}
		\simeq
			\left\{\begin{array}{c}
				\text{\normalfont spatial ordered} \\
				\text{\normalfont locales with~\eqref{axiom:P}}
			\end{array}\right\}
		.
		\]
\end{theorem}

\begin{remark}
	We have seen that any ordered locale $\loc(S)$ coming from an ordered space $S$ with open cones satisfies the axioms~\eqref{axiom:cones give order},~\eqref{axiom:LV} and~\eqref{axiom:wedge}. None of these axioms are actually used in constructing the adjunction above. It is therefore not clear to what extent our notion of ordered locale in \cref{definition:ordered locale} is precisely a point-free axiomatisation of ordered spaces (with open cones). Besides our discussion of alternative definitions throughout the thesis (and below), we leave this question for future work, and for now take \cref{theorem:adjunction ordtopOC bullet and ordloc bullet,theorem:equivalence ordtop and ordloc} as sufficient evidence that the current framework of ordered locales is an interesting one in which to study point-free spacetimes.
\end{remark}

\begin{remark}
	The classic adjunction $\loc\dashv\pt$ from \cref{theorem:adjunction top and loc} is special in that it is \emph{idempotent}, meaning that $\eta_{\pt}$ (equivalently $\loc\circ \eta$, $\epsilon_{\loc}$ or $\pt\circ \epsilon$) defines a natural isomorphism. Explicitly, it means that for every locale $X$ the map $\eta_{\pt(X)}$ is an isomorphism, or equivalently that for every space $S$ the map $\epsilon_{\loc(S)}$ is an isomorphism. The new adjunction in \cref{theorem:adjunction ordtopOC bullet and ordloc bullet} is also idempotent, due to \cref{lemma:counit inverse monotone iff bullet,lemma:T0 ordered iff inverse unit monotone}.
\end{remark}

\section{Upper and lower adjunctions}
\label{section:upper-lower adjunctions}
This work is based on Sections~7 and~8 in~\cite{heunenSchaaf2024OrderedLocales}. The proof of the adjunction in \cref{theorem:adjunction ordtopOC bullet and ordloc bullet} can be split into upper and lower parts, giving two new adjunctions for ordered locales whose morphisms are only upper or lower monotone:
\begin{enumerate}[label = \textbullet]
	\item the category of ordered locales with upper monotone maps: $\OrdLoc^{\Up}$;
	\item the category of ordered locales with lower monotone maps: $\OrdLoc^{\Down}$.
\end{enumerate}

Inspecting the proof of \cref{proposition:loc well defined iff OC}, we see that the upper/lower open cone condition corresponds precisely to $\loc(g)$ being upper/lower monotone for every continuous monotone function $g$. Denote:
\begin{enumerate}[label = \textbullet]
	\item the full subcategory of ordered spaces with open upper cones $\OrdTop_\OC^{\up}$;
	\item the full subcategory of ordered spaces with open lower cones $\OrdTop_\OC^{\down}$.
\end{enumerate}
Then \Cref{corollary:functor ordtopOC to ordloc} can be rewritten as follows.

\begin{corollary}
	There are functors:
	\begin{center}
		\begin{minipage}{.45\textwidth}
			\begin{align*}
				\OrdTop_\OC^{\up} &\xrightarrow{~\loc^{\up}~} \OrdLoc^{\Up}\\
				(S,\leq) & \longmapsto (\loc(S),\LeqU)\\
				g & \longmapsto \loc(g),
			\end{align*}
		\end{minipage}%
		\begin{minipage}{.45\textwidth}
			\begin{align*}
				\OrdTop_\OC^{\down} &\xrightarrow{~\loc^{\down}~} \OrdLoc^{\Down}\\
				(S,\leq) & \longmapsto (\loc(S),\LeqL)\\
				g & \longmapsto \loc(g).
			\end{align*}
		\end{minipage}
	\end{center}
\end{corollary}

Note here that we equip the locale $\loc(S)$ with the upper order $\LeqU$ and lower order $\LeqL$, respectively, and not the Egli-Milner order. These define ordered locales by \cref{proposition:space with EM order is ordered locale}.

Similarly to how $\LeqEM$ is defined in terms of $\LeqU$ and $\LeqL$, there are two new orders on the points of an ordered locale:
\begin{enumerate}[label = \textbullet]
	\item $\mathcal{F}\leqU\mathcal{G}$ if and only if $\Up U\in\mathcal{G}$ for all $U\in \mathcal{F}$;
	\item $\mathcal{F}\leqL\mathcal{G}$ if and only if $\Down V\in\mathcal{F}$ for all $V\in\mathcal{G}$.
\end{enumerate}
We call these the \emph{upper} and \emph{lower order} on points, respectively. Note that $\leqU$ and $\leqL$ are just the orders on points induced through \cref{definition:order on points} by localic orders $\LeqU$ and $\LeqL$, respectively. Also, by \cref{lemma:order on points in terms of localic cones}, the order on the points of an ordered locale from \cref{definition:order on points} simply reads $\mathcal{F}\leq\mathcal{G}$ if and only if $\mathcal{F}\leqU\mathcal{G}$ and $\mathcal{F}\leqL\mathcal{G}$.

\begin{example}\label{example:ordered locale with inclusion order gives specialisation order}
	Recall from \cref{example:inclusion is ordered locale} that any locale $X$ induces an ordered locale $(X,\sqleq)$, where the order is the intrinsic inclusion relation on the opens. For any $U\in\Opens X$ we find $\Up U = X$, the largest open, and $\Down U = U$. Hence the induced orders on the points become:
	\begin{enumerate}[label = \textbullet]
		\item $\mathcal{F}\leqU \mathcal{G}$ for all $\mathcal{F},\mathcal{G}\in\pt(X)$;
		\item $\mathcal{F}\leqL\mathcal{G}$ if and only if $\mathcal{G}\subseteq\mathcal{F}$.
	\end{enumerate}
	In this way, if $X=\loc(S)$ is some topological space $S$, the (opposite of the) \emph{specialisation order} is reobtained: $x\in\overline{\{y\}}$ if and only if $\mathcal{F}_x\subseteq \mathcal{F}_y$ if and only if $\mathcal{F}_y\leqL \mathcal{F}_x$. 
	
	Further, $X\mapsto (X,\sqleq)$ defines a left adjoint to the forgetful functor:
	\[
	\begin{tikzcd}[column sep = 1.5cm]
		\Loc \arrow[r, ""{name= X}, shift left=2] & \OrdLoc^{\Down}. \arrow[l, ""{name=Y}, shift left=2]
		\ar[phantom,from=X, to=Y,  "\vdash" rotate=90]
	\end{tikzcd}
	\]
	(We actually get a functor $\Loc\to \OrdLoc$, but the adjunction only works when restricted to lower monotone maps.)
\end{example}

Looking at the proof of \cref{lemma:pt(f) monotone}, we see that if a map of locales $f$ is upper/lower monotone, then $\pt(f)$ is monotone with respect to the upper/lower order on points. Hence \cref{corollary:functor ordloc to ordtop} can be rewritten as follows.

\begin{corollary}
	There are functors%
	\begin{center}
		\begin{minipage}{.45\textwidth}
			\begin{align*}
				\OrdLoc^{\Up} &\xrightarrow{~\pt^{\up}~} \OrdTop\\
				(X,\Leq) & \longmapsto (\pt(X),\leqU)\\
				f & \longmapsto \pt(f)
			\end{align*}
		\end{minipage}%
		\begin{minipage}{.45\textwidth}
			\begin{align*}
				\OrdLoc^{\Down} &\xrightarrow{~\pt^{\down}~} \OrdTop\\
				(X,\Leq) & \longmapsto (\pt(X),\leqL)\\
				f & \longmapsto \pt(f)
			\end{align*}
		\end{minipage}
	\end{center}
\end{corollary}


\begin{example}
	Starting with an ordered space $(S,\leq)$, we can first take the upper-ordered locale $(\loc(S),\LeqU)$. Using the same proof as \cref{lemma:localic cones in a space}, it is easy to verify that $\Up U = (\up U)^\circ$. On the other hand, we find $\Down U = S$, since $W\LeqU U$ holds in particular for every open subset $W\supseteq U$. Hence applying the functor $\pt^{\down}$ to $(\loc(S),\LeqU)$ gives the space $\pt(\loc(S))$ that is equipped with the following order: $\mathcal{F}\leqL\mathcal{G}$ if and only if for all $V\in\mathcal{G}$ we have $S=\Down V \in \mathcal{F}$, which is always true. Therefore the induced lower order on the points becomes trivial. Similarly, the resulting space under $\pt^{\up}\circ \loc^{\down}$ will have the trivial order. 
\end{example}

The following result generalises \cref{lemma:open cones iff unit monotone}.

\begin{lemma}
	For any ordered space $(S,\leq)$:
	\begin{enumerate}[label = \textbullet]
		\item $\eta_S\colon S\to \pt^{\up}(\loc^{\up} (S))$ is monotone if and only if $S$ has open upper cones;
		\item $\eta_S\colon S\to \pt^{\down}(\loc^{\down} (S))$ is monotone if and only if $S$ has open lower cones.
	\end{enumerate}
\end{lemma}

\begin{theorem}\label{theorem:adjunction upperlower}
	There are adjunctions
	\begin{equation*}
		\begin{tikzcd}[column sep = 1.5cm]
			\OrdTop_{\OC}^{\up\bullet} \arrow[r, "\Opens^{\up}"{name= X}, shift left=2] & \OrdLoc^{\Up\bullet}
			\arrow[l, "\pt^{\up}"{name=Y}, shift left=2]
			\ar[phantom,from=X, to=Y,  "\vdash" rotate=90]
		\end{tikzcd}
		\quad\text{and}\quad
		\begin{tikzcd}[column sep = 1.5cm]
			\OrdTop_{\OC}^{\down\bullet} \arrow[r, "\Opens^{\down}"{name= X}, shift left=2] & \OrdLoc^{\Down\bullet}
			\arrow[l, "\pt^{\down}"{name=Y}, shift left=2]
			\ar[phantom,from=X, to=Y,  "\vdash" rotate=90]
		\end{tikzcd}
	\end{equation*}
	between the full subcategories of $\OrdTop$ of ordered spaces with open upper/lower cones and enough points, and the full subcategories of $\OrdLoc^{\Up}$ and $\OrdLoc^{\Down}$ of ordered locales satisfying, respectively:
	\begin{align}
		U \LeqU V \text{ in }X 
		\qquad \text{implies} \qquad
		\pt(U) \LeqU \pt(V) \text{ in } \loc^{\up}(\pt^\uparrow(X))
		\tag{\ref{axiom:P}$^\uparrow$}, \\
		U \LeqL V \text{ in }X 
		\qquad \text{implies} \qquad
		\pt(U) \LeqL \pt(V) \text{ in } \loc^{\down}(\pt^\downarrow(X)).
		\tag{\ref{axiom:P}$^\downarrow$}
	\end{align}
\end{theorem}

Next, we discuss the fixed points of these adjunctions.

\begin{definition}
	An ordered space $(S,\leq)$ is \emph{$T_\mathrm{U}$-ordered} if $x\nleqslant y$ implies there exists an open neighbourhood $U$ of $x$ such that $y\notin \up U$. Similarly, we say the space is \emph{$T_\mathrm{L}$-ordered} if $x\nleqslant y$ implies there exists an open neighbourhood $V$ of $y$ such that $x\notin \down V$.
\end{definition}

Note that the $T_\mathrm{U}$- and $T_\mathrm{L}$-order axioms both imply the $T_0$-ordered axiom.

\begin{lemma}\label{lemma:unit iso upperlower}
	For a sober ordered space $(S,\leq)$:
	\begin{itemize}
		\item $\eta_S^{-1}\colon \pt^{\up}(\loc^{\up}(S))\to S$ is monotone if and only if $S$ is $T_\mathrm{U}$-ordered;
		\item $\eta_S^{-1}\colon \pt^{\down}(\loc^{\down}(S))\to S$ is monotone if and only if $S$ is $T_\mathrm{L}$-ordered.
	\end{itemize}
\end{lemma}
\begin{proof}
	The map $\eta_S^{-1}$ is monotone in the upper setting if and only if $x\nleqslant y$ implies $\mathcal{F}\nleqslant_\textsc{u}\mathcal{F}_y$. The latter holds if and only if there exists $U\in\mathcal{F}_x$ such that $\Up U\notin \mathcal{F}_y$. Since $S$ has open upper cones and by the proof of \cref{lemma:localic cones in a space} we see $\Up U = \up U$, so we obtain the equivalence between the $T_\mathrm{U}$-ordered condition. The proof in the lower setting is analogous.
\end{proof}

Generalising the proof of \cref{lemma:pt(X) is T0-ordered} gives the following.

\begin{lemma}
	For any ordered locale $(X,\Leq)$, the space $\pt^{\up}(X)$ is $T_\mathrm{U}$-ordered, and the space $\pt^{\down}(X)$ is $T_\mathrm{L}$-ordered.
\end{lemma}

Putting it all together gives two new dualities.

\begin{theorem}\label{theorem:duality upperlower}
	The adjunctions of \cref{theorem:adjunction upperlower} restrict to equivalences between the categories of ordered locales and upper/lower monotone morphisms satisfying the upper/lower version of~\eqref{axiom:P}, and the categories of sober $T_\mathrm{U}$/$T_\mathrm{L}$-ordered spaces with open upper/lower cones.
\end{theorem}

\subsection{Relation to Esakia duality}\label{section:esakia duality}
Priestley and Esakia spaces arise as topological representations of distributive lattices \cite{priestley1972OrderedTopologicalSpaces} and Heyting algebras \cite{esakia1974topologicalKripkemodels}, respectively. This section shows how Esakia duality fits into our adjunctions. We recall basic definitions, but refer to the literature for detail: see \cite{bezhanishvili2010BitopologicalDualityDistributive} for a detailed treatment from a bitopological point of view, \cite{landsman2021LogicQuantumMechanics} for an expository account, and the appendix of \cite{davey2003CoalgebraicViewHeyting} for a detailed proof of Esakia duality. We remark also that these results are different from the \emph{localic} Priestley duality from \cite{townsend1997LocalicPriestleyDuality}.

\begin{definition}
	A \emph{Priestley space} is an ordered space $(S,\leq)$ where $S$ is compact, and the \emph{Priestley separation axiom} holds: if $x\nleqslant y$ there exists a clopen upset $U$ such that $x\in U$ and $y\notin U$. 
\end{definition}

Priestley spaces are also known as \emph{ordered Stone spaces}, since any Priestley space is a Stone space by e.g.~\cite[Lemma~3.2]{bezhanishvili2010BitopologicalDualityDistributive}. Hence, importantly here, every Priestley space is sober. 

\begin{lemma}\label{lemma:priestley spaces are TU and TL ordered}
	Any Priestley space is $T_\mathrm{U}$- and $T_\mathrm{L}$-ordered.
\end{lemma}
\begin{proof}
	If $x\nleqslant y$, then by the Priestley separation axiom there exists a clopen upset $U$ with $x\in U$ but $y\notin U$. Since $U$ is an upset we get $\up U = U$, so this gives the $T_\mathrm{U}$-ordered condition. Further, the complement of $U$ is a clopen downset that contains $y$ but not $x$, giving the $T_\mathrm{L}$-ordered condition.
\end{proof}

Using \cite[Theorem~4.2]{davey2003CoalgebraicViewHeyting} we can adapt the definition of Esakia spaces in terms of open cone conditions as follows.

\begin{definition}\label{definition:esakia spaces}
	A Priestley space is called:
	\begin{enumerate}[label = \textbullet]
		\item an \emph{Esakia space} if it has open lower cones;
		\item a \emph{co-Esakia space} if it has open upper cones;
		\item a \emph{bi-Esakia space} if it has open cones.
	\end{enumerate}
\end{definition}
The definition of morphisms between Esakia spaces are slightly more involved, essentially due to the fact that they have to ensure the induced map between Heyting algebras respects the Heyting implication. For us, it suffices to know that morphisms between Esakia spaces are, in particular, continuous and monotone. Refer to \cite[Section~7]{bezhanishvili2010BitopologicalDualityDistributive} for more details. Denote by $\Esakia$, $\coEsakia$, and $\biEsakia$ the categories of Esakia, co-Esakia, and bi-Esakia spaces, respectively. Similarly, denote by $\Heyt$, $\coHeyt$, and $\biHeyt$ the categories of Heyting, co-Heyting, and bi-Heyting algebras, respectively. \emph{Esakia duality} says that there are equivalences of categories:
\[
\Esakia \cong \Heyt^\op,
\qquad\quad
\coEsakia \cong \coHeyt^\op,
\qquad\quad
\biEsakia \cong \biHeyt^\op.
\]
\Cref{definition:esakia spaces} lets us interpret Esakia spaces as ordered spaces with open cones.

\begin{lemma}\label{lemma:esakia spaces include in ordtop open lower cones}
	There are inclusion functors:
	\[
	\Esakia \hookrightarrow \OrdTop_{\OC}^{\down\bullet},
	\quad
	\coEsakia \hookrightarrow \OrdTop_{\OC}^{\up\bullet},
	\quad
	\biEsakia \hookrightarrow \OrdTop_{\OC}^\bullet.
	\]
\end{lemma}

Moreover, the following proposition shows that these inclusion functors land exactly in the respective fixed points.

\begin{proposition}\label{proposition:esakia spaces are fixed points}
	The following holds:
	\begin{enumerate}[label =\textbullet]
		\item any Esakia space $S$ is isomorphic to $\pt^{\down}(\loc^{\down} (S))$;
		\item any co-Esakia space $S$ is isomorphic to $\pt^{\up}(\loc^{\up}(S))$;
		\item any bi-Esakia space $S$ is isomorphic to $\pt(\loc(S))$.
	\end{enumerate}
\end{proposition}
\begin{proof}
	This follows from the fact that Priestley spaces are sober, and by combining \cref{lemma:unit iso upperlower,lemma:priestley spaces are TU and TL ordered}.
\end{proof}


We now exhibit one sense in which Esakia duality fits into our framework. To obtain an Esakia space from a Heyting algebra, one constructs the space of prime filters, suitably topologised. We denote the resulting functor by ${\pf\colon \Heyt^\op\to\Esakia}$, and similarly for the other varieties of Esakia spaces. Conversely, the lattice of clopen upsets of an Esakia space defines a Heyting algebra, and we denote the resulting functor by $\Clopup\colon\Esakia\to \Heyt^\op$. Consider the following diagrams:
\[
\begin{tikzcd}[column sep = 1cm]%
	\Esakia \arrow[d,hookrightarrow]\arrow[r, "\Clopup"{name= A}, shift left=2]
	& 
	\Heyt^\op  \arrow[l, "\pf"{name=B}, shift left=2]\arrow[d, "\loc^{\down}\circ \pf"]
	\\
	\OrdTop_{\OC}^{\down\bullet} \arrow[r, "\loc^{\down}"{name= X}, shift left=2] 
	&
	\OrdLoc^{\Down\bullet},
	\arrow[l, "\pt^{\down}"{name=Y}, shift left=2]
	\ar[phantom,from=X, to=Y,  "\vdash" rotate=90]
\end{tikzcd}
\qquad
\begin{tikzcd}[column sep = 1cm]%
	\coEsakia \arrow[d,hookrightarrow]\arrow[r, "\Clopup"{name= A}, shift left=2]
	& 
	\coHeyt^\op  \arrow[l, "\pf"{name=B}, shift left=2]\arrow[d, "\loc^{\up}\circ \pf"]
	\\
	\OrdTop_{\OC}^{\up\bullet} \arrow[r, "\loc^{\up}"{name= X}, shift left=2] 
	&
	\OrdLoc^{\Up\bullet},
	\arrow[l, "\pt^{\up}"{name=Y}, shift left=2]
	\ar[phantom,from=X, to=Y,  "\vdash" rotate=90]
\end{tikzcd}
\]%
\[
\begin{tikzcd}[column sep = 1cm]%
	\biEsakia \arrow[d,hookrightarrow]\arrow[r, "\Clopup"{name= A}, shift left=2]
	& 
	\biHeyt^\op  \arrow[l, "\pf"{name=B}, shift left=2]\arrow[d, "\loc\circ \pf"]
	\\
	\OrdTop_{\OC}^{\bullet} \arrow[r, "\loc"{name= X}, shift left=2] 
	&
	\OrdLoc^{\bullet}.
	\arrow[l, "\pt"{name=Y}, shift left=2]
	\ar[phantom,from=X, to=Y,  "\vdash" rotate=90]
\end{tikzcd}
\]

Here $\loc^{\down}\circ\pf\colon\Heyt^\op\to\OrdLoc^{\Down\bullet}$ lands in the fixed points of $\OrdLoc^{\Down\bullet}$, since it is in the image of $\loc^{\down}$. The following shows that the left- and right-directed squares in the above diagrams commute up to isomorphism.
\begin{proposition}
	There are natural isomorphisms:
	\begin{align*}
		\pf & \cong \pt^{\down}\circ \loc^{\down}\circ \pf
		& \loc^{\down}\circ \pf\circ \Clopup &\cong \loc^{\down} \\
		\pf & \cong \pt^{\up}\circ \loc^{\up}\circ \pf
		& \loc^{\up}\circ \pf\circ \Clopup&\cong \loc^{\up} \\
		\pf & \cong \pt\circ \loc \circ \pf
		& \loc \circ \pf\circ \Clopup & \cong \loc
	\end{align*}
\end{proposition}
\begin{proof}
	The latter isomorphisms follow immediately from the fact that $\pf$ and $\Clopup$ are mutually inverse up to isomorphism, and the former from \cref{proposition:esakia spaces are fixed points}.
\end{proof}

\part[Techniques of point-free topology in relativity]{Techniques of point-free topology in relativity\raisebox{.15\baselineskip}{\Large\footnotemark}}\label{chapter:topftir}
\chapter[Secondary structure of ordered locales]{Secondary structure\\of ordered locales}\label{section:secondary structure}
\footnotetext{Cf.~\citetitle{penrose1972TechniquesDifferentialTopology}, \cite{penrose1972TechniquesDifferentialTopology}.}
In this chapter we investigate the secondary structure of ordered locales. These techniques were developed in part to study how several of the rungs on the \emph{causal ladder} (\cref{section:causal ladder}), which are apparently of topological nature, could be captured in purely localic language. This work was conducted joint with Prakash Panangaden.

In particular, we will be able to construct a \emph{convex hull} operator, a \emph{causal complement} operator, and a \emph{causal diamond} operator. Such operations are needed when defining algebraic quantum field theories, which we discuss briefly in \cref{section:aqft on ordered locales}. We use the convex hull to define a new class of \emph{convex} ordered locales, which generalise the ``convex'' ordered spaces, and in particular generalise strongly causal spacetimes. Building on the insights of \cref{section:ordered locales from biframes}, we define the \emph{locale of futures} and \emph{pasts} $X^\triup$ and $X^\tridown$.\vspace*{-1em}

\section{Convexity}\label{section:convexity}
For this section, fix an ordered locale $(X,\Leq)$. We define and study the \emph{convex hull} operator. For ease of reading, instead of continually referring to \cref{lemma:properties of localic cones} for properties of the localic cones, we will from now on take them to be understood.

First, we recall the notion of convexity in the setting of ordered spaces. \begin{definition}\label{definition:convex pointwise}
	A subset $C\subseteq S$ of an ordered space $(S,\leq)$ is called \emph{pointwise convex} if $z\in C$ as soon as $x\leq z\leq y$ for some $x,y\in C$. 
\end{definition}

This can be translated into the language of regions as follows.
\begin{lemma}\label{lemma:convex pointwise in terms of cones}
	The following are equivalent, for a subset $C\subseteq S$ of an ordered space:
	\begin{enumerate}[label = (\alph*)]
		\item$C$ is pointwise convex;
		\item for all $A,B\subseteq C$ we have $\up A\cap \down B\subseteq C$;
		\item $\up C\cap \down C\subseteq C$.
	\end{enumerate}
\end{lemma}
\begin{proof}
	This follows straightforwardly by unpacking the definition of cones.
\end{proof}

\begin{figure}[t]
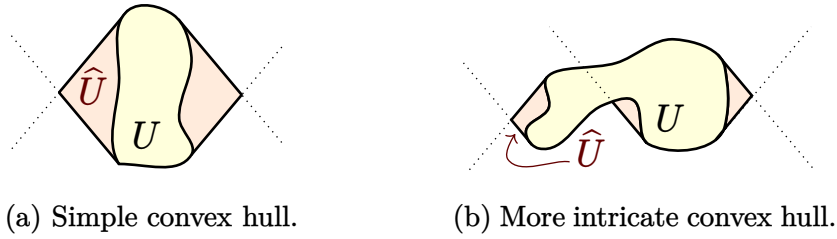
\centering
	\begin{subfigure}[b]{0.35\textwidth}\centering
		\input{tikz/inkscape-export/fig-27-simple}
		\caption{Simple convex hull.}
	\end{subfigure}\hfil
	\begin{subfigure}[b]{0.35\textwidth}\centering
		\input{tikz/inkscape-export/fig-27}
		\caption{More intricate convex hull.}
	\end{subfigure}
	\caption{Illustration of typical convex hulls in a Minkowski-like space.}
	\label{figure:convex hull}
\end{figure}

Hence we see that $C\subseteq S$ is pointwise convex if and only if $C=\up C\cap \down C$ (cf. \cref{figure:convex hull}). The latter condition can now be straightforwardly adapted to the localic setting.

\begin{definition}\label{definition:convex hull}\label{definition:convex in locale}
	The \emph{convex hull operator} of an ordered locale $(X,\Leq)$ is the function
	\[
	\hull{(-)}\colon \Opens X\longrightarrow \Opens X;\qquad U\longmapsto \hull{U}:=\Up U\wedge \Down U.
	\]
	An open region $U$ in an ordered locale is called \emph{convex} if it is equal to its convex hull: $U=\hull{U} =\Up U \wedge \Down U$. See \cref{figure:convex hull} for visual intuition.
\end{definition}

\begin{lemma}\label{lemma:properties of convex hull}
	For any ordered locale $(X,\Leq)$ the convex hull operator satisfies the following properties:
	\begin{enumerate}
		\item[(a)] $U\sqleq \hull{U}$;
		\item[(b)] $\hull{\hull{U}}\sqleq \hull{U}$;
		\item[(c)] If $U\sqleq V$ then $\hull{U}\sqleq \hull{V}$;
		\item[(d)] $\Up U = \Up \hull{U}$ and $\Down U = \Down \hull{U}$;
	\end{enumerate}
	If $X$ satisfies~\eqref{axiom:cones give order}, then furthermore:
	\begin{enumerate}
		\item[(e)] $U\Leq V$ if and only if $\hull{U}\Leq \hull{V}$.
	\end{enumerate}
\end{lemma}
\begin{proof}
	Properties (a) and (c) follow immediately from \cref{lemma:properties of localic cones}. For point~(b), note that monotonicity of the localic cones gives $\Up (U\wedge V)\sqleq \Up U\wedge \Up V$, and similarly for downsets. Use this together with the idempotence of the localic cones to calculate:
	\[
	\Up (\Up U\wedge \Down U) \wedge \Down (\Up U\wedge \Down U)
	\sqleq \Up U \wedge \Up \Down U \wedge \Down \Up U \wedge \Down U
	= \hull{U}\wedge \Up \Down U \wedge \Down \Up U \sqleq \hull{U}.
	\]
	
	For (d), note by (a) and monotonicity we get $\Up U\sqleq \Up \hull{U}$. For the converse, simply use idempotence to calculate $\Up\hull{U} = \Up(\Up U \wedge \Down U)\sqleq \Up U \wedge \Up\Down U\sqleq \Up U$. Point~(e) is an immediate corollary of (d).
\end{proof}

The localic cones themselves determine convex regions.

\begin{lemma}\label{lemma:cones are convex}
	For any $U\in\Opens X$, the localic cones $\Up U$ and $\Down U$ are convex.
\end{lemma}
\begin{proof}
	Simply note that the hull of $\Up U$ is contained in $\Up\Up U \wedge\Down\Up U = \Up U\wedge \Down \Up U\sqleq \Up U$, so equality follows by \cref{lemma:properties of convex hull}(b). The proof for $\Down U$ is similar.
\end{proof}

\begin{lemma}\label{lemma:convex closed under meets}
	Convex regions are closed under meets.
\end{lemma}
\begin{proof}
	This is again an elementary fact about closure operators. The proof is identical to the first paragraph in the proof of \cref{lemma:past sets form a frame}.
\end{proof}

\begin{remark}
	Note that an arbitrary convex region is therefore always of the form $\Up U\wedge \Down V$ for some $U,V\in\Opens X$.
\end{remark}

With a new localic definition of convexity, we investigate to what extent it is related to the pointwise one. There are two cases we need to consider: the relation between convexity in $S$ and $\loc(S)$, and the relation between convexity in $X$ and $\pt(X)$. It turns out that in a space with open cones, localic convexity coincides with pointwise convexity.

\begin{lemma}\label{lemma:convex pointwise iff localic convex}
	Any pointwise convex open subset $U$ of $(S,\leq)$ is convex in~$\loc(S)$. If $(S,\leq)$ has open cones, then $U$ is pointwise convex if and only if it is convex.
\end{lemma}
\begin{proof}
	That pointwise convexity of $U$ implies convexity of $U$ with respect to $\Leq$ follows immediately by combining \cref{lemma:convex pointwise in terms of cones,lemma:localic cones in a space}, and using that interiors respect finite intersections. For the converse, if $S$ has open cones and $U$ is convex then $U=\hull{U}=\Up U\cap \Down U = \up U\cap \down U$, showing $U$ is pointwise convex.
\end{proof}

In a locale that lacks points, we cannot expect convexity in $\pt(X)$ to translate to convexity in $X$. However, if $X$ is spatial the situation is nice.

\begin{lemma}\label{lemma:convex iff pointwise convex}
	If $U$ is a convex open in an ordered locale, then $\pt(U)$ is pointwise convex. If the locale is spatial and satisfies~\eqref{axiom:P}, then $U$ is convex if and only if $\pt(U)$ is pointwise convex.
\end{lemma}
\begin{proof}
	If $U$ is convex then $\hull{U}\sqleq U$, so applying the points functor and using \cref{lemma:cones of points contained in localic cones}:
	\[
	\up \pt(U) \cap \down \pt(U)\subseteq \pt(\Up U)\cap \pt(\Down U)=\pt(\hull{U}) \subseteq \pt(U),
	\]
	showing $\pt(U)$ is pointwise convex by \cref{lemma:convex pointwise in terms of cones}.
	
	For the converse, suppose that $\pt(U)$ is pointwise convex, meaning $\pt(U)=\up\pt(U) \cap\down\pt(U)$. Using the equations from \cref{lemma:characterisation of axiom P} this gives $\pt(U)=\pt(\hull{U})$, so by spatiality we get $U=\hull{U}$. 
\end{proof}

\subsection{The definition of convex ordered locales}
\label{section:convex locales}
In this section we formally introduce a new class of ordered locales that we call \emph{convex}. We informally encountered these during our discussion of biframes in \cref{section:ordered locales from biframes}. First, recall the definition of a convex ordered space.

\begin{definition}\label{definition:convex ordered space}\label{definition:strongly causal space}
	An ordered space $(S,\leq)$ is called \emph{convex} if for every element $x\in S$ and every open neighbourhood $V$ of~$x$, there exists a pointwise convex open neighbourhood $U$ such that $x\in U\subseteq V$.
	
	Equivalently: $S$ is convex if the pointwise convex open subsets of $S$ form a basis for $\Opens S$. Explicitly: $S$ is convex if every open subset $U\in \Opens S$ can be written as a union $\bigcup (F_i \cap P_i)$, where $\up F_i = F_i$ and $\down P_i = P_i$ are open sets.
\end{definition}

\begin{remark}
	Note that a spacetime $M$ is strongly causal in the sense of relativity theory (\cref{section:causal ladder}) if and only if the ordered space $(M,\caus)$ is convex in the sense of \cref{definition:convex ordered space}.
\end{remark}

Recall from \cref{definition:locale basis} that a \emph{base} (or \emph{basis}) for a locale $X$ is a subset $\calB\subseteq \Opens X$ so that every element in $\Opens X$ can be written as a join of elements in $\calB$. This generalises the notion of a base for a topology.

\begin{definition}
	An ordered locale $(X,\Leq)$ is called \emph{convex} if the convex opens form a base of $\Opens X$. Explicitly: $X$ is convex if every $U\in \Opens X$ can be written as a join of the form $\bigvee (\Up U_i \wedge \Down V_i)$.
\end{definition}

In the remainder of this section we briefly study properties of convex locales, and how this condition translates between spaces and locales. First, recall that a spacetime $M$ is \emph{causal} (\cref{section:causal ladder}) if the causality relation $\caus$ is anti-symmetric, i.e. $x\caus y\caus x$ implies $x= y$. Any strongly causal spacetime is causal \cite[Figure~8]{minguzzi2019LorentzianCausalityTheory}. It turns out that this has a localic analogue. For that, we investigate the anti-symmetry of $\Leq$.

\begin{lemma}\label{lemma:anti-symmetry and convex hulls}
	In an ordered locale, if $U\Leq V\Leq U$ then $\hull{U}=\hull{V}$. If~\eqref{axiom:cones give order} holds, then conversely $\hull{U}= \hull{V}$ implies $U\Leq V\Leq U$.
\end{lemma}
\begin{proof}
	If $U\Leq V\Leq U$ then the equality of the hulls follows directly from \cref{lemma:properties of localic cones}(a) and \cref{lemma:properties of convex hull}(b). Conversely, if $\hull{U}=\hull{V}$, then $U\sqleq \hull{U} = \Up V\wedge \Down V\sqleq \Down V$, and similarly $V\sqleq \Up U$, so by~\eqref{axiom:cones give order} we get $U\Leq V$. A symmetric argument gives ${V\Leq U}$, as desired.
\end{proof}

Given this, it makes little sense to talk about the order $\Leq$ on the opens of a locale being anti-symmetric, since in an ordered topological space this means two opens are equal as soon as their convex hulls are equal, which is rarely the case (e.g. certainly not in Minkowski space). Instead, it is more appropriate to study when the induced order on the \emph{points} of the locale becomes anti-symmetric. For convenience, recall from \cref{section:functor ordloc to ordtop} that the order on $\pt(X)$ is given by:
\[
\calF \leq \calG
\quad \text{if and only if} \quad
\forall U \in \calF \colon \Up U \in \calG
\;\;\text{ and }\;\;
\forall V \in \calG \colon \Down V \in \calF.
\]
We define the \emph{hull} of a point $\calF\in\pt(X)$ as:
\[
\hull{\calF}:=\{\hull{U}:U\in\calF\}\subseteq\Opens X.
\]
The following is the point-analogue of \cref{lemma:anti-symmetry and convex hulls}.

\begin{lemma}\label{lemma:anti-symmetry and convex hulls points}
	For points $\calF$ and $\calG$ in an ordered locale we have $\hull{\calF}=\hull{\calG}$ if and only if $\calF\leq \calG\leq\calF$ .
\end{lemma}
\begin{proof}
	From $\calF\leq \calG\leq\calF$ it follows $\hull{\calF}\subseteq \calG$ and $\hull{\calG}\subseteq \calF$, so by idempotence of the hull \cref{lemma:properties of convex hull}(b) we get the desired equality. 
	
	Conversely, if $\hull{\calF}=\hull{\calG}$, take an arbitrary open $U\in\calF$. Then $\hull{U}\in\hull{\calG}$, so there exists $V\in\calG$ such that $\hull{U}=\hull{V}$. But then, using \cref{lemma:properties of convex hull}, we find $V\sqleq \hull{U}\sqleq \Up U$, so since $\calG$ is a filter it follows that $\Up U\in\calG$. It follows similarly that if $V\in\calG$ then $\Down V\in\calF$, so we get $\calF\leq \calG$. A symmetric argument gives $\calG\leq\calF$.
\end{proof}

\begin{proposition}\label{proposition:convex basis implies anti-symmetry}
	If $(X,\Leq)$ is convex, then the induced order on $\pt(X)$ is anti-symmetric.
\end{proposition}
\begin{proof}
	Using \cref{lemma:anti-symmetry and convex hulls points} we see it suffices to prove that $\hull{\calF}=\hull{\calG}$ implies $\calF=\calG$. Hence, suppose that $\hull{\calF}=\hull{\calG}$. Take any open $U\in\calF$. Since the convex opens form a base, we can write $U=\bigvee_{i\in I}\hull{W}_i$, for some $W_i\in\Opens X$. Since $\calF$ is completely prime, there exists some $i\in I$ such that $\hull{W}_i\in\calF$. Since the hulls of $\calF$ and $\calG$ are equal, we can find $V\in\calG$ such that $\hull{V}=\hull{W}_i$. Hence it follows that $V\sqleq \hull{W}_i\sqleq U$, and since $\calG$ is a filter this gives $U\in\calG$. This proves $\calF\subseteq \calG$, and a symmetric argument yields equality of the points.
\end{proof}

\begin{remark}
	We currently do not know of a fully localic analogue of the causality condition for spacetimes. The previous proposition suggests to choose the locales for which $\pt(X)$ is anti-symmetric, however this definition says little about non-spatial locales.
\end{remark}

Next we investigate how the convexity condition translates between spaces and locales. The following result shows that a convex locale relates to convexity of its underlying ordered space.

\begin{proposition}\label{convex locale iff pt convex}
	If $(X,\Leq)$ is convex, then ${(\pt(X),\leq)}$ is convex. If $X$ is spatial and satisfies~\eqref{axiom:P}, it is convex if and only if $(\pt(X),\leq)$ is convex.
\end{proposition}
\begin{proof}
	Let $\calF\in\pt(X)$, and pick an open neighbourhood $\pt(V)$ of $\calF$. Since the convex regions form a basis for $\Opens X$, we can write $V=\bigvee_{i\in I} \hull{U}_i$, for some family of opens $U_i\in\Opens X$. Since $\calF$ is completely prime, there exists some index $i\in I$ such that $\hull{U}_i\in\calF$. By \cref{lemma:convex iff pointwise convex} we get the open convex neighbourhood $\pt(\hull{U}_i)$ of $\calF$, as desired.
	
	Conversely, if $\pt(X)$ is convex and $X$ is spatial, it follows by \cref{lemma:convex iff pointwise convex} that~$X$ has a basis consisting of convex opens.
\end{proof}

Dually, we can start with a space and characterise its convexity in terms of its induced ordered locale. 

\begin{proposition}\label{proposition:convex space iff locale convex}
	If $(S,\leq)$ is convex, then $\loc(S)$ is convex. If $(S,\leq)$ has open cones, then it is convex if and only if $\loc(S)$ is convex.
\end{proposition}
\begin{proof}
	This follows directly by applying~\cref{lemma:convex pointwise iff localic convex}.
\end{proof}

\begin{corollary}\label{corollary:strongly causal iff convex}
	A smooth spacetime $M$ is strongly causal if and only if its ordered locale $\loc(M)$ is convex.
\end{corollary}

Interestingly, if we string \cref{convex locale iff pt convex,proposition:convex space iff locale convex} together, we see that for an ordered space with open cones and enough points:
\[
\parbox{8em}{\centering $(S,\leq)$\\convex}
\xLeftrightarrow{\ref{proposition:convex space iff locale convex}}
\parbox{8em}{\centering $(\loc(S),\Leq)$\\convex}
\xLeftrightarrow{\ref{convex locale iff pt convex}}
\parbox{8em}{\centering $(\pt(\loc(S)),\leq)$\\convex.}
\]
Dually, any spatial ordered locale $(X,\Leq)$ satisfying~\eqref{axiom:P} has:
\[
\parbox{8em}{\centering $(X,\Leq)$\\convex}
\xLeftrightarrow{\ref{convex locale iff pt convex}}
\parbox{8em}{\centering $(\pt(X),\leq)$\\convex}
\xLeftrightarrow{\ref{proposition:convex space iff locale convex}}
\parbox{8em}{\centering $(\loc(\pt(X)),\Leq)$\\convex.}
\]
The second chain of equivalences is unsurprising, since in that setting we know that the counit $\epsilon_X$ is actually an isomorphism of ordered locales (\cref{lemma:counit inverse monotone iff bullet}). The first chain of equivalences actually says something new, since it works even in the case that $S$ is not $T_0$-ordered, so when the unit $\eta_S$ does \emph{not} define an order isomorphism to $\pt(\loc(S))$. We can demonstrate this with an example.

\begin{example}
	In particular, these results tell us that a spacetime $M$ is strongly causal if and only if $\pt(\loc(M))$ is convex. For this example, let $M$ denote $1+1$ Minkowski space, and let $N= M\setminus \{p\}$ be the same space but with one point removed. Then $N$ is another spacetime. In particular it is sober, so there is a homeomorphism $N\cong \pt(\loc(N))$. However, since $N$ is not $T_0$-ordered (recall \cref{example:minkowski space point removed not T0}), this homeomorphism is \emph{not} an order isomorphism. We see that in $\pt(\loc(N))$ any two points on either of the null geodesics passing through $p$ will be causally related, while this is not the case in $N$. Hence the causal order of $\pt(\loc(N))$ is just that of $M$ restricted to $M\setminus\{p\}$, which is different from the causal structure of~$N$ itself. Nevertheless, their strong causalities are equivalent.
\end{example}


\section{Causal complements and diamonds}\label{section:causal complements}
In this section we study a point-free analogue of the notion of a ``causal complement'' operator. The notion of a causal complement has been defined in the context of spacetimes in \cite{casini2002LogicCausallyClosed,cegla1977CausalLogicMinkowski,cegla2005OrthoCausalClosure}.\footnote{We kindly thank Ryszard Paweł Kostecki for making us aware of these references.}

The idea is to assign to any region $U$ another region $U^\perp$ that is the largest one not in causal contact with the original. If we think of the region of causal influence of $U$ as $\Up U$ and $\Down U$, then in an ordered locale, we can make this more precise by saying $U^\perp$ is the largest open region that does not overlap $\Up U$ and $\Down U$:
\[
U^\perp := \bigvee \{V\in\Opens X: V\wedge (\Up U \vee \Down U)=\varnothing\}.
\]
We can use the inherent Heyting operator of the frame $\Opens X$ to write this more neatly as
\[
U^\perp = \neg\left( \Up U \vee \Down U \right) = \neg \Up U \wedge \neg \Down U,
\]
which resonates with our initial motivation, and as a bonus gives us access to the nice properties of Heyting negation \cite[Section III.3]{picado2012FramesLocalesTopology}; see \cref{lemma:properties of heyting negation}. For instance, from easy calculations using \cref{lemma:properties of localic cones} it will follow that this operator indeed behaves like a complement.

We need to know how the Heyting negation interacts with the localic cones. Just as the complement of a down-closed set is upwards-closed in a preorder, and dually, the following result shows that the same is true in a parallel ordered locale.
\begin{lemma}\label{lemma:heyting negation of past and future set}
	For every open $U$ in a parallel ordered locale:
	\[
	\Up\neg\Down U = \neg \Down U
	\qquad\text{and}\qquad
	\Down \neg\Up U = \neg\Up U.
	\]
\end{lemma}
\begin{proof}
	Since localic cones are monads, for the first equation it suffices to prove the inclusion $\Up \neg\Down U \sqleq \neg \Down U$. By definition of Heyting negation this is equivalent to $\Up \neg\Down U \wedge \Down U =\varnothing$. In turn, by parallel orderedness (\cref{proposition:cones are parallel}) this holds if and only if $\neg\Down U \wedge\Down\Down U = \neg\Down U \wedge \Down U = \varnothing$, which is true. The proof for the other equation is dual.
\end{proof}

\begin{definition}
	The \emph{causal complement} operator of an ordered locale $(X,\Leq)$ is defined as
	\[
	(-)^\perp\colon\Opens X\longrightarrow\Opens X; \qquad U \longmapsto U^\perp:= \neg \Up U \wedge \neg \Down U.
	\]
\end{definition}

\begin{lemma}\label{lemma:properties of causal complement}
	The causal complement operator of an ordered locale $(X,\Leq)$ satisfies the following properties:
	\begin{enumerate}
		\item[(a)] $U\wedge U^\perp = \varnothing$;
		\item[(b)] if $U\sqleq V$ then $V^\perp\sqleq U^\perp$;
		\item[(c)] $X^\perp = \varnothing$.
	\end{enumerate}
	If $X$ satisfies~\eqref{axiom:empty} we further have:
	\begin{enumerate}
		\item[(d)] $\varnothing^\perp = X$.
	\end{enumerate}
	If $X$ is parallel ordered we have:
	\begin{enumerate}
		\item[(e)] $U\sqleq U^{\perp\perp}$;
		\item[(f)] $U\sqleq V^\perp$ if and only if $V\sqleq U^\perp$.
	\end{enumerate}
	Lastly, if $X$ satisfies~\eqref{axiom:LV} we have:
	\begin{enumerate}
		\item[(g)] $\left(\bigvee U_i\right)^\perp = \bigwedge U_i^\perp$.
	\end{enumerate}
\end{lemma}
\begin{proof}
	We freely use properties of Heyting negation (\cref{lemma:properties of heyting negation}). For (a), use $U\sqleq \Up U\vee \Down U$ to calculate:
	\[
	U\wedge U^\perp = U \wedge \neg(\Up U \vee \Down U) \sqleq U\wedge \neg U =\varnothing.
	\]
	Property (b) follows immediately from monotonicity of the localic cones and the antimonotonicity of Heyting negation.
	
	For (c), it is easy to see that $\Up X = \Down X = X$, so that $X^\perp = \neg X = \varnothing$. Similarly, for (d), note~\eqref{axiom:empty} implies $\Up\varnothing = \Down \varnothing =\varnothing$, so that we get $\varnothing^\perp = \neg \varnothing = X$.
	
	For (e), we show that $U\wedge \Up U^\perp=\varnothing$. By \cref{proposition:cones are parallel}, this holds if and only if $\Down U \wedge U^\perp =\varnothing$. But this is true by definition. Similarly we get $U\wedge \Down U^\perp =\varnothing$, so we conclude $U\sqleq U^{\perp\perp}$, as desired. Property (f) follows directly by combining this with (b).
	
	Lastly, for (g), using the infinite first De Morgan law:
	\[
	\left(\bigvee U_i\right)^\perp = \neg\left( \Up \bigvee U_i \vee \Down \bigvee U_i\right)
	=\neg \bigvee (\Up U_i\vee \Down U_i)
	= \bigwedge \neg(\Up U_i\vee \Down U_i) = \bigwedge U_i^\perp,
	\]
	where the second equation follows using~\eqref{axiom:LV}.
\end{proof}

\begin{example}\label{example:complement of cones and hull}
	Take $U\in\Opens X$ in an ordered locale. Then:
	\[
	\left(\Up U\right)^\perp = \neg\Down\Up U
	\qquad\text{and}\qquad
	\left(\Down U\right)^\perp = \neg\Up\Down U.
	\]
	Using \cref{lemma:properties of convex hull}(d) we also find $ (\hull{U})^\perp = U^\perp$.
\end{example}

\subsection{Causal diamonds}
\label{section:causal diamonds}
In this section, we very briefly consider an abstract notion of \emph{diamond}, and see that it is in general distinct from the typical notion of diamonds in a spacetime (\cref{example:chronological diamonds are not localic diamonds}). A successful localic generalisation of diamonds would help study aspects from relativity theory that rely explicitly on chronological cones defined on points.

\begin{definition}\label{definition:causal diamonds}
	The \emph{causal diamond operator} of an ordered locale is defined as:
	\[
	(-)^\diamond\colon \Opens X\longrightarrow\Opens X;\qquad U\longmapsto U^\diamond := U^{\perp\perp}.
	\]
	A \emph{diamond} is an open $U\in\Opens X$ such that $U= U^\diamond$.
\end{definition}

\begin{lemma}\label{lemma:properties of diamonds}
	In a parallel ordered locale $(X,\Leq)$ the diamond operator has the following properties:
	\begin{enumerate}
		\item[(a)] $U\sqleq U^\diamond$;
		\item[(b)] $U^{\diamond\diamond}\sqleq U^\diamond$;
		\item[(c)] If $U\sqleq V$ then $U^\diamond \sqleq V^\diamond$.
	\end{enumerate}
\end{lemma}
\begin{proof}
	This follows immediately from the properties of the causal complement in \cref{lemma:properties of causal complement}.
\end{proof}

\begin{lemma}\label{lemma:hulls contained in diamonds}
	In a parallel ordered locale we have that $\hull{U}\sqleq U^\diamond$.
\end{lemma}
\begin{proof}
	It suffices to show that $\hull{U}\wedge \Up U^\perp = \varnothing$, and similarly for downsets. Using \cref{proposition:cones are parallel} this holds if and only if $\Down \hull{U}\wedge U^\perp =\varnothing$, which in turn can be unpacked as follows:
	\[
	\Down \hull{U} \wedge U^\perp = \Down\left( \Up U \wedge \Down U\right) \wedge U^\perp \sqleq \Down\Up U \wedge \Down U \wedge U^\perp.
	\]
	But $\Down U\wedge U^\perp =\varnothing$ by construction, so the result follows.
\end{proof}

\begin{lemma}\label{lemma:diamonds are convex}
	In a parallel ordered locale diamonds are convex.
\end{lemma}
\begin{proof}
	From \cref{lemma:hulls contained in diamonds,lemma:properties of diamonds} we get $\hull{U^\diamond}\sqleq U^{\diamond\diamond}\sqleq U^\diamond$, so the result follows in combination with \cref{lemma:properties of convex hull}(a).
\end{proof}

Hence the diamonds form a special class of convex opens, but it is clear from examples that not every convex region is a diamond. For instance, in Minkowski space the diamond induced by $\Up U$ for non-empty $U$ will be the entire space, while $\Up U$ is generally much smaller. The next example shows that the localic diamonds are generally ``too big.''


\begin{example}\label{example:chronological diamonds are not localic diamonds}
	In a smooth spacetime $M$, the \emph{chronological diamonds} are sets of the form $\langle x,y\rangle := I^+(x)\cap I^-(y)$. These are \emph{not} generally diamonds in the sense of \cref{definition:causal diamonds}. Consider for instance the open subset of Minkowski space $\mathbb{R}^{2}$ defined by $(-1,1)\times \mathbb{R}$. Take $x=(0,-1)$ and $y= (0,1)$, so that $x\chron y$. Their chronological diamond $\langle x,y\rangle$ is the interior of the square in $\mathbb{R}$ induced by the vertices $(1,0)$, $(0,1)$, $(-1,0)$ and $(0,-1)$. Therefore $\langle x,y\rangle^\perp = \varnothing$, and so by \cref{lemma:properties of causal complement} we find $\langle x,y\rangle^\diamond= M\neq\langle x,y\rangle$.	
\end{example}\vspace*{-2ex}

\section{Locales of pasts and futures}\label{section:locale of pasts and futures}
In our discussion on biframes (\cref{section:ordered locales from biframes}) we have already seen that if an ordered locale $(X,\Leq)$ has localic cones that preserve all joins, then $\im(\Up)$ and $\im(\Down)$ are subframes of $\Opens X$. Instead of viewing this (bi)frame-theoretically, we can interpret it localically.

\begin{definition}\label{definition:locale of pasts and futures}
	Let $(X,\Leq)$ be an ordered locale with~\eqref{axiom:LV}. Its \emph{locale of futures} is the locale $X^\triup$ defined by the frame $\Opens X^\triup := \im(\Up)$. Dually, the \emph{locale of pasts} is the locale $X^\tridown$ defined by $\Opens X^\tridown:=\im(\Down)$. They come with canonical surjective locale maps:
		\[
			\begin{tikzcd}[ampersand replacement=\&]
				X \& {X^\triup}
				\arrow["{\eta_\triup}", two heads, from=1-1, to=1-2]
			\end{tikzcd}
			\qquad\text{and}\qquad
			\begin{tikzcd}[ampersand replacement=\&]
				X \& {X^\tridown.}
				\arrow["{\eta_\tridown}", two heads, from=1-1, to=1-2]
			\end{tikzcd}
		\]
	Explicitly, the preimage maps $\eta_\triup^{-1}\colon \im(\Up)\hookrightarrow \Opens X$ and $\eta_\tridown^{-1}\colon \im(\Down)\hookrightarrow \Opens X$ are just the subframe inclusions.
\end{definition}

\begin{example}
	For the real line $(\mathbb{R},\leq)$, the induced locales $\mathbb{R}^\triup$ and $\mathbb{R}^\tridown$ are just the locales of \emph{lower reals} and \emph{upper reals}, respectively (note the switch in ``direction''). A function $g\colon S\to \mathbb{R}\cup \{-\infty,\infty\}$ taking values in the extended real numbers is \emph{lower semicontinuous} precisely when $\loc(g)$ defines a map of locales $\loc(S)\to \mathbb{R}^\triup$ to the \emph{lower} reals. Dually for upper semicontinuity. See e.g.~\cite[\S 2.5]{vickers2007LocalesToposesSpaces} for more on localic reals.
\end{example}

Let us unpack this structure further. The maps of locales $\eta_\triup$ and $\eta_\tridown$ in particular induce direct image maps. By \cref{theorem:join preserving implies left adjoint} we can explicitly calculate
	\[
		\eta^\triup_\ast\colon \Opens X\longrightarrow \im(\Up);
		\qquad
		U\longmapsto \bigvee\left\{\Up F\in\im(\Up): \Up F\sqleq U\right\}.
	\]
Here, to streamline notation, we have written $(\eta_\triup)_\ast = \eta^\triup_\ast$, and shall also write $(\eta_\tridown)_\ast = \eta^\tridown_\ast$. Therefore $\eta^\triup_\ast(U)$ is thought of as the largest upclosed subregion of $U$. Of course, often $\eta^\triup_\ast(U)$ will be empty, but clearly $\eta^\triup_\ast(\Up F) = \Up F$. To understand the behaviour of these direct images, it helps to evaluate them on the points of the locale. Take a completely prime filter $\calF\in \pt(X)$, with corresponding prime element $P_\calF\in\Opens X$ (recall \cref{lemma:prime elements are cpf}). Then:
	\[
		\eta^\triup_\ast(P_\calF) = \bigvee\left\{\Up F\in\im(\Up): \Up F\sqleq P_\calF\right\}
		=
		\bigvee\left\{\Up F\in\im(\Up): \Up F\notin\calF\right\}.
	\]
In an ordered space $S$, with $\calF = \calF_x$ for some $x\in S$, this becomes:
	\[
		\eta^\triup_\ast\left(S\setminus \overline{\{x\}}\right) = 
		\bigcup\left\{\Up F\in\im(\Up): x\notin \Up F\right\}.
	\]
For visual intuition, see \cref{figure:eta on prime}. This is a glimpse of a localic \emph{generalised ideal point}, which we discuss further in \cref{section:generalised ideal points}.

\begin{figure}[t]\centering
	\input{tikz/inkscape-export/fig-07}
	\caption{Illustration of how $\eta_\ast^\triup$ acts on  primes.}
	\label{figure:eta on prime}
\end{figure}

As we know from \cref{section:ordered locales from biframes}, the subframe inclusions $\im(\Up),\im(\Down)\hookrightarrow \Opens X$, which here are just the preimage maps $\eta_\triup^{-1}$ and $\eta_\tridown^{-1}$, both admit left adjoints. They are precisely the localic cones. We thus have adjoint triples:
\[
\begin{tikzcd}[ampersand replacement=\&]
	{\Opens X} \&\& {\Opens X^\triup}
	\arrow[""{name=0, anchor=center, inner sep=0}, "{\eta_\triup^{-1}}"{description, pos=0.3}, hook', from=1-3, to=1-1]
	\arrow[""{name=1, anchor=center, inner sep=0}, "\Up", shift left=4, from=1-1, to=1-3]
	\arrow[""{name=2, anchor=center, inner sep=0}, "{\eta^\triup_\ast}"', shift right=4, from=1-1, to=1-3]
	\arrow["\dashv"{anchor=center, rotate=-90}, draw=none, from=1, to=0]
	\arrow["\dashv"{anchor=center, rotate=-90}, draw=none, from=0, to=2]
\end{tikzcd}
\quad\text{and}\quad
\begin{tikzcd}[ampersand replacement=\&]
	{\Opens X} \&\& {\Opens X^\tridown.}
	\arrow[""{name=0, anchor=center, inner sep=0}, "{\eta_\tridown^{-1}}"{description, pos=0.3}, hook', from=1-3, to=1-1]
	\arrow[""{name=1, anchor=center, inner sep=0}, "\Down", shift left=4, from=1-1, to=1-3]
	\arrow[""{name=2, anchor=center, inner sep=0}, "{\eta^\tridown_\ast}"', shift right=4, from=1-1, to=1-3]
	\arrow["\dashv"{anchor=center, rotate=-90}, draw=none, from=1, to=0]
	\arrow["\dashv"{anchor=center, rotate=-90}, draw=none, from=0, to=2]
\end{tikzcd}
\]

\subsection{Convexity as an internal condition}
In the following we provide an alternative characterisation of convexity in terms of the induced future and past locales. In general, the product $X\times Y$ of locales is described frame-theoretically by the coproduct $\Opens(X\times Y)=\Opens X\oplus \Opens Y$. It is somewhat non-trivial to describe this frame, but we present more details in \cref{section:product locales}, and refer to \cite[\S IV.4]{picado2012FramesLocalesTopology} for a full description. Briefly, this frame is generated by pairs $U\oplus V$, where $U\in \Opens X$ and $V\in \Opens Y$, subject to a suitable equivalence relation. We think of $U\oplus V$ as the basic open ``boxes'' $U\times V$ that are used to define the product topology.

Given an ordered locale $(X,\Leq)$ with~\eqref{axiom:LV}, we obtain the canonical maps $\eta_\triup\colon X\tworightarrow X^\triup$ and $\eta_\tridown\colon X\tworightarrow X^\tridown$, so by universality of products this induces a new map:
	\[
	\begin{tikzcd}[ampersand replacement=\&]
		\& X \\
		{X^\triup} \& {X^\triup\times X^\tridown} \& {X^\tridown.}
		\arrow["{\eta_\triup}"', two heads, from=1-2, to=2-1]
		\arrow["{\eta_\tridown}", two heads, from=1-2, to=2-3]
		\arrow["{\pr_1}", from=2-2, to=2-1]
		\arrow["{\pr_2}"', from=2-2, to=2-3]
		\arrow["{\exists !\eta}", dashed, from=1-2, to=2-2]
	\end{tikzcd}
	\]
Here the product locale $X^\triup\times X^\tridown$ is defined by the frame $\im(\Up)\oplus \im(\Down)$. On basic elements the corresponding frame map is calculated as (see \cref{corollary:pair map preimage on basic opens}):
\[
\eta^{-1}\colon\im(\Up)\oplus \im(\Down)\longrightarrow \Opens X;
\qquad 
\Up U\oplus \Down V\longmapsto \Up U\wedge \Down V.
\]
Hence this map sends the basis elements of $\im(\Up)\oplus \im(\Down)$ precisely to the convex opens of $X$, and every convex open can be written in this way. We therefore have the following observation. 

\begin{proposition}\label{proposition:convex iff sublocale}
	An ordered locale $X$ with~\eqref{axiom:LV} is convex if and only if the map $\eta\colon X\to X^\triup\times X^\tridown$ defines a sublocale.
\end{proposition}
\begin{proof}
	That $\eta$ is a sublocale means that the frame map $\eta^{-1}$ is a surjection (\cref{proposition:extremal monomorphisms in Loc}), i.e.~every $W\in \Opens X$ can be written as
		\[
			\eta^{-1}\left(\bigvee (\Up U_i\oplus \Down V_i)\right) = \bigvee\eta^{-1}(\Up U_i\oplus \Down V_i) = \bigvee (\Up U_i\wedge \Down V_i).
		\]
	which is a join of convex opens.
\end{proof}

Combining with \cref{corollary:strongly causal iff convex}, we thus obtain a localic, structural characterisation of strong causality.

\begin{theorem}\label{theorem:strongly causal iff sublocale}
	A spacetime $M$ is strongly causal if and only if the canonical map $M\to M^\triup\times M^\tridown$ defines a sublocale.
\end{theorem}

\chapter{Towards localic causal boundaries}\label{section:causal boundaries}
In this section we work towards a \emph{causal boundary} construction for ordered locales. These ideas are inspired by the point-wise analogously defined causal boundaries of spacetimes. Our central motivation is the seminal paper \cite{geroch1972IdealPointsSpaceTime}, which contains one of the first attempts at such a construction. Even though this project was initiated over fifty years ago, the precise definition of the assignment $M\mapsto \overline{M}$ of a spacetime to its causal boundary has only recently been somewhat settled, see e.g.~\cite{budic1974CausalBoundariesGeneral,szabados1988CausalBoundaryStrongly,harris1998UniversalityFutureChronological,harris2000TopologyFutureChronological,marolf2003NewRecipeCausala, garcia-parrado2005CausalStructuresCausala,flores2007CausalBoundarySpacetimes,sanchez2009CausalBoundariesHolography,herrera2010CausalBoundaryIts,flores2011FinalDefinitionCausal}.

Common to many of these approaches is the notion of an \emph{ideal point}, first introduced in \cite{geroch1972IdealPointsSpaceTime}. As we highlight, these can be characterised purely using lattice-theoretic techniques, and therefore lends to point-free generalisation. First, ideal points can be characterised point-wise as \emph{ideals} with respect to the chronology relation (\cref{section:chronological ideals}). Second, ideal points can be characterised as the \emph{coprime elements} of the distributive lattices of future and past sets (\cref{example:localic IP are spacetime IP}). Interpreted localically, we prove that ideal points can be recovered as the \emph{points} of the \emph{locales of futures} and \emph{pasts} (\cref{corollary:IPs as primes in spacetime}). A key step in this proof is \cref{lemma:chronological cones are regular opens}, which says that the chronological cones are always \emph{regular opens}. Together with \cref{lemma:chronological cone of double negation} this has the remarkable consequence that the ideal points of a spacetime can be recovered from a strictly \emph{pointless} locale (\cref{corollary:IPs from regular opens}).

%
%

\section{Ideal points in spacetime}\label{section:ideal points}
Ever since Penrose's influential and celebrated singularity theorem \cite{penrose1965GravitationalCollapseSpaceTime}, it has been understood that spacetimes generically exhibit \emph{singularities} (as opposed to only arising from some exceptional class of solutions to the Einstein equations). From this point on, the study of singularities became an interesting one from the perspective of Lorentzian geometry, as opposed to only for physics. The precise history towards formal definition of ``singularity'' is too complex to be outlined here, so we refer to the overview in \cite[Chapter~6]{landsman2021FoundationsGeneralRelativity}.

For our purposes, the basic intuition is that a singular point in spacetime $M$ is modelled by a timelike curve $\gamma$ that does not have a \emph{limit}. The would-be limit of $\gamma$ represents a ``hole'' in the manifold $M$ that we interpret as a singularity. The basic function of the causal boundary $\overline{M}$ is then to formally add in these singular points. Most straightforwardly, we take the timelike curves $\gamma$ as elements in $\overline{M}$ directly to represent their would-be limit points. In particular, for any genuine point $x\in M$ in the spacetime, we can take a timelike curve $\gamma$ converging to $x$ to define the embedding $M\to \overline{M}\colon x\mapsto \gamma$. From this it is clear that we need to make some sort of identification between the curves $\gamma$, since there are generally many distinct curves converging to the same point. One elegant way to achieve this is by taking instead the set $I^-[\gamma]$, defined as the chronological past of the curve $\im(\gamma)\subseteq M$. All timelike curves $\gamma$ converging to $x$ will then determine the same cone $I^-(x)$. Thus, we shall model the \emph{ideal points} of $M$ precisely by the sets $I^-[\gamma]$, where $\gamma$ is any timelike curve.

It was the paper \cite{geroch1972IdealPointsSpaceTime} by Geroch, Kronheimer and Penrose that first introduced this method of constructing causal boundaries, and is the one still most widely used and studied today (see also \cite{garcia-parrado2005CausalStructuresCausala} for an overview). It also happens to describe a method that is most readily generalisable to our point-free language, as we now highlight.

\begin{definition}
	An \emph{indecomposable past set} (or \emph{IP} for short) is a subset $P\subseteq M$ of spacetime such that $P\neq\varnothing$ and $P= I^-(P)$, that is indecomposable in the following sense:
		\[
			P = I^-(A)\cup I^-(B)
			\qquad\text{implies}\qquad
			P = I^-(A)\text{~or~}I^-(B).
		\]
	Dually, an \emph{indecomposable future set} (or \emph{IF} for short) is a subset $F\subseteq M$ such that $F\neq\varnothing$ and $F=I^+(F)$, that is indecomposable in the following sense:
		\[
			F = I^+(A) \cup I^+(B)
			\qquad\text{implies}\qquad
			F = I^+(A) \text{~or~}I^+(B).
		\]
	The set of IPs and IFs of $M$ are denoted $\IP(M)$ and $\IF(M)$, respectively.
\end{definition}

The following remarkable theorem relates the purely lattice-theoretic IPs and IFs to the timelike curves in the spacetime.

\begin{theorem}[{\cite[Theorem~2.1]{geroch1972IdealPointsSpaceTime}}]\label{theorem:GKP Theorem 2.1}
	A subset of a spacetime is of the form~$I^-[\gamma]$ for some timelike curve $\gamma$, if and only if it is an indecomposable past-set. 
\end{theorem}

\begin{example}
	For any $x\in M$ we get $I^-(x)\in\IP(M)$ and $I^+(x)\in \IF(M)$. More explicitly, in two-dimensional Minkowski space $M$ we can take the curve $t\mapsto (0,t)$ that extends towards the future indefinitely, and this gives $M\in \IP(M)$. This ideal point is called \emph{future timelike infinity}. By \cite[\S 2.2]{geroch1972IdealPointsSpaceTime} the set $I^-[\gamma]$ defines an IP even if $\gamma$ is causal, so also the lightlike curve $t\mapsto (t,t)$ defines an IP, which is the open diagonal plane consisting of all points strictly below the diagonal. By defining lightlike curves $t\mapsto (t,t+v)$, for arbitrary $v\in\mathbb{R}$, we obtain infinitely many distinct IPs that do not come from elements of $M$. These types of ideal points are called points at \emph{future lightlike infinity}.
\end{example}

\subsection{Indecomposable past and future regions in locales}\label{section:ideal points in locales}
Note that the definition of IPs and IFs is \emph{purely lattice-theoretic}. Note also that IPs and IFs are necessarily \emph{open}, so everything takes place in the lattice~$\Opens M$. Moreover, since the chronological cones $I^\pm$ are just the localic cones (\cref{example:localic cones in spacetime}), the entire definition is stated in terms of the data of $(\loc(M),\Leq)$. It therefore now straightforward to generalise this to arbitrary ordered locales.

\begin{definition}\label{definition:IPs and IFs in locales}
	Let $(X,\Leq)$ be an ordered locale. An \emph{indecomposable past region} (or \emph{IP} for short) is an element $P\in \Opens X$ such that $\varnothing \neq P = \Down P$, and
		\[
			P = \Down U \vee \Down V
			\qquad\text{implies}\qquad
			P = \Down U \text{~or~} \Down V.
		\]
	An \emph{indecomposable future region} (or \emph{IF} for short) is an element $F\in\Opens X$ such that $\varnothing \neq F=\Up F$, and
		\[
			F= \Up U \vee \Up V
			\qquad\text{implies}\qquad
			F = \Up U \text{~or~} \Up V.
		\]
	We denote the set of IPs and IFs of $X$ by $\IP(X)$ and $\IF(X)$, respectively. Note these are subsets of $\Opens X$. 
\end{definition}


\begin{remarknumbered}\label{remark:coprimes}
	A \emph{coprime} in a distributive lattice $L$ is a prime element (\cref{definition:prime element}) in the opposite lattice $L^\op$ (\cref{definition:opposite category}). Explicitly: it is a non-zero element $d\in L$ such that
		\[
			d= x\vee y
			\qquad\text{implies}\qquad
			d=x \text{~or~} y.
		\]
	Therefore $\IP(X)$ and $\IF(X)$ are exactly the sets of coprimes in $\im(\Down)$ and $\im(\Up)$, respectively. See e.g.~\cite[\S 9.3]{vickers1989TopologyLogic} for how coprimes are used in locale theory.
\end{remarknumbered}

\begin{example}\label{example:localic IP are spacetime IP}
	Of course, per the discussion above, it is clear that for spacetimes the localic definition agrees completely with the one from \cite{geroch1972IdealPointsSpaceTime}:
		\[
			\IP(\loc(M))= \IP(M)
			\qquad\text{and}\qquad
			\IF(\loc(M))= \IF(M).
		\] 
\end{example}

\begin{example}
	On the other hand, it is easy to find examples of ordered locales that have no indecomposable past regions. Take, for instance, any locale $X$ equipped with equality as the causal order. In that case the localic cones are the identity functions, so $P\in\IP(X)$ is just a coprime element of $\Opens X$. These will often not exist: if $S$ is a manifold then $\Opens S$ has no coprimes whatsoever (cf.~the proof of \cref{proposition:points in regular opens Hausdorff space} below).
\end{example}

\begin{figure}[b]\centering
	\input{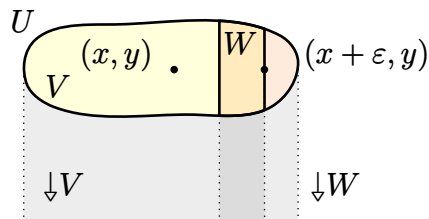}
	\caption{Vertical-$\mathbb{R}^2$ has no IPs.}
	\label{figure:vertical R2 has no IPs}
\end{figure}

\begin{example}\label{example:space without IPs}
	Consider again vertical-$\mathbb{R}^2$ from \cref{example:vertical R2}, seen as an ordered locale. We claim that if $U\in\Opens \mathbb{R}^2$ is not empty, then $\Down U$ cannot be an IP. If $U$ is non-empty we can find some element $(x,y)\in U$. Since $U$ is open there exists some positive real $\epsilon>0$ such that $(x+\epsilon,y)\in U$. Define two opens $V,W\in\Opens \mathbb{R}^2$, the first containing all elements of $U$ whose first component is strictly smaller than $x+\epsilon$, and the second containing all elements of $U$ whose first component is strictly larger than $x+\epsilon/2$, see \cref{figure:vertical R2 has no IPs}. Then $U=V\cup W$. Moreover, by axiom~\eqref{axiom:LV} we find $\Down U = \Down V \cup \Down W$, and by construction the downsets $\Down V$ and $\Down W$ are unequal, proving that $\Down U$ is not an IP. This proves there are (even spatial) ordered locales that have no IPs. Contrast this to spacetimes, where every point induces a principal IP. However, we will see in \cref{example:generalised IPs in vertical R2} that vertical-$\mathbb{R}^2$ still admits \emph{generalised} ideal points.
\end{example}

\subsection{Directed completeness of IPs and IFs}
What lattice structure do the IPs inherit from $\Opens X$? Clearly $\IP(X)\subseteq \Opens X$ is not generally closed under joins with respect to the inclusion order $\sqleq$. In the following we show that $\IP(X)$ is instead \emph{directed} complete. These ideas generalise ones found in \cite{krolak1984MaximalIndecomposableSets}. Recall the definition of a directed set.

\begin{definition}\label{definition:directed set}
	A \emph{directed set} is a non-empty set $D$ equipped with a preorder~$\leq$ such that for every $x,y\in D$ there exists $z\in D$ with $x,y\leq z$. A \emph{directed join} is a join taken over a directed indexing set.
	
	A \emph{directed-complete partial order} (dcpo) is a partially ordered $(P,\leq)$ set where all directed joins exist, i.e.~where every directed subset with respect to $\leq$ has a join in $P$.
\end{definition}


\begin{lemma}\label{lemma:directed set is indecomposable}
	If $D$ is a directed set and $J,K\subseteq D$ are downsets, then $D= J\cup K$ implies $D= J$ or $D=K$.
\end{lemma}
\begin{proof}
	Suppose for the sake of contradiction that $J\not\subseteq K$ and $K\not\subseteq J$. In that case we can find $j\in J\setminus K$ and $k\in K\setminus J$, and since $J$ and $K$ are downsets, we get in particular that $j\notin\down K$ and $k\notin\down J$. But using directedness there exists $x\in D$ such that $j,k\leq x\in J\cup K$, which contradicts this.
\end{proof}

The following is the point-free analogue of \cite[Lemma~1]{krolak1984MaximalIndecomposableSets}. We also generalise their claim from chains to arbitrary directed families.


\begin{lemma}\label{lemma:directed family of IPs has join}
	In an ordered locale $(X,\Leq)$ satisfying~\eqref{axiom:LV}, every directed set in $\IP(X)$ has a join (with respect to the inclusion order).
\end{lemma}
\begin{proof}
	Take a directed family $(P_i)_{i\in D}$ of IPs. We claim that their join $P=\bigvee_{i\in D}P_i$ calculated in $\Opens X$ is also their join in $\IP(X)$. First, it is clear that $P$ is non-empty as long as $D$ is non-empty, and axiom~\eqref{axiom:LV} guarantees $P = \Down P$. Now suppose that $P= \Down U\vee \Down V$ for some opens $U,V\in\Opens X$. Define
		\[
			J = \{j\in D: P_j\sqleq \Down U\}
			\qquad\text{and}\qquad
			K = \{k\in D: P_k\sqleq \Down V\}.
		\]
	The sets $J$ and $K$ are clearly downsets, and since each individial $P_i$ is an IP, we get $D= J\cup K$. Therefore by \cref{lemma:directed set is indecomposable} we get $D= J$ or $K$. This gives, respectively, that $P\sqleq \Down U$ or $P\sqleq \Down V$, proving indecomposability.
\end{proof}

Assuming the Axiom of Choice, we get the following.
\begin{corollary}
	If $\IP(X)$ is non-empty, then it has a maximal element.
\end{corollary}
\begin{proof}
	This follows by the Kuratowski–Zorn lemma.
\end{proof}


\section{Chronological ideals}\label{section:chronological ideals}
In this section we discuss a point-wise order-theoretic characterisation of IPs. The upshot is that they are precisely the ``ideals'' with respect to the chronology relation $\chron$, and this allows us to think of ideal points as a ``directed completion'' of a spacetime. Ideals are usually only defined in a setting as general as that of preorders, see e.g. \cite[Section~VII.2.1]{johnstone1982StoneSpacesa} and \cite[Definition~O-1.3]{gierz2003ContinuousLatticesDomains}. Since the chronology is not a preorder, we therefore need to define a more general notion of ideal with respect to an arbitrary relation.

For this section we work with a set $S$ equipped with a relation $\aleq$. We denote the downset of $A\subseteq S$ under the relation $\aleq$ by $\down A :=\{ x\in S:\exists a\in A: x\aleq a\}$. All dual definitions are left implicit.

\begin{definition}
	Consider a set $S$ equipped with a relation $\aleq$. An \emph{$\aleq$-ideal} is a subset $I\subseteq S$ satisfying the following conditions:
	\begin{enumerate}[label = (J\arabic*)]\setcounter{enumi}{-1}
		\item\label{axiom:J0} $I$ is non-empty;
		\item\label{axiom:J1} if $x,y\in I$, there exists $z\in I$ with $x,y\aleq z$; \hfill(\emph{upwards directedness})
		\item\label{axiom:J2} if $y\aleq x\in I$ then $y\in I$. \hfill(\emph{down closure})
	\end{enumerate}
	Denote by $\Ideals_{\aleq}(S)$ the set of $\aleq$-ideals in $S$, naturally equipped with the subset inclusion order.
\end{definition}

Unlike in a preorder, it is not guaranteed that ideals always exist. For instance it may happen that the downsets $\down x$ are not always $\aleq$-ideals. In the following we give precise conditions on $\aleq$ that makes this the case. This generalises the notion of \emph{past semi-fullness} of causal spaces from \cite[\S 3]{geroch1972IdealPointsSpaceTime}, which in turn generalises the definition of \emph{fullness} from \cite{kronheimer1967StructureCausalSpaces}.

\begin{definition}\label{definition:past semifull}
	A relation $\aleq$ on a set $S$ is called \emph{past semi-full} if:
	\begin{enumerate}[label = (\roman*)]
		\item for all $x\in S$ there exists $y\in S$ such that $y\aleq x$;
		\item if $y_1,y_2\aleq x$, there exists $z\in S$ such that $y_1,y_2\aleq z\aleq x$.
	\end{enumerate}
\end{definition}

\begin{remark}
	This property is also called more neutrally the \emph{finite interpolative property} in \cite[Definition~III-4.15]{gierz2003ContinuousLatticesDomains}. They call a (non-empty) set with a transitive and finitely interpolative relation $\aleq$ an \emph{abstract basis}, and our $\aleq$-ideals are called \emph{rounded ideals}.
\end{remark}

\begin{example}
	The strict order relation $<$ on $\mathbb{Q}$ is past semi-full. More generally, any \emph{dense} partial order similarly induces a past semi-full relation.
\end{example}

\begin{example}
	The chronology relation $\chron$ of any smooth spacetime is past semi-full, as discussed in~\cite{kronheimer1967StructureCausalSpaces}. But note that their causal spaces in general are not assumed to be past semi-full. 
\end{example}


\begin{lemma}\label{lemma:transitive and psf iff principal ideals exist}
	A relation $\aleq$ on $S$ is transitive and past semi-full if and only if $\down x$ is an $\aleq$-ideal for every $x\in S$.
\end{lemma}
\begin{proof}
	We see that the past semi-fullness of $\aleq$ corresponds precisely to the axioms~\ref{axiom:J0} and~\ref{axiom:J1} for the set $\down x$. Similarly, transitivity of $\aleq$ corresponds precisely to axiom~\ref{axiom:J2} for $\down x$.
\end{proof}

Another way to think about an ``ideal point'' is that it is the would-be \emph{join} of a directed family of points in the order. This is the appropriate order-theoretic generalisation of the ``limits of curves'' idea. In general $(S,\aleq)$ will not permit joins of directed subsets. Similar to how we think of $I^-[\gamma]$ as a representation of the would-be limit of $\gamma$, here we think of an $\aleq$-ideal $I$ as the representation of the would-be directed join $\bigvee I$. Thus, what we are looking for is a \emph{directed join completion} of $(S,\aleq)$. Here we prove that the set $\Ideals_{\aleq}(S)$ with the subset inclusion order provides this completion.

To do this, we show that $\Ideals_{\aleq}(S)$ has similar properties to the ideal completion of posets, discussed in e.g.~\cite[\S 9.1]{vickers1989TopologyLogic}. The following results are analogous to their Proposition~9.1.2. Recall the notion of directed join from \cref{definition:directed set}.

%

\begin{proposition}
	For any relation $\aleq$ on $S$, the poset $\Ideals_{\aleq}(S)$ has directed joins.
\end{proposition}
\begin{proof}
	If $(I_d)_{d\in D}$ is a directed family of $\aleq$-ideals in $S$, then clearly their union $\bigcup_{d\in D}I_d$ is non-empty and downwards closed. If $x,y\in \bigcup_{d\in D}I_d$ then by directedness there exists $a\in D$ such that $x,y\in I_a$, and so by~\ref{axiom:J1} there exists $z\in I_a$ with $x,y\aleq z$. Hence the union is a $\aleq$-ideal, and hence must be their join.
\end{proof}

The following follows immediately from \cref{lemma:transitive and psf iff principal ideals exist}.
\begin{lemma}
	If $\aleq$ is transitive and past semi-full on $S$, then there is a monotone function
		\[
			S\longrightarrow \Ideals_{\aleq}(S);
			\qquad x\longmapsto \down x.
		\]
\end{lemma}

\begin{lemma}\label{lemma:ideals are joins of principal ideals}
	Let $\aleq$ be a transitive past semi-full relation on $S$. Then any $\aleq$-ideal $I\in\Ideals_{\aleq}(S)$ can be written as the directed join
	\[
	I = \bigcup_{x\in I} \down x.
	\]
\end{lemma}
\begin{proof}
	By downwards closure~\ref{axiom:J2} it is clear that the right-hand side is contained in $I$. Conversely, let $x\in I$. By~\ref{axiom:J1} there exists $z\in I$ such that $x\aleq z$, so $x\in \down z$ is contained in the right-hand side. That this join is directed follows since~$I$ satisfies~\ref{axiom:J1}.
\end{proof}

%

A function $f\colon P\to Q$ between posets is called \emph{Scott continuous} if it preserves all directed joins. The following is an elementary generalisation of \cite[Proposition~9.1.2]{vickers1989TopologyLogic}, and says that $\Ideals_{\aleq}(S)$ is in some sense the nicest directed join completion of $S$.

\begin{proposition}\label{proposition:dcpo extension universality}
	Let $(S,\aleq)$ be a set with transitive and past semi-full relation, together with a monotone function $f\colon S\to Q$ into a directed complete poset $Q$ such that $f(x) = \bigvee\{f(z):z\aleq x\}$. Then there exists a unique Scott continuous extension:
		\[
			\begin{tikzcd}[ampersand replacement=\&,column sep=2.25em,row sep=scriptsize]
				S \\
				{\Ideals_{\aleq}(S)} \& {Q.}
				\arrow[from=1-1, to=2-1]
				\arrow["f", from=1-1, to=2-2]
				\arrow["{\exists!\overline{f}}"', dashed, from=2-1, to=2-2]
			\end{tikzcd}
		\]
\end{proposition}
\begin{proof}
	Anticipating Scott continuity and using \cref{lemma:ideals are joins of principal ideals}, we would find:
		\[
			\overline{f}(I) = \bigvee_{x\in I} \overline{f}(\down x) = \bigvee_{x\in I}f(x),
		\]
	so we can define $\overline{f}$ uniquely according to the right-hand side. This is a directed join in $Q$ since $I$ satisfies~\ref{axiom:J1} and $f$ is monotone. With this definition and the hypothesis on $f$ we obtain:
		\[
			\overline{f}(\down x) = \bigvee_{z\mspace{1mu}\aleq\mspace{2mu}x} f(z) = f(x),
		\]
	giving commutativity of the diagram. We are left to show $\overline{f}$ is Scott continuous. For that, take a directed family $(I_d)_{d\in D}$ in $\Ideals_{\aleq}(S)$. Then
		\[
			\overline{f}\left(\bigcup_{d\in D}I_d\right) 
			=
			\bigvee\left\{f(x): x\in \bigcup_{d\in D} I_d\right\}
			=
			\bigvee_{d\in D} \overline{f}(I_d),
		\]
	as desired.
\end{proof}

\begin{remark}
	The condition that $f(x)= \bigvee\{f(z):z\aleq x\}$ is not necessary in the setting of posets, since in that case reflexivity gives $x\in \down x$. Here we need it to make the diagram commute. Note that the embedding $S\to \Ideals_{\aleq}(S)$ satisfies this condition, since $\down x = \bigcup_{z\mspace{1mu}\aleq\mspace{2mu}x} \down z$ is a special case of \cref{lemma:ideals are joins of principal ideals}.
\end{remark}


\subsection{The $\chron$-ideals are the indecomposable past sets}
In the following we show that IPs in a spacetime are really nothing but the \mbox{$\chron$-ideals} (read as:~``chronological ideals''). Most of the work has already been done in \cite{geroch1972IdealPointsSpaceTime}, which contains the following result (cf.~\cref{lemma:directed set is indecomposable}).

\begin{lemma}[{\cite[Lemma~3.3]{geroch1972IdealPointsSpaceTime}}]\label{lemma:IP iff finite chronological join}
	A subset $\varnothing\neq P = I^-(P)$ of spacetime is an IP if and only if for every $x,y\in P$ we have $I^+(x)\cap I^+(y)\cap P\neq\varnothing$.
\end{lemma}
\begin{proof}
	We prove the ``if'' direction by contrapositive. If $P$ is not an IP then we can find subsets $A$ and $B$ of spacetime such that $P = I^-(A)\cup I^-(B)$ but $I^-(A)\not\subseteq I^-(B)$ and $I^-(B)\not\subseteq I^-(A)$. In turn, we can find $x\in I^-(A)\setminus I^-(B)$ and $y\in I^-(B)\setminus I^-(A)$. Suppose for the sake of contradiction that there exists $z\in I^+(x)\cap I^+(y)\cap P$. Without loss of generality, $z\in I^-(A)$. But since $y\chron z$ this implies $y\in I^-(A)$, a contradiction. Hence $I^+(x)\cap I^+(y)\cap P=\varnothing$.
	
	The converse direction relies on a trick from the proof of \cite[Theorem~2.1]{geroch1972IdealPointsSpaceTime}. Let $P$ be an IP, and take any $p\in P$. Set-theoretically decomposing~$P$ as the union over its intersection and complement with $I^+(p)$, and then applying the chronological past operator (which preserves joins), gives the following equation:
		\[
			P = I^-(P) = I^-\left[P\cap I^+(p)\right]\cup I^-\left[P\setminus I^+(p)\right].
		\]
	By indecomposability, it follows that $P$ must equal either of the summands. However $p\notin I^-\left[P\setminus I^+(p)\right]$, so it must be the case that $P = I^-\left[P\cap I^+(p)\right]$. The desired property now follows: if $x,y\in P$ then there exists $z\in P\cap I^+(p)$ such that $x,y\chron z$, and hence $z\in I^+(x)\cap I^+(y)\cap P\neq\varnothing$.	
\end{proof}

\begin{lemma}\label{lemma:chronological ideals are open}
	The $\chron$-ideals in a spacetime are open.
\end{lemma}
\begin{proof}
	Let $P\subseteq M$ be a $\chron$-ideal. Axiom~\ref{axiom:J2} just says that $P$ is a past set in the sense of \cref{definition:past set}, i.e.~$I^-(P)\subseteq P$. By \cref{proposition:interior of future/past set} it thus suffices to show $P\subseteq P^\circ = I^-(P)$. But this follows by~\ref{axiom:J1}: if $x\in P$ there exists $z\in P$ such that $x\chron z$, so $x\in I^-(P)$.
\end{proof}

\begin{proposition}
	The IPs of a spacetime $M$ are precisely the $\chron$-ideals:
		\[
			\IP(M) = \Ideals_{\chron}(M).
		\]
\end{proposition}
\begin{proof}
	If $P\subseteq M$ is an IP we know $\varnothing\neq P = I^-(P)$, so axioms~\ref{axiom:J0} and~\ref{axiom:J2} are satisfied. That axiom~\ref{axiom:J1} holds follows by \cref{lemma:IP iff finite chronological join}.
	
	Conversely, let $P$ be an $\chron$-ideal. Then by~\ref{axiom:J0} we get $P\neq\varnothing$, and~\cref{proposition:interior of future/past set,lemma:chronological ideals are open} together imply $I^-(P)= P$. Finally, axiom~\ref{axiom:J1} now implies $P$ is an IP through \cref{lemma:IP iff finite chronological join}.
\end{proof}

\begin{remark}
	Dually, the $\chron$-\emph{filters} will correspond to the indecomposable \emph{future} sets of the spacetime. Of course, more generally, for the principal $\aleq$-filters to exist we need a property dual to \cref{definition:past semifull} called \emph{future semi-fullness}.
\end{remark}

\begin{remark}
	Note that IPs are \emph{not} equal to the ideals with respect to the causality relation $\caus$. For instance, in $1+1$ Minkowski space we get the IP $I^-(x)$ induced by some point. If $y,z\caus x$ are on the boundary of the lightcone, there exists no $p\in I^-(x)$ such that $y,z\caus p$, essentially because $x\not\chron x$. Thus~\ref{axiom:J1} fails.
\end{remark}

The result in \cref{proposition:dcpo extension universality} thus shows that $\IP(M)$ is in some sense the ``universal'' directed join completion of the chronology relation $\chron$. This is order-theoretic justification to think of $\IP(M)$ as a future completion of $M$.

\section{The ideal completion of an ordered locale}
In this intermediary section we discuss an abstract ``compactification'' procedure for ordered locales. There are lattice-theoretic techniques to construct compactifications of frames. The most elementary such construction is as follows. A distributive lattice~$L$ gives rise to a frame $\Ideals(L)$ consisting of the \emph{ideals} of $L$. Recall that an ideal now means a subset $I\subseteq L$ such that:
	\begin{enumerate}[label = (I\arabic*)]\setcounter{enumi}{-1}
		\item $\bot\in I$, or equivalently: $I\neq \varnothing$;
		\item if $x,y\in I$ then $x\vee y\in I$; \hfill(\emph{upwards directedness})
		\item if $y\sqleq x\in I$ then $y\in I$. \hfill(\emph{down closure})
	\end{enumerate}
Ordered by subset inclusion, $\Ideals(L)$ is a frame. Explicitly, the meets are calculated by set-theoretic intersection, and the join of a family $(I_a)_{a\in I}$ of ideals in $L$ is calculated as
	\[
	\bigvee_{a\in A} I_a
	=
	\left\{
	\bigvee F: \text{ finite } F\subseteq\bigcup_{a\in A}I_a
	\right\}.
	\]
Further, the frame $\Ideals(L)$ is \emph{compact} (\cref{proposition:frame of ideals}). In the case that $L$ itself is already a frame, we can take the join $\bigvee I$ of any ideal $I\in\Ideals(L)$, and this defines a map of frames $\Ideals(L)\to L$. The right adjoint of this map is given by the map $x\mapsto \ideal(x)$, where $\ideal(x) = \{y\in L: y\sqleq x\}$ is the principal ideal generated by $x\in L$. For more details we refer to \cref{section:frame of ideals}, and see further \cite[Section~II.2.11]{johnstone1982StoneSpacesa} or \cite[Chapter~VII.4]{picado2012FramesLocalesTopology}. Unfortunately, this compactification procedure is not generally satisfactory: even if $L$ itself is compact the frame $\Ideals(L)$ can be considerably bigger. We here work out the general case, which will hopefully serve as a basis for refinements to better-behaved constructions in future work.

Viewed localically, for every locale $X$ there is a compact locale $\Ideals(X)$ whose frame of opens is defined by $\Opens(\Ideals(X)):= \Ideals(\Opens X)$, the frame of ideals in $\Opens X$, together with a dense sublocale inclusion $X\rightarrowtail \Ideals(X)$ \cite[\S VII.4.1.2]{picado2012FramesLocalesTopology}. In this section we show that this procedure extends to the setting of ordered locales. In particular, for an ordered locale $(X,\Leq)$, we will define the structure of an ordered locale on $\Ideals(X)$ by defining two monads on the frame $\Ideals(\Opens X)$. Namely, for an ideal $J\subseteq \Opens X$, we define
	\[
		\Down J := \{ U\in \Opens X: \exists W\in J: U\sqleq \Down W\}.
	\]
In other words, this is the ideal generated by the set $\{\Down W:W\in J\}$ (\cref{remark:ideal generated by set}).

\begin{lemma}\label{lemma:localic cones on ideals}
	The function $J\mapsto \Down J$ defines a monad on $\Ideals(\Opens X)$.
\end{lemma}
\begin{proof}
	We first need to show that if $J\in \Ideals(\Opens X)$, then $\Down J$ is again an ideal. First, since $U\sqleq \Down U$ we see that $J\subseteq \Down J$, guaranteeing that $\Down J$ is non-empty, and also providing the unit of the monad. Next, if $U_1,U_2\in \Down J$, then we can find $W_1,W_2\in J$ such that $U_i\sqleq\Down W_i$. Taking joins and using \cref{lemma:monotone laxly respects meets and joins} we find $U_1\vee U_2\sqleq \Down (W_1\vee W_2)$, where $W_1\vee W_2\in J$ since $J$ is an ideal. Hence $U_1\vee U_2\in \Down J$. Lastly, it is easy to see that $\Down J$ is downwards closed, since if $V\sqleq U\in\Down J$ then there exists $W\in J$ such that $V\sqleq U\sqleq \Down W$. Hence $\Down J$ is an ideal. 
	
	The function is monotone, since if $I\subseteq J$ are ideals, then for any $U\in \Down I$ we get $W\in I\subseteq J$ such that $U\sqleq \Down W$, and hence $U\in \Down J$. Since we already obtained the unit, we are only left to show idempotence. This amounts to proving $\Down\Down J \subseteq \Down J$. If $U\in \Down \Down J$ then we can find $V\in \Down J$, and in turn $W\in J$ such that $U\sqleq \Down V\sqleq \Down \Down W=\Down W$, giving that $U\in \Down J$. 
\end{proof}

By dualising, we get a second monad $J\mapsto \Up J$ on $\Ideals(\Opens X)$. Together, they allow us to define the structure of an ordered locale on $\Ideals(X)$.

\begin{definition}\label{definition:ordered locale of ideals}
	Let $(X,\Leq)$ be an ordered locale. We define a new ordered locale $\Ideals(X)$ as having the frame of opens $\Ideals(\Opens X)$ and the order $\Leq$ induced by the monads of \cref{lemma:localic cones on ideals}:
	\[
	I\Leq J
	\qquad\text{if and only if}\qquad
	I\subseteq \Down J\text{ and } J\subseteq \Up I.
	\]
\end{definition}

\begin{remark}
	It can help to spell out the condition $I\Leq J$ in terms of the elements of the ideals. The inclusion $I\subseteq \Down J$ translates to:
		\[
			\forall U\in I ~\exists V\in J : U\sqleq \Down V.
		\]
	Compare this to the \cref{definition:order on points} of the order on completely prime filters of an ordered locale.
\end{remark}

The following shows that this order on ideals reduces to the correct order on opens of the original locale.

\begin{lemma}
	If $\ideal(U)$ denotes the principal ideal generated by $U\in\Opens X$, then:
		\[
			\Up \ideal(U) = \ideal (\Up U)
			\qquad\text{and}\qquad
			\Down \ideal(U) = \ideal (\Down U).
		\]
	Thus, in an ordered locale with~\eqref{axiom:cones give order} we have:
		\[
		\ideal(U)\Leq \ideal(V)
		\qquad\text{if and only if}\qquad
		U\Leq V.
		\]
\end{lemma}
\begin{proof}
	If $V\in \Up \ideal(U)$ then we can find $W\sqleq U$ such that $V\sqleq \Up W\sqleq \Up U$, so $V\in \ideal(\Up U)$. Hence $\Up \ideal(U)\subseteq \ideal(\Up U)$. The converse inclusion also holds, since $U\in \ideal(U)$, so $V\sqleq \Up U$ immediately implies $V\in\Up \ideal(U)$.
	
	It now follows that $\ideal(U)\Leq \ideal(V)$ if and only if $\ideal(U)\subseteq \ideal(\Down V)$ and $\ideal(V)\subseteq \ideal(\Up U)$, which in turn is equivalent to $U\sqleq \Down V$ and $V\sqleq \Up U$. Under~\eqref{axiom:cones give order} this holds if and only if $U\Leq V$.
\end{proof}

Recall that if $f\colon X\to Y$ is a map of locales, this induces a map of locales $\Ideals(f)\colon \Ideals(X)\to \Ideals(Y)$ whose map of frames sends an ideal $J\subseteq \Opens Y$ to the ideal in $\Opens X$ generated by $\{f^{-1}(V):V\in J\}$ (\cref{remark:ideal generated by set}).

\begin{lemma}
	If $f\colon X\to Y$ is monotone, then $\Ideals(f)\colon \Ideals(X)\to \Ideals(Y)$ is monotone.
\end{lemma}
\begin{proof}
	We need to show that $\Up \Ideals(f)^{-1}(J)\subseteq \Ideals(f)^{-1}(\Up J)$, for any ideal $J\subseteq \Opens Y$. If $U\in \Up \Ideals(f)^{-1}(J)$ then there exists $W\in \Ideals(f)^{-1}(J)$ such that $U\sqleq \Up W$, and in turn there exists $V\in \Opens Y$ such that $W\sqleq f^{-1}(V)$. Using monotonicity of $f$, this gives $U\sqleq \Up f^{-1}(V) \sqleq f^{-1}(\Up V)$, directly implying $U\in \Ideals(f)^{-1}(\Up J)$. The proof for past cones is dual.
\end{proof}

\begin{corollary}\label{corollary:ideal functor ordloc}
	There is a functor:
	\begin{align*}
		\Ideals(-): \OrdLoc &\longrightarrow \OrdLoc;\\
		(X,\Leq) & \longmapsto (\Ideals(X),\Leq);\\
		f & \longmapsto \Ideals(f).
	\end{align*}
\end{corollary}

This construction is defined on arbitrary ordered locales, but in \cref{section:parallel ordered locales} we identified the parallel ordered locales with~\eqref{axiom:LV} as the physically reasonable ones. In the remainder of this section we prove that $\Ideals(-)$ preserves these properties.

\begin{lemma}
	If $(X,\Leq)$ is an ordered locale whose localic cones preserve finite joins, then $\Ideals(X)$ satisfies axiom~\eqref{axiom:LV}.
\end{lemma}
\begin{proof}
	Since we already know the localic cones of $\Ideals(X)$ are monotone, it suffices to show that $\Down \bigvee_{a\in A} J_a \subseteq \bigvee_{a\in A} \Down J_a$. If $U\in \Down \bigvee_{a\in A} J_a$ then we can find some $W\in \bigvee_{a\in A}J_a$ such that $U\sqleq \Down W$. By the description of the joins in $\Ideals(\Opens X)$ from above, we can find some finite set $F\subseteq \bigcup_{a\in A}J_a$ such that $W=\bigvee F$. Since the localic cones of $X$ preserve finite joins, we hence get $U\sqleq \bigvee_{V\in F}\Down V$, or equivalently: $U= \bigvee_{V\in F} U\wedge \Down V$. Now clearly, if $V\in J_a$, then $U\wedge \Down V\in \Down J_a$, and hence $U$ can be written as a finite join of elements in $\bigcup_{a\in A}\Down J_a$, proving that $U\in \bigvee_{a\in A}\Down J_a$.
\end{proof}

\begin{lemma}
	If $(X,\Leq)$ is an ordered locale with~\eqref{axiom:LV}, then ${i\colon X\rightarrowtail \Ideals(X)}$ is monotone.
\end{lemma}
\begin{proof}
	As mentioned, the map of frames $i^{-1}$ sends an ideal $I\subseteq \Opens X$ to the join $\bigvee I\in \Opens X$. Monotonicity therefore is the condition that $\Up\bigvee I \sqleq \bigvee \Up I$ and $\Down \bigvee I \sqleq \bigvee \Down I$. With~\eqref{axiom:LV} we get $\Up \bigvee I = \bigvee \{\Up U : U \in I\}$, and clearly $\Up U \in \Up I$ for any $U\in I$, so the first inclusion follows. The second inclusion is dual.
\end{proof}

\begin{corollary}\label{corollary:ideal functor ordloc V}
	The functor from \cref{corollary:ideal functor ordloc} restricts to a functor
		\[
			\Ideals(-)\colon \OrdLoc_\vee\longrightarrow \OrdLoc_\vee.
		\]
	Moreover, there is a natural transformation $\id_{\OrdLoc_\vee}\to \Ideals(-)$.
\end{corollary}

\begin{lemma}
	If $(X,\Leq)$ is a parallel ordered locale, then so is $(\Ideals(X),\Leq)$.
\end{lemma}
\begin{proof}
	First, the ideal completion satisfies~\eqref{axiom:empty} since $\Up \varnothing = \varnothing$ in $X$, so that $\Up \{\varnothing\} =\{\varnothing\}$ in $\Ideals(X)$, and dually for down cones. We are thus left to prove~\eqref{axiom:wedge}.
	
	By definition $\Ideals(X)$ satisfies~\eqref{axiom:cones give order}, so by \cref{lemma:wedge iff frobenius} it in turn suffices to prove~\eqref{axiom:frobenius}. We prove that for ideals $I,J\in\Ideals(\Opens X)$ we get
		\[
			I\cap \Down J\subseteq \Down(\Up I\cap J).
		\]
	For that, take $U\in I\cap \Down J$, which means that $U\in I$ and there exists $V\in J$ such that $U\sqleq \Down V$. Using that~\eqref{axiom:frobenius-} holds in $X$, we get that
		\[
			U= U\wedge \Down V\sqleq \Down (\Up U\wedge V),
		\]
	so we see that $U$ is contained in the past of $\Up U \wedge V\in \Up I\cap J$, hence giving $U\in \Down (\Up I \cap J)$. The dual inclusion is proved analogously.
\end{proof}

\begin{corollary}
	The functor from \cref{corollary:ideal functor ordloc} restricts to a functor
	\[
		\Ideals(-)\colon \OrdLoc_{\newparallel}\longrightarrow\OrdLoc_{\newparallel}.
	\]
	on the full subcategory of parallel ordered locales.
\end{corollary}

Combined, any parallel ordered locale $(X,\Leq)$ with~\eqref{axiom:LV} gives rise to a monotone sublocale embedding $i\colon X\rightarrowtail \Ideals(X)$ into a compact parallel ordered locale with~\eqref{axiom:LV}.

\section{Generalised localic ideal points}
\label{section:generalised ideal points}
In this section we introduce a localic, generalised notion of ``ideal point.'' We already generalised IPs and IFs from spacetimes to the setting of ordered locales in \cref{section:ideal points in locales}, but saw in \cref{example:space without IPs} that e.g.~vertical-$\mathbb{R}^2$ does not admit any ideal points in this sense. Moreover, the construction of $\IP(X)$ from $X$ is treated as a \emph{point-set} construction. To stay in the realm of point-free topology we rather want a \emph{locale} whose \emph{points} correspond to the ``ideal points'' of $X$.\footnote{We learned this lesson in localic thinking from \cite{vickers2011MonadValuationLocales}. On locales (or more generally distributive lattices) it is natural to consider so-called \emph{valuations}, which we can think of as abstract measures on spaces. However, rather than studying the \emph{set} of valuations on a given locale, one rather wants a \emph{locale} whose \emph{points} return the valuations. We employ the same philosophy here.} In this section we build towards such a construction.

We have already seen a major ingredient of it in \cref{section:locale of pasts and futures}: recall that any ordered locale $(X,\Leq)$ with~\eqref{axiom:LV} induces locales of futures and pasts
	\[
	\begin{tikzcd}[ampersand replacement=\&]
		X \& {X^\triup}
		\arrow["{\eta_\triup}", two heads, from=1-1, to=1-2]
	\end{tikzcd}
	\qquad\text{and}\qquad
	\begin{tikzcd}[ampersand replacement=\&]
		X \& {X^\tridown.}
		\arrow["{\eta_\tridown}", two heads, from=1-1, to=1-2]
	\end{tikzcd}
	\]
We claim that we can think of $X^\triup$ and $X^\tridown$ as a type of \emph{future} and \emph{past} completions of $X$, respectively. The ultimate motivation for this is that, as we shall see in \cref{corollary:IPs as primes in spacetime}, for a spacetime $M$ we get $\pt(M^\triup)\cong\IP(M)$ and $\pt(M^\tridown)\cong \IF(M)$. The following examples show how $X^\triup$ can be thought of as a future completion.

\begin{example}\label{example:generalised IP in lower reals}
	Recall that $\mathbb{R}^\triup$ is the locale defined by the topology of the intervals $(a,\infty)$. We saw in \cref{remark:completely prime filters vs elements} that $\pt(\mathbb{R}^\triup)\cong\mathbb{R}\cup\{\infty\}$, where the point at infinity is represented by the completely prime filter $\calF_\infty = \Opens \mathbb{R}^\triup\setminus\{\varnothing\}$. With the order $\leq$ on points from \cref{definition:order on points} we get $\calF_x\leq \calF_\infty$ for every $x\in \mathbb{R}$, so that we can indeed think of $\calF_\infty$ as an ideal point at future infinity.
\end{example}

\begin{example}\label{example:generalised IPs in vertical R2}
	Consider vertical-$\mathbb{R}^2$ from \cref{example:vertical R2}, viewed as an ordered locale $\mathbb{V}$. We saw already in \cref{example:space without IPs} that it does not admit any indecomposable past regions: $\IP(\mathbb{V})=\varnothing$. Here we show that it \emph{does} admit ``generalised'' ideal points. As a set, we claim that:
		\[
			\pt'\left(\mathbb{V}^\triup\right) \cong \mathbb{R}^2\sqcup \mathbb{R},
		\]
	where the points in $\mathbb{R}^2$ shall correspond to the ``internal'' points of $\mathbb{V}$, and $\mathbb{R}$ represents the points in the future boundary, thought of as the would-be limits of curves that propagate vertically into positive infinity.
	
	To start, note that the frame of opens of the locale of futures $\mathbb{V}^\triup$ is just the topology on the set~$\mathbb{R}^2$ generated by basic opens of the form $(l,r)\times (b,\infty)$, for $l,r,b\in\mathbb{R}$. Therefore any $x=(x_1,x_2)\in \mathbb{R}^2$ gives a prime $P_x\in\pt'(\mathbb{V})$, and thus a generalised ideal point
		\[
			\eta_\ast^\triup(P_x) = \bigcup\{\Up U\in \Opens\mathbb{R}^2 : x\notin \Up U \}
			=
			\mathbb{R}^2\setminus\left(\{x_1\}\times (-\infty,x_2]\right).
		\]
	Note that this is precisely the \emph{complement} of the would-be IP $\down x\notin \Opens \mathbb{R}^2$ (cf.~\cref{example:space without IPs} and \cref{figure:eta on prime}). This assignment is clearly injective, so gives a full copy of $\mathbb{R}^2$ in the set of points of~$\mathbb{V}^\triup$. It is now straightforward to prove that any $x_1\in \mathbb{R}$ similarly gives rise to a prime element
		\[
			\mathbb{R}^2\setminus \left(\{x_1\}\times \mathbb{R}\right)\in \pt'\left(\mathbb{V}^\triup\right),
		\]
	which represents a point at future infinity with horizontal coordinate $x_1$. This assignment similarly injects a full copy of $\mathbb{R}$ into the points of $\mathbb{V}^\triup$. Combining \cref{example:generalised IP in lower reals} with the fact that $\mathbb{R}$ is sober shows this exhausts all possible points in $\mathbb{V}^\triup$. 
\end{example}

For the rest of this section, we introduce the following terminology.
\begin{definition}\label{definition:generalised ideal points}
	A \emph{generalised} or \emph{localic ideal point} in an ordered locale $(X,\Leq)$ with~\eqref{axiom:LV} is a point of $X^\triup$ or $X^\tridown$. The points of $X^\triup$ are called \emph{future ideal points}, and the points in $X^\tridown$ are called \emph{past ideal points}.
\end{definition}

\subsection{Ideal points in locales with regular cones}
This section builds towards the advertised bijections of $\pt(M^\triup)\cong\IP(M)$ and $\pt(M^\tridown)\cong \IF(M)$. The equations in \cref{lemma:chronological cone of double negation,lemma:chronological cones are regular opens} showcase a remarkable interaction between the localic cones and the Heyting negation in the ordered locale of a spacetime. In particular, the images $\im(I^\pm)$ of the chronological cones live in the \emph{Boolean} double-negation sublocale $\Opens M_{\neg\neg}$ (\cref{section:double-negation sublocale}). This turns out to be the key to relate the boundary constructions from the relativity literature to our definitions, as will soon be explained. In what follows, the Heyting negation is always calculated in the frame $\Opens X$, unless otherwise specified.

For convenience, we recall \cref{lemma:heyting negation of past and future set}.
\begin{lemma}\label{lemma:heyting negation of past and future set 2}
	For every open $U$ in a parallel ordered locale:
	\[
	\Up\neg\Down U = \neg \Down U
	\qquad\text{and}\qquad
	\Down \neg\Up U = \neg\Up U.
	\]
\end{lemma}

This means that in any parallel ordered locale, the Heyting negation determines a monotone function $\neg\colon \im(\Down)\to\im(\Up)^\op$, and dually. In fact:

\begin{lemma}\label{lemma:galois adjunction past future}
	In a parallel ordered locale, Heyting negation determines a Galois adjunction between past and future sets:
		\[
		\begin{tikzcd}[column sep =large]
			\im(\Down) & \im(\Up)^\op.
			\arrow[""{name=0, anchor=center, inner sep=0}, "\neg", shift left=2, from=1-1, to=1-2]
			\arrow[""{name=1, anchor=center, inner sep=0}, "\neg", shift left=2, from=1-2, to=1-1]
			\arrow["\dashv"{anchor=center, rotate=-90}, draw=none, from=0, to=1]
		\end{tikzcd}
		\]
\end{lemma}
\begin{proof}
	This follows since Heyting negation forms a Galois adjunction between $\Opens X$ and $\Opens X^\op$ (\cref{corollary:heyting negation gives adjunction}), which restricts to the desired adjunction by \cref{lemma:heyting negation of past and future set 2}.
\end{proof}

\begin{definition}\label{definition:regular localic cones}
	An ordered locale $X$ has \emph{regular cones} if for all $U\in\Opens X$:
		\[
			\neg\neg\Up U = \Up U = \Up \neg\neg U
			\qquad\text{and}\qquad
			\neg\neg\Down U = \Down U = \Down \neg\neg U .
		\]
\end{definition}

\begin{example}\label{exampel:spacetime has regular cones}
	For any smooth spacetime $M$ the ordered locale $\loc(M)$ has regular cones by \cref{lemma:chronological cone of double negation,lemma:chronological cones are regular opens}:
		\[
		\neg\neg \circ I^\pm = I^\pm = I^\pm\circ\neg\neg.
		\]
\end{example}

\begin{example}\label{example:identity does not have regular cones}
	For any locale $X$, the ordered locale $(X,=)$ has regular cones if and only if $\neg\neg U = U$ for all $U\in \Opens X$, i.e.~precisely when its frame of opens is \emph{Boolean}. Thus the regular cone condition is a generalisation of Boolean locales.
\end{example}

In an ordered locale with regular cones, we see that both $\im(\Up)$ and $\im(\Down)$ are contained in the double-negation sublocale $\Opens X_{\neg\neg}$. If the locale is additionally parallel ordered, we therefore get the following.

\begin{proposition}\label{proposition:heyting negation isomorphism past and future}
	In a parallel ordered locale with regular cones, Heyting negation determines an order isomorphism $\neg\colon \im(\Down)\xrightarrow{\sim} \im(\Up)^\op$.
\end{proposition}
\begin{proof}
	Clearly the conditions of \cref{definition:regular localic cones} ensure that the Galois adjunction from \cref{lemma:galois adjunction past future} involves maps that are mutually inverse. Explicitly, every $\Up U$ can be written as the complement of a downset by writing $\Up U = \neg\neg\Up U = \neg\Down\neg\Up U$, and similarly $\Down U = \neg\Up\neg\Down U$. 
\end{proof}

Recall from \cref{remark:pt'} that $\pt'(X)$ is the space of prime elements of $\Opens X$, which is homeomorphic to the space $\pt(X)$ of completely prime filters in $\Opens X$ via the correspondence in \cref{lemma:prime elements are cpf}. We have already observed that an irreducible past-region $P\in \IP(X)$ is really just a \emph{coprime} element of $\im(\Down)$, and therefore dually a \emph{prime} element of $\im(\Down)^\op$. Therefore Heyting negation will preserve prime elements, so we obtain a prime element ${\neg P \in \im(\Up)= \Opens X^\triup}$, which in turn is just a \emph{point} $\neg P \in \pt'(X^\triup)$. The visual intuition here is represented in \cref{figure:eta on prime}. This is made precise as follows. 

\begin{lemma}\label{lemma:primes preserved by isos}
	An order isomorphism $h\colon L\to M$ between bounded distributive lattices determines a bijection between the prime elements.
\end{lemma}
\begin{proof}
	First, if $h$ is an order isomorphism then $h\dashv h^{-1}\dashv h$, so by \cref{theorem:adjoints preserve limits} both $h$ and its inverse preserve all existing joins and meets, including the top element. Thus if $p\in L$ is prime, we get that $h(p)\neq\top$. Further, if $x\wedge y\sqleq h(p)$, then $h^{-1}(x)\wedge h^{-1}(y)\sqleq p$, so by primeness we get $h^{-1}(x)\sqleq p$ or $h^{-1}(y)\sqleq p$, which implies $x\sqleq h(p)$ or $y\sqleq h(p)$, respectively. Thus $h$ preserves primes, so it is clear that it restricts to a bijection between the sets of prime elements.
\end{proof}

\begin{theorem}\label{corollary:IPs as primes}
	In a parallel ordered locale $X$ with regular cones there are bijections:
		\[
			\pt'(X^\triup)\cong \IP(X)
			\qquad\text{and}\qquad
			\pt'(X^\tridown)\cong \IF(X).
		\]
\end{theorem}
\begin{proof}
	Using \cref{proposition:heyting negation isomorphism past and future,lemma:primes preserved by isos} we see that Heyting negation in $\Opens X$ determines a bijection between the primes of $\im(\Up)$ and the primes of~$\im(\Down)^\op$, i.e.~the coprimes in $\im(\Down)$. The former are precisely the points of $X^\triup$ (\cref{remark:pt'}) and the latter are precisely the IPs of $X$ (\cref{remark:coprimes}). Explicitly, this sends $F\in \pt'(X^\triup)$ to $\neg F\in \IP(X)$, and $P\in \pt'(X^\tridown)$ to $\neg P\in \IF(X)$.
\end{proof}

\begin{corollary}\label{corollary:IPs as primes in spacetime}
	For any smooth spacetime $M$ there are bijections:
		\[
		\pt'(M^\triup)\cong \IP(M)
		\qquad\text{and}\qquad
		\pt'(M^\tridown)\cong \IF(M).
		\]
\end{corollary}
\begin{proof}
	Observe that spacetimes have regular cones by \cref{lemma:chronological cone of double negation,lemma:chronological cones are regular opens}.
\end{proof}

\begin{example}
	The vertically-ordered $\mathbb{R}^2$ from \cref{example:vertical R2} does not have regular cones. For instance, consider the open subset $U=\left(\mathbb{R}\setminus\{0\}\right)\times(0,1)$. The double negation of this open is just $\neg\neg U = \mathbb{R}\times (0,1)$, so
		\[
			\Up U 
			=
			\left(\mathbb{R}\setminus\{0\}\right)\times (0,\infty)
			\neq 
			\mathbb{R}\times (0,\infty)
			=
			\Up \mathbb{R}\times (0,1)
			=
			\neg\neg \Up U = \Up \neg\neg U.
		\]
	And this makes sense: we know from \cref{example:space without IPs} that this space has no IPs and IFs, while we obtained non-trivial generalised ideal points in \cref{example:generalised IPs in vertical R2}. This shows that, in general, the points of $X^\tridown$ have a richer structure than the IFs of~$X$. It just so happens that for spacetimes these notions coincide. Generalised ideal points are therefore a genuine generalisation of the boundary constructions of \cite{geroch1972IdealPointsSpaceTime}.
\end{example}

\begin{example}
	The real line $(\mathbb{R},\leq)$ with the Euclidean topology and standard order clearly has regular cones, but is not a spacetime. We get $\pt'(\mathbb{R}^\triup)\cong \mathbb{R}\sqcup \{\infty\}$.
\end{example}

\subsection{Thoughts on a combined causal boundary}
The results above suggests that we can indeed think of $X^\triup$ and $X^\tridown$ as localic future and past completions of $X$, respectively. However, a full-fledged causal boundary completion should produce a locale~$\overline{X}$ that incorporates both the future and past ideal points simultaneously. For a spacetime, the idea is that a genuine point $x\in M$ is equivalently described by $I^-(x)\in\IP(M)$ and $I^+(x)\in \IF(M)$. In general, the idea is to introduce a pairing of IPs and IFs, and treat the resulting quotient of $\IP(M)\times \IF(M)$ as the causal boundary $\overline{M}$. For spacetimes, there are several approaches to this, appearing already in e.g.~\cite{budic1974CausalBoundariesGeneral,szabados1988CausalBoundaryStrongly}. In setting of (convex) locales, we therefore suggest to look for a suitable locale $\overline{X}$ that either factorises the sublocale~$\eta$ from \cref{section:locale of pasts and futures}, or appears as a quotient locale:
	\[
		\begin{tikzcd}[cramped,column sep=2em]
			X & {\overline{X}} & {X^\triup\times X^\tridown}
			\arrow[tail, from=1-1, to=1-2]
			\arrow[tail, from=1-2, to=1-3]
		\end{tikzcd}
		\qquad\text{or}\qquad
		\begin{tikzcd}[cramped,column sep=2em]
			X & {X^\triup\times X^\tridown} & {\overline{X}.}
			\arrow[tail, from=1-1, to=1-2]
			\arrow[two heads, from=1-2, to=1-3]
		\end{tikzcd}
	\]
In spacetime approaches one of the difficulties is to define a suitable topology on the point-set $\overline{M}$ if identified pairs of IPs and IFs. The benefit of a localic construction would be that the topological data is automatically encoded in $\overline{X}$. 	
	
One of the most recently successful approaches is in \cite{marolf2003NewRecipeCausala}, generalising that of \cite{szabados1988CausalBoundaryStrongly}. See also \cite{flores2011FinalDefinitionCausal}. An important notion here is that of \emph{common future} and \emph{common past:}
	\[
		C^+(U):= \left(\bigcap_{x\in U} I^+(x)\right)^{\!\circ}
		\qquad\text{and}\qquad
		C^-(U):= \left(\bigcap_{x\in U}I^-(x)\right)^{\!\circ}.
	\]
The set $C^+(U)$ is the largest region such that as soon as $W\cap C^+(U)\neq\varnothing$ we get $U\subseteq \Down W$. The \emph{S-relation} $\simS$ is then defined on pairs $(P,F)\in \IP(M)\times \IF(M)$ as:
	\[
		P\simS F
		\qquad\text{if and only if}\qquad
		\begin{array}{l}
			\text{$P$ is a maximal subset of $C^-(F)$, and}\\
			\text{$F$ is a maximal subset of $C^+(P)$.}
		\end{array}
	\]
It was proved in \cite[Proposition~5.1]{szabados1988CausalBoundaryStrongly} that $M$ is strongly causal if and only if $I^-(x)\simS I^+(x)$ for all $x\in M$. Moreover, in strongly causal spacetimes ${I^-(x)\simS F}$ if and only if $F= I^+(x)$. This shows the relation $\simS$ forms a correct pairing of ideal points coming from elements $x\in M$, and this allows us to give a well-behaved embedding $M\to \overline{M}$. The behaviour of this pairing on genuine boundary points is more intricate, see e.g.~\cite[\S 3.1.1]{flores2011FinalDefinitionCausal} for more details.

It is relatively straightforward to adopt these ideas to IPs and IFs in ordered locales. In particular, it is possible to write down a notion of localic common future and past operators that has properties analogous to $C^\pm$ from \cite{budic1974CausalBoundariesGeneral,hawking1974CausallyContinuousSpacetimes}. The main open problem is how to translate these point-wise ideas to localic constructions involving $X^\triup$ and $X^\tridown$, for instance a localic version of the $\simS$ relation to the frame $\im(\Up)\oplus\im(\Down)$, which we leave to future research.
One possibility is to use a type of \emph{Dedekind-MacNeille} completion. The techniques used to construct full causal boundaries in the physics literature strongly resemble the general completion procedure outlined in \cite[\S 3.2]{davey2002IntroductionLatticesOrder}. These could then be generalised to the point-free setting using e.g.~\cite{ball2014ExtendingSemilatticesFrames}. Perhaps more appropriately, we can attempt to generalise the directed join completion procedure from \cref{section:chronological ideals} to incorporate both future and past ideals. Similar ideas appear in \emph{order compactification} techniques in e.g.~\cite{choe1979WallmanTypeOrder} Lastly, we suggest the use of \emph{biframe compactification} techniques from e.g.~\cite{schauerte1992Biframes,bezhanishvili2015ProximityBiframesCompactifications} to construct compactifications of ordered locales.

Lastly, a successful localic boundary construction could also lead to a localic notion of \emph{global hyperbolicity} (\cref{section:causal ladder}) using characterisations in \cite[\S 6]{seifert1971CausalBoundarySpacetimes}.

\subsection{Ideal points in double-negation sublocales}\label{section:IPs in double-negation sublocale}
Recall from \cref{section:double-negation sublocale} that the frame of opens $\Opens X_{\neg\neg}$ consists precisely of the \emph{regular opens:} those $U\in \Opens X$ satisfying $\neg\neg U = U$, and that the sublocale inclusion map $i\colon X_{\neg\neg}\rightarrowtail X$ is defined by the frame map $i^{-1}\colon U\mapsto \neg\neg U$, and that the direct image map $i_\ast$ is just the inclusion function $\Opens X_{\neg\neg}\subseteq \Opens X$. This procedure often produces genuinely pointless locales. Recall that an \emph{isolated point} of a topological space $S$ is an element $x\in S$ such that $\{x\}\in\Opens S$.

\begin{proposition}\label{proposition:points in regular opens Hausdorff space}
	If $S$ is a Hausdorff space, then
		\[
			\pt(S_{\neg\neg})\cong \{\text{\normalfont isolated points in $S$}\}.
		\]
\end{proposition}
\begin{proof}
	The frame of regular opens $\Opens S_{\neg\neg}$ is a Boolean algebra (\cref{lemma:double-negation sublocale is boolean}), so the negation operator defines a bijective correspondence between prime elements and atoms by \cref{proposition:points in boolean locale are atoms}. If $x\in S$ is an isolated point, since singletons are closed in a Hausdorff space, we get $\neg\neg \{x\} = (\overline{\{x\}})^\circ = \{x\}$, which clearly defines an atom in $\Opens S_{\neg\neg}$.
	
	Conversely, let $A\in\Opens S_{\neg\neg}$ be an atom. Suppose for the sake of contradiction that $A$ contains two distinct elements $x\neq y$. By Hausdorffness we can find disjoint opens $U,V\in\Opens X$ with $x\in U$ and $y\in V$. But then $x\in U\cap A\subsetneq A$ is a non-empty set, contradicting that $A$ is an atom. Hence $A$ must be a singleton, defining a unique isolated point of $S$.
\end{proof}

In the following we show how the ideal points of an ordered locale $(X,\Leq)$ with regular cones can equivalently be calculated in the double-negation sublocale. On~$X_{\neg\neg}$ we induce an ordered locale structure $\Leq_i$ via the sublocale inclusion map~$i$, according to the procedure in \cref{example:ordered locale on domain}. Recall this means for $U = \neg\neg U$ and $U'= \neg\neg U'$ in $\Opens X_{\neg\neg}$ we have:
	\[
	U\Leq_i U'
	\qquad\text{if and only if}\qquad
	\begin{array}{l}
		U\sqleq i^{-1}(V)\text{~implies~}U'\sqleq i^{-1}(\Up V) \text{~and}\\
		U'\sqleq i^{-1}(V') \text{~implies~}U\sqleq i^{-1}(\Down V').
	\end{array}
	\]
Note here that the inclusion $\sqleq$ is calculated equivalently in $\Opens X$ or the frame of regular opens. 

\begin{lemma}\label{lemma:order in double-negation sublocale with regular cones}
	Let $(X,\Leq)$ be an ordered locale with regular cones and~\eqref{axiom:cones give order}. Then in the double-negation sublocale $(X_{\neg\neg},\Leq_i)$ we have for all $U,U'\in\Opens X$:
		\[
			\neg\neg U\Leq_i \neg\neg U'
			\qquad\text{if and only if}\qquad
			U\Leq U'.
		\]
\end{lemma}
\begin{proof}
	Suppose that $U,U'\in\Opens X$ are (not necessarily regular) opens such that $\neg\neg U \Leq_i \neg\neg U'$. In the above characterisation of $\Leq_i$, set $V= U$ and $V'=U'$, which then implies 
		\[
			U'\sqleq \neg\neg U'\sqleq i^{-1}(\Up U) = \neg\neg\Up U = \Up U,
		\]
	where the last step follows from the regular cone assumption. Dually, we obtain $U\sqleq \Down U'$, so by~\eqref{axiom:cones give order} we get $U\Leq U'$ in $X$.
	
	For the converse, suppose that $U\Leq U'$, and take $V\in\Opens X$ with ${\neg\neg U\sqleq i^{-1}(V)}$. We need to show that $\neg\neg U'\sqleq i^{-1}(\Up V)=\neg\neg \Up V =\Up V$. But $U\Leq U'$ implies $U'\sqleq \Up U$, so using regular cones we get:
		\[
			\neg\neg U' \sqleq \neg\neg \Up U = \Up \neg\neg U \sqleq \Up \neg\neg V =\Up V.
		\]
	The dual condition is proved similarly, so we get $\neg\neg U \Leq_i \neg\neg U'$. 
\end{proof}

\begin{proposition}\label{corollary:localic cones double-negation sublocale regular}
	If $(X,\Leq)$ has regular cones and~\eqref{axiom:cones give order} then the localic cones in the double-negation sublocale $(X_{\neg\neg},\Leq_i)$ are calculated as:
		\[
			\Upsub{i}\neg\neg U = \Up U
			\qquad\text{and}\qquad
			\Downsub{i}\neg\neg U = \Down U.
		\]
\end{proposition}
\begin{proof}
	We saw in \cref{example:ordered locale on domain} that monotonicity $i\colon X_{\neg\neg}\rightarrowtail X$ already gives the inclusions $\Upsub{i}\neg\neg U \sqleq \Up \neg\neg U = \Up U$. The converse inclusions follow from \cref{lemma:order in double-negation sublocale with regular cones}: if $U\Leq V$ then $\neg\neg U\Leq_i \neg\neg V$, so $V\sqleq \Upsub{i}\neg\neg U$.
\end{proof}

\begin{theorem}
	If $(X,\Leq)$ has regular cones and~\eqref{axiom:cones give order}, then there are isomorphisms of locales
		\[
			X^\triup \cong (X_{\neg\neg})^\triup
			\qquad\text{and}\qquad
			X^\tridown \cong (X_{\neg\neg})^\tridown.
		\]
\end{theorem}
\begin{proof}
	It suffices to show that the map of frames $i^{-1}\colon \Opens X\to \Opens X_{\neg\neg}$ restricts to an isomorphism of frames $\im(\Up)\xrightarrow{\sim}\im(\Upsub{i})$. But this follows by \cref{corollary:localic cones double-negation sublocale regular}, since $i^{-1}(\Up U) = \neg\neg \Up U = \Up U = \Upsub{i}\neg\neg U$ defines a bijection.
\end{proof}

\begin{remark}
	Note that the order matters here: $(X^\triup)_{\neg\neg}$ is something more trivial, since there the double Heyting negation is taken in the locale $X^\triup$. For instance, in Minkowski space $M$ we always get $I^+(U)\cap I^+(V) \neq\varnothing$ as long as $U$ and $V$ are non-empty, so that the Heyting negation of $I^+(U)$ calculated in $M^\triup$ will be empty. In turn, the double negation is just the whole space. Hence $(M^\triup)_{\neg\neg}\cong 1$. 
\end{remark}

\begin{corollary}\label{corollary:IPs from regular opens}
	For any smooth spacetime there are bijections:
		\[
			\IP(M) \cong \IP(M_{\neg\neg})
			\qquad\text{and}\qquad
			\IF(M)\cong\IF(M_{\neg\neg}).
		\]
\end{corollary}

Remarkably, we know from \cref{proposition:points in regular opens Hausdorff space} that, since $M$ is a manifold, $\pt(M_{\neg\neg})=\varnothing$. Yet this corollary says the ideal points recovered from this pointless locale $M_{\neg\neg}$ include all the original ones from $M$, and in particular the original elements $x\in M$ in the guise of $I^+(x)$ and $I^-(x)$. Thus the ideal points procedure recovers all the points of $M$ seemingly ``out of nowhere.'' We take this as evidence that there exist genuinely point-free locales whose causal structure is nevertheless interesting.

\chapter{Causal coverages}
\label{section:causal coverage}
In this section we study how a natural notion of \emph{coverage} arises in an ordered locale. A similar idea was first proposed in \cite[Definition~2.12]{christensen2005CausalSitesQuantum} in the setting of \emph{causal sites} (briefly discussed above in \cref{section:causal sites vs ordered locales}), but has to the best of our knowledge never been developed anywhere else in the literature. Our considerations here are deeply inspired by this definition.

The intuition is as follows. We think of an ordered space $(S,\leq)$ as an arena in which information can ``flow,'' subject to the restraints of the order~$\leq$. For example, in spacetime it is generally an accepted notion that physical information can only propagate along causal curves. Interpreted in our more general setting, we might imagine information being able to flow along continuous monotone curves $\gamma\colon [0,1]\to S$ of the ordered space. Note that at this point we are agnostic as to the type of information that is ``contained'' in the space, but merely study the possible ways of propagation.

We can now ask to what extent the information of a given region $U$ is determined by some other region $A$ in its past. This is made more precise as follows. We say $A$ \emph{covers $U$ from below} if:
	\begin{quote}
		every monotone path $\gamma$ landing in $U$ either intersects $A$, or can be extended into the past to intersect $A$.
	\end{quote}
In other words: all flows of information that land in $U$ have to have come from~$A$. We call this \emph{causal coverage}. In this chapter we develop the point-free analogues of these ideas. One of our main results is a well-behaved causal coverage relation $\CovLeq^\pm$ (\cref{lemma:properties of coverage from locale}) for any parallel locale $(X,\Leq)$ with~\eqref{axiom:LV}, giving rise to a new type of \emph{Grothendieck topology} (\cref{theorem:causal grothendieck topology}). This also fixes a defect of causal sites \cite[Remark~2.13]{christensen2005CausalSitesQuantum}. Out of this causal coverage emerges a natural abstract notion of \emph{domain of dependence} (\cref{section:domains of dependence}). We compare the localic domains of dependence to the ones typically used in relativity theory~\cite{geroch1970DomainofDependence}, and find that they are generally distinct. We conjecture in \cref{section:holes in spacetime} that this difference may have implications for the problem of ``holes'' in spacetime. 

For most of this section, we work in a fixed parallel ordered locale $(X,\Leq)$ that satisfies axiom~\eqref{axiom:LV}. 

\section{Paths in an ordered locale}
\noindent\begin{minipage}{0.64\textwidth}
		To start, we need the point-free analogue of monotone paths $\gamma\colon [0,1]\to S$. Instead of points, the idea is to take totally ordered families $(p_i)_{i\in I}$ of open regions that cover the image $\im(\gamma)\subseteq S$, or at least part of it, such that $p_i\Leq p_j$ whenever $i < j$ (\cref{figure:covered path}). Specifically, we might think of a family $(p_t)_{t\in [0,1]}$ of open neighbourhoods $p_t\ni \gamma(t)$. However, since the unit interval is compact and $\gamma$ is continuous, the image $\im(\gamma)$ is compact, and so $(p_i)_{i\in I}$ will have a finite subcover. For this reason, we shall only consider paths to be defined in terms of \emph{finite} sequences, instead of arbitrary total orders.
\end{minipage}%
\hfill%
\begin{minipage}{.33\textwidth}\centering
	\input{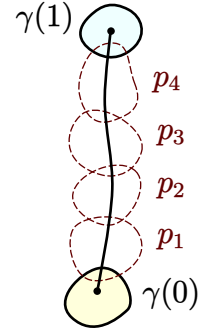}
	\captionof{figure}{Illustration of a path $(p_i)_{i\in I}$ covering $\im(\gamma)$.}
	\label{figure:covered path}
\end{minipage}

\begin{definition}\label{definition:paths}
	A \emph{path} in an ordered locale $(X,\Leq)$ is a finite sequence $p=(p_n)_{n=0}^N$ of non-empty open regions, called \emph{steps}, such that $p_n\Leq p_{n+1}$ for all ${0\leq n <N}$. 
	
	We call $p_0$ and $p_N$ the \emph{starting point} and \emph{endpoint} of $p$, respectively. When the indexing is not specified, the endpoint of $p$ will be denoted $p_\top$, and the starting point $p_\bot$. 
	
	We say $p$ \emph{inhabits} $V\in\Opens X$ if $p_n\sqleq V$ for some index $n$. We say $p$ \emph{lands in} $V$ if its endpoint is contained in $V$.
	
	Given paths $p=(p_n)_{n=0}^N$ and $q=(q_m)_{m=0}^M$ with $p_N = q_0$, we define their \emph{concatenation} as the path $q\cdot p$ whose steps are defined as $(q\cdot p)_k = p_k$ for $0\leq k\leq N$ and $(q\cdot p)_k=q_{k-N}$ for $N\leq k \leq M+N$. 
\end{definition}
\begin{remark}
	To avoid degeneracy, we shall always assume that every step of a path is \emph{non-empty}. In an ordered locale satisfying~\eqref{axiom:empty} this is true as soon as at least one of the steps is non-empty. 
\end{remark}
\begin{definition}\label{definition:refinement of paths}
	We say a path $q$ \emph{refines} a path $p$ if every step of $p$ contains some step of $q$: for every $n$ we can find $m$ such that $q_m\sqleq p_n$. In that case we write~$q \refines p$. 
\end{definition}
\noindent\begin{minipage}{0.64\textwidth}
	\begin{remark}
		Note that $\refines$ is just an instance of the \emph{upper order}, seen in \cref{section:ordered spaces as ordered locales}. That $q\refines p$ does not mean that every $q_m$ has to be contained in some $p_n$. For instance, any path $p$ can be refined by specifying further steps $p_\top \Leq A$ or~$B\Leq p_\bot$. Conversely, a given step~$p_n$ can contain multiple steps of $q$. See \cref{figure:refined path} for intuition. Interpreting a path $p$ as an approximation of a point-wise curve $\gamma$, the fact that $q\refines p$ is interpreted as $q$ being a more accurate approximation of~$\gamma$. Alternatively, we can think of~$q$ as a tighter ``restriction'' of where information can flow, a curve $\gamma$ being an idealised ``infinitely tight'' flow.
	\end{remark}
\end{minipage}%
\hfill%
\begin{minipage}{.33\textwidth}\centering
	\input{tikz/inkscape-export/fig-03}
	\captionof{figure}{Typical $q\refines p$.}
	\label{figure:refined path}
\end{minipage}

\begin{lemma}\label{lemma:properties of refinement relation}
	The refinement relation $\refines$ is a preorder and respects concatenation:
		\begin{enumerate}[label = (\alph*)]
			\item $p\refines p$;
			\item $r\refines q\refines p$ implies $r\refines p$;
			\item if $q\refines p$ and $q'\refines p'$, then $(q'\cdot q)\refines (p'\cdot p)$. 
		\end{enumerate}
\end{lemma}
\begin{proof}
	The first two statements follow since $\refines$ is the upper order induced by the inclusion relation $\sqleq$. For (c), the steps of $p'\cdot p$ corresponding to $p$ are refined by~$q$, and those corresponding to $p'$ are refined by $q'$.
\end{proof}

The following construction shows how to canonically refine a path $p$, given some open subregion in its endpoint. 

\begin{construction}\label{construction:path restriction}
	Consider a path $p=(p_n)_{n=0}^N$ in a parallel ordered locale $(X,\Leq)$, and take any open region $W\sqleq p_N$ in the endpoint. We are going to define a new path: $p|_W$, called the \emph{restriction of $p$ to $W$,} that refines $p$ and has endpoint $W$. Each step of $p|_W$ shall be recursively defined, based on the following procedure using~\eqref{axiom:wedge-}:
	\[
	\begin{tikzcd}[column sep=-.1cm, row sep=-.15cm]
		\cdots & \Leq & {p_{N-2}\wedge\Down(p_{N-1}\wedge\Down W)} & \Leq & {p_{N-1}\wedge\Down W} & \Leq & W \\
		&& \rotatebox[origin=c]{-90}{$\sqleq$} && \rotatebox[origin=c]{-90}{$\sqleq$} && \rotatebox[origin=c]{-90}{$\sqleq$} \\
		\cdots & \Leq & {p_{N-2}} & \Leq & {p_{N-1}} & \Leq & {p_N}.
	\end{tikzcd}
	\]
	More formally, define the endpoint by $(p|_W)_N:= W$, and then use the recursive formula
		\[
			(p|_W)_n := p_n\wedge \Down(p|_W)_{n+1},
		\]
	for $0\leq n < N$. We get $(p|_W)_n \Leq (p|_W)_{n+1}$ by~\eqref{axiom:wedge-} (via \cref{lemma:wedge iff strong wedge}). Since the locale is parallel ordered, it follows that all steps of $p|_W$ are non-empty if and only if $W$ itself is non-empty. Lastly, $p|_W\refines p$ since $(p|_W)_n\sqleq p_n$ by construction. Note that if $W$ is the entire endpoint we get $p|_{p_\top}=p$.
	
	Using axiom~\eqref{axiom:wedge+} we get an analogous construction, where for $V\sqleq p_0$ we get a refinement of $p$ with starting point $V$ that we denote by $p|^V$. 
\end{construction}
\begin{figure}[t]\centering
	\begin{subfigure}[b]{0.55\textwidth}\centering
		\input{tikz/inkscape-export/fig-08}
		\caption{Path restriction in Minkowski-type space.}
		\label{figure:path restriction Minkowski}
	\end{subfigure}\hfil
	\begin{subfigure}[b]{0.4\textwidth}\centering
		\input{tikz/inkscape-export/fig-09}
		\caption{Path restriction in vertical-$\mathbb{R}^2$.}
	\end{subfigure}
	\caption{Examples of \cref{construction:path restriction}.}
	\label{figure:path restriction}
\end{figure}

\begin{remark}
	Note that the construction only works on paths defined in terms of finite sequences, at least without resorting to a higher form of induction.
	
	It is tempting to define $p|_W$ by the simpler formula $(p|_W)_n = p_n\wedge \Down W$. Unfortunately, this does not always define a path. As \cref{figure:path restriction}(a) shows, it can happen that intersecting a path step-wise with $\Down W$ does not preserve the causal order~$\Leq$. Here $p_0\Leq p_1\Leq p_\top\sqgeq W$, where $p_1$ is the disjoint union of the two regions in the middle, but $p_0\not\Leq p_1\wedge \Down W$ since $p_0\not\sqleq \Down (p_1\wedge \Down W)$.

	Similarly, it might seem desirable to consider only paths $p$ whose steps $p_n$ are \emph{connected} (e.g.~as done in the unordered setting of \cite{kennison1989WhatFundamentalGroup}). However, we can see in \cref{figure:path restriction}(b) that \cref{construction:path restriction} does not preserve this property: $p_0$ is connected but $(p|_W)_0$ is not. A more suitable option is to only consider paths consisting of \emph{convex} regions, which \cref{construction:path restriction} does preserve, but for our purposes here we do not find it necessary to do so.
\end{remark}

In the next section we will use the path refinement relation $\refines$ together with path restriction to state the definition of ``causal coverage'' in any ordered locale. First we obtain some elementary properties of these structures.

\begin{lemma}\label{lemma:restriction of paths is functorial}
	Let $p$ be a path, with $W\sqleq V\sqleq p_\top$. Then $p|_W = (p|_V)|_W$.
\end{lemma}
\begin{proof}
	By construction, we get $((p|_V)|_W)_n\sqleq (p|_V)_n$. Using this, we find that
		\[
			((p|_V)|_W)_n 
			=
			p_n \wedge \Down (p|_V)_{n+1} \wedge \Down ((p|_V)|_W)_{n+1}
			=
			p_n \wedge \Down ((p|_V)|_W)_{n+1}.
		\]
	From this equation we see that the result will follow by induction as soon as we can show $((p|_V)|_W)_N = (p|_W)_N$, where $N$ is the top index of $p$. But this equation is just the statement that $W=W$, so we are done.
\end{proof}

The following result shows that the path restriction operation loses no information about a path, as long as the entire endpoint is covered.
\begin{lemma}\label{lemma:join over restrictions}
	Suppose that $p=(p_n)_{n=0}^N$ is a path, and take an open cover $(W_i)_{i\in I}$ of the endpoint $p_N$. Then for all $0\leq n\leq N$:
		\[
			p_n = \bigvee_{i\in I} \left( p|_{W_i} \right)_n.
		\]
\end{lemma}
\begin{proof}
	The statement holds for $n=N$ since $(p|_{W_i})_N=W_i$. We prove the remaining cases via induction. Suppose the equation holds for the index $n+1$. We then calculate:
		\[
			\bigvee_{i\in I}(p|_{W_i})_{n} 
			=
			p_{n}\wedge\bigvee_{i\in I}\Down(p|_{W_i})_{n+1}
			=
			p_{n}\wedge \Down \bigvee_{i\in I}(p|_{W_i})_{n+1}
			= p_{n}\wedge \Down p_{n+1}
			= p_{n}.
		\]
	The first step is infinite distributivity applied to the definition of $(p|_{W_i})_n$, the second is due to~\eqref{axiom:LV}, the third is the induction hypothesis, and the final equality follows by ${p_{n}\Leq p_{n+1}}$.
\end{proof}

Under certain conditions, the path restriction operation preserves refinements.
\begin{lemma}\label{lemma:path restriction preserves refinement}
	Let $q\refines p$ with the property that
		\begin{quote}
			if $q_m\sqleq p_n\neq p_\top$, then $\exists k\geq m: q_k\sqleq p_{n+1}$,
		\end{quote}
	and take non-empty $W\sqleq q_\top \sqleq p_\top$. Then $q|_W\refines p|_{W}$.
\end{lemma}
The condition says that for every step $q_m$ inhabiting $p_n$, there exists some $q_k$ in the future that inhabits the next step $p_{n+1}$.
\begin{proof}
	Write explicitly $p=(p_n)_{n=0}^N$ and $q=(q_m)_{m=0}^M$. We need to prove that for every $n$ there exists an $m$ such that $(q|_W)_m\sqleq (p|_W)_n$. For the top indices we get $(q|_W)_M = W = (p|_W)_N$. Similarly, since $q\refines p$ we get $q_m\sqleq p_{N-1}$ for some~$m$, and since $(q|_W)_m\Leq W$ it follows that ${(q|_W)_m\sqleq p_{N-1}\wedge \Down W = (p|_W)_{N-1}}$. In fact, this inclusion holds for every $m$ such that $q_m$ is contained in $p_{N-1}$.
	
	Pose as an induction hypothesis (case $n+1$) that for all $m$ with $q_m\sqleq p_{n+1}$ we have $(q|_W)_m \sqleq (p|_W)_{n+1}$. We claim the same holds for $n$. If $q_m\sqleq p_n$ then by assumption there exists $k\geq m$ such that ${q_m\Leq q_k\sqleq p_{n+1}}$. We thus get
		\begin{align*}
			(q|_W)_m
			&:=
			q_m \wedge \Down (q|_W)_{m+1}
			\\&\sqleq 
			p_n \wedge \Down (q|_W)_k &&\text{($q_m\sqleq p_n$ and $m\leq k$)}
			\\&\sqleq 
			p_n \wedge \Down (p|_W)_{n+1}, &&\text{(induction hypothesis)}
		\end{align*}
	so it follows that $(q|_W)_m\sqleq (p|_W)_n$. The fact that $q|_W\refines p|_{W}$ now follows directly from~${q\refines p}$.
\end{proof}

\begin{remark}[Speculation]
	In a locale we consider completely prime filters with respect to the inclusion order $\sqleq$ as modelling the would-be infinitesimal points. Here we could attempt a similar construction, where filters with respect to $\refines$ model the image of a would-be point-wise curve.
	
	Similar ideas appear in e.g.~\cite{kennison1989WhatFundamentalGroup}, where they define a fundamental group(oid) of a locale. Within the ordered locale framework, this might lead to a \emph{fundamental category} construction. These ideas could in turn be compared to the study of spacetimes via point-wise fundamental categories \cite{gunther2021CategoricalTopologicalStructure}. The basis would be the \emph{path category} $\Pi_{\Leq}(X)$ of $(X,\Leq)$, which is the category whose objects are opens $U\in\Opens X$, and where arrows $U\to V$ consist of causal paths $p$ that start at $U$ and end at~$V$. This category is $\Ord$-enriched, in the sense that the sets $\Pi_{\Leq}(X)(U,V)$ of arrows are equipped with the refinement preorder $\refines$, which respects the categorical structure by \cref{lemma:properties of refinement relation}. We anticipate a construction $\Pi_{\Leq}(X)\mapsto \Pi_\leq(\pt(X))$, where the right hand side is the point-wise path category of monotone continuous functions of $\pt(X)$, generalising the functor $X\mapsto \pt(X)$.
\end{remark}

\section{Causal coverages of ordered locales}
We are now ready to define the point-free analogue of ``causal coverage'' discussed in the introduction of this section. First we need to refine our definition of refinements.
\begin{definition}\label{definition:locally refines}
	A \emph{local past refinement} of a path $p$ consists of a family $(q^j)_{j\in J}$ of paths such that:
		\begin{enumerate}[label= (\roman*)]
			\item $\bigvee_{j\in J} q^j_\top = p_\top$;
			\item $q^j\refines p|_{q^j_\top}$ for all $j\in J$.
		\end{enumerate}
	In that case we write $(q^j)_{j\in J}\refines p$.
	
	A \emph{(global) past refinement} is a local past refinement where $J$ contains one element. Explicitly, for a path $p$ this consists of a refinement $q\refines p$ such that $q_\top = p_\top$. 
\end{definition}

\begin{remark}
	The definition of a past refinement specifically excludes refinements $q$ that propagate beyond the original path $p$. Moreover, the fact that $\bigvee q^j_\top = p_\top$ guarantees that all information at the endpoint $p_\top$ is reached by the refinements~$q^j$. We use this to define past causal coverages. For future causal coverages, one defines (local) future refinements dually.
\end{remark}

\begin{definition}\label{definition:causal coverage locale}\label{definition:causal coverage from ordered locale}
	Let $(X,\Leq)$ be a parallel ordered locale with~\eqref{axiom:LV}, and pick opens $U,A\in\Opens X$. We say that $A$ \emph{covers $U$ from below} if $A\sqleq \Down U$, and every path landing in $U$ can be \emph{locally refined} to inhabit $A$:
		\begin{quote}
			for every path $p$ landing in $U$ there exists a local past refinement $(q^j)\refines p$ such that each $q^j$ inhabits $A$.
		\end{quote}
	
	An analogous definition can be made for regions that cover $U$ from \emph{above}. Denote by $\CovLeq^\pm(U)\subseteq \Opens X$ the set of all opens that cover $U$ from above/below. We call $\CovLeq^\pm$ the \emph{(future/past) causal coverage} induced by $\Leq$.
\end{definition}

\noindent\begin{minipage}{0.65\textwidth}
	\begin{remark}
		Very explicitly, if $A\in \CovLeq^-(U)$ and $p$ is a path landing in $U$, this means we can find some family~$(q^j)$ of paths, together with an open cover $(W_j)_{j\in J}$ of the endpoint $p_\top$, such that $q^j_\top = W_j$ and $q^j$ refines the path $p|_{W_j}$, and such that for every $j\in J$ there is some index $a_j$ such that $q^j_{a_j}\sqleq A$.
	\end{remark}
\begin{example}
	In \cref{figure:non-globally refinable} we illustrate an example of a path in vertical-$\mathbb{R}^2$ that is not globally past refinable to inhabit $A$, but only locally. Nevertheless, we want to think of $A$ as covering $U$ from below in this picture. This shows the local refinements are necessary for a proper definition of $\Cov^\pm$.
\end{example}	
\end{minipage}%
\hfill%
\begin{minipage}{.3\textwidth}\centering
	\input{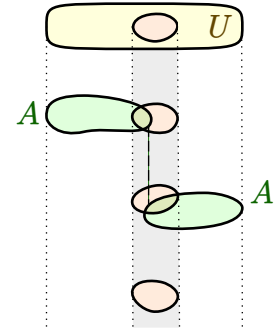}
	\captionof{figure}{Non-globally refinable path.}
	\label{figure:non-globally refinable}
\end{minipage}

\vspace*{5pt}In the following we collect some important properties of $\CovLeq^\pm$. See \cref{figure:properties of coverage from locale} for intuition. In proving these properties we have the definition of a Grothendieck topology (\cref{definition:grothendieck topology}) in mind.

\begin{lemma}\label{lemma:properties of coverage from locale}
	For any parallel ordered locale $(X,\Leq)$ satisfying~\eqref{axiom:LV}, the causal coverages have the following properties:
		\begin{enumerate}
			\item[(a)] $U\in\CovLeq^\pm(U)$;
			\item[(b)] $\Down U \in \CovLeq^-(U)$ and $\Up U \in\CovLeq^+(U)$;
			\item[(c)] if $B\in \CovLeq^\pm(A)$ and $A\in \CovLeq^\pm(U)$, then $B\in \CovLeq^\pm(U)$;
			\item[(d$^-$)] if $A\in \CovLeq^-(U)$ and $W\sqleq U$, then $A\wedge \Down W\in\CovLeq^-(W)$;
			\item[(d$^+$)] if $B\in \CovLeq^+(U)$ and $W\sqleq U$, then $B\wedge \Up W\in\CovLeq^+(W)$;
			\item[(e)] if $A\in \CovLeq^-(U)$ then $A\Leq U$, and if $B\in\CovLeq^+(U)$ then $U\Leq B$;
			\item[(f)] $\CovLeq^\pm(\varnothing)=\{\varnothing\}$.
		\end{enumerate}
\end{lemma}
\begin{proof}
	As usual, we only prove the past versions. Statement (a) is trivial: if $p$ is a path with endpoint in $U$, then $p$ already inhabits $U$. The claim in (b) follows similarly.
	
	For (c), note first that if $B\in \CovLeq^-(A)$ and $A\in\CovLeq^-(U)$ then $B\sqleq \Down A\sqleq \Down \Down U = \Down U$. Now let $p=(p_n)_{n=0}^N$ be a path landing in $U$. Since $A$ covers $U$ we get an open cover $(W_j)_{j\in J}$ of $p_N$, together with refinements $q^j$ of $p|_{W_j}$ with endpoints $q^j_{M_j}=W_j$. These refinements inhabit $A$, so there are indices $a_j$ such that $q^j_{a_j}\sqleq A$. Denote $\hat{q}^j = (q^j_m)_{m=0}^{a_j}$ and $\check{q}^j = (q^j_m)_{m=a_j}^{M_j}$. In turn, since $B$ covers $A$ we get open covers $(V^j_i)_{i\in I_j}$ of $q^j_{a_j}$, together with refinements $r^{ji}$ of $\hat{q}^j|_{V^j_i}$ that inhabit $B$ and have endpoint $V^j_i$. The concatenation $\check{q}^j|^{V^j_i}\cdot r^{ji}$ then forms a refinement of $q^j$, and hence of $p|_{W_j}$, that inhabits $B$. Moreover, using (the dual of) \cref{lemma:join over restrictions} we find
		\[
			\bigvee_{j\in J}\bigvee_{i\in I_j}(\check{q}^j|^{V^j_i})_{M_j}
			=
			\bigvee_{j\in J}q^j_{M_j}
			= 
			\bigvee_{j\in J}W_j
			=
			p_N.
		\]
	This gives the desired open cover and refinements to obtain $B\in\CovLeq^-(U)$.	
	
	For (d$^-$), clearly we have $A\wedge \Down W\sqleq \Down W$. Further, if $p$ is a path landing in~$W$, it also lands in $U$. Thus we can find local refinements $q^j$ that inhabit $A$ and land in the endpoint of $p$. But since the endpoint of $q^j$ is contained in $W$, each of its steps is contained in $\Down W$, so $q^j$ must inhabit $A\wedge \Down W$. 
	
	For (e), if $A\in\CovLeq^-(U)$ we already know $A\sqleq \Down U$. Now consider the path $p$ consisting of the single step:~$U$. Hence there exists an open cover $(U_i)_{i\in I}$ of~$U$, together with paths $q^i$ that inhabit $A$ and have endpoint $U_i$. Hence, for each $i\in I$, there exists some index $a_i$ such that $q^i_{a_i}\sqleq A$, which implies $U_i\sqleq \Up A$. In total, this gives $U=\bigvee_{i\in I} U_i\sqleq \Up A$. That $A\Leq U$ now follows by~\eqref{axiom:cones give order}.
	
	Lastly, for (f), if $A\in\CovLeq^-(\varnothing)$, then $A\sqleq \Down \varnothing=\varnothing$, so $\CovLeq^-(\varnothing)\subseteq \{\varnothing\}$, and the converse inclusion follows by (a).
\end{proof}
\begin{figure}[t]\centering
	\begin{subfigure}[b]{0.4\textwidth}\centering
		\input{tikz/inkscape-export/fig-18-transitive}
		\caption{Transitivity.}
	\end{subfigure}\hfil
	\begin{subfigure}[b]{0.4\textwidth}\centering
		\input{tikz/inkscape-export/fig-18-pullback}
		\caption{Pullback stability.}
	\end{subfigure}
	\caption{Illustrations of the properties (c) and (d$^-$) of $\CovLeq^-$ from \cref{lemma:properties of coverage from locale}.}
	\label{figure:properties of coverage from locale}
\end{figure}

The following gives an analogue of the fundamental axiom~\eqref{axiom:V}.

\begin{lemma}\label{lemma:coverage locale satisfies Cov-V weak}
	In a parallel ordered locale $(X,\Leq)$ with~\eqref{axiom:LV} we have:
		\[
			\forall i\in I: A_i\in\CovLeq^\pm(U_i)
			\quad \text{implies} \quad 
			\bigvee_{i\in I}A_i\in \CovLeq^\pm\left(\bigvee_{i\in I}U_i\right).
		\]
\end{lemma}
\begin{proof}
	Suppose that $A_i\in \CovLeq^-(U_i)$ for all $i\in I$, for some families $(A_i)_{i\in I}$ and $(U_i)_{i\in I}$ of opens. In the case that $\bigvee_{i\in I}U_i$ is empty, the result follows by \cref{lemma:properties of coverage from locale}(f). For the rest of the proof we can therefore assume that the join, and hence the index $I$, are non-empty. In that case, using monotonicity of the localic cones we get $\bigvee_{i\in I}A_i\sqleq \Down\bigvee_{i\in I}U_i$.
	
	We are left to show any path that lands in $\bigvee_{i\in I}U_i$ can be locally refined to inhabit $\bigvee_{i\in I} A_i$. For that, take a path $p=(p_n)_{n=0}^N$ that lands in $\bigvee_{i\in I}U_i$. Denote by $J\subseteq I$ the set of indices $j$ for which $W_j:=p_N\wedge U_j$ is non-empty. Since $p_N$ is non-empty, the set $J$ is non-empty, and $p_N = \bigvee_{j\in J} W_j$. The restrictions $p|_{W_j}$ land in $U_j$. Since $A_j$ covers $U_j$ from below, we get an open cover $(V^j_k)_{k\in K_j}$ of $W_j$, together with refinements $q^{jk}\refines (p|_{W_j})|_{V^j_k}$ that inhabit $A_j$ and have endpoint $V^j_k$. Using \cref{lemma:restriction of paths is functorial} we get $(p|_{W_j})|_{V^j_k} = p|_{V^j_k}$. This shows that the $q^{jk}$'s provide the desired local refinement of $p$. 
%
%
\end{proof}

Put differently, for any family $(U_i)_{i\in I}$ of opens the function
	\[
		\prod_{i\in I}\CovLeq^\pm(U_i)\longrightarrow\Opens X;\qquad (A_i)_{i\in I}\longmapsto \bigvee_{i\in I} A_i
	\]
lands in $\CovLeq^\pm\left(\bigvee_{i\in I} U_i\right)$. The following result sharpens this further:~the image of this function in fact gives the whole set $\CovLeq^\pm\left(\bigvee_{i\in I} U_i\right)$. This means that the causal cover of a join can always be written as the union of causal covers of the constituents.

\begin{lemma}\label{lemma:coverage locale satisfies Cov-V strong}
	In a parallel ordered locale $(X,\Leq)$ with~\eqref{axiom:LV} we have:
	\[\tag{$\Cov$-$\vee$}\label{axiom:Cov-V strong}
		\CovLeq^\pm\left(\bigvee_{i\in I} U_i\right) = \left\{\bigvee_{i\in I} A_i : (A_i)_{i\in I}\in\prod_{i\in I} \CovLeq^\pm(U_i)\right\}.
	\]
\end{lemma}
\begin{proof}
	The inclusion from right-to-left is precisely \cref{lemma:coverage locale satisfies Cov-V weak}. For the converse, take $A\in \CovLeq^-\left(\bigvee_{i\in I} U_i\right)$. Given \cref{lemma:properties of coverage from locale}(f) we can assume without loss of generality that the indexing set is non-empty. We need to show that $A$ can be written as a join of open regions that cover the $U_i$'s from below. By \cref{lemma:properties of coverage from locale}(d$^-$) we get $A\wedge \Down U_i\in \CovLeq^-(U_i)$ for every $i\in I$. Using axiom~\eqref{axiom:LV} we find the desired expression: $A = A\wedge\Down \bigvee_{i\in I} U_i= \bigvee_{i\in I} A\wedge \Down U_i$. The proof for future covers is dual.
\end{proof}


\subsection{Causal coverages as Grothendieck topologies}
In this section we explore the fact that the properties of $\CovLeq^\pm$ showcased in \cref{lemma:properties of coverage from locale} share a remarkable resemblance to the definition of a \emph{Grothendieck topology}. These are very important mathematical structures in category theory, and relate to both geometry and logic~\cite{maclane1994SheavesGeometryLogic}. A Grothendieck topology is a categorical generalisation of the notion of \emph{open covers}, and they are used to study how ``local data'' is amalgamated into ``global data.'' In the causal setting we shall interpret it as dictating how ``past data'' evolves into ``future data.''

For convenience, we first recall the definition of a Grothendieck topology in full generality on a category~$\cat{C}$. To state it, we need the technology of sieves.

\begin{definition}\label{definition:sieve}
	A \emph{sieve} $R$ on an object $C\in\cat{C}$ of a category is a set of arrows with codomain $C$, such that if $D\xrightarrow{g}C\in R$ and $E\xrightarrow{h}D$ is an arbitrary arrow in~$\cat{C}$, then $E\xrightarrow{g}D\xrightarrow{h}C\in R$. 
\end{definition}

For any object $C\in\cat{C}$ we get the \emph{maximal sieve}:
	\[
		t_C:= \left\{\text{all arrows with codomain $C$}\right\}.
	\]
It is easy to see that this in fact defines a sieve. If $R$ is a sieve on $C\in\cat{C}$ and $h\colon D\to C$ is an arrow, the \emph{pullback} of $R$ along $h$ is the sieve on the domain~$D$ defined as:
	\[
		h^\ast(R):=\left\{E\xrightarrow{~f~}D: h\circ f\in R \right\}.
	\]
	
\begin{definition}[{\cite[\S III.2]{maclane1994SheavesGeometryLogic}}]\label{definition:grothendieck topology}
	A \emph{Grothendieck topology} on a category $\cat{C}$ is a function $J$ that assigns to each object $C\in\cat{C}$ a collection $J(C)$ of sieves on $C$, called \emph{covering sieves}, such that the following axioms are satisfied:
	\begin{enumerate}
		\item[(i)] the maximal sieve $t_C$ is in $J(C)$;
		\item[(ii)] if $S\in J(C)$ then $h^\ast(S)\in J(D)$ for any arrow $h\colon D\to C$;
		\item[(iii)] if $S\in J(C)$ and $R$ is a sieve on $C$ such that $h^\ast(R)\in J(D)$ for all $h\colon D\to C$ in $S$, then $R\in J(C)$.
	\end{enumerate}
	A category $\cat{C}$ equipped with a Grothendieck topology $J$ is called a \emph{site}.
\end{definition}

\begin{remark}[Intuition]
	As this definition is quite abstract, it helps to unpack it in the setting where $\cat{C}$ is the frame $\Opens X$ of open subsets of some locale (or space). Thus the objects are the opens of $X$, and an arrow $f\colon V\to U$ is just an inclusion $f\colon V\sqleq U$, which is unique if it exists. In the following we shall be somewhat informal in conflating the arrow $f$ with the open subset $V\sqleq U$ itself. In so doing, the maximal sieve on $U\in\Opens X$ can be identified with the principal ideal generated by~$U$:
		\[
			t_U= \{V\in\Opens X:V\sqleq U\}.
		\]
	In general, a sieve $R$ on $U\in\Opens X$ is a collection of subregions of $U$ that are down-closed under the inclusion relation $\sqleq$. If $h\colon V\sqleq U$ is an arrow and $R$ is a sieve on $U$, then its pullback can be calculated as
		\[
			h^\ast(R) = \{V\wedge W: W\in R\}.
		\]
	With this in mind, we can reinterpret the axioms of a Grothendieck topology in \cref{definition:grothendieck topology} for the \emph{canonical coverage} $J$ on $\Opens X$, which is defined as:
		\[
			R\in J(U)
			\qquad\text{if and only if}\qquad
			\bigvee R = U.
		\]
	In other words, the covering sieves $R\in J(U)$ are precisely the down-closed open covers of $U$. It is straightforward to see that the axioms hold:
	\begin{enumerate}
		\item[(i)] $\bigvee\{V\in\Opens X: V\sqleq U\} = U$;
		\item[(ii)] if $\bigvee S = U$ and $V\sqleq U$ then $V\wedge \bigvee S = V$;
		\item[(iii)] if $\bigvee S = U$ and $R$ is a sieve on $U$ such that $V\wedge \bigvee R = V$ for all $V\in S$, then $\bigvee R = U$. This is true because: $\bigvee S = \bigvee_{V\in S} (V\wedge \bigvee R) = \bigvee S \wedge \bigvee R$.
	\end{enumerate}
	
	As explained more deeply in e.g.~\cite[\S II.2.11]{johnstone1982StoneSpacesa} or \cite{ball2014ExtendingSemilatticesFrames}, the role of the coverage $J$ is to say which sieves $R$ ``cover'' $U$ without the use of an actual join operator $\bigvee$. Thus sites can in some sense be thought of as categorical generalisations of frames. In \emph{formal topology}, the notion of ``coverage'' is taken as fundamental~\cite{sambin2003PointsFormalTopology}.
\end{remark}

Motivated by~\cref{lemma:properties of coverage from locale}, we now want to interpret the causal coverages $\CovLeq^\pm$ of an ordered locale as Grothendieck topologies $J^\pm$. To motivate the definition, we compare this structure to the canonical coverage $J$ on $\Opens X$. There, a sieve ``$R$ covers $U$'' if $R$ forms an open cover for $U$. In the causal setting we want to say ``$R$ \emph{causally} covers $U$'' if $R$ forms an open cover for some $A\in \CovLeq^-(U)$. 

But this poses a problem, since with respect to the inclusion order $\sqleq$ the sieves~$R$ on $U$ by definition only contain opens $V\in R$ that are subregions $V\sqleq U$, while the interesting causal covering behaviour of $\CovLeq^\pm$ happens when the regions $V$ are disjoint from $U$. We thus need some means to say when a sieve $R$ is covering for $U$ while $R$ itself is not strictly a sieve on $U$ itself. For this, we propose the following generalisation of Grothendieck topologies.

\begin{definition}\label{definition:modified grothendieck topology}
	Let $\cat{C}$ be a category, together with an endofunctor $T\colon \cat{C}\to \cat{C}$ and a natural transformation $\mu\colon T^2\to T$. A \emph{$(T,\mu)$-Grothendieck topology} on $\cat{C}$ is a function $J$ that assigns to each object $C\in\cat{C}$ a collection $J(C)$ of sieves on~$T(C)$, called \emph{covering sieves}, such that the following axioms are satisfied:
	\begin{enumerate}
		\item[(i)] the maximal sieve $t_{T(C)}$ is in $J(C)$;
		\item[(ii)] if $S\in J(C)$ then $T(h)^\ast(S)\in J(D)$ for any arrow $h\colon D\to C$;
		\item[(iii)] if $S\in J(C)$ and $R$ is a sieve on $T(C)$ such that $(\mu_C\circ T(h))^\ast(R)\in J(D)$ for all $h\colon D\to T(C)$ in $S$, then $R\in J(C)$.
	\end{enumerate}
\end{definition}

Setting $T=\id_\cat{C}$ and $\mu = \id_T$ of course returns the ordinary notion of Grothendieck topology on $\cat{C}$, so this is a genuine generalisation. The natural transformation $\mu$ is needed only to make sense of axiom~(iii), where for $h\colon D\to T(C)$ in a covering sieve $S\in J(C)$ we need to be able to construct the sieve $(\mu_C\circ T(h))^\ast(R)$ on $T(D)$.

We are not aware of these types of structures being studied in the literature. It begs a lot of questions: what are sheaves in this setting, and do they form toposes just as in the ordinary setting? Are there any other natural examples, besides causal coverages? This we (unfortunately have to) leave for future research. 

The following result shows that the causal coverage $\CovLeq^-$ of an ordered locale gives rise to a canonical $\Down$-Grothendieck topology, and also explains the intuition of \cref{definition:modified grothendieck topology}.

\begin{figure}[t]
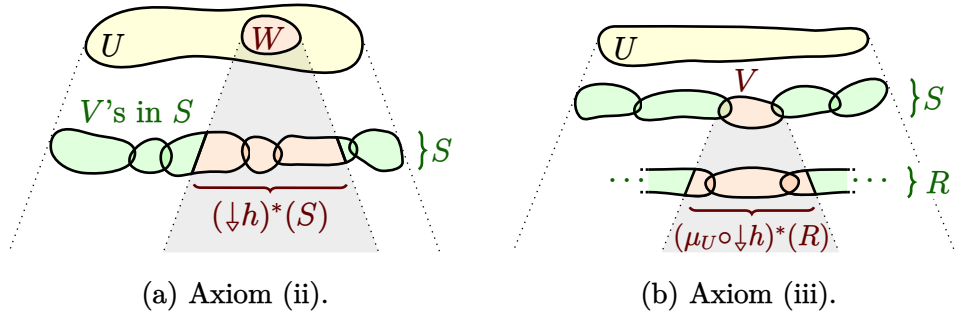
\centering
	\begin{subfigure}[b]{0.45\textwidth}\centering
		\input{tikz/inkscape-export/fig-29-ii}
		\caption{Axiom (ii).}
	\end{subfigure}
	\begin{subfigure}[b]{0.45\textwidth}\centering
		\input{tikz/inkscape-export/fig-29-iii}
		\caption{Axiom (iii).}
	\end{subfigure}
	\caption{Illustrations of the axioms of a $\Down$-Grothendieck topology for $J^-$.}
	\label{figure:properties of causal grothendieck topology}
\end{figure}

\begin{theorem}\label{theorem:causal grothendieck topology}
	Let $(X,\Leq)$ be a parallel ordered locale with~\eqref{axiom:LV}. The causal coverage $\CovLeq^-$ defines a $\Down$-Grothendieck topology~$J^-$ on~$\Opens X$ via
		\[
			R\in J^-(U)
			\qquad\text{if and only if}\qquad
			\bigvee R\in \CovLeq^-(U).
		\]
\end{theorem}
\begin{proof}
	We check the axioms (see \cref{figure:properties of causal grothendieck topology} for intuition):
	\begin{enumerate}
	\item[(i)] that $t_{\Down U}\in J^-(U)$ just means $\bigvee \{V\in\Opens X: V\sqleq \Down U\} = \Down U \in \CovLeq^-(U)$, which holds by \cref{lemma:properties of coverage from locale}(b);
	
	\item[(ii)] let $S\in J^-(U)$ and take an arrow $h\colon W\sqleq U$. We need to show that $(\Down h)^\ast(S)\in J^-(W)$. As above, we identify arrows $V\sqleq \Down U$ in $S$ with their codomains, which allows us to write the pullback as
		\[
			(\Down h)^\ast(S) = \{V\wedge \Down W: V\in S\}.
		\]
	Indeed, if $V\in S$ then $V\wedge \Down W\in (\Down h)^\ast(S)$ since $V\wedge \Down W\sqleq V\in S$, and~$S$ is downwards closed with respect to the inclusion relation. Conversely, if $V\in (\Down h)^\ast(S)$ then the composition $V\sqleq \Down W\sqleq \Down U$ is in $S$, and in particular we get $V=V\wedge \Down W$.
	Thus taking the join over $(\Down h)^\ast(S)$ and using infinite distributivity we see
		\[
			(\Down h)^\ast(S) \in J^-(W)
			\qquad\text{if and only if}\qquad
			\Down W\wedge \bigvee S\in \CovLeq^-(W),
		\]
	and the right hand side holds by \cref{lemma:properties of coverage from locale}(d$^-$).
	
	\item[(iii)] let $S\in J^-(U)$, and let $R$ be a sieve on $\Down U$ with the property that for every $h\colon V\sqleq \Down U$ in $S$ we get $(\mu_U\circ \Down h)^\ast(R)\in J^-(V)$. In this case the natural map $\mu$ is just the family of equations  $\mu_U\colon \Down\Down U = \Down U$. Similar to the above, we can calculate the pullback as
		\[
			(\mu_U\circ \Down h)^\ast(R) = \{W\wedge \Down V: W\in R\}.
		\]
	Translating into the language of causal coverages, we get ${\bigvee S \in \CovLeq^-(U)}$, and $R$ is a sieve such that $\Down V\wedge \bigvee R\in \CovLeq^-(V)$ for every $V\in S$. From~\eqref{axiom:Cov-V strong} and using~\eqref{axiom:LV} with infinite distributivity it then follows that
		\[
			\Down \bigvee S \wedge \bigvee R = \bigvee_{V\in S}\left(\Down V\wedge \bigvee R\right)\in \CovLeq^-\left(\bigvee S\right).
		\]
	By \cref{lemma:properties of coverage from locale}(c) it thus follows that $\Down \bigvee S\wedge \bigvee R\in \CovLeq^-(U)$. Since $R$ is a sieve on $\Down U$ we have $\bigvee R\sqleq \Down U$, and since it contains a region that covers~$U$ from below we must in fact have $\bigvee R\in \CovLeq^-(U)$. \qedhere
	\end{enumerate}
\end{proof}

\begin{remark}
	One might think $J^-$ could define a Grothendieck topology (in the original sense of \cref{definition:grothendieck topology}) on the so-called \emph{Kleisli category} $\Kleisli(\Down)$ induced by the past localic cone monad. Explicitly, this is the category whose objects are opens $U\in \Opens X$, and where there is a unique arrow $V\to U$ precisely if $V\sqleq \Down U$. We denote these arrows by $V\sqleqdown U$. That $U\sqleqdown U$ follows since $U\sqleq \Down U$, and that $W\sqleqdown V\sqleqdown U$ implies $W\sqleqdown U$ follows since $W\sqleq \Down V\sqleq \Down \Down U = \Down U$. Therefore $\Kleisli(\Down)$ is just the set $\Opens X$ equipped with the preorder $\sqleqdown$.
	
	In this preorder, a sieve $R$ on $U$ is indeed some collection of open regions $V\sqleq \Down U$, which at first glance seems like what we were looking for. However, defining $J^-$ as above, we see that axiom~(ii) of a Grothendieck topology can fail. Namely, if $S\in J^-(U)$ is a covering sieve and $h\colon W\sqleqdown U$, this does not force $W\sqleq U$. It can happen that $W$ is in the past of and disjoint from $\bigvee S$, so that $h^\ast(S) = \{V\wedge \Down W:V\in S\}\notin J^-(W)$, since $\Down W\wedge \bigvee S=\varnothing\notin \CovLeq^-(W)$. In order to interpret $\CovLeq^\pm$ as Grothendieck topologies we are therefore seemingly forced to consider the generalised notion of Grothendieck topology in \cref{definition:modified grothendieck topology}.
\end{remark}
%

\begin{remark}
	The fact that the localic cones are monads suggests a refinement of the \cref{definition:modified grothendieck topology} of $(T,\mu)$-Grothendieck topologies on $\cat{C}$ to the setting where $T$ is in fact a monad. This gives the additional data of a natural unit map $\eta\colon \id_\cat{C}\to T$, allowing us to state the additional axiom:
		\begin{enumerate}
			\item[(i$'$)] $(\eta_C)_\ast(t_C)\in J(C)$.
		\end{enumerate}
	Here $t_C$ is the maximal sieve on $C\in\cat{C}$, and we have the \emph{pushforward} sieve on~$T(C)$:
		\[
			(\eta_C)_\ast(t_C):=\{ \eta_C\circ f:f\in t_C\}.
		\]
	
	In the setting of causal coverages, where $J = J^-$ from \cref{theorem:causal grothendieck topology}, we get unit $\eta_U\colon U\sqleq \Down U$, and axiom~(i$'$) translates to $\bigvee (\eta_U)_\ast(t_U) = U\in \CovLeq^-(U)$, which holds by \cref{lemma:properties of coverage from locale}(a). We can interpret this as a compatibility condition between $J^-$ and the canonical Grothendieck topology induced by the frame structure of $\Opens X$.
\end{remark}

%

\section{Causal coverages in spacetimes}
In the following we calculate the behaviour of $\CovLeq^\pm$ in spacetimes, and relate it to the intuition using monotone curves outlined at the start of this chapter. For spacetimes we get two distinct notions of coverages, depending on whether one works with the causal or chronological curves.

\begin{definition}\label{definition:spacetime causal coverages}
	In a smooth spacetime $M$, we define the following notion of coverage. For $A,U\in \Opens M$ we say \emph{$A$ covers $U$ from below with (causal) curves} and write $A\in\Covcaus^-(U)$ if $A\subseteq J^-(U)$ and:
		\begin{quote}
			for every fd causal curve $\gamma\colon [a,b]\to M$ with $\gamma(b)\in U$, either ${\im(\gamma)\cap A \neq \varnothing}$ or $J^-(\gamma(a))\cap A\neq\varnothing$.
		\end{quote}
	There is an analogous chronological definition. We say \emph{$A$ covers $U$ from below with timelike curves} and write $A\in\Covchron^-(U)$ if $A\subseteq I^-(U)$ and:
		\begin{quote}
			for every fd timelike curve $\gamma\colon [a,b]\to M$ with $\gamma(b)\in U$, either ${\im(\gamma)\cap A\neq \varnothing}$ or $I^-(\gamma(a))\cap A\neq\varnothing$.
		\end{quote}
\end{definition}
\noindent\begin{minipage}{0.64\textwidth}
	\begin{remarknumbered}\label{remark:causal cover contained in chronological cover}
		We note that since all timelike curves are causal, it follows immediately that \[\Covcaus^-(U)\subseteq \Covchron^-(U).\]
		The converse inclusion does not generally hold. In \cref{figure:non-intersecting lightlike curve} we have a typical $A\in \Cov^-(U)$, together with a causal curve $\gamma$ that is lightlike as it crosses~$A$. Then $A\setminus\im(\gamma)\in \Covchron^-(U)$, but of course ${A\setminus\im(\gamma)\notin\Covcaus^-(U)}$, since $\gamma$ itself provides a counterexample.
	\end{remarknumbered}
	\vspace{1ex}
\end{minipage}%
\hfill%
\begin{minipage}{.34\textwidth}\centering
	\input{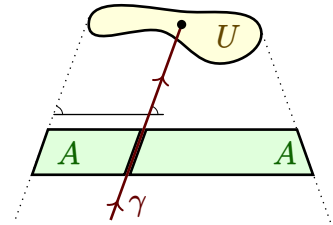}
	\captionof{figure}{Non-intersecting lightlike curve.}
	\label{figure:non-intersecting lightlike curve}
\end{minipage}

In this definition we made use of ``bounded'' causal and timelike curves, defined on closed intervals. In relativity theory it is common to make use of causal curves defined on half-open intervals $(a,b]$ that are \emph{inextendible} into the past. We think of these as causal curves that (when in reverse) propagate indefinitely into the past. For instance, these are used to define \emph{domains of dependence} \cite[Definition~3.1]{minguzzi2019LorentzianCausalityTheory}, which we will meet later in this chapter (\cref{section:domains of dependence}), and \emph{Cauchy surfaces} \cite[Definition~3.35]{minguzzi2019LorentzianCausalityTheory}. To make this rigorous, we adopt \cite[Definition~5.17]{landsman2021FoundationsGeneralRelativity}.

\begin{definition}\label{definition:past inextendible}
	A curve $\gamma\colon (a,b]\to M$ in a smooth spacetime is called \emph{past extendible} if the limit $\lim_{t\searrow a}\gamma(t)$ exists. Otherwise $\gamma$ is called \emph{past inextendible}.
\end{definition}


\begin{lemma}\label{lemma:past inextendible intersects then bounded intersects}
	In a smooth spacetime $M$, we have $A\in\Covcaus^-(U)$ as soon as $A\subseteq J^-(U)$ and:
		\begin{quote}
			every past inextendible causal curve landing in $U$ intersects $A$.
		\end{quote}
\end{lemma}
\begin{proof}
	Let $\gamma\colon [a,b]\to M$ be a causal curve landing in $U$. By e.g.~\cite[Theorem~2.5.7]{chrusciel2020ElementsCausality} we can extend $\gamma$ to a past inextendible causal curve $\overline{\gamma}\colon (c,b]\to M$ such that $\overline{\gamma}|_{[a,b]} = \gamma$. By assumption $\overline{\gamma}^{-1}(A)$ is open and non-empty, so since~$\mathbb{R}$ is regular there exists some closed interval $[n,m]\subseteq \overline{\gamma}^{-1}(A)$. If $n\geq a$ then $\gamma$ already intersects $A$. If $n\leq a$ then $\overline{\gamma}(n)\in A \cap J^-(\gamma(a))$. Since $\gamma$ was arbitrary this gives $A\in\Covcaus^-(U)$.
%
%
\end{proof}
\noindent\begin{minipage}{0.64\textwidth}
	\begin{remarknumbered}\label{remark:difference causal coverage bounded or inextendible}
		An analogous result holds for the chronological coverage. The converse of \cref{lemma:past inextendible intersects then bounded intersects} is not generally true. In \cref{figure:non-intersecting inextendible curve} we have a past inextendible curve $\gamma$ in Minkowski space with one point removed that does not intersect~$A$. Note however that $A\in\Covcaus^-(U)$ does hold, since for every restriction $\gamma|_{[c,b]}$ to a bounded curve we do get $J^-(\gamma(c))\cap A\neq\varnothing$.
	\end{remarknumbered}
\end{minipage}%
\hfill%
\begin{minipage}{.34\textwidth}\centering
		\input{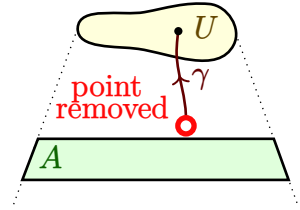}
		\captionof{figure}{Non-intersecting inextendible curve.}
		\label{figure:non-intersecting inextendible curve}
		\vspace{1ex}
\end{minipage}

\vspace*{8pt}In the following we compare the localic and curve-wise definitions of coverages in a spacetime. First we need a small lemma.
\begin{lemma}\label{lemma:path restriction preserves points}
	Let $(S,\leq)$ be an ordered space with open cones, and consider a path $(p_n)_{n=0}^N$ and a chain $x_0\leq \cdots \leq x_N$ with $x_n\in p_n$ for all $n$. If $x_N\in W\subseteq U$ for $W\in\Opens S$, then $x_n\in (p|_W)_n$ for all $n$.
\end{lemma}
\begin{proof}
	Since $x_N\in W = (p|_W)_N$, for an induction proof it suffices to show that if $x_{n+1}\in (p|_W)_{n+1}$ then $x_n \in (p|_W)_n$. But $x_n\leq x_{n+1}$, so using the open cone condition:
		\[
			x_n \in p_n \cap \down x_{n+1}
			\subseteq
			p_n \cap \Down (p|_W)_{n+1}=: (p|_W)_n.\qedhere
		\] 
\end{proof}

\begin{figure}[b]\centering
	\begin{subfigure}[b]{0.4\textwidth}\centering
		\input{tikz/inkscape-export/fig-10}
		\caption{Construction of $W$.}
	\end{subfigure}\hfil
	\begin{subfigure}[b]{0.4\textwidth}\centering
		\input{tikz/inkscape-export/fig-11-curved}
		\caption{Counterexample of converse.}
	\end{subfigure}
	\caption{Illustrations for the proof of \cref{proposition: localic coverage contains chronological coverage}.}
	\label{figure:causal coverage in spacetime}
\end{figure}
\begin{proposition}\label{proposition: localic coverage contains chronological coverage}
	In any smooth spacetime: $\Covchron^-(U)\subseteq \CovLeq^-(U)$.
\end{proposition}
\begin{proof}
	Take a curve-wise covering region $A\in \Covchron^-(U)$, and take a localic path ${p=(p_n)_{n=0}^N}$ landing in~$U$. By definition, each step $p_n$ is non-empty. Starting with $x_N\in p_N\subseteq U$, this allows us to inductively construct a chain of chronologically related points:
	\[
	\begin{tikzcd}[cramped,row sep=tiny]
		{x_0} & {x_1} & \cdots & {x_{N-1}} & {x_N} \\
		{p_0} & {p_1} & \cdots & {p_{N-1}} & {p_N.}
		\arrow["\chron"{marking, allow upside down}, draw=none, from=1-1, to=1-2]
		\arrow["\chron"{marking, allow upside down}, draw=none, from=1-2, to=1-3]
		\arrow["\chron"{marking, allow upside down}, draw=none, from=1-3, to=1-4]
		\arrow["\chron"{marking, allow upside down}, draw=none, from=1-4, to=1-5]
		\arrow["\in"{marking, allow upside down}, draw=none, from=1-1, to=2-1]
		\arrow["\in"{marking, allow upside down}, draw=none, from=1-2, to=2-2]
		\arrow["\Leq"{marking, allow upside down}, draw=none, from=2-1, to=2-2]
		\arrow["\Leq"{marking, allow upside down}, draw=none, from=2-2, to=2-3]
		\arrow["\in"{marking, allow upside down}, draw=none, from=1-5, to=2-5]
		\arrow["\in"{marking, allow upside down}, draw=none, from=1-4, to=2-4]
		\arrow["\Leq"{marking, allow upside down}, draw=none, from=2-3, to=2-4]
		\arrow["\Leq"{marking, allow upside down}, draw=none, from=2-4, to=2-5]
	\end{tikzcd}
	\]
	By definition of the chronology relation, from this chain we obtain a timelike curve $\gamma \colon [a,b]\to M$ that starts at $x_0$, passes each $x_n$, and ends in $x_N$. By hypothesis, there either exists $t\in [a,b]$ such that $\gamma(t)\in A$, or $I^-(\gamma(a))\cap A \neq\varnothing$. The proof in the latter case is analogous and simpler, so we leave it to the reader. We can therefore assume there exists an index~$k$ such that $x_k\chron \gamma(t)\chron x_{k+1}$ with $\gamma(t)\in A$, and this gives a non-empty open
	\[
	\gamma(t)\in A\cap \Up p_k\cap \Down p_{k+1}=: W.
	\]
	See \cref{figure:causal coverage in spacetime}(a) for intuition. Since $p_k\Leq \Up p_k\cap \Down p_{k+1}\Leq p_{k+1}$, by parallel orderedness (via \cref{lemma:wedge iff strong wedge}) we obtain causal relations
	\[
	p_k\cap \Down W\Leq W\Leq p_{k+1}\cap \Up W.
	\]
	We will use this fact to refine $p$ into a path that inhabits $W\subseteq A$. This refinement is constructed as follows. First, split the original path in two by defining ${\hat{p}= (p_n)_{n=0}^k}$ and $\check{p}= (p_n)_{n=k+1}^N$, and take the refinements $\hat{p}|_{p_k\cap \Down W}$ and $\check{p}|^{p_{k+1}\cap \Up W}$. The desired path is defined via the concatenation of these refinements, interjected with the single-step path $W$ in the middle:
	\[
	q^{x_N}:= \check{p}|^{p_{k+1}\cap \Up W}\cdot W \cdot \hat{p}|_{p_k\cap \Down W}.
	\]
	This is well-defined by the previous equation. Using \cref{lemma:properties of refinement relation}(c) and \cref{lemma:path restriction preserves refinement} we see that $q^{x_N}\refines (\check{p}\cdot \hat{p})|_{q^{x_N}_\top}=p|_{q^{x_N}_\top}$,
	and by construction it inhabits $A$. Finally, since in particular we have $W\ni\gamma(t)\caus x_{k+1}\in p_{k+1}$ it follows by (the dual of) \cref{lemma:path restriction preserves points} that ${x_N\in q^{x_N}}$. Repeating this construction for arbitrary $x\in p_N$ gives a family $(q^x)_{x\in p_N}$ of refinements of $p$ that inhabit $A$, such that ${\bigvee_{x\in p_N}q^x_\top = p_N}$. In other words, $(q^x)_{x\in p_N}\refines p$ is a local past refinement, and this proves $A\in\CovLeq^-(U)$. 
\end{proof}

This proof does not rely on any particular property of spacetimes other than the fact that the chronology relation $\chron$ is related to the existence of timelike curves. We therefore expect a more general result to hold in the setting of ordered spaces induced by e.g.~\emph{d-spaces} \cite{grandis2009DirectedAlgebraicTopology}. In particular, $\Covcaus^\pm$ and $\Covchron^\pm$ will have similar properties as to $\CovLeq^\pm$.
	
\begin{remarknumbered}\label{remark:difference causal coverage localic or curve}
	The converse of \cref{proposition: localic coverage contains chronological coverage} is unfortunately not generally true. \cref{figure:causal coverage in spacetime}(b) provides a simple counterexample. Here we have a past inextendible timelike curve $\gamma$ in Minkowski space terminating in~$U$. The covering region is constructed by starting with a typical covering region $A\in\CovLeq^-(U)$, and cutting out the image $\im(\gamma)$. Since no localic path $p$ can pass through this infinitesimal hole, in the localic setting we still get that~$A\setminus \im(\gamma)\in\CovLeq^-(U)$.	
	From a philosophical point of view it can be argued that the localic definition is more appropriate, since in the definition of $\Covcaus^-$ or $\Covchron^-$ the only paths that can dodge $A\setminus\im(\gamma)$ are those where the information it carries is compacted into an infinitesimal space, which seems physically unreasonable.
	
	It is currently an open question under what conditions on $A$ or $U$ we get converse inclusions, if at all. We conjecture that it holds at least when $A$ is of the form $I^-(A)$, which excludes the counterexample in \cref{figure:causal coverage in spacetime}(b) since the chronological past ``closes'' any infinitesimal holes removed by curves.
\end{remarknumbered}

\begin{conjecture}
	In a smooth spacetime:
		\[
			\CovLeq^-(U)\cap \im(I^-)
			\subseteq \Covchron^-(U)\cap \im(I^-)
			\subseteq
			\Covcaus^-(U).
		\]
\end{conjecture}
%

The conjecture cannot be sharpened to e.g.~regular opens or convex opens. A counterexample can be constructed by ``sharpening'' the middle edges of the region $A$ in \cref{figure:causal coverage in spacetime}(b) to a point. More explicitly, construct $A$ as the disjoint union of two chronological diamonds that ``touch'' at $\im(\gamma)$. Both of these form regular convex opens, but still allow the curve $\gamma$ to pass.

\section{Abstract causal coverages}
In this section we take a more abstract approach, and view $\CovLeq^\pm$ as defining the structure of a \emph{coverage} on $\Opens X$. The following definition is an abstraction of the properties in \cref{lemma:properties of coverage from locale,lemma:coverage locale satisfies Cov-V strong}, and is similar to the philosophy of viewing it as a Grothendieck topology. The challenge here is to state axioms for abstract causal coverages $\Cov^\pm$ without referring to the pre-existing structure of $\Leq$ or the localic cones $\Up$ and $\Down$. We emphasise that the following definition is tentative.

\begin{definition}\label{definition:causal coverage abstract}
	Let $X$ be a locale. A \emph{causal coverage} on $X$ consists of two functions
		\[
			\Cov^\pm\colon \Opens X\longrightarrow\Powerset(\Opens X)
		\]
	satisfying the following axioms:
			\begin{enumerate}[label = (C\arabic*)]
			\item\label{axiom:Cov-unit} $U\in\Cov^\pm(U)$;
			\item\label{axiom:Cov-joins} axiom~\eqref{axiom:Cov-V strong} holds;
			\item\label{axiom:Cov-transitivity} if $B\in \Cov^\pm(A)$ and $A\in\Cov^\pm(U)$ then $B\in\Cov^\pm(U)$;
			\item\label{axiom:Cov-ideal} if $A,B\in \Cov^\pm(U)$ and $A\sqleq C\sqleq B$, then $C\in\Cov^\pm(U)$;
			\item\label{axiom:Cov-flip} if $A\in\Cov^\pm(U)$ then there exists $W\in\Cov^\mp(A)$ such that $U\sqleq W$.
		\end{enumerate}
	A \emph{causal site}\footnote{Now not to be confused with the homonymous but distinct notion in \cite{christensen2005CausalSitesQuantum}.} $(X,\Cov^\pm)$ is a locale $X$ equipped with a causal coverage $\Cov^\pm$.
\end{definition}

These axioms should be directly compared to the properties in \cref{lemma:properties of coverage from locale,lemma:coverage locale satisfies Cov-V weak,lemma:coverage locale satisfies Cov-V strong}. Indeed, we have the following.

\begin{proposition}
	If $(X,\Leq)$ is a parallel ordered locale satisfying~\eqref{axiom:LV}, then $\CovLeq^\pm$ from \cref{definition:causal coverage from ordered locale} define a causal coverage on $X$.
\end{proposition}
\begin{proof}
	Axioms~\ref{axiom:Cov-unit}--\ref{axiom:Cov-transitivity} follow directly from \cref{lemma:properties of coverage from locale,lemma:coverage locale satisfies Cov-V strong}. Axiom~\ref{axiom:Cov-ideal} holds since if $A,B\in\CovLeq^-(U)$ and $A\sqleq C\sqleq B$, then $C\sqleq B\sqleq \Down U$, and any path landing in $U$ can be refined to inhabit $A\sqleq C$. Thus $C\in \CovLeq^-(U)$. Lastly, for axiom~\ref{axiom:Cov-flip}, if $A\in\CovLeq^-(U)$ then by \cref{lemma:properties of coverage from locale}(e) we get $A\Leq U$, so $U\sqleq \Up A\in\CovLeq^+(A)$.
\end{proof}

In the following we show how we can recover the structure of an ordered locale from a causal site.
\begin{lemma}\label{lemma:causal coverage determines monads}
	If $\Cov^\pm\colon \Opens X\to\Powerset(\Opens X)$ satisfies axioms~\ref{axiom:Cov-unit}--\ref{axiom:Cov-transitivity} of a causal coverage, then the following functions define join-preserving monads on $\Opens X$:
		\[
			L^\pm\colon \Opens X\longrightarrow \Opens X;\qquad U\longmapsto \bigvee \Cov^\pm(U).
		\]
\end{lemma}
\begin{proof}
	Axiom~\ref{axiom:Cov-unit} immediately gives the unit $U\sqleq L^\pm (U)$. Axiom~\eqref{axiom:Cov-V strong} gives that $L^\pm (U)\in \Cov^\pm (U)$. In turn, it follows that $L^\pm \circ L^\pm (U)\in \Cov^\pm(L^\pm (U))$, so by axiom~\ref{axiom:Cov-transitivity} we get that $L^\pm \circ L^\pm (U)\in \Cov^\pm (U)$. This provides the multiplication $L^\pm \circ L^\pm (U)\sqleq L^\pm (U)$ of the monad. Lastly, if $U\sqleq V$ and $A\in \Cov^\pm(U)$, using~\ref{axiom:Cov-unit} and~\eqref{axiom:Cov-V strong} we get $A\sqleq A\vee V\in \Cov^\pm(U\vee V) = \Cov^\pm(V)$, which proves that $L^\pm$ are monotone.
	
	Finally, to prove $L^\pm$ preserve all joins, take opens $(U_i)_{i\in I}$. Since $L^\pm$ are monotone, we only need to show $L^\pm\left(\bigvee_{i\in I} U_i\right)\sqleq \bigvee_{i\in I} L^\pm (U_i)$. For that, take $A\in \Cov^\pm (\bigvee_{i\in I} U_i)$. By~\eqref{axiom:Cov-V strong} we can find $A_i\in \Cov^\pm (U_i)$ for each $i\in I$ such that $A=\bigvee_{i\in I} A_i\sqleq \bigvee_{i\in I} L^\pm (U_i)$, so we are done.
\end{proof}

The open region $L^\pm (V)$ is also called the \emph{future/past region of influence} of $V$. Since they are join-preserving monads, we think of $L^\pm$ as a different incarnation of the localic cones of an ordered locale.

We now show that they also relate to $\Cov^\pm$ analogously as how $\Up$ and $\Down$ relate to $\CovLeq^\pm$. First, $U\in \Cov^\pm(U)$ is baked in as~\ref{axiom:Cov-unit}, and $L^\pm(U)\in\Cov^\pm(U)$ follows by~\eqref{axiom:Cov-V strong}.

The following is the abstract analogue of \cref{lemma:properties of coverage from locale}(d).

\begin{lemma}\label{lemma:Cov-pullbacks}
	Let $\Cov^\pm$ be a causal coverage on $X$. If $A\in \Cov^\pm(U)$ and $W\sqleq U$, then $A\wedge L^\pm (W)\in \Cov^\pm(W)$.
\end{lemma}
\begin{proof}
	Note that $U= W\vee U$, so $A\in \Cov^\pm(W\vee U)$, and hence~\eqref{axiom:Cov-V strong} gives regions $B\in\Cov^\pm(W)$ and $C\in \Cov^\pm(W)$ such that $A= B\vee C$. Now observe that $B\sqleq A\wedge L^\pm (W) \sqleq L^\pm (W)\in \Cov^\pm(W)$, so the result follows by axiom~\ref{axiom:Cov-ideal}.
\end{proof}

The following shows that if the causal coverage comes from an ordered locale, this is literally the case.

\begin{proposition}\label{proposition:ordered locale from canonical coverage}
	If $(X,\Leq)$ is a parallel ordered locale satisfying~\eqref{axiom:LV}, then $\Leq$ is recovered by the monads induced by the causal coverage $\CovLeq^\pm$. 
\end{proposition}
\begin{proof}
	By \cref{lemma:properties of coverage from locale}(b) it follows that $L^-(U)= \bigvee\CovLeq^-(U) = \Down U$, and similarly $L^+ (U) = \Up U$. Hence the preorder induced by $L^\pm$ on $\Opens X$ is equal to $\Leq$ by~\eqref{axiom:cones give order}.
\end{proof}

We currently do not know if the converse to \cref{proposition:ordered locale from canonical coverage} also holds: if we start with an abstract causal coverage $\Cov^\pm$ on $X$ and induce the ordered locale $(X,\Leq)$ from the localic cones $L^\pm$, does $\CovLeq^\pm$ have to equal $\Cov^\pm$? We conjecture that there is an adjunction between ordered locales and ordered sites, analogous to the adjunction between frames and sites in \cite{ball2014ExtendingSemilatticesFrames}. This might allow us to build a connection to the theory of \emph{formal topology} \cite{sambin2003PointsFormalTopology}.

%
%

\section{Domains of dependence}
\label{section:domains of dependence}
In this section we study an abstract version of \emph{domains of dependence} \cite{geroch1970DomainofDependence} in the setting of causal coverages. The idea is that, given a region $A\in\Opens X$, the future domain of dependence $D^+(A)$ is the largest region that is covered below by $A$. Thinking of information flows, it is the largest region where all information is determined by $A$.
 
\begin{figure}[b]
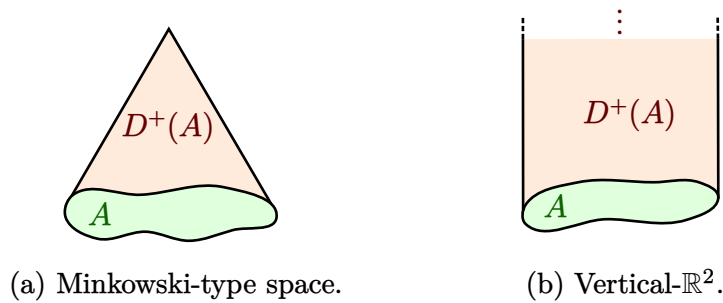
\centering
	\begin{subfigure}[b]{0.4\textwidth}\centering
		\input{tikz/inkscape-export/fig-20}
		\caption{Minkowski-type space.}
	\end{subfigure}
	\begin{subfigure}[b]{0.4\textwidth}\centering
		\input{tikz/inkscape-export/fig-20-vertical}
		\caption{Vertical-$\mathbb{R}^2$.}
	\end{subfigure}
	\caption{Illustration of typical future domains of dependence.}
	\label{figure:typical D+}
\end{figure}

\begin{definition}\label{definition:domain of dependence localic}
	Let $\Cov^\pm$ be a causal coverage on a locale $X$. The \emph{future/past domain of dependence} of $A\in\Opens X$ is defined as
		\[
			D^\pm (A) :=\bigvee\left\{V\in\Opens X: A\in\Cov^\mp(V)\right\}.
		\]
	By~\eqref{axiom:Cov-V strong}, the future/past domain of dependence $D^\pm(A)$ is uniquely characterised as the largest open region that is covered from below/above by $A$. See \cref{figure:typical D+} for intuition. In particular we have
		\[
			A\in\Cov^\mp \left(D^\pm(A) \right).
		\]
\end{definition}
\begin{lemma}
	If $\Cov^\pm$ is a causal coverage on $X$, the domains of dependence define monads $D^\pm\colon\Opens X\to \Opens X$.
\end{lemma}
\begin{proof}
	First, axiom~\ref{axiom:Cov-unit} of a causal coverage gives $A\in \Cov^\mp(A)$, from which it follows immediately that $A\sqleq D^\pm (A)$. Next, if $D^\pm (A)\in\Cov^\mp(V)$, using ${A\in\Cov^\mp(D^\pm A)}$, we get by~\ref{axiom:Cov-transitivity} that $A\in \Cov^\mp(V)$, providing the multiplication ${D^\pm \circ D^\pm (A)\sqleq D^\pm (A)}$. 
	
	We are left to show monotonicity, for which we need the following equation:
	\[
	D^\pm (A) = \bigvee\left\{V\in\Opens X: A\wedge L^\mp (V)\in\Cov^\mp(V)\right\}.
	\]
	To prove this, denote the expression on the right hand side by $\tilde{D}^\pm (A)$. It is clear that we get an inclusion $D^\pm (A)\sqleq \tilde{D}^\pm (A)$. For the converse, it suffices to show that $A\in\Cov^\mp (\tilde{D}^\pm (A))$. To obtain this, note that~\eqref{axiom:Cov-V strong} gives
	\[
	A\wedge \bigvee\left\{V\in\Opens X: A\wedge L^\mp (V)\in\Cov^\mp(V)\right\} \in \Cov^\mp(\tilde{D}^\pm (A)).
	\]
	Further note that this region fits into a chain of inclusions $A\sqleq D^\pm (A)\sqleq \tilde{D}^\pm (A)\sqleq L^\mp \circ \tilde{D}^\pm (A)\in \Cov^\mp(\tilde{D}^\pm (A))$. Therefore, by axiom~\ref{axiom:Cov-ideal} of a causal coverage, we get $A\in\Cov^\mp (\tilde{D}^\pm (A))$, proving the equation.
	
	Now take opens $A\sqleq B$, and pick $V\in\Opens X$ such that $A\in\Cov^\mp (V)$. We get a chain of inclusions $A\sqleq B\wedge L^\mp (V)\sqleq L^\mp (V)$, so through axiom~\ref{axiom:Cov-ideal} it follows $B\wedge L^\mp (V)\in \Cov^\mp(V)$. This shows $V\sqleq \tilde{D}^\pm (B)$, and using the equation above it follows that $D^\pm (A)\sqleq D^\pm (B)$.
\end{proof}

Note that $D^\pm$ are generally not join-preserving. Typically $D^\pm(A\vee B)$ is much larger than $D^\pm(A)\vee D^\pm(B)$, cf.~the examples in \cref{figure:differences domains of dependence}.

The following shows that a causal coverage $\Cov^\pm$ is fully determined by its regions of influence and domain of dependence operators.

\begin{lemma}\label{lemma:causal coverage determined by L and D}
	Let $\Cov^\pm$ be a causal coverage on a locale $X$. Then:
		\[
			A\in \Cov^\pm(V)
			\qquad\text{if and only if}\qquad
			A\sqleq L^\pm (V)\text{ and }V\sqleq D^\mp (A).
		\]
\end{lemma}
\begin{proof}
	The ``only if'' direction follows by construction. Conversely, suppose that $A\sqleq L^\pm (V)$ and $V\sqleq D^\mp (A)$. Since $A\in \Cov^\pm(D^\mp (A))$ we get by \cref{lemma:Cov-pullbacks} some $B\in \Cov^\mp(V)$ such that $B\sqleq A$. From~\eqref{axiom:Cov-V strong} we know $L^\pm (V)\in \Cov^\pm (V)$, so now~\ref{axiom:Cov-ideal} gives $A\in\Cov^\pm(V)$, as desired. 
\end{proof}

\begin{lemma}
	We have $L^\pm\circ D^\pm = L^\pm$, and hence $D^\pm \sqleq L^\pm$. 
\end{lemma}
\begin{proof}
	Since $D^\pm$ are monads, the inclusions $L^\pm (A)\sqleq L^\pm\circ D^\pm(A)$ are clear. For the converse, recall that by \cref{lemma:causal coverage determines monads} the maps $L^\pm$ preserve all joins. Hence $L^\pm \circ D^\pm (A) = \bigvee\{L^\pm (V) : A\in \Cov^\mp (V)\}$. Now, if $A\in \Cov^\mp(V)$ we get by axiom~\ref{axiom:Cov-flip} an open $W\in \Cov^\pm(A)$ such that $V\sqleq W$. Thus $L^\pm (V)\sqleq L^\pm (W)\sqleq L^\pm \circ L^\pm (A) = L^\pm (A)$, and we may conclude $L^\pm \circ D^\pm (A) = L^\pm (A)$. Using the units of $L^\pm$ it immediately follows that $D^\pm (A)\sqleq L^\pm \circ D^\pm (A) = L^\pm (A)$. 
\end{proof}

\subsection{Domains of dependence in spacetime}\label{section:domains of dependence in spactime}
In spacetime there is already a well-established notion of domain of dependence, making use of inextendible curves. In this section we compare these to the localic versions, and see that they behave quite differently. Here we adopt \cite[Definition~3.1]{minguzzi2019LorentzianCausalityTheory}.

\begin{definition}\label{definition:domain of dependence in spacetime}
	For a subset $A\subseteq M$ of a smooth spacetime, the \emph{future causal domain of dependence} of $A$ is the set
		\[
			\Dcaus^+(A) := \left\{x\in M:~\parbox{.4\textwidth}{every past inextendible \emph{causal} curve through $x$ intersects $A$}~\right\}.
		\]
	Analogously, the \emph{future chronological domain of dependence} of $A$ is the set\footnote{Beware that \cite{minguzzi2019LorentzianCausalityTheory} instead uses the notation $D^+$ and $\tilde{D}^+$ for the causal and chronological domains of dependence, respectively.}
		\[
			\Dchron^+(A) := \left\{x\in M:~\parbox{.4\textwidth}{every past inextendible \emph{timelike} curve through $x$ intersects $A$}~\right\}.
		\]
\end{definition}

\begin{remark}
	Note that domains of dependence are usually mainly considered on \emph{achronal} subsets: those $S\subseteq M$ for which $I^-(S)\cap I^+(S)=\varnothing$ \cite[p.~105]{landsman2021FoundationsGeneralRelativity}. These are spacelike hypersurfaces on which initial data of partial differential equations are typically defined. In an ordered locale the condition $\Down U\wedge \Up U = \varnothing$ implies $U=\varnothing$, so the notion becomes trivial. This is one reason why we think a proper notion of causal order $\Leq$ on the full lattice $\Sl(X)$ of sublocales will be fruitful: it will allow us to consider non-trivial achronal ``sets'' in an ordered locale. 
\end{remark}

\begin{definition}
	For a smooth spacetime $M$ we define the following notions of domains of dependence for $A\subseteq M$:
		\begin{align*}
			\DCaus^+(A) &:= \bigcup\left\{U\in\Opens M: A\in\Covcaus^-(U) \right\},\\
			\DChron^+(A) &:= \bigcup\left\{U\in\Opens M: A\in\Covchron^-(U) \right\},\\
			\DLeq^+(A) &:= \bigcup\left\{U\in\Opens M: A\in\CovLeq^-(U) \right\}.
		\end{align*}
\end{definition}

Note that $\DLeq^+$ is just the localic domain of dependence from \cref{definition:domain of dependence localic} determined by the canonical causal coverage $\CovLeq^-$, and is hence a monad. We will now show how these five definitions are related but distinct.

\begin{figure}[h]
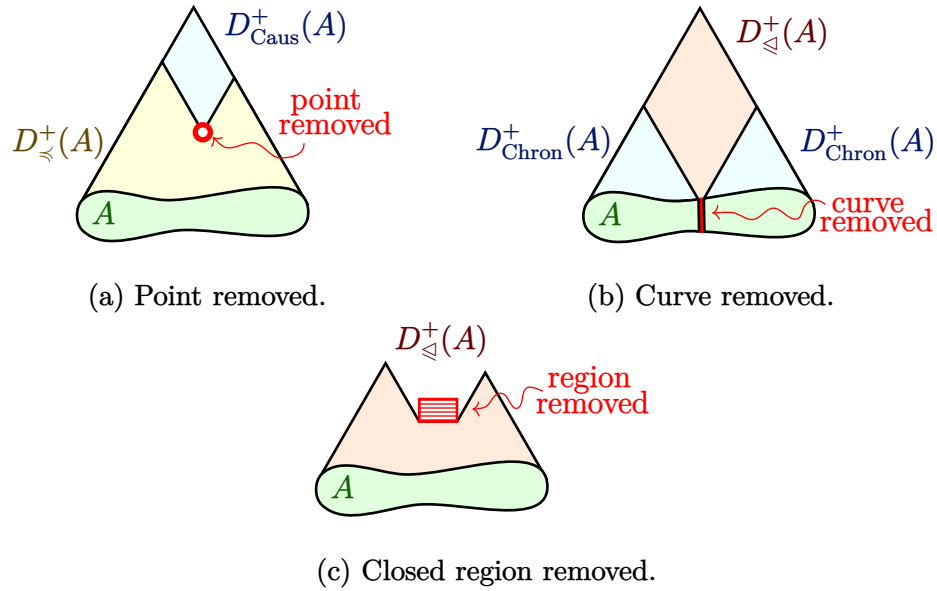
\centering
	\begin{subfigure}[b]{0.4\textwidth}\centering
		\input{tikz/inkscape-export/fig-17-point}
		\caption{Point removed.}
	\end{subfigure}
	\begin{subfigure}[b]{0.5\textwidth}\centering
		\input{tikz/inkscape-export/fig-17-curve}
		\caption{Curve removed.}
	\end{subfigure}
	\begin{subfigure}[b]{0.7\textwidth}\centering
		\input{tikz/inkscape-export/fig-17-region}
		\caption{Closed region removed.}
	\end{subfigure}
	\caption{Illustration of differences between localic and curve-wise domains of dependence in Minkowski spaces with either a point or a larger closed region removed. In (b) a curve is removed from $A$ (not the spacetime).}
	\label{figure:differences domains of dependence}
\end{figure}

\begin{enumerate}[label = \textbullet]
	\item \emph{Bounded curves vs.~localic paths.}
	By \cref{proposition: localic coverage contains chronological coverage} we get for every ${U\in\Opens M}$ that $\Covcaus^-(U)\subseteq \Covchron^-(U)\subseteq \CovLeq^-(U)$,
	from which it immediately follows that
		\[
			\DCaus^+\subseteq \DChron^+\subseteq\DLeq^+.
		\]
	Due to the counterexample in \cref{remark:difference causal coverage localic or curve} that $\Covchron^-(U)\not\subseteq \CovLeq^-(U)$, see \cref{figure:causal coverage in spacetime}(b), the second inclusion will be strict. An example of what this can look like is in \cref{figure:differences domains of dependence}(b).
	
	\item \emph{Inextendible vs.~bounded curves.} Similarly, using \cref{lemma:past inextendible intersects then bounded intersects} and its chronological analogue we obtain inclusions
		\[
			\Dcaus^+\subseteq \DCaus^+
			\qquad\text{and}\qquad
			\Dchron^+\subseteq \DChron^+.
		\]
	But, again, these inclusions are strict by virtue of the counterexample in \cref{remark:difference causal coverage bounded or inextendible} and \cref{figure:non-intersecting inextendible curve}. This is showcased in \cref{figure:differences domains of dependence}(a).
	
	\item \emph{Chronological vs.~causal.} First, since all timelike curves are causal curves, we get an inclusion $\Dcaus^+\subseteq \Dchron^+$ of the traditional domains of dependence. It is well-known that this inclusion is strict in general, cf.~\cite[Remark~3.11]{minguzzi2019LorentzianCausalityTheory}. Though we conjecture that they are equal when restricted to open regions.
	
	The counterexample in \cref{remark:causal cover contained in chronological cover} that ${\Covchron^-(U)\not\subseteq \Covcaus^-(U)}$ (see \cref{figure:non-intersecting lightlike curve}) also implies that the inclusion $\DCaus^+\subseteq \DChron^+$ is strict. To see how, one can draw a situation analogous to \cref{figure:differences domains of dependence}(b) where the curve removed is lightlike. In that case $\DCaus^+(A)$ equals the two smaller triangles as drawn, but $\DChron^+(A)$ will equal $\DLeq^+(A)$.
\end{enumerate}

In conclusion, the domains $\Dcaus^+$ and $\Dchron^+$ defined using inextendible curves differ from the domains $\DCaus^+$ and $\DChron^+$ defined in terms of bounded curves when we remove a point from the spacetime. The domains $\DCaus^+$ and $\DChron^+$ differ from the localic domain $\DLeq^+$ as soon as we remove a curve from $A$. In summary, we have:
\[
\begin{tikzcd}[cramped,sep=small]
	{\Dcaus^+} & {\Dchron^+} \\
	{\DCaus^+} & {\DChron^+} & {\DLeq^+.}
	\arrow["\subsetneq"{marking, allow upside down}, draw=none, from=1-1, to=1-2]
	\arrow["\subsetneq"{marking, allow upside down}, draw=none, from=1-1, to=2-1]
	\arrow["\subsetneq"{marking, allow upside down}, draw=none, from=2-1, to=2-2]
	\arrow["\subsetneq"{marking, allow upside down}, draw=none, from=1-2, to=2-2]
	\arrow["\subsetneq"{marking, allow upside down}, draw=none, from=2-2, to=2-3]
\end{tikzcd}
\]
All the presented examples where the non-localic domains of dependence differ from $\DLeq^+$ involve the removal of some infinitesimal point or curve. In \cref{figure:differences domains of dependence}(c) we see an example where the localic domain $\DLeq^+$ indeed behaves like one would expect from the traditional definition of domain of dependence, given that one replaces points by larger closed regions.

In \cref{section:domains of dependence future work} we discuss some potential future work relating to \emph{sheaves} for causal coverages, and \emph{holes} in spacetime.

\vspace*{\fill}

\chapter{Future work}\label{chapter:applications}
In this concluding chapter we suggest some outlines for future research.


\section{Order on sublocales}\label{section:order on sublocales}
Our main \cref{definition:ordered locale} of ordered locales is an axiomatisation of the Egli-Milner order on the lattice $\Opens S$ of open subsets of an ordered space. But, like topological spaces, we have seen that locales admit a more general notion of ``subspace'' that are not necessarily open: \emph{sublocales} (\cref{section:sublocales}). As $\Opens X$ is a generalisation of $\Opens S$, so we can think of the lattice of sublocales $\Sl(X)$ as a generalisation of the powerset $\Powerset(S)$. An alternative to \cref{definition:ordered locale} would thus be to axiomatise the Egli-Milner order on the richer structure $\Sl(X)$, instead of just restricting it to the open regions~$\Opens X$. This also gives us access, for instance, to possible relations between the causal order and closed and compact sublocales. It connects also to the story of internally ordered locales in \cref{section:internal preorders}, via \cref{construction:ordered locale from internal preorder}.

One immediate benefit of such a definition is that it makes the study of a \emph{localic causal ladder} (cf.~\cref{section:causal ladder}) easier: for instance the \emph{global hyperbolicity} condition can then be rephrased as saying that $\up K\wedge \down L$ is a compact sublocale, for any compact sublocales $K,L\in\Sl(X)$. This condition is equivalent to global hyperbolicity in spacetimes by~\cite[Lemma~5.29]{landsman2021FoundationsGeneralRelativity}, and note importantly that it works even if the locale $X$ has no points. Similarly, we can define \emph{causally simple} locales as those where $\up K$ and $\down K$ are closed sublocales for every compact $K\in\Sl(X)$, cf.~\cite[Theorem~4.12]{minguzzi2019LorentzianCausalityTheory}.

\begin{example}\label{example:simple signalling}
	Continuing speculatively, here is a superficial example of ``localic superluminal signalling'' in an ordered space. As a set, take $S:=\{x,y,z\}$, and define a topology by $\Opens S:= \{\varnothing, U, S\}$, where the sole non-trivial open is $U:=\{x,z\}$. Generate an order on $S$ by setting $x\leq y$, but $x\not\leq z$. We interpret this as $x$ not being able to signal to $z$. Equivalently, this means that $\up x \cap\{z\}=\varnothing$. However, it can be shown that in the lattice of sublocales $\varnothing\neq \tilde{x}= \tilde{x}\wedge \tilde{z} \subseteq \widetilde{\up x}\wedge \tilde{z}$, where $\tilde{A}$ denotes the sublocale induced by $A\subseteq S$ (\cref{section:sublocales}), which we could interpret as a type of topological signalling from $x$ to $z$ that did not exist in the original space. 
\end{example}

\section{AQFTs on ordered locales}\label{section:aqft on ordered locales}
\emph{Algebraic quantum field theory} (AQFT) is a mathematical formalisation of quantum fields that is based on the notion of a \emph{local observable}, first studied by Haag and Kastler \cite{haag1964AlgebraicApproachQuantum}. The following exposition is based on the introductory text \cite{fewster2020AlgebraicQuantumField}; see \cite{halvorson2007AlgebraicQuantumField} for a technical development.

The basic idea is that it consists of an assignment $\calA$ of regions $U\subseteq M$ in spacetime to a certain type of algebra $\calA(U)$ of physical observables at $U$. Unfortunately, there is no single agreed upon definition of an AQFT. Instead there are slight differences in axioms and presentations, depending also on whether one is working solely on Minkowski space or more general spacetimes. To make the point, we paraphrase the definition of an AQFT over Minkowski space $M$ from \cite[\S 4.1]{fewster2020AlgebraicQuantumField}.
\begin{definition}\label{definition:aqft minkowski}
	An AQFT over $M$ is consists of a \emph{local net of observables}:
		\begin{enumerate}[label = (A\arabic*)]
			\item There is a functor
				\[
					\calA\colon \hull{\Opens M}\longrightarrow \cat{Alg}
				\]
			from the lattice of convex open subsets of $M$ to some category of algebras. For a convex open subset $U\subseteq M$ we call $\calA(U)$ the \emph{local algebra (of observables)} at $U$.
			\item If $U\subseteq V$, then the arrow $\calA(U\subseteq V)$ is a subalgebra inclusion.
			\item If $U$ and $V$ are spacelike separated then the commutator $[\calA(U),\calA(V)]$ of their local algebras is trivial in the encompassing algebra $\calA(M)$.
			\setcounter{enumi}{4}
			\item For any $U$ we have $\calA(U) = \calA\left(D^+(U)\cup D^-(U)\right)$.
		\end{enumerate}
\end{definition}

Here, for simplicity, we have skipped axiom~(A4) relating to Poincaré invariance. The main point is that this definition does not depend on the points $x\in M$ of the spacetime, but can be phrased purely in terms of algebras defined on open subsets $U\subseteq M$. This gives a localic flavour. Moreover, the causal notions used in this definition: convexity, spacelike separation, and domains of dependence, all exist in an \emph{ordered} locale. In fact, in some sense an ordered locale $(X,\Leq)$ provides the minimal amount of ingredients to state the definition of an AQFT:
	\begin{enumerate}[label = \textbullet]
		\item a base lattice $\Opens X$ of open regions $U$;
		\item a notion of convexity (\cref{section:convexity}): $U = \Up U \wedge \Down U$;
		\item a notion of spacelike separation (\cref{section:causal complements}): $U\sqleq V^\bot$ iff $V\sqleq U^\bot$;
		\item a notion of domain of dependence (\cref{section:domains of dependence}): $D^\pm(U)$.
	\end{enumerate}
The \cref{definition:aqft minkowski} of AQFTs on Minkowski spacetime can therefore be generalised directly to the setting of an abstract ordered locale $(X,\Leq)$.

\begin{definition}
	A \emph{local net of observables} on a parallel ordered locale $(X,\Leq)$ with~\eqref{axiom:LV} consists of the following data and axioms:
		\begin{enumerate}[label = (A\arabic*)]
			\item There is a functor
			\[
			\calA\colon \hull{\Opens X}\longrightarrow \cat{Alg}
			\]
			from the lattice of convex regions of $X$ to some category of algebras.
			\item If $U\sqleq  V$, then the arrow $\calA(U\sqleq V)$ is a subalgebra inclusion.
			\item If $U\sqleq V^\bot$, or equivalently $V\sqleq U^\bot$, then the commutator $[\calA(U),\calA(V)]$ of their local algebras is trivial in the encompassing algebra $\calA(X)$.
			\setcounter{enumi}{4}
			\item For any $U$ we have $\calA(U) = \calA\left(D^+(U)\vee D^-(U)\right)$.
		\end{enumerate}
\end{definition}

\begin{example}
	Consider the ordered locale induced by vertical-$\mathbb{R}^2$, equipped with a local net of observables $\calA$. From \cref{figure:typical D+} we see that the domains of dependence can be calculated as $D^+(U) = \Up U$ and $D^-(U) = \Down U$. Axiom~(A5) then gives the restriction: $\calA(U) = \calA(\Up U\vee \Down U)$. Thus, intuitively, $\calA$ is simply a family of `static' algebras indexed by the spatial axis $\mathbb{R}$ that do not interact or evolve throughout time. To make this more precise, for ${x=(x_1,x_2)\in\mathbb{R}^2}$ consider the \emph{costalk} $\calA_x$ of the local net $\calA$ at $x$:
			\[
				\calA_x := \lim_{U\ni x} \calA(U),
			\]
	provided $\cat{Alg}$ admits the required limits (\cref{definition:limits}). By axiom~(A5), we see that if $x,y\in \mathbb{R}^2$ are on the same vertical axis, $x_1=y_1$, then there is a canonical isomorphism $\calA_x\cong \calA_y$.
\end{example}

More abstractly, \cite{grant-stuart2023SpacetimesCategoriesDisjointness} introduces the notion of a \emph{disjointless relation} on a category, in order to provide an abstract framework in which to study AQFTs. The disjointness relation generalises spacelike separation. We conjecture that the causal complement of an ordered locale leads to a natural choice for a disjointness relation on the category $\OrdLoc$ of ordered locales. 

\section{Atomistic ordered locales}\label{section:localic causal sets}
\emph{Causal set theory} was proposed in \cite{bombelli1987SpacetimeCausalSet} (for a more recent review article see e.g.~\cite{surya2019CausalSetApproacha}), and models spacetime as a ``discrete'' set of events related by a causal order. The discreteness is supposed to model the conjectured \emph{indivisibility} of physical spacetime at the Planck scale $\ell_\mathrm{Planck} := \sqrt{G \hbar c^{-3}}$, which is of the order of~$10^{-35}$ metres. The indivisibility is modelled by a finiteness axiom.

\begin{definition}\label{definition:causal set}
	A \emph{causal set} is a partial order $(C,\caus)$ that is \emph{locally finite}: for every $x,y\in C$ the \emph{interval} $\up x\cap \down y$ is a finite subset of $C$.
\end{definition}
Thinking of causal sets as discrete spaces, where finiteness is equivalent to compactness, we can view $(C,\caus)$ as a ``globally hyperbolic'' discrete spacetime.

The motivation behind causal sets is very close to the motivation behind ordered locales, taking the causal order $\caus$ of a spacetime as fundamental, the only difference being the substance that carries this causal order: spacetime points or topology, respectively. Can we devise a localic analogue of causal sets? We provide some speculation in this section. More refined connections can perhaps be made using topological techniques in \cite{isham1990IntroductionGeneralTopology,sorkin1991FinitarySubstituteContinuous}.

Instead of starting with a set of spacetime \emph{events} $C$ with a causal order, one suggestion is to directly model the discrete nature of physical spacetime by postulating the existence of indivisible \emph{``Planck regions''} in an ordered locale. One way to do this is to make use of lattice-theoretic atoms. We recall \cref{definition:atom}; see also \cite[\S 5.2]{davey2002IntroductionLatticesOrder}.

\begin{definition}
	An \emph{atom} in a bounded lattice $L$ is a non-empty element $a\in L$ such that if $x\sqleq a$ then either $x=a$ or $x=\bot$. Say $L$ is \emph{atomic} if every non-empty element contains an atom. It is called \emph{atomistic} if every element can be written as a join of atoms.
\end{definition}

\begin{remark}
	That is: atoms are the minimal non-empty elements. Many topologies do not admit them. For instance, if $S$ is a $T_0$-space, then the atoms in $\Opens S$ correspond precisely to the open sets of the form~$\{x\}$, for $x\in S$. Singletons, when open, are clearly atoms. Conversely, if $A\in\Opens S$ is an atom, suppose there exist $x,y\in A$ such that $x\neq y$. By $T_0$-separation we get, without loss of generality, some open $U\in \Opens S$ such that $x\in U\not\ni y$, so $\varnothing\subsetneq U\cap A\subsetneq A$, contradicting that $A$ is an atom. Thus $A$ is a singleton. In particular, (spacetime) manifolds do not admit atoms (cf.~\cref{proposition:points in regular opens Hausdorff space}).
\end{remark}

For Boolean algebras, being atomic and atomistic is the same thing, but not for general lattices. A starting point for a localic analogue of causal sets could therefore be those ordered locales $(X,\Leq)$ whose underlying frame $\Opens X$ is atomistic. This turns out to be a drastic assumption.

\begin{lemma}
	Any atomistic frame is a Boolean algebra.
\end{lemma}
\begin{proof}
	This proof is from \cite{apostolos2018representationTheormeCompleteAtomicHeytingAlg}. Take $x\in L$ in an atomistic frame. Let $\calA\subseteq L$ denote the atomic elements. So we get $x=\bigvee B$ for some $B\subseteq\calA$. If $y\in \calA\setminus B$ then $x\wedge y =\bigvee\{b\wedge y: b\in B\}=\varnothing$. Thus $\bigvee(\calA\setminus B)\sqleq \neg x$, and so $\top = \bigvee B \vee \bigvee (\calA\setminus B) \sqleq x\vee \neg x$, showing the \emph{law of the excluded middle} holds. Thus $L$ is Boolean by \cite[\S AI.7.4]{picado2012FramesLocalesTopology}.
\end{proof}

Hence atomistic frames are just complete atomic Boolean algebras, which are the subject of the well-known \emph{Tarski duality} \cite{tarski1935ZurGrundlegungBoole}, stating the following.

\begin{theorem}[{\cite[\S II.5.4.2]{picado2012FramesLocalesTopology}}]
	A complete Boolean algebra is atomic if and only if it is isomorphic to a powerset.
\end{theorem}

Therefore assuming $\Opens X$ is atomistic forces discreteness: $\Opens X\cong \Powerset(S)$ for some set $S$. Of course, $S$ can just be taken as the set of atoms of $\Opens X$. The order $\Leq$ on $\Opens X$ reduces to a preorder $\leq$ on $S$, and with a suitable global hyperbolicity condition, we obtain a causal set $(S,\leq)$. But we also obtain nothing more: every ordered locale whose frame is atomistic is of this form. All topological structure is lost. 

For future work, laxer atomicity conditions on $\Opens X$ could be considered. For instance, there are alternative non-equivalent definitions for atom-like elements in a general lattice, such as \emph{(completely) join-irreducible} elements, \emph{compact elements}, and \emph{finite} elements \cite[Definition~7.15]{davey2002IntroductionLatticesOrder}. For instance, we could consider \emph{algebraic lattices}, which are complete lattices $L$ whose elements are generated by their compact elements. We leave investigating the connection between these ideas and causal set theory to future work.

%
%

\section{Categorification of ordered locales}
In this speculative section we discuss some ideas related to categorical generalisations of ordered locales. In the introduction \cref{section:tensor topology} we mentioned that ordered locales were in fact developed with possible categorifications to monoidal categories in mind. In this section we provide speculative remarks on an intermediary step towards categorification.

Locales categorify to \emph{toposes} \cite{maclane1994SheavesGeometryLogic}. With an eye towards the (far) future: what could be an appropriate categorification of ordered locales to this setting? The discussion in \cref{section:locale of pasts and futures} possibly sheds some light on this. Combining all structure, an ordered locale $(X,\Leq)$ with~\eqref{axiom:LV} induces subframes $\im(\Up)$ and $\im(\Down)$ of $\Opens X$, which in localic language can be interpreted as epimorphisms $\eta_\triup \colon X\tworightarrow X^\triup$ and $\eta_\tridown \colon X\tworightarrow X^\tridown$. But these morphisms are special in the sense that they induce adjoint triples:
		\[
		\begin{tikzcd}[ampersand replacement=\&]
			{\Opens X} \&\& {\Opens X^\triup}
			\arrow[""{name=0, anchor=center, inner sep=0}, "{\eta_\triup^{-1}}"{description, pos=0.3}, hook', from=1-3, to=1-1]
			\arrow[""{name=1, anchor=center, inner sep=0}, "\Up", shift left=4, from=1-1, to=1-3]
			\arrow[""{name=2, anchor=center, inner sep=0}, "{\eta^\triup_\ast}"', shift right=4, from=1-1, to=1-3]
			\arrow["\dashv"{anchor=center, rotate=-90}, draw=none, from=1, to=0]
			\arrow["\dashv"{anchor=center, rotate=-90}, draw=none, from=0, to=2]
		\end{tikzcd}
		\quad\text{and}\quad
		\begin{tikzcd}[ampersand replacement=\&]
			{\Opens X} \&\& {\Opens X^\tridown.}
			\arrow[""{name=0, anchor=center, inner sep=0}, "{\eta_\tridown^{-1}}"{description, pos=0.3}, hook', from=1-3, to=1-1]
			\arrow[""{name=1, anchor=center, inner sep=0}, "\Down", shift left=4, from=1-1, to=1-3]
			\arrow[""{name=2, anchor=center, inner sep=0}, "{\eta^\tridown_\ast}"', shift right=4, from=1-1, to=1-3]
			\arrow["\dashv"{anchor=center, rotate=-90}, draw=none, from=1, to=0]
			\arrow["\dashv"{anchor=center, rotate=-90}, draw=none, from=0, to=2]
		\end{tikzcd}
		\]
	
	The point is that these notions have topos-theoretic analogues that have been well-studied. One option is therefore to think of an ``ordered topos'' as a topos $\calE$, together with surjective essential geometric morphisms $\eta_\triup\colon \calE \to \calE^\triup$ and ${\eta_\tridown\colon \calE\to \calE^\tridown}$, essentialness here meaning that there are adjoint triples
	\[
	\begin{tikzcd}[ampersand replacement=\&]
		\calE \&\& {\calE^\triup}
		\arrow[""{name=0, anchor=center, inner sep=0}, "{\eta_\triup^*}"{description, pos=0.3}, from=1-3, to=1-1]
		\arrow[""{name=1, anchor=center, inner sep=0}, "{\eta^\triup_!}", shift left=4, from=1-1, to=1-3]
		\arrow[""{name=2, anchor=center, inner sep=0}, "{\eta^\triup\ast}"', shift right=4, from=1-1, to=1-3]
		\arrow["\dashv"{anchor=center, rotate=-90}, draw=none, from=1, to=0]
		\arrow["\dashv"{anchor=center, rotate=-90}, draw=none, from=0, to=2]
	\end{tikzcd}
	\quad\text{and}\quad
	\begin{tikzcd}[ampersand replacement=\&]
		\calE \&\& {\calE^\tridown,}
		\arrow[""{name=0, anchor=center, inner sep=0}, "{\eta_\tridown^*}"{description, pos=0.3}, from=1-3, to=1-1]
		\arrow[""{name=1, anchor=center, inner sep=0}, "{\eta^\tridown_!}", shift left=4, from=1-1, to=1-3]
		\arrow[""{name=2, anchor=center, inner sep=0}, "{\eta^\tridown_\ast}"', shift right=4, from=1-1, to=1-3]
		\arrow["\dashv"{anchor=center, rotate=-90}, draw=none, from=1, to=0]
		\arrow["\dashv"{anchor=center, rotate=-90}, draw=none, from=0, to=2]
	\end{tikzcd}
	\]
	and surjectiveness being the appropriate generalisation of epimorphisms of locales, see e.g.~\cite[p.~367]{maclane1994SheavesGeometryLogic}. These sorts of structures suggests a possible connection to Lawvere's work on \emph{cohesion} \cite{lawvere2007AxiomaticCohesion}. See also \cite{johnstone2011RemarksPunctualLocal}. Very informally, the idea behind a \emph{cohesive topos} $\calE$ is that we can think of its objects as a type of ``locally connected spaces,'' abstractly interpreted.
	
	We point out however that the traditional definition of cohesion will not apply here directly: this requires the left adjoints $\eta^\triup_!$ and $\eta^\tridown_!$ to preserve finite products, which is simply not the case for the localic cones $\Up$ and $\Down$. Similarly, it seems $\eta_\triup$ and $\eta_\tridown$ will not be \emph{locally connected}, since the maps $\eta^{-1}_\triup$ and $\eta^{-1}_\tridown$ appear to only laxly preserve Heyting implication, not strictly. We leave it to future work to devise a notion of \emph{``temporal cohesion.''}
\section{Causal Heyting implication}\label{section:caucal heyting implication}
In this section we take a slight logical turn. Our constructions here are inspired by similar ones in \cite{akbartabatabai2021ImplicationSpacetime}. We have already mentioned that any frame $L$ has an intrinsic Heyting implication, which is right adjoint to the finite meet operator. For an open region $U\in \Opens X$ in a locale:
\[
(U\wedge -)\dashv (U\to -).
\]
The existence of the Heyting implication is equivalent to the fact that $U\wedge -$ is a join-preserving map (\cref{theorem:join preserving implies left adjoint}), which is of course just the infinite distributivity law. For spatial intuition, recall from \cref{example:heyting implication in space} that in a topological space $S$ this takes the form $U\to V = (V\cup (S\setminus U))^\circ$. But here we shall interpret things logically: elements $U\in \Opens X$ are propositions, and the inclusion relation $\sqleq$ is interpreted as an \emph{inference} relation. Thus we read $U\sqleq V$ as ``from premise $U$ it can be concluded that $V$.'' In logic one often writes instead ${U\vdash V}$ (which is here not to be confused with the symbol used for adjunctions). The lattice operations of $\Opens X$ are interpreted as the logical connectives in the natural way: meets $\wedge$ corresponds to conjunction (``and''), and joins $\bigvee$ to infinitary disjunction~(``or''). The Heyting implication gives an internal notion of ``implication.'' Thus we interpret $U\to V$ as the proposition ``$U$ implies $V$.'' It is uniquely characterised via
\[
U\wedge W\sqleq V
\qquad\text{if and only if}\qquad
W\sqleq U\to V,
\]
which in particular implies that the following inclusion holds:
\[
U\wedge (U\to V)\sqleq V.
\]
This is of course just the internal \emph{modus ponens} law in $\Opens X$.

How is the order $\Leq$ related to this structure? One way is as follows. Suppose now that we have an ordered locale $(X,\Leq)$ satisfying~\eqref{axiom:LV}. Following \cite{akbartabatabai2021ImplicationSpacetime}, we interpret the localic cones as modal operators, where $\Down U$ reads as ``$U$ holds at some point in the past,'' and $\Up U$ as ``$U$ holds at some point in the future.'' Since we assume the localic cones preserve joins, together with the infinite distributive law we obtain new join-preserving maps $U\wedge \Up(-)$ and $U\wedge \Down(-)$, and therefore get new right adjoints:
\[
(U\wedge \Up(-))\dashv (U\to_\triup -)
\qquad\text{and}\qquad
(U\wedge \Down(-))\dashv (U\to_\tridown -).
\]
To spell this out explicitly, the new ``past Heyting implication'' $\to_\tridown$ is characterised uniquely as follows:
\[
\Down W\sqleq U\to V
\qquad\text{iff}\qquad
U\wedge \Down W\sqleq V
\qquad\text{iff}\qquad
W\sqleq U\to_\tridown V.
\]
The corresponding \emph{past modus ponens} law is obtained via:
\[
(U\to_\tridown V)\sqleq (U\to_\tridown V)
\qquad\text{if and only if}\qquad
U\wedge\Down(U\to_\tridown V)\sqleq V,
\]
which we interpret as: ``if $U$ holds, and it is true in the past that $U\to_\tridown V$, then $V$ holds.''

We can recover the causal Heyting implication in terms of the ordinary one via the following lemma.
\begin{lemma}\label{lemma:past heyting implication is direct image}
	For any ordered locale satisfying~\eqref{axiom:LV}:
	\[
	U\to_\triup V = \eta^\triup_\ast(U\to V)
	\qquad\text{and}\qquad
	U\to_\tridown V = \eta^\tridown_\ast(U\to V).
	\]
\end{lemma}
\begin{proof}
	First use~\eqref{axiom:LV} to evaluate $\eta^\tridown_\ast$:
	\begin{align*}
		\eta^\tridown_\ast(U\to V)
		&= \bigvee\left\{ \Down P\in \im(\Down): \Down P\sqleq U\to V\right\}
		\\&= 
		\Down\bigvee\left\{W\in\Opens X: \Down W\sqleq U\to V
		\right\}
		\\&=
		\Down(U\to_\tridown V).
	\end{align*}
	We are thus only left to show $\Down(U\to_\tridown V)\sqleq U\to_\tridown V$. But this follows by idempotence of the cones:
	\begin{align*}
		\Down(U\to_\tridown V)\sqleq U\to_\tridown V 
		&\quad\text{iff}\quad 
		U\wedge \Down \Down(U\to_\tridown V)\sqleq V
		\\&\quad\text{iff}\quad 
		U\wedge \Down(U\to_\tridown V)\sqleq V
		\\&\quad\text{iff}\quad 
		U\to_\tridown V\sqleq U\to_\tridown V,
	\end{align*}
	which is true.
\end{proof}


\begin{remark}
	Note that \cite{akbartabatabai2021ImplicationSpacetime} uses a join-preserving operation $\tridown\colon\Opens X\to \Opens X$ that is \emph{not} idempotent. In that case \cref{lemma:past heyting implication is direct image} weakens to $\eta_\ast^\tridown(U\to V) = \tridown(U\to_\tridown V)$, whereas $U\to_\tridown V$ is generally strictly smaller. They argue that $\tridown(U\to_\tridown V)$ and $U\to_\tridown V$ should actually be distinguished. In particular they argue against the correctness of the modus ponens law of the latter, which reads $U\wedge (U\to_\tridown V)\sqleq V$, because here ``[$U\to_\tridown V$] is just constructed and it can not be applicable at the moment,'' \cite[p.~166]{akbartabatabai2021ImplicationSpacetime}. Their Definition 10 actually defines ``spacetime'' just as a locale $X$ equipped with a certain type of join-preserving modal operator $\tridown\colon \Opens X\to \Opens X$. Note for instance that they assume $\tridown U\sqleq U$ (``if $U$ holds at some point in the past then $U$ holds''), which is markedly different from our $U\sqleq \Down U$ (``if $U$ holds then $U$ holds at some point in the past''). In future work, we could investigate the relationship between ordered locales $(X,\Leq)$ and these more general structures. We note especially that \cite{akbartabatabai2021ImplicationSpacetime} works in the more general setting of \emph{quantales}, which are related to linear logic, and therefore might provide new insights into possible categorifications of ordered locales to tensor topology.
\end{remark}

\section{Domains of dependence}\label{section:domains of dependence future work}
In this section we outline some ideas relating to a potential definition for sheaves of causal coverages, and discuss how the localic domains of dependence could be used to investigate holes in spacetime.

\subsection{Deterministic sheaves}\label{section:deterministic sheaves}
The incarnation of causal coverages as a type of Grothendieck topology (\cref{theorem:causal grothendieck topology}) raises the question of \emph{sheaves} \cite{maclane1994SheavesGeometryLogic}. Recall that a \emph{presheaf} on a locale $X$ is a functor $F\colon \Opens X^\op\to \Set$. We think of the set $F(U)$ as the collection of ``local data'' living on $U$. If $U\sqleq V$ then we get a \emph{restriction function} $(-)|_U\colon F(V)\to F(U)$. A presheaf $F$ is called a \emph{sheaf} when a family of compatible local data amalgamates uniquely to a piece of ``global data.'' An important example is the sheaf of real-valued functions on a topological space $S$:
	\[
		C(-,\mathbb{R})\colon \Opens S^\op\longrightarrow \Set,
	\]
the restriction function literally takes a real-valued function $f\colon V\to \mathbb{R}$ to its restriction $f|_U$. This is a sheaf since any locally compatible family of continuous real-valued functions can be glued uniquely to a continuous function on the union of their domains.

In summary, a sheaf tells us how ``local data'' is glued uniquely into ``global data.'' Extrapolating this type of behaviour, we see that for an appropriate notion of sheaf with respect to the causal coverage relation on an ordered locale, we additionally get that ``past data'' evolves deterministically into ``future data.'' Taking domains of dependence as fundamental, one possible definition could be as follows. For the ordinary definition of a sheaf see \cite[\S III.4]{maclane1994SheavesGeometryLogic}.

\begin{definition}
	Let $(\cat{C},J)$ be a site, equipped with a monad $D\colon \cat{C}\to\cat{C}$. A \emph{$D$-sheaf} on $\cat{C}$ is a presheaf $F\colon \cat{C}^\op\to\Set$ satisfying the following property:
	\begin{quote}
		if $(s_g)_{g\in R}$ is a compatible family of sections in $F$ indexed by a covering sieve $R\in J(C)$ on an object $C\in\cat{C}$, then there exists a unique amalgamation $s\in F(D(C))$ such that $s_g= F(\eta_C\circ g)(s)$ for all $g\in R$. 
	\end{quote}
\end{definition}

\begin{lemma}
	Any separated $D$-sheaf is a sheaf.
\end{lemma}
\begin{proof}
	Suppose $F\colon \cat{C}^\op \to \Set$ is a $D$-sheaf. Then a compatible family of sections $(s_g)_{g\in R}$ in $F$ indexed by $R\in J(C)$ has a unique amalgamated section $s\in F(D(C))$. We claim that $s|_C = F(\eta_C)(s)$ is the unique amalgamation of $s$ over $C$ that restricts to the $s_g$. Suppose that $t\in F(C)$ is another such section. Then $t|_{\dom(g)} = s_g = (s|_C)|_{\dom(g)}$ for every $g\in R$, and since $F$ is separated it follows that $s|_C = t$. 
\end{proof}

The interpretation is that an ordinary sheaf $F$ on $\cat{C}$ only allows you to glue local data on objects $C\in\cat{C}$, whereas the $D$-sheaf condition says that if you have a compatible family indexed by a covering sieve, then you can uniquely extend it to $D(C)$. For a locale this intuition becomes more concrete.

\begin{example}
	If $D\colon\Opens X\to \Opens X$ is a monad on a locale, the $D$-sheaf condition says the following. Consider a presheaf $F\colon \Opens X^\op \to \Set$, an open cover $(U_i)_{i\in I}$ of $U\in\Opens X$, and a compatible family $(s_i)_{i\in I}$ of section $s_i\in F(U_i)$. Then $F$ is a $D$-sheaf if there exists a unique amalgamation $s\in F(D(U))$ such that $s_i = s|_{U_i}$. Note that the restriction of $s$ to $U_i$ makes sense using the unit of $D$, which gives $U_i\sqleq U\sqleq D(U)$.
\end{example}

\begin{example}
	Consider Minkowski space $\mathbb{R}^n$, together with the future domain of dependence monad $D^+$. Then consider the presheaf $F\colon (\Opens \mathbb{R}^n)^\op \to \Set$ defined by
	\[
	F(U):= \{\text{solutions to the wave equation on $U$}\},
	\]
	and where $F(U\subseteq V)$ is just the restriction map of functions. This defines a $D^+$-sheaf, since initial data on $U$ uniquely determines a solution to the wave equation on $D^+(U)$ \cite[\S 2.4]{evans2022PartialDifferentialEquations}.
\end{example}

\subsection{Holes in spacetime}\label{section:holes in spacetime}
As mentioned in the introductory \cref{section:mathematical structure of spacetimes}, there are philosophical problems with the contemporary definition of spacetimes. One of which is the problem of \emph{hole-freeness} \cite{krasnikov2009EvenMinkowskiSpace,manchak2009SpacetimeHolefree}. An accessible overview of this problem is in \cite{robertsNotesOnHoles}.

Adopting the definition on \cite[p.~3]{krasnikov2009EvenMinkowskiSpace}, a spacetime $M$ is called \emph{hole-free} if for any achronal hypersurface $S\subseteq M$ and any isometric embedding $\pi$ of an open neighbourhood $U$ of the future domain of dependence $\Dchron^+(S)$ into another spacetime $M'$, we have ${\pi(\Dchron^+(S)) = \Dchron^+(\pi(S))}$. It is easy to see that the examples in \cref{figure:differences domains of dependence}(a) and (c) provide spacetimes that are not hole-free, since they can be embedded into Minkowski space, wherein the domains of dependence suddenly~``grow.''

However, the paper \cite{krasnikov2009EvenMinkowskiSpace} shows that even Minkowski space is not hole-free in this sense. Does using the \emph{localic} domain of dependence $\DLeq^+$ solve this problem? Concretely: is Minkowski space hole-free with respect to $\DLeq^+$?

%
%
%

\newpage\section{Questions and smaller ideas}
To conclude, we propose some questions and smaller ideas for future research.

\begin{enumerate}[label = \textbullet]
	\item The axiom~\eqref{axiom:P} is a somewhat \emph{ad hoc} solution to ensure that $\pt(X)$ has open cones. This condition is unsatisfactory from a localic perspective, since it is stated in terms of points. Similarly, the need to restrict to spaces with enough points in \cref{lemma:enough points and OC implies P} is unsatisfactory. See also the first bullet point in \cite[\S 9]{heunenSchaaf2024OrderedLocales}. Can these problems be remedied?	
	
	\item Does the forgetful functor $\OrdLoc\to \Loc$ define a topological functor in the sense of \cite[\S 21]{adamek1990AbstractConcreteCategories}? This would involve the structure of the lattice $\OL(X)$ of orders $\Leq$ from \cref{section:lattice of ordered locale orders}. The corresponding point-wise fact is that the forgetful functor $\Ord\to\Set$ of preordered sets defines a topological functor.
	
	\item In \cite{martin2012SpacetimeGeometryCausal} it is shown that the spacetime metric can be recovered from the causal order \emph{plus} a ``measurement.'' Similarly, can the structure of a spacetime $M$ be recovered from the ordered locale $\loc(M)$ together with a \emph{valuation}?
	
	\item We have seen in \cref{section:ordered locales from biframes} that ordered locales with~\eqref{axiom:LV} correspond to a certain class of biframes. How is parallel orderedness characterised in the biframe setting? 
	
	\item Relatedly, for the purpose of biframe compactifications \cite{schauerte1992Biframes}, it is interesting to ask when the biframe $(\Opens X,\im(\Up),\im(\Down))$ induced by an ordered locale is \emph{regular}. Note that if $X$ is parallel ordered, the biframe pseudocomplement reduces simply to Heyting negation:
		\[
			\bigvee\{\Up V\in\im(\Up): \Down U \wedge \Up V =\varnothing \} = \neg\Down U.
		\]
	This implies that the induced \emph{rather below} relations similarly reduce to the one from the original frame $\Opens X$. For example, the biframe induced by some ordered space with open cones is regular as soon as
		\[
			\Up U = \bigcup\{\Up V\in\im(\Up): \overline{\Up V}\subseteq \Up U\},
		\]
	and dually for down cones. A sufficient condition for this to be true is that $\overline{V}\subseteq U$ implies $\overline{\Up V}\subseteq \Up U$. In particular we ask: when is the biframe $(\Opens M,\im(I^+),\im(I^-))$ induced by a smooth spacetime regular? Since $M$ is a smooth manifold, it is a completely regular space, and hence $\Opens M$ is a completely regular frame. Therefore, the regularity of the induced biframe will come down to a compatibility condition between the causal structure and the topology of $M$. Does this correspond to a known causality condition from the causal ladder? As above, it would be sufficient if $\overline{V}\subseteq U$ implied $\overline{I^\pm(V)}\subseteq I^\pm(U)$. It is easily seen that this holds true in Minkowski space, but fails in Minkowski space with a point removed. In general, what is a good notion of compactification for ordered locales?
	
	\item As briefly alluded to in \cref{section:definition ordered locales}, it is possible to define locales \emph{internal} to a topos~\cite[Section~C1.6]{johnstone2002SketchesElephantTopos2}. For example, if $Y$ is locale, then a locale $X$ internal to the topos $\Sh(Y)$ of sheaves over $Y$ is equivalently a map of locales $X\to Y$, which is often interpreted as a type of bundle. Since $\Set \simeq \Sh(1)$, a locale internal to the topos $\Set$ is just a bundle $X\to 1$, which returns the ordinary definition. What is an ordered locale internal to a topos? In particular, how does the causal structure manifest in terms of the bundle picture $X\to Y$?
	
	\item How can causal order be interpreted on the \emph{measurable locales} from~\cite{pavlov2022GelfandtypeDualityCommutative}? Does this connect to causality of probability distributions in~\cite[\S 6]{eckstein2022CausalityTimeOrder}?
	
	\item How does our locale theoretic view of spacetimes compare to the domain theoretic treatment in \cite{martin2006DomainSpacetimeIntervals,martin2012SpacetimeGeometryCausal,mazibuko2023CausalStructureSpacetime}?
	
	\item Ordered spaces generalise to \emph{locally} ordered spaces, also called \emph{streams} \cite{krishnan2009ConvenientCategoryLocally,haucourt2012StreamsDSpacesTheir}. These are described by \emph{circulations}, which are cosheaves on the frame of opens of a space taking values in the category $\Ord$ of preordered sets. This definition straightforwardly generalises to the setting of locales. Does the adjunction between ordered spaces and ordered locales lift to these categories of cosheaves? 
	
	\item Can ordered locales provide a setting in which to study causality and contextuality, as in \cite{abramsky2011SheafTheoreticStructureNonLocality,gogioso2021SheafTheoreticStructureDefinite}?
	
	\item \emph{Quantales} are a \emph{noncommutative} generalisation of locales, with close ties to the theory of C$^\ast$-algebras \cite{coniglio2000NonCommutativeTopologyQuantales,mulvey2001QuantisationPoints}. Do ordered locales generalise to this setting, and can the resulting spaces be interpreted as \emph{quantum spacetimes?}
	
	\item Can causal coverages be realised as a genuine Grothendieck topology on some category? Alternatively, what is the theory of the newly introduced Grothendieck topologies from \cref{definition:modified grothendieck topology}? \enlargethispage{\baselineskip}What is the appropriate notion of \emph{sheaf} in this setting?
\end{enumerate}


\printbibliography[heading=bibintoc]
\appendix
\part*{Appendix}\label{chapter:appendix}
\markboth{Appendix}{}\addcontentsline{toc}{chapter}{Appendix}
\chapter{Basic category theory}\label{section:categories}
In this section we provide the very basic idea behind \emph{categories}. This is mostly to serve as a guideline, and refer to the classic textbook \cite{maclane1998CategoriesWorkingMathematician} for a more detailed and precise development of the theory. We also recommend the lecture notes~\cite{vanOosten2016basic}.

\begin{definition}\label{definition:category}
	A \emph{category} $\cat{C}$ consists of collections of:
		\begin{enumerate}[label = (\roman*)]
			\item \emph{objects} $C\in\cat{C}$;
			\item \emph{arrows} or \emph{morphisms} $f\in \cat{C}(A,B)$ between the objects $A,B\in\cat{C}$, also written $f\colon A\to B$ or $A\xrightarrow{f}B$;
		\end{enumerate}
	such that:
		\begin{enumerate}[label = (\roman*)]\setcounter{enumi}{2}
			\item every object $A\in\cat{C}$ admits an \emph{identity arrow} $\id_A\in\cat{C}(A,A)$;
			\item if $f\colon A\to B$ and $g\colon B\to C$ are any morphisms, we get its \emph{composition} $g\circ f\colon A\to C$;
			\item the identity arrows form units for the composition operator: if $f\colon A\to B$ is any arrow then
				\[
					\id_B \circ f =f = f\circ \id_A;
				\]
			\item the composition operation is associative, meaning for any triple of composable arrows:
				\[
					f\circ (g\circ h) = (f\circ g)\circ h.
				\]
		\end{enumerate}
\end{definition}

\begin{definition}\label{definition:epi mono iso}
	There are several classes of well-behaved arrows in a category $\cat{C}$ that we can consider. An \emph{isomorphism} is an arrow $f\colon C\to D$ that admits an \emph{inverse} $f^{-1}\colon D\to C$, which is the unique arrow satisfying $f\circ f^{-1}=\id_D$ and $f^{-1}\circ f = \id_C$.

	An \emph{epimorphism} is an arrow $e\colon C\to D$ such that for any parallel arrows $f,g\colon D\to E$, the condition $f\circ e = g\circ e$ implies $f=g$.
	
	Dually, a \emph{monomorphism} is an arrow $m\colon C\to D$ such that for any parallel arrows $f,g\colon E\to C$ we have that $m\circ f = m\circ g$ implies $f=g$.
\end{definition}

\begin{lemma}\label{lemma:cancellation property epi mono}
	If $f\circ g$ is an epimorphism, then so is $f$. If $f\circ g$ is a monomorphism, then so is $g$.
\end{lemma}

\begin{definition}\label{definition:functors}
	A \emph{functor} $F\colon \cat{C}\to \cat{D}$ between categories consists of an assignment $C\mapsto F(C)$ of objects, also often written just $FC$, together with an assignment
		\[
			\cat{C}(C,D)\longrightarrow \cat{D}(FC,FD);
			\qquad F\longmapsto F(f)
		\]
	on arrows, such that:
		\begin{enumerate}[label = \textbullet]
			\item $F$ preserves identity arrows: $F(\id_C)= \id_{FC}$ for all $C\in\cat{C}$;
			\item $F$ preserves composition: $F(g\circ f) = F(g)\circ F(f)$.
		\end{enumerate}
\end{definition}

Informally, a \emph{diagram} in a category $\cat{C}$ is a collection of objects and arrows between them. A diagram is said to \emph{commute} if any pair of directed paths with the same start and endpoint compose to the same arrow. Formally, a diagram is a functor $\cat{I}\to\cat{C}$, where the \emph{index category} $\cat{I}$ determines the \emph{shape} of the diagram.

\begin{definition}\label{definition:natural transformations}
	Consider two parallel functors $F,G\colon \cat{C}\to \cat{D}$. A \emph{natural transformation} $\nu\colon F\to G$ from $F$ to $G$ consists of a family $(\nu_C\colon FC\to GC)_{C\in\cat{C}}$ of arrows in $\cat{D}$ such that for every arrow $f\colon C\to D$ in $\cat{C}$ the following square commutes:
		\[
			\begin{tikzcd}
				FC & GC\\
				FD & GD.
				\arrow["F(f)"', from= 1-1, to= 2-1]
				\arrow["\nu_C", from = 1-1, to = 1-2]
				\arrow["G(f)", from = 1-2, to = 2-2]
				\arrow["\nu_D"', from = 2-1, to = 2-2]
			\end{tikzcd}
		\]
\end{definition}

\begin{definition}\label{definition:monads}
	A \emph{monad} on a category $\cat{C}$ is a triple $(T,\mu,\eta)$ consisting of an endofunctor $T\colon \cat{C}\to \cat{C}$, together with natural transformations $\mu\colon T\circ T\to T$ called the \emph{multiplication}, and $\eta\colon \id_\cat{C}\to T$ called the \emph{unit}, such that the following two diagrams commute:
		\[
			\begin{tikzcd}[ampersand replacement=\&]
				T \& {T\circ T} \\
				{T\circ T} \& T
				\arrow["{\eta T}", from=1-1, to=1-2]
				\arrow["{T \eta}"', from=1-1, to=2-1]
				\arrow[from=1-1, to=2-2,equal]
				\arrow["\mu", from=1-2, to=2-2]
				\arrow["\mu"', from=2-1, to=2-2]
			\end{tikzcd}
			\qquad\text{and}\qquad
			\begin{tikzcd}[ampersand replacement=\&]
				{T\circ T\circ T} \& {T\circ T} \\
				{T\circ T} \& {T.}
				\arrow["T\mu", from=1-1, to=1-2]
				\arrow["{\mu T}"', from=1-1, to=2-1]
				\arrow["\mu", from=1-2, to=2-2]
				\arrow["\mu"', from=2-1, to=2-2]
			\end{tikzcd}
		\]
\end{definition}

\begin{definition}\label{definition:opposite category}
	If $\cat{C}$ is a category, we define its \emph{opposite} as the category $\cat{C}^\op$ containing the same objects, but where the morphisms are formally reversed. Thus an arrow $g^\op \colon C\to D$ in $\cat{C}^\op$ consists of an arrow $g\colon D\to C$ in $\cat{C}$, and composition is defined via $f^\op \circ g^\op := (g\circ f)^\op$.
\end{definition}

\begin{definition}\label{definition:adjunction}
	Consider functors $F\colon \cat{C}\to \cat{D}$ and $G\colon \cat{D}\to \cat{C}$. We say $F$ is \emph{left adjoint} to $G$, and $G$ is \emph{right adjoint} to $F$ if there is a natural bijection
		\[
			\cat{D}(FC,D)\cong \cat{C}(C,GD).
		\]
	In that case, we write $F\dashv G$.
	
	The \emph{unit} of an adjunction $F\dashv G$ is the natural transformation ${\eta\colon \id_\cat{C}\to G\circ F}$ whose components are defined as the unique arrow $\eta_C\colon C\to GFC$ corresponding to the identity arrow $\id_{FC}\in \cat{D}(FC,FC)$ in the above natural bijection. Dually, the \emph{counit} of $F\dashv G$ is the natural transformation $\epsilon\colon F\circ G\to \id_{\cat{D}}$, where $\epsilon_D\colon FGD\to D$ is the unique arrow corresponding to the identity arrow $\id_{GD}$ in the natural bijection. With this notation, the natural bijection can actually be rewritten explicitly as:
		\begin{align*}
			\cat{D}(FC,D)  &\xrightarrow{~ \sim ~}  \cat{C}(C,GD)\\
			f&\xmapsto{\hspace{3ex}} G(f)\circ \eta_C\\
			\epsilon_D\circ F(g)&{~\reflectbox{\ensuremath{\xmapsto{\hspace{3ex}}}}} ~g	.	
		\end{align*}
	In fact, adjunctions can be equivalently characterised purely in terms of $\eta$ and $\epsilon$ satisfying the so-called \emph{triangle identities} and \emph{zig-zag laws} \cite{nlab2024adjunction}, but we do not need to do so here.
\end{definition}

\begin{definition}\label{definition:initial and terminal objects}
	Let $\cat{C}$ be a category. A \emph{terminal object} is an object $1\in\cat{C}$ such that for every $A\in\cat{C}$ there exists a unique arrow $A\to 1$.
	
	Dually, an \emph{initial object} is an object $0\in \cat{C}$ such that for every $A\in\cat{C}$ there exist a unique arrow $0\to A$.
\end{definition}

\begin{definition}\label{definition:limits}
	For any object $A\in\cat{C}$ in a category, the \emph{constant functor} $\Delta_A\colon \cat{I}\to\cat{C}$ sends every object of $\cat{I}$ to $A$, and every arrow of $\cat{I}$ to $\id_A$.
	
	Now consider a diagram $F\colon \cat{I}\to \cat{C}$ in a category. A \emph{cone} on $F$ consists of an object $A\in\cat{C}$ called the \emph{vertex}, together with a natural transformation $\nu\colon \Delta_A\to F$. Explicitly, this consists of a family of arrows $(\nu_i\colon A\to F(i))_{i\in \cat{I}}$, such that for every arrow $k\colon i\to j$ the following triangle commutes:
		\[
			\begin{tikzcd}[ampersand replacement=\&]
				\& A \\
				{F(i)} \&\& {F(j).}
				\arrow["{\nu_i}"', from=1-2, to=2-1]
				\arrow["{\nu_j}", from=1-2, to=2-3]
				\arrow["{F(k)}"', from=2-1, to=2-3]
			\end{tikzcd}
		\]
	A \emph{morphism of cones} $(A,\nu)\to (B,\lambda)$ on the functor $F$ consists of an arrow $f\colon A\to B$ in $\cat{C}$ such that $\lambda_i\circ f = \nu_i$ for every $i\in\cat{I}$. This gives the category of cones on $F$.
	
	A \emph{limiting cone} for $F\colon \cat{I}\to\cat{C}$ is a terminal object in the category of cones over $F$. Explicitly, it consists of a cone $(P,\lambda)$ on $F$, such that for every other cone $(A,\nu)$ on $F$ there exists a unique arrow $\alpha\colon A\to P$ such that $\lambda_i \circ \alpha = \mu_i$ for all $i\in\cat{I}$:
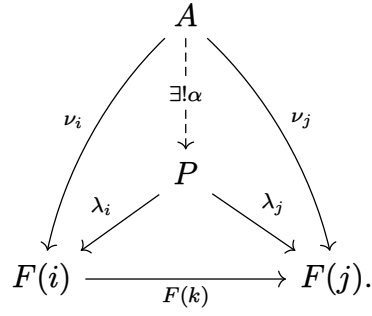
		\[
			\begin{tikzcd}[ampersand replacement=\&]
				\& A \\
				\\
				\& P \\
				{F(i)} \&\& {F(j).}
				\arrow["{\exists!\alpha}"{description}, dashed, from=1-2, to=3-2]
				\arrow["{\nu_i}"', curve={height=12pt}, from=1-2, to=4-1]
				\arrow["{\nu_j}", curve={height=-12pt}, from=1-2, to=4-3]
				\arrow["{\lambda_i}"', from=3-2, to=4-1]
				\arrow["{\lambda_j}", from=3-2, to=4-3]
				\arrow["{F(k)}"', from=4-1, to=4-3]
			\end{tikzcd}
		\]
	The vertex $P$, together with the data of the arrows $\lambda_i$, is called the \emph{limit} of $F$, and sometimes denoted $\lim F$ or $\lim_{i\in\cat{I}}F(i)$. When the limit of $F$ exists, it is unique up to unique isomorphism.
\end{definition}

\begin{remark}
	The \emph{colimit} of $F\colon\cat{I}\to\cat{C}$ is defined dually. In diagrams, the colimit of $F$ is the initial cocone $(Q,\iota)$:
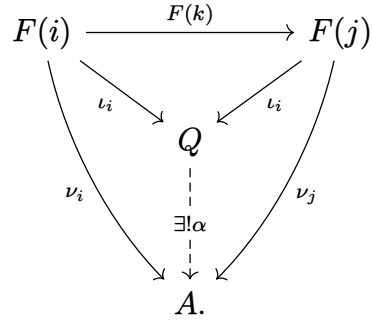
		\[
			\begin{tikzcd}[ampersand replacement=\&]
				{F(i)} \&\& {F(j)} \\
				\& Q \\
				\\
				\& {A.}
				\arrow["{F(k)}", from=1-1, to=1-3]
				\arrow["{\iota_i}"', from=1-1, to=2-2]
				\arrow["{\nu_i}"', curve={height=12pt}, from=1-1, to=4-2]
				\arrow["{\iota_i}", from=1-3, to=2-2]
				\arrow["{\nu_j}", curve={height=-12pt}, from=1-3, to=4-2]
				\arrow["{\exists!\alpha}"{description}, dashed, from=2-2, to=4-2]
			\end{tikzcd}
		\]
\end{remark}

\chapter{Lattice theory}\label{section:lattice theory}
In this section we give an overview of the relevant lattice theoretic background for this thesis. We draw from many sources, amongst which mainly \cite{johnstone1982StoneSpacesa,vickers1989TopologyLogic,gierz2003ContinuousLatticesDomains,picado2012FramesLocalesTopology}. We have already seen the definition of distributive lattices, complete lattices, and frames in \cref{section:locales and spaces}. Here are some elementary lemmas.


\begin{lemma}\label{lemma:join respects inclusion}
	Let $L$ be a complete lattice, and consider an arbitrary family $(x_i)_{i\in I}$ of elements in $L$. If $x_i\sqleq y$ for all $i\in I$, then $\bigvee_{i\in I}x_i\sqleq y$.
\end{lemma}
\begin{proof}
	This follows immediately from the definition of a supremum.
\end{proof}


\begin{lemma}\label{lemma:join respects parametrised inclusion}
	In a complete lattice $L$, if $x_i\sqleq y_i$ for all $i\in I$, then ${\bigvee_{i\in I} x_i\sqleq \bigvee_{i\in I}y_i}$.
\end{lemma}
\begin{proof}
	That $x_j\sqleq y_j$ implies $x_j\sqleq \bigvee_{i\in I} y_i$ for all $j\in J$, so the result follows by the previous \cref{lemma:join respects inclusion}.
\end{proof}

\begin{lemma}\label{lemma:monotone laxly respects meets and joins}
	Let $h\colon L\to M$ be an order-preserving function between complete lattices. Then for any family $(x_i)_{i\in I}$ in $L$ we have
	\[
	\bigvee_{i\in I} h(x_i)\sqleq h\left(\bigvee_{i\in I} x_i\right)
	\qquad\text{and}\qquad
	h\left(\bigwedge_{i\in I} x_i\right)\sqleq \bigwedge_{i\in I}h(x_i).
	\]
\end{lemma}
\begin{proof}
	It is clear that these inclusions hold whenever the indexing set $I$ is empty, so we may assume without loss of generality that it is not. In that case, for every $j\in I$ we have $x_j\sqleq \bigvee_{i\in I} x_i$, so since $h$ preserves order we get $h(x_j)\sqleq h\left(\bigvee_{i\in I} x_i\right)$. Since this holds for arbitrary $j\in I$, the first inclusion follows by \cref{lemma:join respects inclusion}. The second inclusion involving meets is proved dually. 
\end{proof}

\begin{theorem}[{\cite[\S AI.4.3]{picado2012FramesLocalesTopology}}]\label{theorem:lattice has joins iff it has meets}
	A lattice is complete if and only if it admits all joins, if and only if it admits all meets.
\end{theorem}
\section{Galois adjunctions}\label{section:galois adjunctions}
In this section we collect the needed results on \emph{Galois adjunctions}, also known as \emph{Galois connections}. The material in this section is taken mostly from the appendix \cite[\S AI.5]{picado2012FramesLocalesTopology}. The labels attached to the results in this section refer to sections in their Appendix I.

\begin{definition}\label{definition:galois adjunction}
	A \emph{Galois adjunction} between posets $P$ and $Q$ constists of monotone functions $f\colon P\to Q$ and $g\colon Q\to P$ satisfying the following property:
	\[
	f(x)\leq y
	\qquad\text{if and only if}\qquad
	x\leq g(y).
	\]
	In that case we write:
	\[
	\begin{tikzcd}[cramped,column sep=large]
		Q & P,
		\arrow[""{name=0, anchor=center, inner sep=0}, "g"', shift right=2, from=1-1, to=1-2]
		\arrow[""{name=1, anchor=center, inner sep=0}, "f"', shift right=2, from=1-2, to=1-1]
		\arrow["\dashv"{anchor=center, rotate=-90}, draw=none, from=1, to=0]
	\end{tikzcd}
	\]
	or more compactly: $f\dashv g$. The function $f$ is called the \emph{left adjoint} of $g$, and $g$ is the \emph{right adjoint} of $f$.
\end{definition}

\begin{remark}
	A Galois adjunction is just an adjunction between functors in the categorical sense (\cref{definition:adjunction}), restricted to the setting of monotone functions between posets. Often we shall simply use the word ``adjunction,'' unless we need emphasis for the fact the domains are posets.
\end{remark}

\begin{proposition}[\S 5.3]\label{proposition:galois adjunction iff unit and counit}
	There is an adjunction $f\dashv g$ if and only if
	\[
	f\circ g(y)\leq y
	\qquad\text{and}\qquad
	x\leq g\circ f(x).
	\]
\end{proposition}

\begin{corollary}\label{corollary:adjoint fgf=f}
	If $f\dashv g$ then $f\circ g\circ f = f$ and $g\circ f\circ g = g$.
\end{corollary}

\begin{theorem}[\S 5.4]\label{theorem:adjoints preserve limits}
	If $f\dashv g$ then $f$ preserves all joins, and $g$ preserves all meets.
\end{theorem}

The following shows that if we restrict to the setting of complete lattices, the preservation of joins or meets is actually equivalent to the existence of an adjoint.

\begin{theorem}[\S 5.5]\label{theorem:join preserving implies left adjoint}
	A monotone map between complete lattices admits a left/right adjoint if and only if preserves all meets/joins. Explicitly, if $f\dashv g$ we obtain:
	\[
	f(x) = \bigwedge\left\{y : x\leq g(y)\right\}
	\quad\text{and}\quad
	g(y) = \bigvee\{x: f(x)\leq y\}.
	\]
\end{theorem}


Lastly, we connect the following result from \cite[Exercise~7.13]{davey2002IntroductionLatticesOrder}.
\begin{proposition}\label{proposition:galois adjunction implies surjection iff injection}
	If $f\dashv g$, then $f$ is surjective if and only if $g$ is injective.
\end{proposition}

\section{Heyting algebras}\label{section:heyting algebras}
Heyting algebras arise naturally in frame theory. Namely, the infinite distributivity law of a frame $L$ says precisely that for every $x\in L$ the function $x\wedge - \colon L\to L$ preserves all joins. By \cref{theorem:join preserving implies left adjoint} this means there exists an adjunction $(x\wedge -)\dashv (x\to -)$, where the right adjoint is calculated at $y\in L$ as:
\[
x\to y = \bigvee\{z\in L: x\wedge z\sqleq y\}.
\]
This is the largest element in $L$ such that when intersected with $x$ is contained in $y$. The operation $\to$ is called the \emph{Heyting implication} of $L$, and fits into the more general framework of Heyting algebras. For more details we refer to \cite[\S AI.7]{picado2012FramesLocalesTopology}.

\begin{definition}\label{definition:heyting algebra}
	A \emph{Heyting algebra} is a lattice $L$ with bottom element, together with a binary operation $L\times L\to L: (x,y)\mapsto x\to y$, called the \emph{Heyting implication}, satisfying
	\[
	x\wedge y \sqleq z
	\qquad\text{if and only if}\qquad
	x\sqleq y\to z.
	\]
	
	In a Heyting algebra $L$ the \emph{Heyting negation} (or \emph{pseudocomplement}) of $x\in L$ is defined as $\neg x:= x\to \bot$.
\end{definition}

\begin{remark}
	From \cref{theorem:adjoints preserve limits} it follows that \emph{complete} Heyting algebras (i.e.~those for which the underlying lattice $L$ is complete) are just frames. In fact, a complete lattice admits a Heyting implication if and only if it is a frame. So while at the object-level these structures are identical, what distinguishes them are their categories: morphisms of frames are not required to respect the Heyting implication.
\end{remark}

\begin{example}\label{example:heyting implication in space}
	It is instructive to calculate the Heyting implication in the frame of opens $\Opens S$ of a topological space. First, note that in the powerset $\Powerset(S)$ we have by elementary set theory that $A\cap B\subseteq C$ if and only if $B\subseteq C\cup (S\setminus A)$. Using this, we calculate:
	\begin{align*}
		U\to V &= \bigcup\{W\in \Opens S: U\cap W\subseteq V\}
		\\&=
		\bigcup\{W\in\Opens S: W\subseteq V\cup (S\setminus U)\}
		\\&=
		\left(V\cup (S\setminus U)\right)^\circ.
	\end{align*}
	In particular, the Heyting negation of an open subset is calculated as the interior of its complement:
	\[
	\neg U = \left(S\setminus U\right)^\circ.
	\]
\end{example}

\begin{definition}\label{definition:boolean algebra}
	A \emph{Boolean algebra} is a Heyting algebra $L$ where $\neg\neg x = x$ for all $x\in L$; see \cite[Proposition~AI.7.4.3]{picado2012FramesLocalesTopology}. The lattice operations and Heyting implication are interpreted straightforwardly as the classical logical connectives.
\end{definition}

%

We list some properties of the Heyting negation~\cite[\S AI.7]{picado2012FramesLocalesTopology}.

\begin{lemma}\label{lemma:properties of heyting negation}
	The Heyting negation in a Heyting algebra $L$ satisfies:
	\begin{enumerate}[label = (\alph*)]
		\item $x\wedge y = \bot$ if and only if $x\sqleq \neg y$;
		\item if $x\sqleq y$ then $\neg y \sqleq \neg x$;
		\item $x\sqleq \neg\neg x$;
		\item $\neg\neg\neg x = \neg x$;
		\item $\neg\neg (x\wedge y) = \neg\neg x \wedge \neg\neg y$;
		\item $\neg\left(\bigvee_{i\in I} x_i\right) = \bigwedge_{i\in I}\neg x_i$, whenever the join exists. \hfill \emph{(first De Morgan law)}
	\end{enumerate}
\end{lemma}

\begin{corollary}\label{lemma:double negation is monad}
	For any Heyting algebra $L$, the double negation map $x\mapsto \neg\neg x$ defines a monad on $L$.
\end{corollary}

\begin{corollary}\label{corollary:heyting negation gives adjunction}
	For any Heyting algebra $L$, there is a Galois adjunction
		\[
			\begin{tikzcd}[column sep =large]
				L & L^\op.
				\arrow[""{name=0, anchor=center, inner sep=0}, "\neg", shift left=2, from=1-1, to=1-2]
				\arrow[""{name=1, anchor=center, inner sep=0}, "\neg", shift left=2, from=1-2, to=1-1]
				\arrow["\dashv"{anchor=center, rotate=-90}, draw=none, from=0, to=1]
			\end{tikzcd}
		\]
\end{corollary}
\begin{proof}
	This follows simply by \cref{lemma:properties of heyting negation}(a), which shows $x\sqleq \neg y$ is determined by the symmetric condition $y\wedge x= x\wedge y = \bot$. 
\end{proof}

\section{The frame of ideals of a distributive lattice}\label{section:frame of ideals}
To any distributive lattice $L$ we can associate a frame $\Ideals(L)$ of its ideals. The exposition in this section closely follows \cite[Section~VII.4]{picado2012FramesLocalesTopology}. 

\begin{definition}\label{definition:ideals}
	An \emph{ideal} in a bounded distributive lattice $L$ is a subset $I\subseteq L$ satisfying the following conditions:
		\begin{enumerate}[label = (I\arabic*)]\setcounter{enumi}{-1}
			\item\label{axiom:I0} $\bot\in I$;
			\item\label{axiom:I1} if $x,y\in I$ then $x\vee y\in I$; \hfill(\emph{upwards directedness})
			\item\label{axiom:I2} if $y\sqleq x\in I$ then $y\in I$. \hfill(\emph{down closure})
		\end{enumerate}
	Any element $x\in L$ generates a \emph{principal ideal} $\ideal(x):=\{y\in L:y\sqleq x\}$.
\end{definition}

\begin{remark}
	Note that for distributive lattices, axiom~\ref{axiom:I2} is equivalent to saying $x\wedge y\in I$ whenever $x\in I$. This makes the analogy to ideals in ring theory more apparent.
\end{remark}

The definition of an ideal makes sense more generally in the setting of non-bounded distributive lattices, where~\ref{axiom:I0} can be replaced by the condition that $I$ is non-empty. For our purposes, the generality provided by distributive lattices will suffice.

The following result and proof is from \cite[Proposition~VII.4.1.1]{picado2012FramesLocalesTopology}. Recall that a frame $L$ is \emph{compact} if whenever $\bigvee_{i\in I} x_i = \top$, there exist a finite subset $F\subseteq I$ such that $\bigvee_{i\in F} x_i = \top$. 

\begin{proposition}\label{proposition:frame of ideals}
	For any bounded distributive lattice $L$, the set of ideals $\Ideals(L)$, ordered by inclusion, is a compact frame.
\end{proposition}
\begin{proof}
	We first prove that $\Ideals(L)$ defines a frame. To start, note that $I=L$ and $I=\{\bot\}$ define ideals in $L$, forming the top and bottom element of $\Ideals(L)$, respectively. More generally, it is easy to verify that the set-theoretic intersection of any family $(I_a)_{a\in A}$ of ideals in $L$ again forms an ideal. This defines their meet in $\Ideals(L)$. We claim that
		\[
			\bigvee_{a\in A} I_a
			=
			\left\{
				\bigvee F: \text{ finite } F\subseteq\bigcup_{a\in A}I_a
			\right\}
		\]
	defines their join. First we show this is an ideal: axiom~\ref{axiom:I0} follows immediately by setting $F$ to be the empty set, and axiom~\ref{axiom:I1} follows since the join over two finite collections can be described as the join of another finite collection. For~\ref{axiom:I2}, if $y\sqleq \bigvee F$ for some finite $F\subseteq \bigcup_{a\in A}I_a$, then the distributivity of $L$ gives $y = \bigvee\{y\wedge z:z\in F\}$, where again $\{y\wedge z:z\in F\}\subseteq \bigcup_{a\in A}I_a$. This proves $\bigvee_{a\in A}I_a\in\Ideals(L)$.
	
	To see that it is an upper bound, note that by setting $F$ to be a singleton we get $I_b\subseteq \bigvee_{a\in A}I_a$ for every $b\in A$. If $J\in\Ideals(L)$ is another upper bound for this family, meaning $\bigcup_{a\in A}I_a\subseteq J$, then axiom~\ref{axiom:I1} guarantees $\bigvee_{a\in A}I_a\subseteq J$. Hence it is the least upper bound.
	
	Lastly, we show that $\Ideals(L)$ is compact. Suppose that $\bigvee_{a\in A}I_a = L$, the top element of $\Ideals(L)$. Since $\top\in L$, this means that we can find a finite set ${F\subseteq \bigcup_{a\in A}I_a}$ such that $\top =\bigvee F$. In particular there exists a finite subset $B\subseteq A$ such that ${F\subseteq \bigcup_{b\in B}I_b}$, and hence $\top\in \bigvee_{b\in B}I_b$. But the latter implies $\bigvee_{b\in B}I_b=L$ by downwards closure, so we are done.
\end{proof}

\begin{remarknumbered}\label{remark:ideal generated by set}
	If $A\subseteq L$ is some subset, not necessarily an ideal, using the complete lattice structure of $\Ideals(L)$ we can define the ideal \emph{generated} by $A$ as the smallest ideal containing $A$:
		\[
			\langle A\rangle := \bigcap\{I\in\Ideals (L): A\subseteq I\}.
		\]
	For instance: $\ideal(x)=\langle\{x\}\rangle$.
\end{remarknumbered}

\begin{lemma}\label{lemma:principal ideals are a sublattice}
	For any $x,y\in L$ in a distributive lattice:
	\[
	\ideal(x\vee y) = \ideal(x)\vee \ideal(y)
	\qquad\text{and}\qquad
	\ideal(x\wedge y) = \ideal(x)\cap \ideal(y).
	\]
\end{lemma}
\begin{proof}
	The first follows by the description of joins of ideals as in the proof of \cref{proposition:frame of ideals}, and the second follows by noting that $z\sqleq x\wedge y$ if and only if $z\sqleq x$ and $z\sqleq y$.
\end{proof}

\begin{remark}
	This means that the principal ideals of $L$ form a sublattice of $\Ideals(L)$.
\end{remark}

\begin{lemma}
	If $J\subseteq L$ is an ideal, then $J= \bigvee_{x\in J}\ideal(x)$ in $\Ideals(L)$.
\end{lemma}

In the case where $L$ is already a frame, the situation is particularly nice.

\begin{lemma}\label{lemma:ideal adjunction}
	For any frame $L$, there is a Galois adjunction%
		\begin{tikzcd}
			{\Ideals(L)} & L
			\arrow[""{name=0, anchor=center, inner sep=0}, "v", shift left=2, from=1-1, to=1-2]
			\arrow[""{name=1, anchor=center, inner sep=0}, "\alpha", shift left=2, from=1-2, to=1-1]
			\arrow["\dashv"{anchor=center, rotate=-90}, draw=none, from=0, to=1]
		\end{tikzcd}
	defined by the functions $v\colon I\mapsto \bigvee I$ and $\alpha\colon x\mapsto \ideal(x)$.
\end{lemma}
\begin{proof}
	It is easy to see that both $v$ and $\alpha$ are monotone. Moreover $v\circ\alpha(x) = \bigvee\ideal(x) = x$, and $I\subseteq \alpha\circ v(I) = \ideal(\bigvee I)$. 
\end{proof}

\begin{corollary}
	For every frame $L$, there is a map of frames $\Ideals(L)\to L$ defined by $I\mapsto \bigvee I$.
\end{corollary}
\begin{proof}
	By \cref{lemma:ideal adjunction} the map $v$ preserves all suprema, since it admits a right adjoint. We are only left to show it preserves finite meets, for which now, since $v$ is monotone, it suffices to show $v(I)\wedge v(J)\sqleq v(I\cap J)$. First, using infinite distributivity we can write
		\[
			v(I)\wedge v(J)
			=
			\bigvee I\wedge \bigvee J 
			=
			\bigvee\{x\wedge y:x\in I,y\in J\}.
		\]
	By axiom~\ref{axiom:I2}, if $x\in I$ and $y\in J$ then $x\wedge y\in I\cap J$. This immediately implies $v(I)\wedge v(J)\sqleq v(I\cap J)$, as required.
\end{proof}

In the language of locales, this map of frames can be interpreted as an embedding $X\to \Ideals(X)$ of a locale into a type of \emph{compactification}.


Recall that an element $x$ in a lattice is called \emph{inaccessible by directed joins} if for every directed set $D$ such that $\bigvee D = x$ there exists $d\in D$ such that $x=d$.

\begin{proposition}[{\cite[Theorem~8]{birkhoff1948RepresentationsLatticesSets}}]
	An ideal $J$ of a lattice $L$ is principal if and only if it is inaccessible by directed joins in $\Ideals(L)$.
\end{proposition}
\begin{proof}
	We first show that if $x\in L$, then the principal ideal $\ideal(x)$ is inaccessible by directed joins. Suppose that $(J_d)_{d\in D}$ is a directed family in $\Ideals(L)$ such that $\bigvee_{d\in D}J_d= \ideal(x)$. In particular $x\in \ideal(x)$, so by the description of the join of ideals in the proof of \cref{proposition:frame of ideals} we can find finitely many $y_{d_i}\in J_{d_i}$ such that $\bigvee_{i=0}^n y_{d_i} = x$. Since $(J_d)_{d\in D}$ is directed, there furthermore exists $e\in D$ such that $\bigvee_{i=0}^n J_{d_i}\subseteq J_e$. Hence $x\in J_e$, from which it is now straightforwardly follows that $J_e=\ideal(x)$.
	
	For the converse, note first that $J = \bigvee_{x\in J}\ideal(x)$. By \cref{lemma:principal ideals are a sublattice} and axiom~\ref{axiom:I1} this join is directed, so if $J$ is inaccessible by directed joins there exists $x\in J$ such that $J =\ideal(x)$.
\end{proof}

\chapter{Locale theory}\label{section:locale theory}
In this section we provide further technical details about locales that were deferred from the introductory \cref{section:locales and spaces} and \cref{section:adjunction top and loc}.

A map of locales $f\colon X\to Y$ is solely defined in tems of a map of frames $f^{-1}\colon \Opens Y\to \Opens X$, called the \emph{preimage map} of $f$. This map preserves all joins, so by \cref{theorem:join preserving implies left adjoint} admits a right adjoint $f^{-1}\dashv f_\ast\colon \Opens X\to \Opens Y$ called the \emph{direct image} map of $f$.

\begin{lemma}\label{lemma:locale map right adjoint}
	If $f\colon X\to Y$ is a map of locales, then for all $U\in \Opens X$:
	\[
	f_\ast(U) = \bigvee\left\{ V\in \Opens Y: f^{-1}(V)\sqleq U\right\}.
	\]
\end{lemma}
\begin{proof}
	This is just \cref{theorem:join preserving implies left adjoint}.
\end{proof}
Thus $f_\ast(U)$ is the largest open region in $Y$ whose preimage is contained in~$U$. It is sometimes called the \emph{direct image} of $U$ under $f$. Dually, $f^{-1}(V)$ is the smallest open in $X$ whose direct image is contained in $V$. For continuous functions $g\colon S\to T$ between topological spaces the direct image $g_\ast(U)$ of $U\in\Opens S$ is \emph{not} to be confused with the set-theoretic image (which is not always open). Being a right adjoint, $f_\ast$ preserves all meets, but it does not in general preserve arbitrary joins, and is hence not a frame map. The theory of locales can equivalently be developed in terms of these types of maps, as is done in \cite{picado2012FramesLocalesTopology}. Our exposition sticks closer to the preimage maps, but we make free use of both.

\begin{definition}\label{definition:locale basis}
	A \emph{basis} for a locale $X$ is a subset $\mathcal{B}\subseteq \Opens X$ such that every open region of $X$ can be written as a join of elements in $\mathcal{B}$.
\end{definition}

If the codomain of $f$ is generated by some basis, then the formula for $f_\ast$ can be simplified.
\begin{lemma}
	\label{lemma:locale map right adjoint using basis}
	Let $f\colon X\to Y$ be a map of locales, and let $\mathcal{B}$ be a basis for $Y$. Then for all $U\in \Opens X$:
	\[
	f_\ast(U) = \bigvee\left\{B\in\mathcal{B}:f^{-1}(B)\sqleq U\right\}.
	\]
\end{lemma}
\begin{proof}
	The inclusion from right to left is obvious. Conversely, by \cref{lemma:locale map right adjoint} it suffices to show that any $V\in\Opens Y$ with $f^{-1}(V)\sqleq U$ is contained in the right hand side. But this follows directly when writing $V$ as a join of elements in $\mathcal{B}$.
\end{proof}

\section{Points of a locale}\label{section:points of locales}
In this section we describe the three most common ways to talk about the notion of ``point'' in a locale. Our main reference is \cite[Section~II.3]{picado2012FramesLocalesTopology}.

The first, and arguably most abstract, description will be our starting point. The idea is to use the existence of the ``one-point'' space $1$, the terminal object, as the archetypal model of what it means to be a point. In $\Top$, we of course find that $1=\{\ast\}$ is a singleton set, equipped with its unique discrete topology. The frame of opens of this space is therefore: $\Opens 1 = \{\varnothing \subseteq \{\ast\}  \}$. This two-element frame, which is a Boolean algebra, is also called the \emph{frame of truth values}. 

\begin{lemma}
	The frame of truth values $\{0<1\}$ is an initial object in $\Frm$.
\end{lemma}
\begin{proof}
	If $L$ is a frame, then there is a map of frames $\{0<1\}\to L$ that sends $0$ and $1$ to the bottom and top elements of $L$, respectively. Since maps of frames preserve top and bottom elements, this map is unique.
\end{proof}

Dually, the locale $1$ whose frame of opens is $\Opens 1 :=\{0<1\}$ is a terminal object in $\Loc$. 

\begin{definition}\label{definition:point via map}
	A \emph{point} in a locale $X$ is a map of locales $1\to X$.
\end{definition}

Recall the definition of a \emph{completely prime filter} $\calF\subseteq \Opens X$ from \cref{section:adjunction top and loc}.

\begin{lemma}
	There is a bijective correspondence between points $1\to X$ of a locale, and completely prime filters in $\Opens X$.
\end{lemma}
\begin{proof}
	If $p\colon 1\to X$ is a map of locales, then $\mathcal{F}:=\{U\in \Opens X: p^{-1}(U)=1\}$ is a completely prime filter that fully determines the frame map $p^{-1}$.
\end{proof}

We also saw in \cref{lemma:prime elements are cpf} that completely prime filters $\calF\subseteq\Opens X$ correspond bijectively to prime elements $P\in\Opens X$. In terms of maps $p\colon 1\to X$, the prime element $P_\calF$ corresponds to the open $\bigvee\{U\in\Opens X: p^{-1}(U)=0\}$.


A \emph{Boolean locale} is a locale $X$ for which the underlying frame $\Opens X$ is a Boolean algebra (\cref{definition:boolean algebra}) with respect to its canonical Heyting algebra structure.

Recall the definition of an \emph{atom} from \cite[\S 5.2]{davey2002IntroductionLatticesOrder}.

\begin{definition}\label{definition:atom}
	An \emph{atom} in a distributive lattice $L$ with bottom element is a non-empty $\bot\neq a\in L$ such that $x\sqleq a$ implies $x=\bot$ or $x=a$. 
\end{definition}

\begin{proposition}[{\cite[\S II.5.4]{picado2012FramesLocalesTopology}}]\label{proposition:points in boolean locale are atoms}
	The points of a Boolean locale $X$ are precisely the atoms in $\Opens X$.
\end{proposition}

\section{Sublocales}\label{section:sublocales}
The material in this section is from \cite[Chapter~IX]{maclane1994SheavesGeometryLogic} and \cite[Chapter~III]{picado2012FramesLocalesTopology}.

The notion of a sublocale is analogous to the notion of subspace of a topological space. Category theoretically, subobjects are modelled by (equivalence classes of) monomorphisms. Recall that in $\Top$, the monomorphisms and epimorphisms are precisely the injective and surjective continuous functions, respectively. However, injective continuous functions $j\colon A\hookrightarrow S$ between topological spaces do not quite match with what we usually think of as subspaces. For instance, $A$ is not necessarily homeomorphic to its image under $j$ in $S$.

Rather, we want $A$ to carry the coarsest topology for which $j$ is continuous. Hence $A$ must carry the \emph{initial topology} induced by $j$, and its open subsets are precisely those of the form $j^{-1}(U)$, where $U\in \Opens S$. In the case that $A$ is literally a subset of $S$, we see that the induced topology on $A$ is just the collection of sets $\{A\cap U:U\in \Opens S\}$ (\cref{example:subset topology}). Note that a function $g\colon Z\to A$ is then continuous if and only if $j\circ g\colon Z\to S$ is continuous. 

\begin{definition}\label{definition:extremal monomorphism}
	An \emph{extremal monomorphism} is a monomorphism $m$ with the following property: if $e$ is an epimorphism and $m=n\circ e$, then $e$ is an isomorphism. We denote extremal monomorphisms by the arrows $\rightarrowtail$.
	
	Dually, an \emph{extremal epimorphism} is an epimorphism $e$ with the property that: if $m$ is a monomorphism and $e=m\circ f$, then $m$ is an isomorphism.
\end{definition}

\begin{lemma}
	The extremal monomorphisms in $\Top$ are precisely the subspaces.
\end{lemma}
\begin{proof}
	First, suppose $m\colon A\hookrightarrow S$ is a subspace inclusion. Take now a factorisation $m = n\circ e$ where $e\colon A\tworightarrow E$ is an epimorphism, i.e.~a surjective continuous function. By \cref{lemma:cancellation property epi mono} we get that $e$ is a monomorphism, i.e.~an injection. Hence $e$ is a bijection, so to show it is a homeomorphism it suffices to show it is open. But since $m$ is a subspace inclusion, the opens of the domain $A$ are of the form $m^{-1}(U)$ for $U\in\Opens S$, and so $e(m^{-1}(U)) = n^{-1}(U)$, which is indeed open.
	
	Conversely, suppose that $m\colon A\rightarrowtail S$ is an extremal monomorphism. Then we can endow the subset $m(A)\subseteq S$ with the subset topology, and this makes the functions $e\colon A\tworightarrow m(A)$ and $n\colon m(A)\hookrightarrow S$ continuous. Thus $e$ is an isomorphism, which proves $A$ has the initial topology generated by $m$. 
\end{proof}

Next we classify the epimorphisms and extremal monomorphisms in $\Loc$.

\begin{lemma}[{\cite[Lemma~III.1.1.1]{picado2012FramesLocalesTopology}}]
	\label{lemma:epimorphisms in Loc}
	The epimorphisms in $\Loc$ are exactly the maps of locales $f\colon X\to Y$ for which the map of frames ${f^{-1}\colon \Opens Y\to \Opens X}$ is injective.
\end{lemma}

\begin{lemma}[{\cite[Corollary~III.1.2]{picado2012FramesLocalesTopology}}]\label{proposition:extremal monomorphisms in Loc}
	The extremal monomorphisms in $\Loc$ are precisely those arrows $f\colon X\to Y$ for which the underlying map of frames $f^{-1}\colon \Opens Y\to \Opens X$ is a surjection.
\end{lemma}
%

\begin{remark}
	We should be careful to distinguish sub\emph{locales} from the notion of a sub\emph{frame}. The former are the appropriate subobjects in the category $\Loc$, while the latter are subobjects in the category $\Frm$. 
\end{remark}

Since the extremal monomorphisms determine the correct notion of subspace in $\Top$, we can now take extremal monomorphisms in $\Loc$ to represent the correct notion of subspaces of locales.

\begin{definition}
	A \emph{sublocale} of $X\in\Loc$ is an equivalence class of extremal monomorphisms $j\colon A\rightarrowtail X$, up to subobject inclusion (see \cref{definition:M-subobject} for details).
\end{definition}

Intuitively, therefore, a sublocale $j\colon A\rightarrowtail X$ is a locale $X$ whose frame of opens $\Opens A$ is fully determined as the image of the surjective frame map $j^{-1}$.

There are several alternative, equivalent ways to present sublocales more concretely. These are worked out in detail in \cite[\S III.2]{picado2012FramesLocalesTopology}. 

\begin{definition}
	A \emph{sublocale set} of a locale $X$ is a subset $S\subseteq \Opens X$ that is closed under all meets, and if $V\in S$ then $U\to V\in S$ for every $U\in\Opens X$.
\end{definition}

Sublocale sets $S$ of $X$ are precisely the subsets of $\Opens X$ such that $S$ is a frame with the induced order, and such that the inclusion function $S\hookrightarrow \Opens X$ defines (the direct image map of) a map of locales. Any extremal monomorphism $j\colon A\rightarrowtail X$ induces a unique sublocale set $j_\ast(\Opens A)\subseteq \Opens X$.

\begin{definition}\label{definition:coframe of sublocales}
	Denote by $\Sl(X)$ the collection of sublocale sets of $X$, ordered by subset inclusion. This defines a \emph{coframe}, i.e.~it is the opposite of a frame \cite[Theorem~III3.2.1]{picado2012FramesLocalesTopology}.
\end{definition}

\begin{definition}\label{definition:open and closed sublocales}
	Let $X$ be a locale, and fix $U\in\Opens X$. The \emph{open} and \emph{closed sublocales} induced by $U$ are described by the following sublocale sets, respectively:
		\[
			\mathfrak{o}(U):= \{U\to V:V\in\Opens X\}
			\qquad\text{and}\qquad
			\mathfrak{c}(U):= \{V\in\Opens X:U\sqleq V\}.
		\]
\end{definition}

\subsection{Double-negation sublocale}\label{section:double-negation sublocale}
For any locale $X$, the double Heyting negation operator $\Opens X\to \Opens X:U\mapsto \neg\neg U$ determines a finite meet-preserving monad. It can be shown that the fixed-points of this monad, i.e.~those opens with $\neg\neg U = U$, form a sublocale set of $\Opens X$. In terms of extremal monomorphisms, we get the map $i\colon X_{\neg\neg}\rightarrowtail X$ defined by $i^{-1}\colon U\mapsto \neg\neg U$ and where $i_\ast$ is the inclusion map $\Opens X_{\neg\neg}\subseteq \Opens X$. The sublocale $X_{\neg\neg}$ is called the \emph{double-negation sublocale} of $X$. It is quite special: it is the smallest dense Boolean sublocale.

\begin{lemma}\label{lemma:double-negation sublocale is boolean}
	The double-negation sublocale $X_{\neg\neg}$ is a Boolean locale.
\end{lemma}
\begin{proof}
	This follows simply by \cref{lemma:properties of heyting negation}(d); cf.~\cite[\S III.10.4]{picado2012FramesLocalesTopology}.
\end{proof}

\begin{definition}\label{definition:dense sublocales}
	A sublocale (set) $S\in\Sl(X)$ is called \emph{dense} if $\varnothing \in S$.
\end{definition}

Note $X_{\neg\neg}$ is dense in $X$, since $\neg\neg \varnothing = \neg X = \varnothing$. \emph{Isbell's density theorem} says that every dense sublocale must contain the double-negation sublocale. 

\begin{theorem}[{\cite[\S III.8.3]{picado2012FramesLocalesTopology}}]
	Every locale $X$ contains a smallest dense sublocale $X_{\neg\neg}$.
\end{theorem}

\section{Product locales}\label{section:product locales}
Since $\Loc = \Frm^\op$, the product $X\times Y$ of locales is determined by the \emph{co}product of their respective frame of opens:
	\[
		\Opens\mspace{1mu}(X\times Y) \cong \Opens X\oplus \Opens Y.
	\]
Calculating the coproduct of frames is a somewhat more involved procedure. For details, we refer to \cite[Chapter~IV]{picado2012FramesLocalesTopology}. For binary coproducts, the most intuitive way to calculate the coproduct $\Opens X\oplus \Opens Y$ is perhaps following \cite{dowker1977sumsCategoryFrames}. Namely, $\Opens X\oplus \Opens Y$ is the frame freely generated by the elements $U\oplus V$, for $U\in\Opens X$ and $V\in\Opens Y$, identified according to the following equations:
	\[
		\bigwedge_{i=1}^n (U_i\oplus V_i) = \left(\bigwedge_{i=1}^n U_i\right)\oplus \left(\bigwedge_{i=1}^n V_i\right);
	\]
	\[
		\bigvee_{i\in I} (U\oplus V_i) = U\oplus \left(\bigvee_{i\in I} V_i\right)
		\quad\text{and}\quad
		\bigvee_{i\in I} (U_i\oplus V) = \left(\bigvee_{i\in I} U_i\right)\oplus V.
	\]	
The latter two equations imply that $\varnothing \oplus V =\varnothing = U\oplus \varnothing$. In fact, $U\oplus V = \varnothing$ if and only if either $U$ or $V$ is empty. The intuition is that $U\oplus V$ is the frame-theoretic analogue of the basic open boxes $U\times V$ in the product topology on~$S\times S$.

\begin{lemma}[{\cite[Proposition~IV.5.2(4)]{picado2012FramesLocalesTopology}}]
	\label{lemma:inclusions in coproducts of frames}
	If $\oplus_{i\in I} x_i\neq\varnothing$ in the coproduct $\bigoplus_{i\in I} L_i$ of frames, then
	\[
	\oplus_{i\in I}x_i \sqleq \oplus_{i\in I} y_i
	\qquad\text{if and only if}\qquad
	\forall i\in I: x_i\sqleq y_i.
	\]
\end{lemma}

\begin{lemma}[{\cite[\S IV.5.5.1]{picado2012FramesLocalesTopology}
	}]
	\label{lemma:pair map preimage}
	If $f\colon X\to Y$ and $g\colon X\to Z$ are maps of locales, then the unique map of locales $(f,g) \colon X\to Y\times Z$ with $\pr_1\circ (f,g) = f$ and $\pr_2\circ (f,g) = g$ is defined by the preimage map
		\[
			(f,g)^{-1}\colon \Opens Y\oplus \Opens Z\longrightarrow \Opens X;
			\qquad
			D\longmapsto 
			\bigvee\left\{f^{-1}(U)\wedge g^{-1}(V):U\oplus V\sqleq D\right\}.
		\]
\end{lemma}
\begin{corollary}\label{corollary:pair map preimage on basic opens}
	On basic opens $V\oplus W\in \Opens Y\oplus \Opens Z$ we have:
		\[
			(f,g)^{-1}\left(V\oplus W\right)= f^{-1}(V)\wedge g^{-1}(W).
		\]
\end{corollary}
\begin{proof}
	Clearly $V\oplus W\sqleq V\oplus W$, which by the formula in \cref{lemma:pair map preimage} gives $f^{-1}(V)\wedge g^{-1}(W)\sqleq (f,g)^{-1}(V\oplus W)$. For the converse, let $A\oplus B\sqleq V\oplus W$ be arbitrary. If $A\oplus B= \varnothing$, then either $A$ or $B$ is empty, and since both $f^{-1}$ and $g^{-1}$ preserve the bottom element, this implies $f^{-1}(A)\wedge g^{-1}(B)=\varnothing$. Thus $(f,g)^{-1}$ can be calculated taking the join only over non-empty $A\oplus B$. In that case, by \cref{lemma:inclusions in coproducts of frames} we get $A\sqleq V$ and $B\sqleq W$, so monotonicity of $f^{-1}$ and $g^{-1}$ give $f^{-1}(A)\wedge g^{-1}(B)\sqleq f^{-1}(V)\wedge g^{-1}(W)$. This gives the desired inclusion ${(f,g)^{-1}(V\oplus W) \sqleq f^{-1}(V)\wedge g^{-1}(W)}$.
\end{proof}


\begin{proposition}\label{proposition:pair map right adjoint}
	If $P\in \Opens X$ is a prime element, then:
		\[
			(f,g)_\ast(P) = \left(f_\ast(P)\oplus Z\right)\vee \left(Y\oplus g_\ast(P)\right).
		\]
\end{proposition}
\begin{proof}
	First, by \cref{lemma:locale map right adjoint using basis,corollary:pair map preimage on basic opens} we can write
		\[
			(f,g)_\ast(P)
			=
			\bigvee\left\{U\oplus V : f^{-1}(U)\wedge g^{-1}(V)\sqleq P\right\}.
		\]
	But since $P$ is prime, this condition is satisfied if and only if either $f^{-1}(U)\sqleq P$ or $g^{-1}(V)\sqleq P$. Hence the join splits into two:
		\begin{align*}
			(f,g)_\ast(P)
			&=
			\bigvee\left\{ U\oplus Z : f^{-1}(U)\sqleq P\right\}
			\vee
			\bigvee\left\{Y\oplus V: g^{-1}(V)\sqleq P\right\}
			\\&=
			\left(\bigvee\left\{ U : f^{-1}(U)\sqleq P\right\}\right)\oplus Z
			\vee
			Y\oplus\left(\bigvee\left\{V: g^{-1}(V)\sqleq P\right\}\right),
		\end{align*}
	which is equal to the desired expression by \cref{lemma:locale map right adjoint}.
\end{proof}

Since $\loc\dashv \pt$, the functor $\pt\colon \Loc\to \Top$ preserves products, and hence we see
\[
\pt(X\times Y) \cong \pt(X)\times \pt(Y).
\]
See also \cite{picado2015NotesProductLocales} for a way to explicitly showcase this homeomorphism. For our discussion in \cref{section:causal boundaries}, the use of prime elements is more favourable, so we restate this isomorphism here.

\begin{proposition}
	For locales $X$ and $Y$, there is a homeomorphism
	\[
		\pt'(X)\times \pt'(Y)\longrightarrow \pt'(X\times Y);
		\qquad
		(P,Q)\longmapsto (P\oplus Y)\vee (X\oplus Q)
	\]
	with inverse $D\longmapsto \left((\pr_1)_\ast (D), (\pr_2)_\ast (D) \right)$.
\end{proposition}


\chapter{Internal preorders}
\label{section:internal preorders}
In category theory, \emph{subobjects} of a given object $X\in\cat{D}$ provide an abstract, internal analogue of the various types of subobjects one encounters in mathematics, such as subsets, subgroups, subspaces, subalgebras, sublattices, etc. Subobjects are represented by monomorphisms $A\hookrightarrow X$. Given that the categorical structure of $\cat{D}$ is rich enough, it is possible to use these ideas to define an internal notion of a \emph{preorder} on ${X\in\cat{D}}$. Namely, reflexivity and transitivity can be stated solely as categorical properties of a subobject $R\hookrightarrow X\times X$, thought of as the graph of the relation.

In some cases the monomorphisms of $\cat{D}$ do not quite reproduce the desired notion of subobject. For instance, in the category of topological spaces and continuous functions $\Top$, the subobjects of $X$ correspond merely to sub\emph{sets} of $X$, equipped with an arbitrary topology, as long as the inclusion map is continuous. Conventionally, the accepted notion of subspace is that of a subset equipped with the \emph{smallest} topology that makes the inclusion map continuous. A similar problem holds in the category $\Loc$ of locales. Fortunately, it is possible to recover the correct notion of subobject by restricting the class of monomorphisms to something more specific. For instance, in $\Top$ and $\Loc$ this can be achieved by taking only the \emph{extremal subobjects} \cite[Chapter~III]{picado2012FramesLocalesTopology}.

In the abstract setting this can be treated by specifying a \emph{factorisation system} on $\cat{D}$: a pair $(\calE,\calM)$ of classes of morphisms that tell us how maps factorise via their image. An $\calM$-subobject of $X$ is then represented by a morphism $a\colon A\rightarrowtail X$ with $a\in\calM$. This leads to a notion of $\calM$-preorders $R\rightarrowtail X\times X$ on $X$, and we shall construct a new category $\ord(\cat{D})$ whose objects $(X,R)$ consist of $X\in\cat{D}$ equipped with a $\calM$-preorder $R$, together with maps $f\colon (X,R)\to (Y,Q)$ that satisfy an appropriate monotonicity condition. 

The point of this section is to lift an adjunction $\begin{tikzcd}[ampersand replacement=\&]
	{\cat{C}} \& {\cat{D}}
	\arrow[""{name=0, anchor=center, inner sep=0}, "F", shift left=2, from=1-1, to=1-2]
	\arrow[""{name=1, anchor=center, inner sep=0}, "G", shift left=2, from=1-2, to=1-1]
	\arrow["\dashv"{anchor=center, rotate=-90}, draw=none, from=0, to=1]
\end{tikzcd}$ to their induced categories of internally preordered objects and monotone morphisms:
	\[
		\begin{tikzcd}[column sep=large]
			\OrdC & \OrdD \\
			{\cat{C}} & {\cat{D}.}
			\arrow[hook, from=2-1, to=1-1]
			\arrow[hook, from=2-2, to=1-2]
			\arrow[""{name=0, anchor=center, inner sep=0}, "G", shift left=2, from=2-2, to=2-1]
			\arrow[""{name=1, anchor=center, inner sep=0}, "{\overline{F}}", shift left=2, from=1-1, to=1-2]
			\arrow[""{name=2, anchor=center, inner sep=0}, "{\overline{G}}", shift left=2, from=1-2, to=1-1]
			\arrow[""{name=3, anchor=center, inner sep=0}, "F", shift left=2, from=2-1, to=2-2]
			\arrow["\dashv"{anchor=center, rotate=-90}, draw=none, from=1, to=2]
			\arrow["\dashv"{anchor=center, rotate=-90}, draw=none, from=3, to=0]
		\end{tikzcd}
	\]
Throughout this section, we consider the data of the adjunction $F\dashv G$ as given, together with the factorisation systems $(\calE_\cat{C},\calM_{\cat{C}})$ on $\cat{C}$ and $(\calE_\cat{D},\calM_{\cat{D}})$ on $\cat{D}$. We sometimes omit the subscripts on $\calE$ and $\calM$, when context is clear. As we progress, we shall be adding assumptions on this data, as and where needed. All results hold in the context of these assumptions.

\section{Orthogonal factorisation systems}
First we recall some definitions and results from the literature, in particular from \cite[\S 4]{borceux1994HandbookCategoricalAlgebra}, \cite[Chapter~IX]{herrlich2007CategoryTheoryIntroduction} and \cite{riehl2008FactorizationSystems}. The idea of a factorisation system $(\calE,\calM)$ is that it provides an abstract, canonical image factorisation of arrows
	\[
		\begin{tikzcd}[ampersand replacement=\&]
			A \& {\im(f)} \& {B.}
			\arrow["\in\calE"' {rotate = -30,pos=.6}, two heads, from=1-1, to=1-2]
			\arrow["f", curve={height=-24pt}, from=1-1, to=1-3]
			\arrow["\in\calM"' {rotate = -30,pos=.7}, tail, from=1-2, to=1-3]
		\end{tikzcd}
	\]
This is made precise as follows. We adopt \cite[Definition~1.9]{riehl2008FactorizationSystems}.

\begin{definition}\label{definition:factorisation system}
	Let $\cat{D}$ be a category. An \emph{orthogonal factorisation system} on~$\cat{D}$ is a pair $(\calE,\calM)$ of classes of morphisms such that:
		\begin{enumerate}[label = (\roman*)]
			\item $\calE$ and $\calM$ are closed under composition;
			\item every arrow $f$ in $\cat{D}$ has a factorisation $f= m\circ e$ with $e\in \calE$ and $m\in \calM$;
			\item the factorisation is unique in the sense that if the outside square commutes
				\[
					\begin{tikzcd}[ampersand replacement=\&]
						A \& C \\
						B \& D
						\arrow[from=1-1, to=1-2]
						\arrow["\calE\ni"', two heads, from=1-1, to=2-1]
						\arrow["\in\calM", tail, from=1-2, to=2-2]
						\arrow["{\exists!}"{description}, dashed, from=2-1, to=1-2]
						\arrow[from=2-1, to=2-2]
					\end{tikzcd}
				\]
			then the dashed arrow exists and is unique completing the diagram.
		\end{enumerate}
	The arrows in the \emph{left class} $\calE$ will be written as $\tworightarrow$, while the arrows in the \emph{right class} $\calM$ will be written as $\rightarrowtail$.
\end{definition}

\begin{example}\label{example:epi-mono factorisation}
	On the category $\Set$ of sets and functions we have the familiar \emph{epi-mono} factorisation system, where $\calE$ and $\calM$ are the classes of surjective and injective functions, respectively (\cite[Example~1.10]{riehl2008FactorizationSystems}). See \cite[\S 4.4.4]{borceux1994HandbookCategoricalAlgebra} for more details on epi-mono factorisations.
\end{example}

The following two results are from \cite[Lemma~1.13]{riehl2008FactorizationSystems}.

\begin{lemma}\label{lemma:factorisation classes contain isos}
	For any orthogonal factorisation system $(\calE,\calM)$ the intersection $\calE\cap \calM$ is precisely the class of isomorphisms.
\end{lemma}

\begin{lemma}\label{lemma:cancellation property factorisation system}
	If $(\calE,\calM)$ is an orthogonal factorisation system on $\cat{D}$ then we have the following \emph{cancellation properties:}
	\begin{align*}
		\text{$g\circ f\in \calE$ and $f\in\calE$}
		\quad&\text{implies}\quad g\in\calE;\\
		\text{$g\circ f\in\calM$ and $g\in \calM$}\quad&\text{implies}\quad f\in\calM.
	\end{align*}
\end{lemma}

\begin{remark}
	Note that axiom~(iii) indeed implies that factorisations are unique up to isomorphism:
	\[
	\begin{tikzcd}[ampersand replacement=\&]
		A \& Y \\
		X \& {B.}
		\arrow[two heads, from=1-1, to=1-2]
		\arrow[two heads, from=1-1, to=2-1]
		\arrow["{\exists!}", shift left, dashed, from=1-2, to=2-1]
		\arrow[tail, from=1-2, to=2-2]
		\arrow["{\exists!}", shift left, dashed, from=2-1, to=1-2]
		\arrow[tail, from=2-1, to=2-2]
	\end{tikzcd}
	\]
	\cref{lemma:cancellation property factorisation system} implies that both dashed arrows are in $\calE\cap \calM$, so by \cref{lemma:factorisation classes contain isos} they are isomorphisms. Uniqueness implies they are mutually inverse.
\end{remark}

Orthogonal factorisation systems also have the pleasant property that $\calM$ is closed under pullbacks. This follows from e.g.~\cite[p.~11]{riehl2008FactorizationSystems}, or from the more explicit proof in \cite[Proposition~3.1]{nlab2024weakFactorizationSystem}.

\begin{lemma}\label{lemma:factorisation right class closed under pullback}
	Let $(a_i\colon A_i\rightarrowtail X)_{i\in I}$ be a family of arrows in $\calM$ admitting a wide pullback $P$ with cone $(n_i\colon P\to A_i)_{i\in I}$. Then $n_i\in \calM$ for all $i\in I$.
\end{lemma}

\section{Internal preorders for factorisation systems}
\begin{definition}\label{definition:M-subobject}
	An \emph{$\calM$-subobject} of $X\in\cat{D}$ is an arrow $a\colon A\rightarrowtail X$ in $\calM$. We denote the set of $\calM$-subobjects of $X$ by $\MSub(X)$. This becomes a preorder where $A\subseteq B$ whenever there is a factorisation:
		\[
			\begin{tikzcd}
				B & X \\
				{A.}
				\arrow["a"', tail, from=2-1, to=1-2,shorten <=-4pt]
				\arrow["b", tail, from=1-1, to=1-2]
				\arrow["m", dashed, from=2-1, to=1-1]
			\end{tikzcd}
		\]
	Note that by the cancellation property in \cref{lemma:cancellation property factorisation system} the arrow $m$ is automatically in $\calM$. In the case that $A\subseteq B$ and $B\subseteq A$, the arrow $m$ is necessarily the unique isomorphism that factorises $a$ through $b$. Quotienting by this equivalence relation, $\MSub(X)$ becomes a partial order. Doing so leads to the pleasant result in \cref{lemma:M-complete}.
\end{definition}

\begin{example}
	For the epi-mono factorisation system on $\Set$ from \cref{example:epi-mono factorisation}, an $\calM$-subobject of a set $X$ is just a monomorphism $A\hookrightarrow X$, which can in turn be identified with a subset of $X$ (after taking equivalence classes of subobjects).
	
	For $\Top$ and $\Loc$ we shall pick $\calM$ to be the \emph{extremal monomorphisms}. An $\calM$-subobject of a topological space $X$ is a subset $A\subseteq X$ equipped with the subset topology, i.e.~the smallest topology making the inclusion function $A\hookrightarrow X$ continuous. An $\calM$-subobject of a locale $X$ corresponds precisely to a \emph{sublocale} (\cref{section:sublocales}).
\end{example}

At this point we introduce our first assumptions on $F\dashv G$. 
	\begin{enumerate}[label = (A\arabic*)]
		\item\label{assumption:well-powered} $\cat{C}$ and $\cat{D}$ are \emph{$\calM$-well-powered}, meaning that each $\MSub(X)$ is a set;
		\item\label{assumption:complete} $\cat{C}$ and $\cat{D}$ have all small limits. This allows us to define relations (which needs products) and their composition (which needs pullbacks), and also gives the next pleasant result.
	\end{enumerate}

\begin{lemma}\label{lemma:M-complete}
	If $\cat{D}$ satisfies \ref{assumption:well-powered}-\ref{assumption:complete}, then $\MSub(X)$ is a complete lattice.
\end{lemma}
\begin{proof}
	Here we show explicitly how for a family $(A_i)_{i\in I}$ in $\MSub(X)$ the meet $\bigwedge_{i\in I}A_i$ is constructed. First, if the index $I$ is empty then the meet $\bigwedge\varnothing$ is just $X$ itself, represented by the identity map $\id_X\in \calM$ (\cref{lemma:factorisation classes contain isos}). If $I$ is non-empty, use the existence of all limits to construct the \emph{wide pullback}:
		\[
			\begin{tikzcd}[ampersand replacement=\&,column sep=1ex, row sep=3ex]
				\& {\bigwedge_{i\in I}A_i} \\
				{A_i} \&\& {A_j} \\
				\& {X,}
				\arrow["{n_i}"', from=1-2, to=2-1]
				\arrow["{n_j}", from=1-2, to=2-3]
				\arrow["a", from=1-2, to=3-2]
				\arrow["{a_i}"', tail, from=2-1, to=3-2]
				\arrow["{a_j}", tail, from=2-3, to=3-2]
			\end{tikzcd}
		\]
	where by  construction we get the map $a$ equalling $a_i\circ n_i$ for every $i\in I$. From \cref{lemma:factorisation right class closed under pullback} we get $n_i\in\calM$, and so $a\in\calM$ by \cref{lemma:cancellation property factorisation system}. That this arrow defines (i.e.~represents) the meet of $(A_i)_{i\in I}$ in $\MSub(X)$ follows by the universal property of the wide pullback.
\end{proof}

Since we assume $\cat{D}$ is complete~\ref{assumption:complete}, we get in particular binary products ${X\times Y}$, the subobjects of which we can think about as abstract graphs of relations.

\begin{definition}
	An \emph{$\calM$-relation} from $X$ to $Y$ is an $\calM$-subobject $r\colon R\rightarrowtail X\times Y$. We write $r=(s,t)$, where $s,t\colon R\to X$ are the \emph{source} and \emph{target maps}. They are sometimes also denoted $r=(r_1,r_2)$, especially when multiple relations are involved.
\end{definition}

We always obtain the largest $\calM$-relation $\id_{X\times Y}$ from $X$ to $Y$. Similarly, we want to allow that:
\begin{enumerate}[label = (A\arabic*)]\setcounter{enumi}{2}
	\item\label{assumption:diagonal in M}
	the canonical diagonal maps $X\to X\times X$ are in $\calM$ for every object $X\in\cat{D}$.
\end{enumerate}
The resulting $\calM$-relation on $X$ is denoted by $\Delta_X\in\MSub(X\times X)$, or sometimes just $\Delta$ when it is clear what object we are working with.

\begin{definition}\label{definition:internal monotone maps}
	Consider pairs $(X,R)$ and $(Y,Q)$, where $X,Y\in\cat{D}$ are objects and $R$ and $Q$ are $\calM$-relations on $X$ and $Y$, respectively. We say that an arrow $f\colon X\to Y$ is \emph{(internally) monotone} if there exists an arrow completing the following diagram:
	\[
	\begin{tikzcd}[ampersand replacement=\&]
		R \& Q \\
		{X\times X} \& {Y\times Y.}
		\arrow["{\exists \overline{f}}"', dashed, from=1-1, to=1-2]
		\arrow["r"', tail, from=1-1, to=2-1]
		\arrow["q", tail, from=1-2, to=2-2]
		\arrow["{f\times f}"', from=2-1, to=2-2]
	\end{tikzcd}
	\]
\end{definition}

\begin{remark}
	For our applications to spaces and locales, the arrows in the class $\calM$ will all be monomorphisms, in which case it follows that the map $\overline{f}$ is unique.
\end{remark}

\begin{lemma}\label{lemma:monotone order modification}
	If $f\colon (X,R)\to (Y,Q)$ is monotone, $P\subseteq R$ and $Q\subseteq T$, then $f\colon (X,P)\to (Y,T)$ is monotone.
\end{lemma}
\begin{proof}
	This follows immediately from the following diagram:
	\[
	\begin{tikzcd}
		P & R & Q & T \\
		& {X\times X} & {Y\times Y.}
		\arrow[dashed, tail, from=1-3, to=1-4]
		\arrow["t", tail, from=1-4, to=2-3]
		\arrow["q"', tail, from=1-3, to=2-3]
		\arrow["{\overline{f}}", dashed, from=1-2, to=1-3]
		\arrow["r"', tail, from=1-2, to=2-2]
		\arrow["{f\times f}"', from=2-2, to=2-3]
		\arrow["p"', tail, from=1-1, to=2-2]
		\arrow[dashed, tail, from=1-1, to=1-2]
	\end{tikzcd}
	\qedhere
	\]
\end{proof}

\begin{example}
	In $\Set$ with the epi-mono factorisation system, a function $f$ is internally monotone if and only if it is monotone in the usual sense:
		\[
			xRy \qquad\text{implies}\qquad f(x)Qf(y).
		\]
	Indeed, using this we can define the function $\overline{f}(x,y):= (f(x),f(y))$, which is clearly the unique one completing the square above.
\end{example}

\begin{construction}
	Consider two $\calM$-relations $r\colon R\to X\times Y$ and ${q\colon Q\to Y\times Z}$. We define their \emph{composition} $R\circ Q$, a new $\calM$-relation from $X$ to $Z$, as follows. First we take the pullback
		\[
			\begin{tikzcd}
				{R\times_Y^{r_2,q_1}Q} & Q \\
				R & {Y.}
				\arrow["{\pr_1}"', from=1-1, to=2-1]
				\arrow["{\pr_2}", from=1-1, to=1-2]
				\arrow["{q_1}", from=1-2, to=2-2]
				\arrow["{r_2}"', from=2-1, to=2-2]
				\arrow["\lrcorner"{anchor=center, pos=0.05}, draw=none, from=1-1, to=2-2]
			\end{tikzcd}
		\]
	The new $\calM$-relation then occurs as the $(\calE,\calM)$-factorisation:
		\[
			\begin{tikzcd}[column sep=small,row sep=scriptsize]
				{R\circ Q} && {X\times Z} \\
				{R\times_Y^{r_2,q_1}Q.}
				\arrow["{(r_1\circ \pr_1,q_2\circ\pr_2)}"', from=2-1, to=1-3]
				\arrow[two heads, from=2-1, to=1-1]
				\arrow[tail, from=1-1, to=1-3]
			\end{tikzcd}
		\]
\end{construction}

\begin{example}\label{example:internal relation composition in Set}
	As this construction appears quite abstract, it helps to unpack the definitions in the case ${\cat{D}=\Set}$. An $\calM$-relation between sets $X$ and $Y$ can then just be identified with a subset $R\subseteq X\times Y$, returning the ordinary notion of relation. For relations $R\subseteq X\times Y$ and ${Q\subseteq Y\times Z}$ the source and target maps are just the projections, so we get a canonical representation of the pullback
		\[
			R\times_Y^{r_2,q_1}Q
			\cong
			\{\left((x,y),(y',z)\right)\in R\times Q: y = y' \},
		\]
	and we find that the composition $R\circ Q$ can indeed be identified with the ordinary composition of relations:
		\[
			R\circ Q
			\cong
			\{(x,z)\in X\times Z:\exists y\colon (x,y)\in R, (y,z)\in Q\}.
		\]
\end{example}

\begin{remarknumbered}\label{remark:non-associative relational composition}
	Unlike in $\Set$, it is not always the case that the composition of $\calM$-relations defines an associative operation. To be precise, it was proved in \cite[Theorem~2.5]{klein1970relations} that composition of relations relative to a factorisation system $(\calE,\calM)$ is associative if and only if the arrows in the left class $\calE$ are stable under pullback. In the case of locales, $\calE$ will be the class of epimorphisms. Due to a result of Banaschewski~\cite{banaschewski1990PushingOutFrames} we know that epimorphisms are generally not stable under pullbacks in $\Loc$, so unfortunately the relational composition of sublocale relations will not be associative in general. Roughly, the class of locales where epimorphisms are stable under pullback are those where the dual of the coframe of sublocales is Boolean, a restriction that in topological language corresponds to being \emph{scattered}. We therefore have to account for the more general case where composition of $\calM$-relations is not associative. Nevertheless, as the next two lemmas show, the familiar calculus of relations goes through in the abstract setting.
\end{remarknumbered}

\begin{lemma}\label{lemma:relation composition respects inclusion}
	Consider relations $R\subseteq S$ in $\MSub(X\times Y)$ and $Q\subseteq T$ in $\MSub(Y\times Z)$. Then $R\circ Q\subseteq S\circ T$ in $\MSub(X\times Z)$. 
\end{lemma}
\begin{proof}
	We have $\calM$-subobject inclusion arrows $n\colon R\rightarrowtail S$ and $m\colon Q\rightarrowtail T$, which allows us to construct the following diagram:
	\[
	\begin{tikzcd}[column sep=5ex,row sep=2ex]
		{R\times_Y^{r_2,q_1}Q} \\
		& {S\times_Y^{s_2,t_1}T} & {T} \\
		& {S} & Y.
		\arrow["{t_1}", from=2-3, to=3-3]
		\arrow["{s_2}"', from=3-2, to=3-3]
		\arrow["{\pr_1}"', from=2-2, to=3-2]
		\arrow["{\pr_2}", from=2-2, to=2-3]
		\arrow["{n\circ\pr_1}"', curve={height=18pt}, from=1-1, to=3-2,to path={ --([yshift=-3ex,xshift=0ex]\tikztostart.south) [pos=2]\tikztonodes |- (\tikztotarget) },rounded corners=12pt]
		\arrow["{m\circ \pr_2}", curve={height=-18pt}, from=1-1, to=2-3,rounded corners=12pt,to path={ --([yshift=0ex,xshift=3ex]\tikztostart.east) [pos=4]\tikztonodes  -|  (\tikztotarget) }]
		\arrow["{\exists! \alpha}"', dashed, from=1-1, to=2-2]
		\arrow["\lrcorner"{anchor=center, pos=0.05}, draw=none, from=2-2, to=3-3]
	\end{tikzcd}
	\]
	The newly found arrow $\alpha$ now allows us to employ the orthogonal factorisation system, facilitating the desired inclusion:
	\[
	\begin{tikzcd}[column sep=0cm,cramped]
		& {S\times_Y^{s_2,t_1}T} \\
		{R\times_Y^{r_2,q_1}Q} && {S\circ T} \\
		{R\circ Q} && {Y\times Y.}
		\arrow["\alpha", from=2-1, to=1-2]
		\arrow[two heads, from=1-2, to=2-3]
		\arrow[tail, from=2-3, to=3-3]
		\arrow[tail, from=3-1, to=3-3]
		\arrow[two heads, from=2-1, to=3-1]
		\arrow["{\exists!}"{description}, dashed, from=3-1, to=2-3]
	\end{tikzcd}
	\qedhere
	\]
\end{proof}

\begin{lemma}\label{lemma:diagonal relations are identities}
	For any relation $R\in\MSub(X\times Y)$, we have isomorphisms of subobjects $\Delta_X\circ R\cong R\cong R\circ \Delta_Y$.
\end{lemma}
\begin{proof}
	By construction of the pullback, the following outside square commutes:
	\[
	\begin{tikzcd}
		{X\times_X^{\id,r_1}R} & R \\
		{\Delta_X\circ R} & {X\times Y,}
		\arrow[two heads, from=1-1, to=2-1]
		\arrow[tail, from=1-2, to=2-2]
		\arrow["{\pr_2}", from=1-1, to=1-2]
		\arrow[tail, from=2-1, to=2-2]
		\arrow["{\exists!}"{description}, dashed, from=2-1, to=1-2]
	\end{tikzcd}
	\]
	and hence the orthogonal factorisation system gives the inclusion $\Delta_X\circ R\subseteq R$. For the converse inclusion, consider the following diagram:
	\[
	\begin{tikzcd}[column sep=5ex,row sep=2ex]
		{R} \\
		& {X\times_X^{\id,r_1}R} & {R} \\
		& {X} & X.
		\arrow["{r_1}", from=2-3, to=3-3]
		\arrow[equal, from=3-2, to=3-3]
		\arrow["{\pr_1}"', from=2-2, to=3-2]
		\arrow["{\pr_2}", from=2-2, to=2-3]
		\arrow["{r_1}"', curve={height=18pt}, from=1-1, to=3-2,to path={ --([yshift=-3ex,xshift=0ex]\tikztostart.south) [pos=2]\tikztonodes |- (\tikztotarget) },rounded corners=12pt]
		\arrow[equal, curve={height=-18pt}, from=1-1, to=2-3,rounded corners=12pt,to path={ --([yshift=0ex,xshift=3ex]\tikztostart.east) [pos=4]\tikztonodes  -|  (\tikztotarget) }]
		\arrow["{\exists! \beta}"', dashed, from=1-1, to=2-2]
		\arrow["\lrcorner"{anchor=center, pos=0.05}, draw=none, from=2-2, to=3-3]
	\end{tikzcd}
	\]
	The outside obviously commutes, and the resulting arrow $\beta$ allows us to employ the orthogonal factorisation system, giving the desired inclusion:
	\[
	\begin{tikzcd}[column sep=.2cm,cramped, row sep=small]
		& {X\times_X^{\id,r_1}R} \\
		R && {\Delta_X\circ R} \\
		R && {X\times Y.}
		\arrow[equals,from=2-1, to=3-1]
		\arrow[tail, from=3-1, to=3-3]
		\arrow[tail, from=2-3, to=3-3]
		\arrow["\beta",from=2-1, to=1-2]
		\arrow[two heads, from=1-2, to=2-3]
		\arrow["{\exists!}"{description}, dashed, from=3-1, to=2-3]
	\end{tikzcd}
	\]
	Together, this gives $\Delta_X\circ R\cong R$. The isomorphism $R\circ \Delta_Y\cong R$ is proved analogously. They clearly preserve the subobject inclusions.
\end{proof}

\begin{definition}\label{definition:internal preorders}
	Fix an $\calM$-relation $R\in\MSub(X\times X)$ on $X$. We say $R$ is:
		\begin{enumerate}[label = (\roman*)]
			\item \emph{reflexive} if $\Delta\subseteq R$;
			\item \emph{transitive} if $R\circ R\subseteq R$.
		\end{enumerate}
	An \emph{$\calM$-preorder} on $X\in\cat{D}$ is a reflexive and transitive $\calM$-relation $R$ on $X$.
\end{definition}

\begin{example}
	It is again helpful to unpack these definitions in the setting of $\Set$ with its epi-mono factorisation system. The diagonal map $X\to X\times X$ of a set~$X$ is the function $x\mapsto (x,x)$, so that reflexivity of $R\in\MSub(X\times X)$ indeed just says that $(x,x)\in R$ for all $x\in R$. Similarly, since by \cref{example:internal relation composition in Set} the composition $R\circ R$ of internal relations is just the ordinary notion of relation composition, it is clear the internal notion of transitivity is equivalent to the ordinary one.
\end{example}

The following results elucidate the reflexivity and transitivity in terms of arrows in $\cat{D}$.

\begin{lemma}\label{lemma:reflexive in terms of factorisation}
	An $\calM$-relation $R$ is reflexive if and only if there exists an arow $\rho\colon X\to R$ such that $s\circ \rho = \id_X = t\circ \rho$. 
\end{lemma}
\begin{proof}
	By definition, $R$ is reflexive if there exists a factorisation
		\[
			\begin{tikzcd}
				R & {X\times X} \\
				X
				\arrow["r", tail, from=1-1, to=1-2]
				\arrow["\delta"', tail, from=2-1, to=1-2]
				\arrow["\rho", dashed, from=2-1, to=1-1]
			\end{tikzcd}
		\]
	where $\delta$ is the diagonal map of $X$. By composing with projection maps it is clear that $r\circ \rho = \delta$ if and only if $s\circ \rho =\id_X =t\circ \rho$. 
\end{proof}

\begin{corollary}\label{corollary:reflexive relations have split epi source}
	The source and target maps of a reflexive $\calM$-relation are split epimorphisms.
\end{corollary}

\begin{lemma}\label{lemma:transitive in terms of factorisation}
	An $\calM$-relation $R$ is transitive if and only if there is a factorisation
		\[
			\begin{tikzcd}
				R & {X\times X} \\
				{R\times_X^{t,s}R.}
				\arrow["r", tail, from=1-1, to=1-2]
				\arrow["\tau", dashed, from=2-1, to=1-1]
				\arrow["{(s\circ\pr_1,t\circ \pr_2)}"', from=2-1, to=1-2]
			\end{tikzcd}			
		\]
\end{lemma}
\begin{proof}
	The existence of such $\tau$ guarantees that $R\circ R\subseteq R$ by definition of an orthogonal factorisation system:
		\[
			\begin{tikzcd}
				R & {X\times X} \\
				{R\times_X^{t,s}R} & {R\circ R.}
				\arrow[two heads, from=2-1, to=2-2]
				\arrow["\tau", from=2-1, to=1-1]
				\arrow["r", tail, from=1-1, to=1-2]
				\arrow[tail, from=2-2, to=1-2]
				\arrow["{\exists!}"{description}, dashed, from=2-2, to=1-1]
			\end{tikzcd}
		\]
	For the converse, $\tau$ can be defined through the bottom left triangle.
\end{proof}


\begin{definition}\label{definition:category of internal preorders}
	We define the category $\ord(\cat{D})$ whose objects are pairs $(X,R)$ consisting of $X\in\cat{D}$, together with an $\calM$-relation $R\in\MSub(X\times X)$. Morphisms $f\colon (X,R)\to (Y,Q)$ are maps $f\colon X\to Y$ that are internally monotone in the sense of \cref{definition:internal monotone maps}.
\end{definition}

\section{Transitive closure}
Since left adjoint functors $F$ do not in general preserve limits, it is not automatic that $F(R)$ is an internal preorder whenever $R$ is. What can fail is transitivity. To fix this, we need an internal transitive closure procedure. This can be achieved using the complete lattice structure from \cref{lemma:M-complete}.

\begin{definition}\label{definition:transitive closure}
	The \emph{transitive closure} of an $\calM$-relation $R\in\MSub(X\times X)$ is defined as the meet over all transitive $\calM$-relations containing it:
		\[
			\widehat{R} := \bigwedge\left\{Q\in \MSub(X\times X):R\subseteq Q\text{~and~}Q\circ Q\subseteq Q \right\}.
		\]
\end{definition}

\begin{lemma}\label{lemma:meet of transitive relations is transitive}
	The meet of transitive $\calM$-relations is transitive.
\end{lemma}
\begin{proof}
	The empty meet produces the largest relation on $X$, which is trivially transitive. So we take a non-empty family $(Q_i)_{i\in I}$ of transitive $\calM$-relations on $X$. Write $q^i = (q^i_1,q^i_2)$ for their inclusions into $X\times X$. Let $q\colon Q\rightarrowtail X\times X$ denote their meet, obtained via wide pullback (see the proof of \cref{lemma:M-complete}). Explicitly, for every $i\in I$ we get an arrow $n_i\colon Q\rightarrowtail Q_i$, and $q=(q_1,q_2)$ is given by $q^i\circ n_i$, which does not depend on the index.
	
	We need to prove that $Q$ is transitive. For this we will use \cref{lemma:transitive in terms of factorisation}, so the goal is to construct an arrow $\tau\colon Q\times_X^{q_2,q_1}Q\to Q$ with the property that $q\circ \tau = (q_1\circ \pr_1,q_2\circ \pr_2)$, where $\pr_1$ and $\pr_2$ are the pullback projections of $Q\times_X^{q_2,q_1}Q$. Since the individual $Q_i$ are transitive, by the same lemma we get for each component $i\in I$ precisely such a factorisation:
		\[
			\begin{tikzcd}
				{Q_i} & {X\times X} \\
				{Q_i\times_X^{q^i_2,q^i_1}Q_i.}
				\arrow["{q^i}", tail, from=1-1, to=1-2]
				\arrow["{\tau_i}", dashed, from=2-1, to=1-1]
				\arrow["{(q_1^i\circ\pr_1^i,q_2^i\circ\pr_2^i)}"', from=2-1, to=1-2]
			\end{tikzcd}
		\]
	Here $\pr_1^i$ and $\pr_2^i$ denote the pullback projections of $Q_i\times_X^{q^i_2,q^i_1}Q_i$.
	We also get for each $i\in I$:
		\[
			\begin{tikzcd}[column sep=5ex,row sep=2ex]
				{Q\times_X^{q_2,q_1}Q} \\
				& {Q_i\times_X^{q_2^i,q_1^i}Q_i} & {Q_i} \\
				& {Q_i} & X,
				\arrow["{q^i_1}", from=2-3, to=3-3]
				\arrow["{q_2^i}"', from=3-2, to=3-3]
				\arrow["{\pr_1^i}"', from=2-2, to=3-2]
				\arrow["{\pr_2^i}", from=2-2, to=2-3]
				\arrow["{n_i\circ\pr_1}"', curve={height=18pt}, from=1-1, to=3-2,to path={ --([yshift=-3ex,xshift=0ex]\tikztostart.south) [pos=2]\tikztonodes |- (\tikztotarget) },rounded corners=12pt]
				\arrow["{n_i\circ \pr_2}", curve={height=-18pt}, from=1-1, to=2-3,rounded corners=12pt,to path={ --([yshift=0ex,xshift=3ex]\tikztostart.east) [pos=4]\tikztonodes  -|  (\tikztotarget) }]
				\arrow["{\exists! m_i}"', dashed, from=1-1, to=2-2]
				\arrow["\lrcorner"{anchor=center, pos=0.05}, draw=none, from=2-2, to=3-3]
			\end{tikzcd}
		\]
	since the outter square commutes by construction of the pullback $Q\times_X^{q_2,q_1}Q$. Using the data of the maps $m_i$ and $\tau_i$, we construct a new cone on the family~$(Q_i)_{i\in I}$:
		\[
			\begin{tikzcd}[column sep=-2ex, row sep=2ex]
				& {Q\times_X^{q_2,q_1}Q} \\
				{Q_i\times_X^{q^i_2,q^i_1}Q_i} && {Q_j\times_X^{q^j_2,q^j_1} Q_j} \\
				& Q \\
				{Q_i} && {Q_j} \\
				& {X\times X.}
				\arrow["{q^i}"', tail, from=4-1, to=5-2]
				\arrow["{q^j}", tail, from=4-3, to=5-2]
				\arrow["{n_i}"', tail, from=3-2, to=4-1]
				\arrow["{n_j}", tail, from=3-2, to=4-3]
				\arrow["{\tau_i}"', from=2-1, to=4-1]
				\arrow["{\tau_j}", from=2-3, to=4-3]
				\arrow["{m_i}"', from=1-2, to=2-1]
				\arrow["{m_j}", from=1-2, to=2-3]
				\arrow["{\exists !\tau}"{description}, dashed, from=1-2, to=3-2]
			\end{tikzcd}
		\]
	Here the outside commutes since the expression
		\begin{align*}
			q^i\circ \tau_i\circ m_i
			&=
			\left(q^i_1\circ \pr_1^i,q^i_2\circ \pr_2^i\right) \circ m_i
			\\&=
			\left(q^i_1\circ n_i\circ\pr_1,q^i_2\circ n_i\circ \pr_2\right)
			\\&=
			\left(q_1\circ\pr_1,q_2\circ\pr_2\right)
		\end{align*}
	does not depend on the index $i\in I$. This gives the arrow $\tau$ as in the diagram. We claim this is the desired arrow. To see that the condition of \cref{lemma:transitive in terms of factorisation} is fulfilled, note that by construction: $q\circ \tau = q^i\circ \tau_i\circ m_i$, which by the previous calculation in turn equals $(q_1\circ\pr_1,q_2\circ\pr_2)$, as desired. Hence $Q$ is transitive.
\end{proof}
	
\begin{proposition}
	Let $R\in\MSub(X\times X)$ be any $\calM$-relation on $X$. Then:
		\begin{enumerate}[label = (\alph*)]
			\item $R\subseteq \widehat{R}$ is the smallest transitive relation containing $R$;
			\item if $R$ is reflexive then $\widehat{R}$ is a preorder;
			\item if $R$ is transitive then $\widehat{R}\cong R$.
		\end{enumerate}
\end{proposition}
\begin{proof}
	For (a), $\widehat{R}$ is transitive by \cref{lemma:meet of transitive relations is transitive}, and it is by construction the smallest such containing~$R$. Claim (b) follows immediately since in that case $\Delta\subseteq R\subseteq \widehat{R}$. Lastly, for (c), if $R$ is already transitive then $\widehat{R}\subseteq R$, so combined with (a) this gives the desired isomorphism $\widehat{R}\cong R$. 
\end{proof}

\begin{lemma}
	Consider a family $(Q_i)_{i\in I}$ of $\calM$-relations on $Y$ that all make $f\colon (X,R)\to (Y,Q_i)$ monotone. Then $f\colon (X,R)\to (Y,\bigwedge_{i\in I}Q_i)$ is monotone.
\end{lemma}
\begin{proof}
	Since $f$ is monotone for each individual $Q_i$, there exist ${\overline{f}_i\colon R\to Q_i}$ with $q^i\circ \overline{f}_i = f\times f\circ r$. Use these to construct a new cone:
		
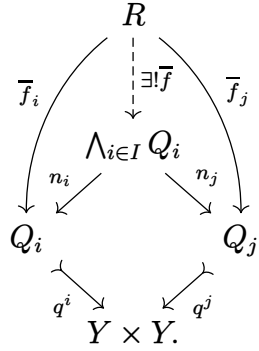
\[
			\begin{tikzcd}[ampersand replacement=\&, column sep=1ex, row sep=3ex]
				\& R \\
				\\
				\& \bigwedge_{i\in I} Q_i\\
				{Q_i} \&\& {Q_j} \\
				\& {Y\times Y.}
				\arrow["{\exists!\overline{f}}", dashed, from=1-2, to=3-2]
				\arrow["{\overline{f}_i}"', curve={height=16pt}, from=1-2, to=4-1]
				\arrow["{\overline{f}_j}", curve={height=-16pt}, from=1-2, to=4-3]
				\arrow["{n_i}"', from=3-2, to=4-1]
				\arrow["{n_j}", from=3-2, to=4-3]
				\arrow["{q^i}"', tail, from=4-1, to=5-2]
				\arrow["{q^j}", tail, from=4-3, to=5-2]
			\end{tikzcd}
		\]
	The outside square commutes since $q^i\circ \overline{f}_i = (f\times f)\circ r$ is independent on the index. Thus we obtain an arrow $\overline{f}\colon R\to Q$, and it is moreover straightforward to verify that $q\circ \overline{f} = q^i\circ n_i \circ \overline{f} = q^i\circ \overline{f}_i = (f\times f)\circ r$, concluding the proof.
\end{proof}

\begin{corollary}\label{corollary:pushforward of internal relation}
	For any arrow $f\colon X\to Y$ and relation $R\in\MSub(X\times X)$, there is a smallest $\calM$-relation $f_\ast (R)$ on $Y$ that makes $f$ monotone. Moreover, ${f\colon (X,R)\to (Y,Q)}$ is monotone if and only if $f_\ast(R)\subseteq Q$.
\end{corollary}

\subsection{Transitive extension}
\cref{definition:transitive closure} will allow us to define the object component of the functor $\overline{F}\colon \ord(\cat{C})\to \ord(\cat{D})$. For the arrow component we will need to be able to lift monotone maps $f\colon (X,R)\to (Y,Q)$ in $\cat{D}$ to monotone maps $f\colon (X,\widehat{R})\to (Y,\widehat{Q})$ between their transitive extensions. Unfortunately, there is no straightforward way to lift $\overline{f}\colon R\to Q$ to the transitive extensions in general.

However, this sort of lifting is certainly possible in $\Set$, since there the transitive closure of $R$ can equivalently be described by the union $\bigcup_{n\in\mathbb{N}}R^n$ of $n$-ary composites $R^n:=R\circ \cdots \circ R$ of the relation with itself. The resulting relation is denoted $R^\infty$, and its graph consists of pairs of elements that are linked by finite chains in $R$:
	\[
		x R^\infty y
		\qquad\text{if and only if}\qquad
		x R x_1 \cdots x_{n-1} R y,
	\]
for some $x_1,\ldots x_{n-1}\in X$. Thus if $f$ is monotone with respect to $R$ and $Q$, it is also monotone with respect to $R^\infty$ and $Q^\infty$:
	\[
		x_0 R x_1 \cdots x_{n-1} R x_n
		\qquad
		\text{implies}
		\qquad
		f(x_0)Q f(x_1)\cdots f(x_{n-1})Q f(x_n).
	\]
It turns out that this idea goes through in the abstract, provided that the category and its factorisation system are sufficiently nice. The crux is that $R^\infty$ is defined in terms of a \emph{colimit}, whereas $\widehat{R}$ is defined in terms of a \emph{limit}. The idea is thus to abstractly define $R^\infty$ as (the image of) a certain coproduct over $n$-ary compositions $R^n$. Recall from \cref{remark:non-associative relational composition} that we need to account for the non-associative case, so the expression $R^n$ needs to be treated a bit more carefully.

To take care of the parenthesisation, we make use of \emph{full binary trees}. These are finite rooted trees where every node has precisely zero or two children. The set of full binary trees is denoted by $\FBT$, and can be constructed recursively as follows. We admit the tree $\mathbf{1}\in\FBT$ consisting only of a single vertex, and if $S,T\in\FBT$ then we obtain a new full binary tree $S\oplus T$ by gluing them at a new root \cite[\S 5.3]{rosen2012DiscreteMathematicsIts}. A \emph{leaf} of $T$ is a vertex that admits no children, and we denote by $|T|$ the number of leaves of $T$. In particular, $|\mathbf{1}|=1$, and we obtain a tree $\mathbf{2}:= \mathbf{1}\oplus \mathbf{1}$ with $|\mathbf{2}|=2$. These are the unique full binary trees with one and two leaves, respectively. In general $|S\oplus T| = |S|+ |T|$. The point is now that the valid $n$-ary compositions of a non-associative binary operation are described precisely by the full binary trees~$T$ with $|T|=n$. 

For the next construction to work, we therefore assume:
\begin{enumerate}[label = (A\arabic*)]\setcounter{enumi}{3}
	\item\label{assumption:coproducts}
	$\cat{D}$ has all small coproducts.
\end{enumerate}

\begin{construction}\label{construction:transitive extension}
	Let $R$ be an $\calM$-relation on $X\in\cat{D}$. We recursively define a new $\calM$-relation $r^T\colon R^T\rightarrowtail X\times X$ on $X$ for every full binary tree $T\in\FBT$ as follows. First, we set $R^\mathbf{1} := R$, and for $S,T\in\FBT$ we set $R^{S\oplus T}:= R^S\circ R^T$. In particular $R^\mathbf{2} = R\circ R$, and in general $R^T$ corresponds to the appropriate parenthesisation of the composition of $|T|$-fold copies of $R$.
	
	Since we assume~\ref{assumption:coproducts}, we obtain a new $\calM$-relation $R^\infty$, called the \emph{transitive extension} of $R$, as follows:
		\[
			\begin{tikzcd}[cramped]
				{R^S} & {\displaystyle\coprod_{T\in\FBT}R^T} & {R^T} \\
				& {R^\infty} \\
				& {X\times X.}
				\arrow["{i^S}", from=1-1, to=1-2]
				\arrow["{i^T}"', from=1-3, to=1-2]
				\arrow["{r^S}"', curve={height=24pt}, tail, from=1-1, to=3-2]
				\arrow["{r^T}", curve={height=-24pt}, tail, from=1-3, to=3-2]
				\arrow["{r^\infty}"', tail, from=2-2, to=3-2]
				\arrow["e"', two heads, from=1-2, to=2-2]
				\arrow["{\exists!\Omega}"{description}, curve={height=-24pt}, dashed, from=1-2, to=3-2]
			\end{tikzcd}
		\]
	That is, $R^\infty$ is the $(\calE,\calM)$-image of the map $\Omega$ obtained via the universal property of the coproduct.
\end{construction}

\begin{remark}
	It can be seen that the diagram above indeed is the appropriate categorical generalisation of set-theoretic unions, see e.g.~\cite[p.~148]{goldblatt2006TopoiCategorialAnalysis}.
\end{remark}

The rest of this section is dedicated to showing that, in suitably nice categories~$\cat{D}$, the transitive extension $R^\infty$ actually coincides with the transitive closure~$\widehat{R}$. The assumptions are quite technical:
\begin{enumerate}[label = (A\arabic*)]\setcounter{enumi}{4}
	\item \label{assumption:coproducts distribute over pullbacks} coproducts distribute over pullbacks in $\cat{D}$;
	\item \label{assumption:beta in E} in $\cat{D}$, if $e\colon Q\tworightarrow R$ is in $\calE$ and $R\in \MSub(X\times X)$ then the canonical map
	\[
	Q\times_X^{r_2\circ e, r_1\circ e} Q \xrightarrow{~e\times_X e~} R\times_X^{r_2,r_1}R
	\]
	is in $\calE$.
\end{enumerate}

\begin{lemma}\label{lemma:transitive extension contains original}
	For any $R\in\MSub(X\times X)$, we have $R\subseteq R^\infty$.
\end{lemma}
\begin{proof}
	By construction, we obtain a coproduct inclusion $i^\mathbf{1}\colon R\to \coprod_{T\in\FBT}R^T$ that satisfies $r^\infty \circ e\circ i^\mathbf{1} = r$. This allows us to employ the orthogonal factorisation system, giving the desired inclusion:
		\[
			\begin{tikzcd}[cramped]
				R & {\displaystyle\coprod_{T\in\FBT}R^T} & {R^\infty} \\
				R && {X\times X.}
				\arrow[equal, from=1-1, to=2-1]
				\arrow["e", two heads, from=1-2, to=1-3]
				\arrow["{r^\infty}", tail, from=1-3, to=2-3]
				\arrow["r"', tail, from=2-1, to=2-3]
				\arrow["{\exists!}"{description}, dashed, from=2-1, to=1-3]
				\arrow["{i^\mathbf{1}}", from=1-1, to=1-2]
			\end{tikzcd}
		\]
	Note here that $\id_R\in\calE$ by \cref{lemma:factorisation classes contain isos}.
\end{proof}


\begin{lemma}\label{lemma:transitive expansion is transitive}
	$R^\infty$ is transitive for any~$R\in\MSub(X\times X)$.
\end{lemma}
\begin{proof}
	Consider the following diagram:
	\[
	\begin{tikzcd}[sep=1.5ex]
		{R^K\times_XR^L} & {\displaystyle\coprod_{S\in\FBT}\coprod_{T\in\FBT} R^S\times_XR^T} & {R^M\times_XR^N} \\
		& {\displaystyle\left(\coprod_{S\in\FBT}R^S\right)\times_X\left(\coprod_{T\in\FBT}R^T\right)} \\
		{R^K\circ R^L} && {R^M\circ R^N} \\
		{R^{K\oplus L}} & {\displaystyle\coprod_{T\in\FBT}R^T} & {R^{M\oplus N}} \\
		\\
		& {R^\infty} \\
		\\
		& {X\times X.}
		\arrow[from=1-1, to=1-2]
		\arrow[from=1-3, to=1-2]
		\arrow["\cong"{marking, allow upside down}, draw=none, from=1-2, to=2-2]
		\arrow["{r^\infty}", tail, from=6-2, to=8-2]
		\arrow["e", two heads, from=4-2, to=6-2]
		\arrow["{\exists!\alpha}", dashed, from=2-2, to=4-2]
		\arrow["{i^{K\oplus L}}", from=4-1, to=4-2]
		\arrow["{i^{M\oplus N}}"', from=4-3, to=4-2]
		\arrow["{r^{K\oplus L}}"', curve={height=18pt}, tail, from=4-1, to=8-2]
		\arrow["{r^{M\oplus N}}", curve={height=-18pt}, tail, from=4-3, to=8-2]
		\arrow["{e^{K,L}}"', two heads, from=1-1, to=3-1]
		\arrow[equals, from=3-1, to=4-1]
		\arrow["{e^{M,N}}", two heads, from=1-3, to=3-3]
		\arrow[equals, from=3-3, to=4-3]
	\end{tikzcd}
	\]
	We here use the distributivity of coproducts over pullbacks~\ref{assumption:coproducts distribute over pullbacks} to obtain the isomorphism
		\[
			\left(\coprod_{S\in\FBT}R^S\right)\times_X\left(\coprod_{T\in\FBT}R^T\right)
			\cong
			\coprod_{S\in \FBT}\left(R^S\times_X \left(\coprod_{T\in\FBT}R^T\right)\right)
			\cong
			\coprod_{S\in \FBT}\coprod_{T\in\FBT}\left(R^S\times_X R^T\right),
		\]
	and the arrow $\alpha$ arises from the universal property of the coproduct $\coprod_{T\in\FBT}R^T$. Similarly, in the following diagram the outermost square commutes by construction of the pullback of the coproducts, so by universality of the pullbacks of $R^\infty$ we get the existence of an arrow $\beta$:
	\[
		\begin{tikzcd}[cramped]
			{\displaystyle\left(\coprod_{S\in \FBT} R^S\right)\times_X \left(\coprod_{T\in\FBT}R^T\right)} && {\displaystyle\coprod_{T\in\FBT}R^T} \\
			& {R^\infty\times_XR^\infty} & {R^\infty} \\
			{\displaystyle\coprod_{S\in\FBT}R^S} & {R^\infty} & {X.}
			\arrow["{\pr_1}"', from=2-2, to=3-2]
			\arrow["{t^\infty}"', from=3-2, to=3-3]
			\arrow["{s^\infty}", from=2-3, to=3-3]
			\arrow["{\pr_2}", from=2-2, to=2-3]
			\arrow["{\pr_1}", from=1-1, to=3-1]
			\arrow["{\pr_2}"', from=1-1, to=1-3]
			\arrow["e"', two heads, from=1-3, to=2-3]
			\arrow["e"', two heads, from=3-1, to=3-2]
			\arrow["\lrcorner"{anchor=center, pos=0.125}, draw=none, from=2-2, to=3-3]
			\arrow["{\exists!\beta}"{description}, dashed, from=1-1, to=2-2]
		\end{tikzcd}
	\]
	Using~\ref{assumption:beta in E} we get $\beta\in \calE$, so the orthogonal factorisation provides the desired inclusion:
	\[
	\begin{tikzcd}[cramped]
		{\displaystyle\left(\coprod_{S\in \FBT} R^S\right)\times_X \left(\coprod_{T\in\FBT}R^T\right)} & {\displaystyle\coprod_{T\in\FBT}R^T} \\
		{R^\infty\times_X R^\infty} \\
		{R^\infty\circ R^\infty} & {R^\infty} \\
		{X\times X} & {X\times X.}
		\arrow["{\beta}"', from=1-1, to=2-1,two heads]
		\arrow[two heads, from=2-1, to=3-1]
		\arrow["\alpha", from=1-1, to=1-2]
		\arrow[two heads, from=1-2, to=3-2]
		\arrow[tail, from=3-2, to=4-2]
		\arrow[tail, from=3-1, to=4-1]
		\arrow["{\exists!}"{description}, dashed, from=3-1, to=3-2]
		\arrow[equals, from=4-1, to=4-2]
	\end{tikzcd}
	\]
\end{proof}

\begin{lemma}\label{lemma:transitive expansion of transitive relation}
	If $R$ is a transitive relation, then $R^\infty \cong R$. 
\end{lemma}
\begin{proof}
	First, we know from \cref{lemma:transitive extension contains original} that $R\subseteq R^\infty$, so it suffices to prove the converse inclusion. For that, from transitivity we get via induction that for every $T\in\FBT$ a $\calM$-subobject inclusion map $n^T\colon R^T\rightarrowtail R$, meaning $r\circ n^T = r^T$. Hence using the universal property of the coproduct we obtain a new arrow:
		\[
			\begin{tikzcd}[cramped]
				{R^S} & {\displaystyle\coprod_{T\in\FBT}R^T} & {R^T} \\
				\\
				& {R.}
				\arrow["{n^S}"', curve={height=18pt}, from=1-1, to=3-2]
				\arrow["{n^T}", curve={height=-18pt}, from=1-3, to=3-2]
				\arrow["{i^S}", from=1-1, to=1-2]
				\arrow["{i^T}"', from=1-3, to=1-2]
				\arrow["{\exists!n}"{description}, dashed, from=1-2, to=3-2]
			\end{tikzcd}
		\]
	To finish the proof, it therefore suffices to show that the outside square in the following diagram commutes:
		\[
			\begin{tikzcd}[cramped]
				{\displaystyle\coprod_{T\in\FBT}R^T} & R \\
				{R^\infty} & {X\times X,}
				\arrow["n", from=1-1, to=1-2]
				\arrow["e"', two heads, from=1-1, to=2-1]
				\arrow["{r^\infty}"', tail, from=2-1, to=2-2]
				\arrow["r", tail, from=1-2, to=2-2]
				\arrow["{\exists!}"{description}, dashed, from=2-1, to=1-2]
			\end{tikzcd}
		\]
	after which the orthogonal factorisation system gives the desired inclusion via the dashed diagonal arrow. To do this, note that by \cref{construction:transitive extension}, the arrow $\Omega =r^\infty\circ e$ is the unique one such that $\Omega\circ i^ T = r^T$, for every $T\in\FBT$. However, by construction of the subobject inclusions $n^T$, we see that $r\circ n \circ i^T = r\circ n^T = r^T$, from which it follows that $\Omega = r\circ n$. This concludes the proof. 
\end{proof}

This entire section essentially stands in service of the following lemma, which says that a monotone map can be lifted to a monotone map between transitive extensions. This lemma is crucial to make $\overline{F}$ well defined on morphisms.

\begin{lemma}\label{lemma:monotone maps lift to transitive extensions}
	Let $f\colon (X,R)\to (Y,Q)$ be a monotone map between ordered objects in $\cat{D}$. Then $f\colon (X,R^\infty)\to (Y,Q^\infty)$ is also monotone.
\end{lemma}
\begin{proof}
	If $f\colon (X,R)\to (Y,Q)$ is monotone, there exists an arrow $\overline{f}\colon R\to Q$ such that $q\circ \overline{f}= (f\times f )\circ r$. Via induction, the arrow $\overline{f}$ extends to arbitrarily parenthesised relations: for every $T\in\FBT$ we get
	\[
	\begin{tikzcd}
		{R^T} & {Q^T} \\
		{X\times X} & {Y\times Y.}
		\arrow["{q^T}", tail, from=1-2, to=2-2]
		\arrow["{r^T}"', tail, from=1-1, to=2-1]
		\arrow["{f\times f}"', from=2-1, to=2-2]
		\arrow["{\exists\overline{f}^T}", dashed, from=1-1, to=1-2]
	\end{tikzcd}
	\] 
	We use these to construct a cocone, and hence the universality of the coproduct produces a new arrow:
	\[
	\begin{tikzcd}
		{R^S} & {\displaystyle\coprod_{T\in\FBT}R^T} & {R^T} \\
		{Q^S} & {\displaystyle\coprod_{T\in\FBT}Q^T} & {Q^T.}
		\arrow["{\overline{f}^S}"', from=1-1, to=2-1]
		\arrow["{\overline{f}^T}", from=1-3, to=2-3]
		\arrow["{i^S}", from=1-1, to=1-2]
		\arrow["{i^T}"', from=1-3, to=1-2]
		\arrow["{j^S}"', from=2-1, to=2-2]
		\arrow["{j^T}", from=2-3, to=2-2]
		\arrow["{\exists!\alpha}", dashed, from=1-2, to=2-2]
	\end{tikzcd}
	\]
	To finish the proof, it now suffices to show that the outside rectangle in the following diagram commutes:
	\[
	\begin{tikzcd}
		{\displaystyle\coprod_{T\in\FBT}R^T} & {\displaystyle\coprod_{T\in\FBT}Q^T} & {Q^\infty} \\
		{R^\infty} & {X\times X} & {Y\times Y,}
		\arrow["\alpha", from=1-1, to=1-2]
		\arrow["{q^\infty}", tail, from=1-3, to=2-3]
		\arrow["d", two heads, from=1-2, to=1-3]
		\arrow["{f\times f}"', from=2-2, to=2-3]
		\arrow["{r^\infty}"', tail, from=2-1, to=2-2]
		\arrow["e"', two heads, from=1-1, to=2-1]
		\arrow["{\exists!\overline{f}^\infty}"', dashed, from=2-1, to=1-3]
	\end{tikzcd}
	\]
	after which the orthogonal factorisation system gives the depicted dashed arrow, proving the monotonicity of $f$ with respect to the transitive extensions. To prove commutativity of the outer rectangle, we use the fact that the coproduct injections $(i^T)_{T\in\FBT}$ are jointly epic. Keeping this in mind, we calculate:
	\[
	q^\infty \circ d\circ \alpha\circ i^T
	=
	q^\infty \circ d\circ j^T\circ \overline{f}^T
	=
	q^T\circ \overline{f}^T
	=
	(f\times f) \circ r^T
	= (f\times f) \circ r^\infty \circ e \circ i^T.
	\]
	This uses the defining property of $\alpha$, the definition of the coproduct injections of $\coprod_{T\in\FBT}Q^T$, the construction of $\overline{f}^T$, and lastly the definition of the coproduct injections of $\coprod_{T\in\FBT}R^T$, respectively. Since this calculation goes through irrespective of $T$, the commutativity of the outer rectangle follows, and this concludes the proof.	
\end{proof}

\begin{corollary}\label{corollary:transitive extension preserves inclusions}
	If $R\subseteq Q$ then $R^\infty \subseteq Q^\infty$.
\end{corollary}
\begin{proof}
	This is a special case of \cref{lemma:monotone maps lift to transitive extensions} where $f=\id_X$, and noting that $\id_X\colon (X,R)\to (X,Q)$ is monotone if and only if $R\subseteq Q$.
\end{proof}

\begin{lemma}\label{lemma:transitive expansion is transitive closure}
	For any $\calM$-relation $R$ we have $\widehat{R}\cong R^\infty$. 
\end{lemma}
\begin{proof}
	First, by \cref{lemma:transitive expansion is transitive}, the relation $R^\infty$ is transitive, so since $\widehat{R}$ is the smallest transitive relation containing $R$ we get $\widehat{R}\subseteq R^\infty$. For the converse, recall that $\widehat{R}$ is the wide pullback over all transitive relations $R_i$ containing $R$. By \cref{lemma:transitive expansion of transitive relation,corollary:transitive extension preserves inclusions} we get $R^\infty \subseteq R_i^\infty \cong R_i$. So we get inclusion arrows $n^i \colon R^\infty \rightarrowtail R_i$ that allow us to construct a cone:
	
		
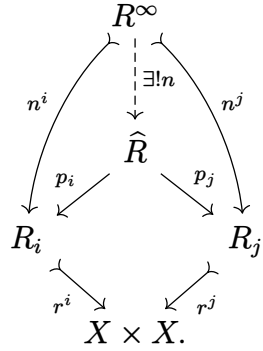
\[
			\begin{tikzcd}[column sep=1ex,row sep=3ex]
				& {R^\infty} \\
				\\
				& {\widehat{R}} \\
				{R_i} && {R_j} \\
				& {X\times X.}
				\arrow["{n^i}"', curve={height=12pt}, tail, from=1-2, to=4-1]
				\arrow["{n^j}", curve={height=-12pt}, tail, from=1-2, to=4-3]
				\arrow["{r^i}"', tail, from=4-1, to=5-2]
				\arrow["{r^j}", tail, from=4-3, to=5-2]
				\arrow["{p_i}"', from=3-2, to=4-1]
				\arrow["{p_j}", from=3-2, to=4-3]
				\arrow["{\exists! n}", dashed, from=1-2, to=3-2]
			\end{tikzcd}
		\]
	The outside clearly commutes, since $r^i\circ n^i = r^\infty$ is independent on the index, and thus we obtain the depicted arrow $n\colon R^\infty \to \widehat{R}$. We claim that this facilitates the desired inclusion. For this, simply note that $\widehat{r}\circ n = r^i\circ p_i \circ n = r^i\circ n^i = r^\infty$, as desired.
\end{proof}

\begin{corollary}\label{corollary:monotone map lifts to monotone map transitive closures}
	In a category $\cat{D}$ satisfying assumptions \ref{assumption:well-powered}-\ref{assumption:beta in E}, a monotone map ${f\colon (X,R)\to (Y,Q)}$ between ordered objects lifts to a monotone map $f\colon (X,\widehat{R})\to (Y,\widehat{Q})$.
\end{corollary}


\section{Lifting an adjunction}
\label{section:internal ordered adjunction}
We are now finally ready to put together all the pieces. The basic idea is that on objects $(S,Q)\in \ord(\cat{C})$ the functor $\overline{F}$ defines a preorder~$\overline{F}(Q)$ on $F(S)$, and dually for $(X,R)\in\ord(\cat{D})$ the functor $\overline{G}$ defines a preorder $\overline{G}(R)$ on $G(X)$.

That $G$ is a right adjoint implies that it preserves all limits, so in particular it preserves products and pullbacks. Explicitly, this means that for every $X\in\cat{D}$ the canonical map~$\alpha$ in the diagram
\[
\begin{tikzcd}[ampersand replacement=\&]
	{G(X\times X)} \& {GX\times GX} \\
	\& GX
	\arrow["{\exists!\alpha}", from=1-1, to=1-2]
	\arrow["{G(\pr_i)}"', from=1-1, to=2-2]
	\arrow["{\pr_i}", from=1-2, to=2-2]
\end{tikzcd}
\]
is an isomorphism. Similarly, by functoriality of $F$ we obtain canonical maps
	\[
		\begin{tikzcd}[ampersand replacement=\&]
			{F(S\times S)} \& {FS\times FS} \\
			\& FS
			\arrow["{\exists!\mu_S}", from=1-1, to=1-2]
			\arrow["{F(\pr_i)}"', from=1-1, to=2-2]
			\arrow["{\pr_i}", from=1-2, to=2-2]
		\end{tikzcd}
	\]
for every $S\in\cat{C}$, but these are not necessarily isomorphisms. This leads to the following new assumptions:

\begin{enumerate}[label = (A\arabic*)]\setcounter{enumi}{6}
	\item\label{assumption:functors preserve M} if $m\in\calM_\cat{C}$ then $F(m)\in\calM_\cat{D}$, and if $n\in\calM_\cat{D}$ then $G(n)\in\calM_\cat{C}$;
	\item \label{assumption:canonical mu in M} $\mu_S\in\calM_\cat{C}$ for every $S\in\cat{C}$.
\end{enumerate}

We can now describe the object parts of the functors $\overline{F}$ and $\overline{G}$. 

\begin{construction}\label{construction:object parts of ordered functors}
	Given $(X,R)\in\ord(\cat{D})$, explicitly given by ${r\colon R\rightarrowtail X\times X}$, by~\ref{assumption:functors preserve M} we obtain an $\calM_\cat{C}$-subobject $G(r)\colon GR\rightarrowtail G(X\times X)$, which under the canonical isomorphism $\alpha$ (which is in $\calM_\cat{D}$ by \cref{lemma:factorisation classes contain isos}) can be identified with an $\calM_\cat{C}$-relation on $GX$. Since $G$ preserves all limits, if $R$ is an $\calM_\cat{D}$-preorder it follows using \cref{lemma:reflexive in terms of factorisation,lemma:transitive in terms of factorisation} that $GR$ is an $\calM_\cat{C}$-preorder on $GX$. This will constitute the object component of the functor $\overline{G}$. 
	
	The construction for $\overline{F}$ is analogous, but requires more care since $F$ does not preserve limits. Take an internal preordered object $(S,Q)\in\ord(\cat{C})$, with $q\colon Q\rightarrowtail S\times S$. Using~\ref{assumption:functors preserve M} and~\ref{assumption:canonical mu in M} we get an $\calM_\cat{D}$-relation
		\[
			\begin{tikzcd}[ampersand replacement=\&]
				FQ \& {F(S\times S)} \& {FS\times FS}
				\arrow["{F(q)}", tail, from=1-1, to=1-2]
				\arrow["{\mu_S}", tail, from=1-2, to=1-3]
			\end{tikzcd}
		\]
	on $FS$. We claim that if $Q$ is reflexive, so is $FQ$. Using \cref{lemma:reflexive in terms of factorisation}, that $Q$ is reflexive means there exists $\rho\colon S\to Q$ such that $q\circ \rho = \delta_S$. The canonical nature of $\mu_S$ implies that $\mu_S\circ F(\delta_S)= \delta_{FS}$, so the reflexivity of $FQ$ follows:
		\[
			\begin{tikzcd}[ampersand replacement=\&]
				FQ \& {F(S\times S)} \& {FS\times FS.} \\
				FS
				\arrow["{F(q)}", tail, from=1-1, to=1-2]
				\arrow["{\mu_S}", tail, from=1-2, to=1-3]
				\arrow["{F(\rho)}", from=2-1, to=1-1]
				\arrow["{F(\delta_S)}"', from=2-1, to=1-2]
				\arrow["{\delta_{FS}}"', curve={height=18pt}, from=2-1, to=1-3]
			\end{tikzcd}
		\]
	Since we cannot assume $F$ preserves pullbacks, transitivity unfortunately does not lift from $Q$ to $FQ$. Instead we define:
		\[
			\overline{F}(S,Q):= \left(FS,\widehat{FQ}\right),
		\]
	where $\widehat{FQ}$ denotes the transitive closure of the $\calM_\cat{D}$-relation $\mu_S\circ F(q)$ above. If~$Q$ is a $\calM_\cat{C}$-preorder, $\widehat{FQ}$ will be a $\calM_{\cat{D}}$-preorder.
\end{construction}

\begin{proposition}\label{proposition:ordered functors}
	There are functors:
	
	\begin{minipage}{.4\textwidth}
		\begin{align*}
			\overline{F}: \ord(\cat{C}) &\longrightarrow \ord(\cat{D});\\
			(S,Q) & \longmapsto (FS,\widehat{FQ});\\
			g & \longmapsto F(g)
		\end{align*}
	\end{minipage}%
	\hfill and\hfill
	\begin{minipage}{.4\textwidth}
		\begin{align*}
			\overline{G}: \ord(\cat{D}) &\longrightarrow \ord(\cat{C});\\
			(X,R) & \longmapsto (GX,GR);\\
			f & \longmapsto G(f).
		\end{align*}
	\end{minipage}
\end{proposition}
\begin{proof}
	By \cref{construction:object parts of ordered functors} we know that these assignments are well-defined on objects. Clearly $\overline{F}$ and $\overline{G}$ will inherit functoriality directly from $F$ and $G$, respectively, so we are left to show that they are indeed well-defined on morphisms. This amounts to showing they preserve monotonicity. If ${g\colon (S,R)\to (T,Q)}$ is a monotone map in $\ord(\cat{C})$ then the functoriality of $F$ and the properties of the canonical map $\mu_S$ tell us that
		\[
			\begin{tikzcd}[ampersand replacement=\&]
				R \& Q \\
				{S\times S} \& {T\times T}
				\arrow["{\exists \overline{g}}", dashed, from=1-1, to=1-2]
				\arrow["r"', tail, from=1-1, to=2-1]
				\arrow["q", tail, from=1-2, to=2-2]
				\arrow["{g\times g}"', from=2-1, to=2-2]
			\end{tikzcd}
			\qquad\text{implies}\qquad 
			\begin{tikzcd}[ampersand replacement=\&]
				FR \& FQ \\
				{F(S\times S)} \& {F(T\times T)} \\
				{FS\times FS} \& {FT\times FT}
				\arrow["{F(\overline{g})}", dashed, from=1-1, to=1-2]
				\arrow["{F(r)}"', tail, from=1-1, to=2-1]
				\arrow["{F(q)}", tail, from=1-2, to=2-2]
				\arrow["{F(g\times g)}"', from=2-1, to=2-2]
				\arrow["{\mu_S}"', tail, from=2-1, to=3-1]
				\arrow["{\mu_T}", tail, from=2-2, to=3-2]
				\arrow["{F(g)\times F(g)}"', from=3-1, to=3-2]
			\end{tikzcd}
		\]
	commutes, showing that $F(g)\colon (FS,FR)\to (FT,FQ)$ is monotone. Using \cref{corollary:monotone map lifts to monotone map transitive closures}, this map lifts to a monotone map between the transitive extensions, which is the desired form. The proof for $\overline{G}$ is analogous (and simpler).
\end{proof}

\subsection{The unit and counits are monotone}\label{section:units and counits are internally monotone}
In this section we show that the unit $\eta_S\colon S\to GFS$ and counit $\epsilon_X\colon FGX\to X$ of the adjunction $F\dashv G$ lift to monotone maps.

\begin{lemma}\label{lemma:functor preserves subobject inclusions}
	If $A\subseteq B$ in $\MSub(X)$ implies $GA\subseteq GB$ in $\MSub(GX)$.
\end{lemma}
\begin{proof}
	That $A\subseteq B$ means there exists an arrow $m\colon A\to B$ such that ${m\circ b = a}$. The arrow ${G(m)\colon G(A)\to G(B)}$ defines an inclusion of $\calM$-subobjects by~\ref{assumption:functors preserve M} and functoriality.
\end{proof}

\begin{proposition}\label{proposition:unit is internally monotone}
	The unit map $\eta_S\colon (S,R)\to(GFS,G(\widehat{FR}))$ is monotone for any $(S,R)\in \OrdC$.
\end{proposition}
\begin{proof}
	Since there is an inclusion $FR\subseteq \widehat{FR}$ in $\MSub(FS\times FS)$, and $G$ preserves this inclusion (\cref{lemma:functor preserves subobject inclusions}), by \cref{lemma:monotone order modification} it suffices to prove that the map $\eta_S\colon (S,R)\to (GFS,GFR)$ is monotone. To prove this, in turn, we need to show that the following diagram commutes:
		\[
			\begin{tikzcd}[column sep=large,row sep=small]
				R & GFR \\
				\\
				{S\times S} & {GF(S\times S)} \\
				\\
				& {G(FS\times FS)} \\
				{S\times S} & {GFS\times GFS.}
				\arrow["r"', tail, from=1-1, to=3-1]
				\arrow["{\eta_R}", from=1-1, to=1-2]
				\arrow["\cong"{marking, allow upside down}, draw=none, from=5-2, to=6-2]
				\arrow["{GF(r)}", from=1-2, to=3-2]
				\arrow["{G(\mu_S)}", from=3-2, to=5-2]
				\arrow["{\eta_S\times \eta_S}"', from=6-1, to=6-2]
				\arrow[equal, from=3-1, to=6-1]
				\arrow["{\eta_{S\times S}}"', from=3-1, to=3-2]
			\end{tikzcd}
		\]
	The top square does so precisely because $\eta$ is natural. To prove that the bottom square commutes, we use that $\mu$ and $G$ preserve products. Denoting the isomorphism in the diagram by $\alpha$, the preservation of products means specifically that the following diagram commutes:
		\[
			\begin{tikzcd}
				{GF(S\times S)} & {G(FS\times FS)} \\
				GFS & {GFS\times GFS.}
				\arrow["{GF(\pr_i)}"', from=1-1, to=2-1]
				\arrow["{G(\mu_S)}", from=1-1, to=1-2]
				\arrow["\alpha", from=1-2, to=2-2]
				\arrow["{\pr_i}", from=2-2, to=2-1]
				\arrow["{G(\pr_i)}"{description}, from=1-2, to=2-1]
			\end{tikzcd}
		\]
	Using this, we find
		\[
			\pr_i \circ \alpha\circ G(\mu_S) \circ \eta_{S\times S}
			=
			GF(\pr_i)\circ \eta_{S\times S}
			=
			\eta_S\circ \pr_i
			=
			\pr_i\circ (\eta_S\times \eta_S),			
		\]
	where the second-to-last step is again due to naturality of $\eta$. That the desired square commutes now follows since projections are jointly monomorphic.
\end{proof}

\begin{lemma}\label{lemma:counit is internally monotone}
	The counit map $\epsilon_X\colon (FGX,FGR)\to (X,R)$ is monotone for any $\calM_\cat{D}$-relation $R$ on $X$.
\end{lemma}
\begin{proof}
	To prove monotonicity, we need to show the following diagram commutes:
		\[
			\begin{tikzcd}
				FGR & R \\
				{FG(X\times X)} & {X\times X} \\
				{F(GX\times GX)} \\
				{FGX\times FGX} & {X\times X.}
				\arrow["{\epsilon_R}", from=1-1, to=1-2]
				\arrow["{FG(r)}"', from=1-1, to=2-1]
				\arrow["{F(\alpha)}"', from=2-1, to=3-1]
				\arrow["{\mu_{GX}}"', from=3-1, to=4-1]
				\arrow["{\epsilon_X\times\epsilon_X}"', from=4-1, to=4-2]
				\arrow[equals,from=2-2, to=4-2]
				\arrow["r", tail, from=1-2, to=2-2]
				\arrow["{\epsilon_{X\times X}}", from=2-1, to=2-2]
			\end{tikzcd}
		\]
	Here, as in the proof of \cref{proposition:unit is internally monotone}, we denote by $\alpha\colon G(X\times X)\to GX\times GX$ the natural product preservation isomorphism. The top square in the diagram commutes by naturality of the counit~$\epsilon$. To prove that the bottom square commutes, note that the preservation of products by $\alpha$ and $\mu$ implies the existence of a commutative diagram:
		\[
			\begin{tikzcd}
				{FG(X\times X)} & {F(GX\times GX)} \\
				FGX & {FGX\times FGX.}
				\arrow["{FG(\pr_i)}"', from=1-1, to=2-1]
				\arrow["{F(\alpha)}", from=1-1, to=1-2]
				\arrow["{\mu_{GX}}", from=1-2, to=2-2]
				\arrow["{\pr_i}", from=2-2, to=2-1]
				\arrow["{F(\pr_i)}"{description}, from=1-2, to=2-1]
			\end{tikzcd}
		\]
	Using this, we calculate:
		\[
			\pr_i \circ (\epsilon_X\times \epsilon_X)\circ \mu_{GX}\circ F(\alpha)
			=
			\epsilon_X\circ\pr_i \circ \mu_{GX}\circ F(\alpha)
			=
			\epsilon_X\circ FG(\pr_i)
			=
			\pr_i\circ \epsilon_{X\times X},
		\]
	where the first step follows by definition of $\epsilon_X\times \epsilon_X$, and the last step follows by naturality of the counit $\epsilon$. Since the projection maps are jointly monomorphic, the commutativity of the first diagram follows.
\end{proof}

\begin{proposition}\label{proposition:counit is internally monotone}
	The counit map $\epsilon_X\colon (FGX, \widehat{FGR})\to (X,R)$ is monotone for any $(X,R)\in \OrdD$.
\end{proposition}
\begin{proof}
	This now follows directly by using \cref{lemma:monotone maps lift to transitive extensions} to lift the monotone map from \cref{lemma:counit is internally monotone} to a monotone map $\epsilon_X\colon (FGX,(FGR)^\infty)\to (X,R^\infty)$, and then using \cref{lemma:transitive expansion is transitive closure,lemma:transitive expansion of transitive relation} to reduce $(FGR)^\infty\cong \widehat{FGR}=:\overline{F}(GR)$ and $R^\infty \cong R$, which is the desired form.
\end{proof}

\subsection{The adjunction}
We are now ready to prove that the functors from \cref{proposition:ordered functors} define an adjunction $\overline{F}\dashv \overline{G}$. With the results from the previous section, the proof is relatively straightforward (cf.~the proof of \cref{theorem:adjunction ordtopOC bullet and ordloc bullet}).

\begin{theorem}\label{theorem:internal ordered adjunction}
	Under assumptions~\ref{assumption:well-powered}-\ref{assumption:canonical mu in M} on
	$\begin{tikzcd}[ampersand replacement=\&]
		{\cat{C}} \& {\cat{D}}
		\arrow[""{name=0, anchor=center, inner sep=0}, "F", shift left=2, from=1-1, to=1-2]
		\arrow[""{name=1, anchor=center, inner sep=0}, "G", shift left=2, from=1-2, to=1-1]
		\arrow["\dashv"{anchor=center, rotate=-90}, draw=none, from=0, to=1]
	\end{tikzcd}$, there exists an adjunction between the induced functors of \cref{proposition:ordered functors}:
		\[
			\begin{tikzcd}[ampersand replacement=\&]
				{\ord(\cat{C})} \& {\ord(\cat{D}).}
				\arrow[""{name=0, anchor=center, inner sep=0}, "{\overline{F}}", shift left=2, from=1-1, to=1-2]
				\arrow[""{name=1, anchor=center, inner sep=0}, "{\overline{G}}", shift left=2, from=1-2, to=1-1]
				\arrow["\dashv"{anchor=center, rotate=-90}, draw=none, from=0, to=1]
			\end{tikzcd}
		\]
\end{theorem}
\begin{proof}
	From $F\dashv G$ we get natural bijections $\cat{D}(FS,X)\cong \cat{C}(S,GX)$ that by \cref{proposition:unit is internally monotone,proposition:counit is internally monotone} restrict to natural bijections
		\begin{align*}
			\ord(\cat{D})((FS,\widehat{FQ}),(X,R))  &\xrightarrow{~ \sim ~}  \ord(\cat{C})((S,Q),(GX,GR))\\
			f&\xmapsto{\hspace{3ex}} G(f)\circ \eta_S\\ \epsilon_X\circ F(g)&{~\reflectbox{\ensuremath{\xmapsto{\hspace{3ex}}}}} ~g	.			\qedhere	
		\end{align*}
\end{proof}

\begin{remark}
	We conjecture that an alternative sufficient assumption on $\cat{D}$ and its factorisation system is that the induced forgetful functor $\ord(\cat{D})\to\cat{D}$ is \emph{topological} (in the sense of \cite{adamek1990AbstractConcreteCategories}). This should allow us to perform the ``monotone lifting'' procedure from \cref{corollary:monotone map lifts to monotone map transitive closures}. We leave the relation between $\ord(\cat{D})\to\cat{D}$ being a topological functor and our assumptions~\ref{assumption:well-powered}-\ref{assumption:canonical mu in M} to future research.
\end{remark}

\section{Internally ordered locales}
In this section we plug in $\cat{C}=\Top$ and $\cat{D}=\Loc$, and build towards a lifted adjunction. This boils down to a proof that $\loc\dashv\pt$ satisfies~{\ref{assumption:well-powered}-\ref{assumption:canonical mu in M}}. Unfortunately, at this stage, the proof is still incomplete.

For convenience, we provide a summarised list:
\begin{enumerate}[label= (A\arabic*)]
	\item $\cat{C}$ and $\cat{D}$ are $\calM$-well-powered;
	\item $\cat{C}$ and $\cat{D}$ have all small limits;
	\item $\delta_S\in \calM_{\cat{C}}$ for every $S\in\cat{C}$ and $\delta_X\in\calM_{\cat{D}}$ for every $X\in\cat{D}$;
	\item $\cat{D}$ has all small coproducts;
	\item in $\cat{D}$, coproducts distribute over pullbacks;
	\item in $\cat{D}$, if $e\colon Q\tworightarrow R$ is in $\calE_\cat{D}$ and $R$ is an $\calM_{\cat{D}}$-relation on $X$, then $e\times_X e\in\calE_\cat{D}$;
	\item $F$ and $G$ preserve the right classes of the factorisation systems;
	\item $\mu_S\in\calM_{\cat{D}}$ for every $S\in\cat{C}$.
\end{enumerate}

Here we equip $\Top$ and $\Loc$ with the epi-extremal monomorphism factorisation systems.

\begin{remark}
	Note that $\Frm$ is not \emph{co-well-powered} \cite[\S IV.6.6]{picado2012FramesLocalesTopology}, meaning that $\Loc$ is not well-powered with respect to the class of all monomorphisms. Switching to extremal monomorphisms solves this issue.
\end{remark}

\begin{lemma}\label{lemma:loc preserves extremal monos}
	If $g\colon S\to T$ is an extremal monomorphism in $\Top$, then $\loc(g)$ is an extremal monomorphism in $\Loc$.
\end{lemma}
\begin{proof}
	That $g$ is an extremal monomorphism implies $g$ is injective and the topology on $S$ is the initial topology generated by it: ${\Opens S = \{g^{-1}(V):V\in\Opens T\}}$, so clearly $\loc(g)$ has a surjective frame map.
\end{proof}

\begin{lemma}\label{lemma:pt preserves extremal monos}
	If $f\colon X\to Y$ is an extremal monomorphism in $\Loc$, then $\pt(f)$ is an extremal monomorphism in $\Top$.
\end{lemma}
\begin{proof}
	That $f$ is an extremal monomorphism means $f^{-1}$ is a surjection of frames. This clearly implies
		\begin{align*}
			\Opens \pt(X) &= \{\pt(U):U\in\Opens X\}
			\\&=\{ \pt(f^{-1}(V)):V\in\Opens Y\}
			\\&= \{ \pt(f)^{-1}(\pt(V)):V\in\Opens Y\},
		\end{align*}
	where in the last step we use \cref{lemma:pt(f) continuity}. Thus $\pt(X)$ has the initial topology induced by $\pt(f)$, and we are only left to show this is an injection. For that, take $\calF,\calG\in\pt(X)$ with $\pt(f)(\calF)=\pt(f)(\calG)$. The latter means $f^{-1}(V)\in \calF$ iff $f^{-1}(V)\in \calG$, and given the equation above this implies $\calF = \calG$, as desired.
\end{proof}

Unfortunately, at this point we do not know if assumption~\ref{assumption:beta in E} holds in $\Loc$. The result certainly holds for spaces, since the product of surjective functions is surjective. We suggest that the proof (or the construction of a counterexample) could involve the explicit description of the product of maps of locales in \cite[Corollary~IV.5.5.3]{picado2012FramesLocalesTopology}. We leave this technical matter to future research, but in the meantime discuss some of the potential consequences.

\begin{conjecture}\label{conjecture:adjunction loc pt satisfies all assumptions}
	The adjunction $\begin{tikzcd}[ampersand replacement=\&]
		{\Top} \& {\Loc}
		\arrow[""{name=0, anchor=center, inner sep=0}, "\loc", shift left=2, from=1-1, to=1-2]
		\arrow[""{name=1, anchor=center, inner sep=0}, "\pt", shift left=2, from=1-2, to=1-1]
		\arrow["\dashv"{anchor=center, rotate=-90}, draw=none, from=0, to=1]
	\end{tikzcd}$ with the epi-extremal monomorphism factorisation systems satisfies~\ref{assumption:well-powered}-\ref{assumption:canonical mu in M}.
\end{conjecture}
\begin{proof}[Proof sketch]
	We proceed in order:
	\begin{enumerate}[label= (A\arabic*)]
		\item in $\Top$ we have $\MSub(S)\cong \Powerset(S)$, and in $\Loc$ we have $\MSub(X)\cong \Sl(X)$, which are both sets (\cite[\S III.2]{picado2012FramesLocalesTopology});
		
		\item it is well-known that $\Top$ is a complete category, and $\Loc$ is complete by \cite[Corollary~IV.4.3.5]{picado2012FramesLocalesTopology};
		
		\item the diagonal function $\delta_S$ is an extremal monomorphism for any $S\in\Top$, since $\delta_S^{-1}(U\times V) = U\cap V\in\Opens S$. Similarly, the diagonal map $\delta_X$ for $X\in\Loc$ has an explicit expression (\cref{lemma:pair map preimage})
			\[
				\delta_X^{-1}(U\oplus V) = U\wedge V
			\]
		which clearly defines a surjection. Thus $\delta_X$ defines a sublocale.
		
		\item $\Loc$ is cocomplete by \cite[Corollary~IV.4.3.5]{picado2012FramesLocalesTopology};
		
		\item that coproducts distribute over pullbacks in $\Loc$ follows from the explicit description of coproducts in \cite[Section~2.2.6]{vickers1999TopicalCategoriesDomains};\footnote{We thank Steve Vickers for pointing out a sketch for the extensivity of $\Loc$ in private communication.}
		
		\setcounter{enumi}{6}
		\item this is \cref{lemma:loc preserves extremal monos,lemma:pt preserves extremal monos};
		
		\item this follows by \cite[Proposition~IV.5.4.1]{picado2012FramesLocalesTopology}.\qedhere
	\end{enumerate}
\end{proof}

\begin{corollary}
	If \cref{conjecture:adjunction loc pt satisfies all assumptions} holds, then there is an adjunction:
		\[
			\begin{tikzcd}[ampersand replacement=\&]
				{\ord(\Top)} \& {\ord(\Loc).}
				\arrow[""{name=0, anchor=center, inner sep=0}, "\loc", shift left=2, from=1-1, to=1-2]
				\arrow[""{name=1, anchor=center, inner sep=0}, "\pt", shift left=2, from=1-2, to=1-1]
				\arrow["\dashv"{anchor=center, rotate=-90}, draw=none, from=0, to=1]
			\end{tikzcd}
		\]
\end{corollary}

\begin{remark}
	Note here that $\ord(\Top)$ is just the category $\OrdTop$ of ordered topological spaces (with no relation between order and topology) and monotone continuous functions from \cref{section:ordered spaces}. The category $\ord(\Loc)$ is something entirely new: it is the category of locales $X$ equipped with a \emph{sublocale} $R\rightarrowtail X\times X$ that is internally reflexive and transitive. This generalises the notion of locales equipped with a \emph{closed} reflexive and transitive sublocale relation in \cite[Chapter~5]{townsend1996preframeTechniquesConstructiveLocale}.
\end{remark}


\subsection{Internally ordered locales with open cones}
To finish this section, and the thesis, we discuss a possible relation between the \emph{internally} ordered locales $(X,R)\in\ord(\Loc)$, and the \emph{externally} ordered locales $(X,\Leq)\in\OrdLoc$ from \cref{definition:ordered locale}. The motivation of the latter definition of ordered locales using an auxiliary preorder~$\Leq$ on the frame of opens might seem somewhat ad hoc from a structural perspective. For instance, it is not clear where the axiom~\eqref{axiom:V} comes from: in essence it is a condition about sup-lattices, rather than of frames. Similarly, it is unsatisfying that we need to impose~\eqref{axiom:LV}, while this is automatic in the case of spaces. Here, we conjecture that externally ordered locales can be recovered from a certain class of internally ordered locales, thereby giving a more structural justification for \cref{definition:ordered locale}.

Even in general, for an internally preordered object $(X,R)\in\ord(\cat{D})$ we get a notion of \emph{internal cones}. These are constructed as follows. For this construction to work, note that we assume $\cat{D}$ is complete, and so admits a terminal object $1$. Thus, an $\calM$-subobject $A\rightarrowtail X$ can equivalently be viewed as an $\calM$\emph{-relation} from $X$ to $1$, or the other way around. Composing this relation with $R$ therefore gives new $\calM$-subobjects:
	\[
		\up A:= A\circ R\rightarrowtail 1\times X\cong X
		\quad\text{and}\quad
		\down A:= R\circ A\rightarrowtail X\times 1\cong X.
	\]

\begin{definition}
	\label{definition:internal cones}
	Let $\cat{D}$ be a category with finite limits and orthogonal factorisation system $(\calE,\calM)$. The \emph{internal cones} of $X\in\cat{D}$ equipped with a relation $R\in\MSub(X\times X)$ are the functions%
	\begin{center}
		\begin{minipage}{.45\textwidth}
			\begin{align*}
				\up(-)\colon \MSub(X) & \longrightarrow \MSub(X);\\
				A & \longmapsto A\circ R,
			\end{align*}
		\end{minipage}%
		\begin{minipage}{.45\textwidth}
			\begin{align*}
				\down(-)\colon \MSub(X) & \longrightarrow \MSub(X);\\
				A & \longmapsto R\circ A.
			\end{align*}
		\end{minipage}
	\end{center}
\end{definition}

\begin{remark}
	As a sanity check, we can calculate the composition $A\circ R$ in $\Set$, where $R\subseteq S\times S$ is the graph of some preorder $\leq$ on $S$. If $a\colon A\hookrightarrow S$ denotes the inclusion, then:
		\[
			A\times_S^{a,s}R\cong \{(a,(x,y)): a\in A, (x,y)\in R: a=x\},
		\]
	and so the image of this set under $(!\circ \pr_1,t\circ \pr_2)$, where $!\colon A\to 1$ is the terminal map, is
		\[
			A\circ R \cong \{y\in S: \exists a\in A: a\leq y\},
		\]
	which is of course just the usual upper cone $\up A$.
\end{remark}

\begin{lemma}\label{lemma:properties of internal cones}
	Let $\cat{D}$ be a category with finite limits and orthogonal factorisation system $(\calE,\calM)$, and let $R\in\MSub(X\times X)$ be a relation. Then:
		\begin{enumerate}[label = (\alph*)]
			\item if $A\subseteq B$ in $\MSub(X)$ then $\up A\subseteq \up B$ and $\down A\subseteq \down B$.
		\end{enumerate}
	Moreover, if $R$ is an $\calM$-preorder:
		\begin{enumerate}[label = (\alph*)]\setcounter{enumi}{1}
			\item $A\subseteq \up A$ and $A\subseteq \down A$;
			\item $\up\up A\subseteq \up A$ and $\down\down A\subseteq \down A$.
		\end{enumerate}
\end{lemma}
\begin{proof}
	We only sketch the ideas; the details are similar to those in the proofs found earlier in this section. For~(a), that $A\subseteq B$ implies there exists $m\colon A\to B$ such that $b\circ m = a$, and by universality of the pullback this lifts to a map $n\colon A\times_X^{a,s}R\to B\times_X^{b,s}R$, which with the factorisation system gives the desired inclusion $A\circ R\subseteq B\circ R$.
	
	For~(b), that $R$ is reflexive implies there exists an inclusion $m\colon X\to R$ such that $r\circ m = \delta\colon X\to X\times X$, the diagonal map. This ensures the existence of a canonical map $n\colon A\cong A\times_X^{a,\id_X}X \to A\times_X^{a,s}R$, which with the factorisation system provides the inclusion $A\subseteq A\circ R$.
	
	Lastly, for~(c), that $R$ is transitive implies there is a map ${\tau\colon R\times_X^{t,s}R\to R}$. Using associativity of pullbacks this gives a map $(A\times_X R)\times_X R\to A\times_X^{s,t}R$ that provides the desired inclusion $\up (A\circ R)\subseteq A\circ R$.
\end{proof}

These definitions of course go through in $\Loc$, so that for every internally ordered locale $(X,R)\in\ord(\Loc)$ we get internal cones that determine functions $\Sl(X)\to \Sl(X)$. This makes it possible to state a localic analogue of the open cone condition (cf.~\cref{definition:open cones}).

\begin{definition}
	An internally ordered locale $(X,R)\in\ord(\Loc)$ has \emph{open cones} if the internal cones preserve open sublocales.
\end{definition}

To close this section, we sketch a conjecture that relates the internal and external definitions of ordered locales, which uses the open cone assumption.

\begin{lemma}\label{lemma:open cones iff graph projections open}
	An ordered space $(S,\leq)$ has open cones if and only if the source and target projections $s,t\colon \graph(\leq)\to S$ are open maps.
\end{lemma}
\begin{proof}
	As an extremal subobject, the topology of $R:=\graph(\leq)$ is generated by the basis consisting of subsets of the form $U\times V\cap R$, where $U,V\in\Opens S$. It is then easy to see the direct images may be calculated as follows:
	\[
	s\left(U\times V\cap R\right)= U\cap \down V
	\qquad\text{and}\qquad
	t\left(U\times V\cap R\right) = \up U \cap V.
	\]
	From this, the equivalence of the two conditions immediately follows.
\end{proof}

A map of locales $f\colon X\to Y$ is called \emph{open} if there exists a left adjoint $f_!\dashv f^{-1}$ satisfying the \emph{Frobenius reciprocity condition} \cite{maclane1994SheavesGeometryLogic}:
	\[
		f_!(U\wedge f^{-1}(V)) = f_!(U)\wedge V.
	\]
Equivalently, $f$ is open precisely when the direct image of every open sublocale of $X$ is an open sublocale in $Y$ \cite[\S III.7]{picado2012FramesLocalesTopology}. The previous lemma motivates the following.

\begin{conjecture}\label{lemma:internally open cones iff projections open}
	An internally ordered locale $(X,R)$ has open internal cones if and only if the source and target maps $s,t\colon R\to X$ are open localic maps, and in that case
		\[
		\up \mathfrak{o}(U) = \mathfrak{o}(t_!\circ s^{-1}(U))
		\qquad\text{and}\qquad
		\down \mathfrak{o}(U) = \mathfrak{o}(s_!\circ t^{-1}(U)).
		\]
\end{conjecture}

\begin{remark}
	Recall the definition of open sublocales from \cref{definition:open and closed sublocales}. It is easy to see that if $R\hookrightarrow X\times X$ is the graph of a preorder $\leq$ on a space $S$, then $\up U = t(U\times S\cap R) = t(s^{-1}(U))$, and dually. We expect this to go through in the localic setting, but leave the proof for future research.
\end{remark}

Using this, we now describe a translation between internally and externally ordered locales.

\begin{construction}\label{construction:ordered locale from internal preorder}
	Fix an internally ordered locale $(X,R)\in\ord(\Loc)$. We define a preorder $\Leq$ on the lattice $\Sl(X)$ of sublocales of $X$ via
		\[
			A\Leq B
			\qquad\text{if and only if}\qquad
			A\subseteq \down B\text{~and~} B\subseteq \up A,
		\]
	where $\up$ and $\down$ are the internal cones on sublocales from \cref{definition:internal cones}. This is a preorder by \cref{lemma:properties of internal cones}.
	
	Suppose now that $(X,R)$ has open cones, which by \cref{lemma:internally open cones iff projections open} implies that~$\Leq$ restricts to a preorder on $\Opens X$. To be precise, we obtain an ordered locale $(X,\Leq)$ in the sense of \cref{definition:ordered locale} by setting:
		\[
			U\Leq V
			\qquad\text{if and only if}\qquad
			\mathfrak{o}(U)\subseteq \down \mathfrak{o}(V)\text{~and~}
			\mathfrak{o}(V)\subseteq \up \mathfrak{o}(U).
		\]
	In fact, by \cref{lemma:internally open cones iff projections open} we see that this restriction can be described by the functions $t_!\circ s^{-1}$ and $s_!\circ t^{-1}$:
		\[
			U\Leq V
		\qquad\text{if and only if}\qquad
		U\sqleq s_!\circ t^{-1}(V)
		\text{~and~}
		V\sqleq t_!\circ s^{-1}(U).
		\]
	Since $s_!$, $t_!$, $s^{-1}$ and $t^{-1}$ all preserve arbitrary joins (\cref{theorem:adjoints preserve limits}), $(X,\Leq)$ becomes an ordered locale with~\eqref{axiom:LV}.
\end{construction}

\begin{construction}\label{construction:internally ordered locale from EM ordered locale}
	Suppose now conversely that we have an ``externally'' ordered locale $(X,\Leq)\in\OrdLoc$. We sketch how an internal preorder ${R\rightarrowtail X\times X}$ might be obtained from $\Leq$. In the spatial case, observe that if $(S,\leq)$ has open cones, and $R$ represents the graph, then for $U,V,A,B\in\Opens S$ we have
		\[
			U\times V\cap R = A\times B \cap R
			\qquad\text{if and only if}\qquad
			\parbox{.35\textwidth}{$\Up U \cap V\subseteq B$, $U\cap \Down V\subseteq A$,\\
			$\Up A\cap B\subseteq V$, $A\cap \Down B\subseteq U$.}
		\]
	To define $R$ we therefore propose to generate a congruence relation on $\Opens X\oplus \Opens X$, uniquely characterising a sublocale of $X\times X$ (\cite[\S III.5]{picado2012FramesLocalesTopology}), as follows:
		\[
			U\oplus V \sim A\oplus B
			\qquad\text{if and only if}\qquad
			\parbox{.35\textwidth}{$\Up U \wedge V\sqleq  B$, $U\wedge \Down V\sqleq A$,\\
				$\Up A\wedge B\sqleq V$, $A\wedge \Down B\sqleq U$.}
		\]
	It is straightforward to verify that this defines an equivalence relation that preserves binary meets on the basis elements of $\Opens X\oplus \Opens X$. We conjecture that this extends to a congruence relation on the whole coproduct frame, and that the reflexivity and transitivity of $\Leq$ transfer to internal reflexivity and transitivity of the induced sublocale $R\rightarrowtail X\times X$.
\end{construction}

\begin{conjecture}\label{conjecture:internally ordered locales with OC are externally ordered locales}
	\cref{construction:ordered locale from internal preorder,construction:internally ordered locale from EM ordered locale} provide an equivalence of categories between the full subcategory of $\ord(\Loc)$ of internally ordered locales with open cones, and the full subcategory $\OrdLoc_\vee$ of ordered locales in the sense of \cref{definition:ordered locale} with~\eqref{axiom:LV}.
\end{conjecture}

\end{document}